\documentclass[11pt]{article}
\usepackage{graphicx}
\usepackage{amssymb,amsfonts,amsthm}
\usepackage{enumerate}
\usepackage[T2B]{fontenc}
\usepackage[cp1251]{inputenc}

\usepackage{enumitem}
\usepackage{amsmath,amsthm,amssymb}
\usepackage{amscd}
\usepackage{amsfonts}

\usepackage{xcolor}

\newcommand{\la}{\langle}
\newcommand{\ra}{\rangle}

\newtheorem{theorem}{Theorem}[section]
\newtheorem{lemma}[theorem]{Lemma}

\theoremstyle{definition}
\newtheorem{df}[theorem]{Definition}

\newtheorem{rk}[theorem]{Remark}
\newtheorem{prob}[theorem]{Problem}

\newcommand{\card}{\mathrm{card}}
\newcommand{\area}{\mathrm{Area}}

\newcommand{\s}{\mathbf{s}}
\newcommand{\ms}{\mathbf{S}}

\newcommand{\me}{\medskip}

\newcommand{\Lab}{{\mathrm{Lab}}}

\newcommand{\tool}{\stackrel{\ell}{\too} }

\newcommand{\rhh}{\mathbf{LR}}
\newcommand{\rhr}{\mathbf{RL}}
\newcommand{\mmm}{\mathbf{M}}
\newcommand{\mtt}{\mathbf{\Theta}}
\renewcommand{\lll}{{\mathcal L}}

\newcommand{\sbf}{\mathbf{S}}

\newcommand{\Z}{\mathbb Z}

\newcommand{\bmm}{\mathbf{\overline{\mmm}}}

\newcommand{\ttt}{{\mathcal T}}
\renewcommand{\tt}{{\tilde t}}

\newcommand{\ddd}{{\mathcal D}}

\newcommand{\bb}{{\mathcal B}}
\newcommand{\topp}{\mathbf{top}}
\newcommand{\ttopp}{\mathbf{ttop}}
\newcommand{\tbott}{\mathbf{tbot}}
\newcommand{\bott}{\mathbf{bot}}
\newcommand{\vk}{van Kampen }

\newcommand{\ccc}{{\mathcal C}}

\newcommand{\iv}{^{-1}}

\newcommand{\tosg}[1]{\stackrel{#1}{\to}}

\newcommand{\too}{\to }

\newcommand{\rrr}{{\mathcal R} }

\newcommand{\qq}{{\mathcal Q} }

\newcommand{\sss}{{\mathbf S} }

\usepackage{imakeidx}
\makeindex[name=g,title=Subject index]
\usepackage[columns=2]{idxlayout}

\usepackage[pagebackref, breaklinks]{hyperref}

    \renewcommand*{\backref}[1]{}
    \renewcommand*{\backrefalt}[4]{%
    \ifcase #1 %
        (Not cited).
    \or
        (Cited on page~#2).%
    \else
        (Cited on pages~#2).
    \fi}

\usepackage{fancyhdr,lipsum}
\begin{document}
\renewcommand{\theequation}{\thesection.\arabic{equation}}
\bigskip

\title{Conjugacy problem in groups with quadratic Dehn function}

 \author{A.Yu. Olshanskii, M.V. Sapir\thanks{Both authors were supported in part by the NSF grants DMS-1500180 and DMS-1901976. The first author was also supported by RFFI grant 15-01-05823}}

\date{}
\maketitle

\renewcommand\rightmark{[Short Title]}

\begin{abstract}
We construct a finitely presented group with quadratic Dehn function and undecidable conjugacy problem. This solves E. Rips' problem formulated in 1994.
\end{abstract}


\medskip


\setcounter{tocdepth}{2}

\tableofcontents

\section{Introduction}

Every group given by a presentation $G=\la X \mid R \ra$ is a factor group $F/N$ of the free group $F=F(X)$
with the set of free generators $X$ over the normal closure $N=\la\la R\ra\ra^F$ of the set of relators $R$. Therefore every word $w$ over
the alphabet $X^{\pm 1}$ vanishing in $G$ represents an element of $N$, and so in $F$, $w$ is a product $v_1\dots v_m$
of factors $v_i=u_ir_i^{\pm 1}u_i^{-1}$ which are  conjugate to the relators $r_i\in R$ or their inverses.

The minimal number of factors $m=m(w)$ is called the {\it area of the word $w$} with respect to the presentation
$G=\la X | R \ra $. M. Gromov \cite{GrHyp} introduced this concept and term in geometric group theory, because $m$ is equal to the
minimal number of $2$-cells (counting with multiplicities) used in a $0$-homotopy of the path $\bf p$ labeled by $w$ in the
Cayley complex of the presentation of $G$ (or the $0$-homotopy of a
singular disk with boundary $\bf p$).

In other words, given equality $w=1$ in $G$, one can construct a van Kampen diagram, that is a finite, connected graph
on Euclidean plane with $m$ bounded regions, where every edge has label from $X^{\pm 1}$, the boundary path of
every region (= $2$-cell) is therefore labelled,  the label of it belongs in $R^{\pm 1}$, and the boundary of  the whole
map is labelled by $w$. (See more details for this visual definition of area and van Kampen diagram in Section \ref{md}.)



The \index[g]{Dehn function of a finitely presented group}
Dehn function of a finitely presented group $G=\la X\mid R\ra$ is the smallest function $f(n)$ such that for every word $w$ of length at most $n$ in the alphabet $X^{\pm 1}$, which is equal to $1$ in $G$, the area of $w$ is at most $f(n)$.

It is well known \cite{GrHyp, Gr} that the Dehn functions of different finite presentations of the same finitely presented group are equivalent, where we call two functions $f(n), g(n)$ equivalent if for some constants $A,B,C,D\ge 1$ and every $n\ge 0$, we have

$$ \frac 1{A} f\left(\frac {n}B\right)-Cn-D < g(n) < Af(Bn)+Cn+D.$$

As usual, we do not distinguish equivalent functions.

The Dehn function is an important invariant of a group for the following reasons.

{\bf A}) It almost obviously follows from the definition that if $G$ is the fundamental group of a compact Riemannian manifold $M$ then the Dehn function of $G$ is equivalent to the smallest isoperimetric function of the universal cover $\tilde M$.

{\bf B}) From the Computer Science point of view, the Dehn function of a group $G$ is equivalent to the time function of a non-deterministic Turing machine "solving" the word problem in $G$ (see \cite[Introduction]{SBR} for details). Moreover as was  shown in \cite{BORS}:

 \begin{quote}
 A not necessarily finitely presented finitely generated group has word problem in $\mathbf{NP} $ if and only if it is a subgroup of a finitely presented group with at most polynomial Dehn function (a similar result holds for other computational complexity classes \cite{BORS}).
\end{quote}

 Papers \cite{SBR, O18, S17} provide, modulo the famous conjecture $\mathbf{P=NP}$, a complete description of all real numbers $\alpha$ such that $n^\alpha$ is equivalent to the Dehn function of a finitely presented group.

{\bf C}) From the geometric point of view the Dehn function measures the "curvature" of the group: linear Dehn functions correspond to negative curvature, quadratic Dehn function correspond to non-positive curvature, etc.

More precisely, a finitely presented group is hyperbolic if and only if it has a subquadratic (hence linear) Dehn function \cite{GrHyp, Bow, Ol91}. In particular, the conjugacy problem in such groups is decidable \cite{GrHyp}.

It is also known that groups with quadratic Dehn functions exhibit certain "non-generic" non-positive curvature behavior as far as geometric and algorithmic properties are concerned. For example their asymptotic cones are simply connected \cite{Pap}. The conjugacy problem is the second of Dehn's basic algorithmic problems. The word problem in groups with quadratic (or any recursive) Dehn function is  decidable.   Thus the following question is very natural:

\begin{prob}[Rips]\label{p:11} Does every finitely presented group with quadratic Dehn function have decidable conjugacy problem?
\end{prob}

Rips explicitly asked this question in his talk at the DYMACS Workshop: Geometric Group Theory in 1994. At that time the largest known class of groups with at most quadratic Dehn functions was the class of bi-automatic groups \cite{Eps}. Note that the decidability of the conjugacy problem for automatic groups is still not known; it is not even known if every automatic group is bi-automatic.

Since 1994, many important groups were proved to have quadratic Dehn function. For each of these groups it was proved (or was already known) that the conjuacy problem is decidable. Here is (an
incomplete) list of these groups.

\begin{itemize}
\item $SL_n(\Z), n\ge 5$; quadratic Dehn function: \cite{Young}; solvability of the conjugacy problem: \cite{Sar, Grunewald},

\item High rank integral Heisenberg groups $H_{n}$: \cite{Al,OS} and \cite{Black},

\item Many metabelian (non-nilpotent) groups, some of them containing Baumslag-Solitar groups $BS(1,p)$, and some groups that are obtained by using the Baumslag-Re\-mes\-len\-ni\-kov construction: \cite{Drutu,CT} and \cite{Noskov},

\item The R. Thompson group $F$:    \cite{Guba} and \cite{GubSa},

\item groups acting geometrically on CAT(0) spaces: \cite{Kokarev} and \cite{BridsonHaefliger},
\item free-by-cyclic groups: \cite{BrGr} and \cite{OS06}.
\end{itemize}

The decidability of conjugacy problem was proved in a completely different way in each of these cases and it is natural to ask whether every group with quadratic Dehn function has decidable conjugacy problem and there is a uniform proof of that fact. That made Problem \ref{p:11} even more intriguing.
In fact Rips had a "quasi-proof" showing that the answer should be positive. That "quasi-proof" first appeared in \cite{OS06}.
Basically the idea is the following (see details in \cite{OS06}).

  Suppose the conjugacy problem is undecidable in a finitely presented group $G=\la X\mid R\ra$.
  This implies that for arbitrary recursive function $f$, there are  infinitely many pairs of words $(u,v)$ in the alphabet $X^{\pm 1}$, such that $v=zuz^{-1}$ in $G$, but the length $||z||$ of the shortest word $z$ with this property exceeds $f(n)+n$ for $n=||u||+||v||$. Therefore every minimal
  area annular diagram $\Delta$ with boundary labels $u$ and $v$ has no path of length $\le f(n)$ connecting the two components of the boundary.
  Let $\bf q$ be a simple path connecting the boundaries of $\Delta$, $t=||\bf q||$. Then there are simple closed paths ${\bf p}_1,...,{\bf p}_m$ of $\Delta$ surrounding the hole such that ${\bf p}_i,..., {\bf p}_j$ do not intersect if $i\ne j$ and $m> c_1 t$ for some constant $c_1$. The area of $\Delta$ is at least a constant times $\sum ||{\bf p}_i||$. If "many" lengths $||{\bf p}_i||$ are less than $c\log t$ where $c=\frac 1{2|X|}$, then two of the paths ${\bf p}_i$, ${\bf p}_j$ ($i\ne j$) have the same labels. That allows us to identify ${\bf p}_i, {\bf p}_j$ and remove the annular subdiagram of $\Delta$ bounded by ${\bf p}_i, {\bf p}_j$, decreasing the area of $\Delta$, a contradiction. Therefore "many" lengths $||{\bf p}_i||$ are at least $c_2\log t$ for some constant $c_2$. Hence the area of $\Delta$ is at least $c_3t\log t$ for some constant $c_3$. If we cut $\Delta$ along the path ${\bf q}$, we obtain a disk van Kampen diagram $\Delta'$ with boundary path subdivided into four parts ${\bf q}_1{\bf p}_1{\bf q}_2\iv {\bf p}_2\iv$ where $||{\bf p}_1||, ||{\bf p}_2||< n$ and the labels of ${\bf q}_1$ and ${\bf q}_2$ coincide with the label of ${\bf q}$. The area of $\Delta'$ is at least $c_3 t\log t$. Since the labels of ${\bf q}_1, {\bf q}_2$ are the same, we can glue $t/n$ copies of $\Delta'$ together to obtain a van Kampen diagram $\Delta''$ with perimeter bounded from above by a linear function in $t$ and area bounded below by $c_3 t^2\log t/n$ since $t$ is bounded from below by  given recursive function in $n$, $n$ is insignificant compared to $t$. The diagram $\Delta''$ can be assumed reduced. So we found a reduced van Kampen diagram of perimeter $\sim t$ and area $\sim t^2\log t$.
Hence the Dehn function cannot be smaller than $n^2\log n$.

The incorrectness of this "quasi-proof" is in the last phrase. Indeed, there may be a smaller area van Kampen diagram with the same boundary label as $\Delta''$. Still there is a lot of flexibility in choosing $\Delta$ and the path $\bf q$ in it. It looks like it would require infinite number of defining relations to ensure that all the boundary paths of various diagrams $\Delta''$ have fillings with much fewer cells than $\Delta''$.  In particular, if $G$ satisfies some mild form of asphericity, the proof should work. We conjectured that this is true for all finitely presented groups \cite{OS06}. In \cite{OS06}, we confirmed this conjecture for a wide class of multiple HNN extensions of free groups. We also constructed in \cite{OS06} a multiple HNN extension of a free group with undecidable conjugacy problem and the minimal possible Dehn function $n^2\log n$.

Nevertheless, in this paper, we give a negative answer to Rips' question (and hence disprove our conjecture as well):

\begin{theorem}\label{t:main} There exists a finitely presented group with undecidable conjugacy
problem and quadratic Dehn function.
\end{theorem}

As in several of our previous papers (\cite{SBR, BORS, OS06, O18} and others) the construction is based on an $S$-machine (we call it $\mmm_5$) which can be viewed as a computing device with undecidable halting problem or as a group which is a multiple HNN extension of a free group. $S$-machines were first introduced by Sapir in \cite{SBR} (see Section \ref{SM} below for the definition used here and  \cite{S} for various other definitions).

In order to describe some ideas of our proof in more details, let us start with a simple example of an $S$-machine $\sbf$. (That $S$-machine first appeared in \cite{OS6}. The corresponding group was the first example of a group with polynomial Dehn function, linear isodiametric function and non-simply connected asymptotic cones answering a question of C. Dru\c tu.) It is a rewriting system \cite{Sbook} with alphabet $\{a,q, a\iv, q\iv\}$ and two "same" rules
$\theta_i\colon q\to aq$ and their inverses $\theta_i\iv\colon q\to a\iv q$, $i=1,2$. The rewriting system works with group words in $\{a,q\}$. And applying a rule $\theta_i^{\pm 1}$ means replacing every letter $q^{\epsilon}$ (where $\epsilon=\pm 1$) by $(a^{\pm 1})q^\epsilon$ and then reducing the word. The $S$-machine $\sbf$ can also be viewed as a multiple HNN extension of the free group $\la a,q\ra$: $$\la a,q, \theta_1, \theta_2\mid q^{\theta_i}=aq, a^{\theta_i}=a, i=1,2\ra.$$
(Note that this is far from the only way to interpret $S$-machines as groups. We are using a different interpretation in this paper, and the most complicated one so far was used in \cite{OS04}. But the main principle is still the same.)

\begin{figure}[ht]
\unitlength .47 mm 
\linethickness{0.3pt}
\ifx\plotpoint\undefined\newsavebox{\plotpoint}\fi 
\begin{picture}(472.707,141)(0,0)
\put(63.75,9.5){\makebox(0,0)[cc]{$q\iv$}}
\put(63.75,41){\makebox(0,0)[cc]{$q\iv$}}
\put(167.75,40.75){\makebox(0,0)[cc]{$q$}}
\put(167.75,73.75){\makebox(0,0)[cc]{$q$}}
\put(104.75,73.75){\makebox(0,0)[cc]{$a$}}
\put(63.75,74){\makebox(0,0)[cc]{$q\iv$}}
\put(63.75,106.25){\makebox(0,0)[cc]{$q\iv$}}
\put(63.75,138){\makebox(0,0)[cc]{$q\iv$}}
\put(105,9.5){\makebox(0,0)[cc]{$a$}}
\put(163.25,9.5){\makebox(0,0)[cc]{$q$}}
\put(221.5,9.5){\makebox(0,0)[cc]{$a$}}
\put(243.5,40.75){\makebox(0,0)[cc]{$a$}}
\put(198.75,40.5){\makebox(0,0)[cc]{$a$}}
\put(279.75,9.5){\makebox(0,0)[cc]{$q$}}
\put(279.75,41){\makebox(0,0)[cc]{$q$}}
\multiput(43,38.75)(16.8333333,-.0333333){15}{\line(1,0){16.8333333}}
\multiput(43,71)(16.8333333,-.0333333){15}{\line(1,0){16.8333333}}
\multiput(42.5,103.75)(16.8333333,-.0333333){15}{\line(1,0){16.8333333}}
\multiput(42.5,136.5)(16.8333333,-.0333333){15}{\line(1,0){16.8333333}}
\multiput(43,6)(16.8333333,-.0333333){15}{\line(1,0){16.8333333}}
\put(43,6.25){\circle{1.414}}
\put(43,71){\circle{1.414}}
\put(43,38.25){\circle{1.414}}
\put(43,103.25){\circle{1.414}}
\put(43,135.25){\circle{1.414}}
\put(76.25,39.25){\circle{1.414}}
\put(76.25,71.25){\circle{1.414}}
\put(76.25,136.25){\circle{1.414}}
\put(76.25,103.25){\circle{1.414}}
\put(137.75,122.75){\circle{1.414}}

\put(137.25,91.75){\circle{1.414}}
\put(226.5,38.75){\circle{1.414}}
\put(226.5,71){\circle{1.414}}
\put(191.5,38.75){\circle{1.414}}
\put(191.5,71){\circle{1.414}}
\put(261,38.75){\circle{1.414}}
\put(261,71){\circle{1.414}}
\put(261.75,103){\circle{1.414}}
\put(223,103.75){\circle{1.414}}
\put(258,135.75){\circle{1.414}}
\put(209.5,70){\circle{1.414}}
\put(76.25,5.75){\circle{1.414}}
\put(137.25,5.75){\circle{1.414}}
\put(137.25,38){\circle{1.414}}
\put(137.25,70.25){\circle{1.414}}
\put(137.25,103.25){\circle{1.414}}
\put(137.25,136){\circle{1.414}}
\put(192.5,5.75){\circle{1.414}}
\put(255.5,6){\circle{1.414}}
\put(42.5,6){\vector(0,1){32.75}}
\put(42.5,38.25){\vector(0,1){32.75}}
\put(42,71){\vector(0,1){32.75}}
\put(42,103.75){\vector(0,1){32.75}}
\put(76,6.25){\vector(0,1){32.75}}
\put(76,39.25){\vector(0,1){32.75}}
\put(192.25,6.5){\vector(0,1){32.75}}
\put(294.75,5.75){\vector(0,1){32.75}}
\put(294.75,38){\vector(0,1){32.75}}
\put(294.25,70.75){\vector(0,1){32.75}}
\put(294.25,103.5){\vector(0,1){32.75}}
\multiput(76.25,6)(-.03125,2.3125){8}{\line(0,1){2.3125}}
\multiput(76.25,39)(-.03125,2.3125){8}{\line(0,1){2.3125}}
\multiput(137,5.75)(-.06256436663,.03372811535){971}{\line(-1,0){.06256436663}}
\multiput(255.25,6.25)(-.03371428571,.03714285714){875}{\line(0,1){.03714285714}}
\put(76,15.25){\vector(0,1){7.5}}
\put(76,48.25){\vector(0,1){7.5}}
\multiput(75.75,31.75)(.0333333,.4){15}{\line(0,1){.4}}
\multiput(75.75,64.75)(.0333333,.4){15}{\line(0,1){.4}}
\multiput(215.25,38.25)(1.1875,.03125){8}{\line(1,0){1.1875}}
\multiput(215.25,70.5)(1.1875,.03125){8}{\line(1,0){1.1875}}
\multiput(246.75,38.25)(1.21875,.03125){8}{\line(1,0){1.21875}}
\multiput(246.75,70.5)(1.21875,.03125){8}{\line(1,0){1.21875}}
\put(281.25,38.25){\vector(1,0){10}}
\put(281.25,70.5){\vector(1,0){10}}
\put(280.75,103.25){\vector(1,0){10}}
\put(280.75,136){\vector(1,0){10}}

\put(232,32.25){\vector(-1,1){1.4142}}

\put(191.5,38.25){\vector(0,1){32.5}}
\put(191.5,71.25){\vector(0,1){32.5}}
\put(191.5,104){\vector(0,1){32.5}}

\multiput(43.25,38.75)(-.0333333,-.0333333){15}{\line(0,-1){.0333333}}
\put(174.5,70.75){\vector(1,0){7.5}}
\multiput(261,38.5)(-.03578838174,.03371369295){964}{\line(-1,0){.03578838174}}
\multiput(226,38.75)(-.0337022133,.0643863179){497}{\line(0,1){.0643863179}}
\multiput(226.25,71.25)(.03768577495,.03370488323){942}{\line(1,0){.03768577495}}
\multiput(208.5,70.75)(.0336879433,.0774231678){423}{\line(0,1){.0774231678}}
\multiput(222.75,103.5)(.03604771784,.03371369295){964}{\line(1,0){.03604771784}}
\put(36.75,23.5){\makebox(0,0)[cc]{$\theta_1$}}
\put(36.75,55.75){\makebox(0,0)[cc]{$\theta_2$}}
\put(36.75, 86.75){\makebox(0,0)[cc]{$\theta_1\iv$}}
\put(36.75,121){\makebox(0,0)[cc]{$\theta_2\iv$}}
\put(81.75,46.5){\makebox(0,0)[cc]{$\theta_2$}}
\put(114.25,53){\makebox(0,0)[cc]{$\theta_2$}}
\put(81.75,11.75){\makebox(0,0)[cc]{$\theta_1$}}
\put(112.75,23.5){\makebox(0,0)[cc]{$\theta_1$}}
\put(198.25,122){\makebox(0,0)[cc]{$\theta_2\iv$}}
\put(198.25,86.75){\makebox(0,0)[cc]{$\theta_1\iv$}}
\put(198.25,50.75){\makebox(0,0)[cc]{$\theta_2$}}
\put(198.5,19.25){\makebox(0,0)[cc]{$\theta_1$}}
\put(251.25,18.5){\makebox(0,0)[cc]{$\theta_1$}}
\put(300.75,18){\makebox(0,0)[cc]{$\theta_1$}}
\put(250.25,53.5){\makebox(0,0)[cc]{$\theta_2$}}
\put(301,53.5){\makebox(0,0)[cc]{$\theta_2$}}
\put(223.75,53.25){\makebox(0,0)[cc]{$\theta_2$}}
\put(226.5,85.5){\makebox(0,0)[cc]{$\theta_1\iv$}}
\put(251.25,85.5){\makebox(0,0)[cc]{$\theta_1\iv$}}
\put(301,85.75){\makebox(0,0)[cc]{$\theta_1\iv$}}
\put(304.25,122.75){\makebox(0,0)[cc]{$\theta_2\iv$}}
\multiput(137,38.25)(-.06282183316,.03372811535){971}{\line(-1,0){.06282183316}}
\put(117.5,70.75){\vector(1,0){16.75}}
\put(112.25,38.5){\vector(1,0){15}}
\put(104.25,41.5){\makebox(0,0)[cc]{$a$}}
\put(79.25,62){\makebox(0,0)[cc]{$a$}}
\put(79.25,29){\makebox(0,0)[cc]{$a$}}
\put(137,103.75){\vector(0,-1){12}}
\put(137,70.5){\vector(0,1){21.5}}
\multiput(75.75,70.75)(.0632780083,.03371369295){964}{\line(1,0){.0632780083}}
\multiput(75.5,103.25)(.06288032454,.03372210953){986}{\line(1,0){.06288032454}}
\put(137.25,103.75){\vector(0,1){19.5}}
\put(137,136.75){\vector(0,-1){13.5}}
\put(50.5,136.5){\vector(1,0){17.5}}

\put(50.5,38.75){\vector(1,0){17.5}}
\put(50.5,70.75){\vector(1,0){17.5}}
\put(50.5,103.75){\vector(1,0){17.5}}
\put(89,136.5){\vector(1,0){15}}
\put(102.25,140){\makebox(0,0)[cc]{$a$}}
\put(138.75,129.5){\makebox(0,0)[cc]{$a$}}
\put(114.75,107.25){\makebox(0,0)[cc]{$a$}}
\put(119.5,103.5){\vector(1,0){10.75}}
\put(146.25,115){\makebox(0,0)[cc]{$\theta_2\iv$}}

\put(109.25,119.25){\makebox(0,0)[cc]{$\theta_2\iv$}}

\put(100,116.5){\vector(3,2){.7}}
\put(236.75,116.5){\vector(1,1){.7}}

\put(100,83.75){\vector(3,2){.7}}

\put(240,83.75){\vector(3,2){.7}}

\put(214.25,83.75){\vector(2,3){.7}}

\put(214.25,60.25){\vector(-2,3){.7}}

\put(237.75,60.25){\vector(-2,3){.7}}

\put(100,58.50){\vector(-3,2){.7}}

\put(100,26.00){\vector(-3,2){.7}}

\put(89,136.5){\vector(1,0){15}}

\put(118,86){\makebox(0,0)[cc]{$\theta_1\iv$}}
\put(146.25,80.75){\makebox(0,0)[cc]{$\theta_1\iv$}}
\put(139.75,98.5){\makebox(0,0)[cc]{$a$}}
\put(165.25,107.5){\makebox(0,0)[cc]{$q$}}
\put(168.5,141){\makebox(0,0)[cc]{$q$}}
\put(253.5,123.25){\makebox(0,0)[cc]{$\theta_2\iv$}}
\put(207.75,107.5){\makebox(0,0)[cc]{$a$}}
\put(246,107.25){\makebox(0,0)[cc]{$a$}}
\put(280.75,106.5){\makebox(0,0)[cc]{$q$}}
\put(199.5,73.5){\makebox(0,0)[cc]{$a$}}
\put(219.75,74.75){\makebox(0,0)[cc]{$a$}}
\put(244.25,74){\makebox(0,0)[cc]{$a$}}
\put(280,74.5){\makebox(0,0)[cc]{$q$}}
\put(227,139.25){\makebox(0,0)[cc]{$a$}}
\put(278.5,139.25){\makebox(0,0)[cc]{$q$}}
\put(63.92,6.243){\vector(1,0){.07}}\multiput(57.082,5.946)(.759771,.033034){9}{\line(1,0){.759771}}
\put(98.407,5.946){\vector(1,0){9.216}}

\put(182.054,5.886){\vector(1,0){.07}}\multiput(176.168,5.676)(.840896,.030032){7}{\line(1,0){.840896}}

\put(182.054,38.25){\vector(1,0){.07}}\multiput(176.168,5.676)(.840896,.030032){7}{\line(1,0){.840896}}

\put(182.054,103.25){\vector(1,0){.07}}\multiput(176.168,5.676)(.840896,.030032){7}{\line(1,0){.840896}}
\put(182.054,136.25){\vector(1,0){.07}}\multiput(176.168,5.676)(.840896,.030032){7}{\line(1,0){.840896}}
\put(247.434,5.676){\vector(1,0){5.256}}
\multiput(286.75,5.25)(.71875,.03125){8}{\line(1,0){.71875}}
\put(292.5,5.5){\vector(1,0){.07}}\multiput(286.75,5.25)(.71875,.03125){8}{\line(1,0){.71875}}

\put(294.25,135.35){\circle{1.414}}
\put(294.25,103.75){\circle{1.414}}
\put(294.25,70.75){\circle{1.414}}

\put(294.25,38.5){\circle{1.414}}
\put(294.25,5.75){\circle{1.414}}
\put(192,135.35){\circle{1.414}}
\put(192,103.75){\circle{1.414}}

\put(235.5,38.25){\vector(1,0){5.5}}

\put(200.5,38.25){\vector(1,0){5.5}}

\put(235.5,70.75){\vector(1,0){5.5}}
\put(200.5,70.75){\vector(1,0){5.5}}
\put(217.5,70.75){\vector(1,0){5.5}}

\put(200.5,103.60){\vector(1,0){5.5}}
\put(235.5,103.60){\vector(1,0){5.5}}

\put(235.5,136.25){\vector(1,0){5.5}}
\end{picture}

\begin{center}
\caption{The trapezium corresponding to a computation of $\sbf$.}
\end{center}
\label{p:1}
\end{figure}
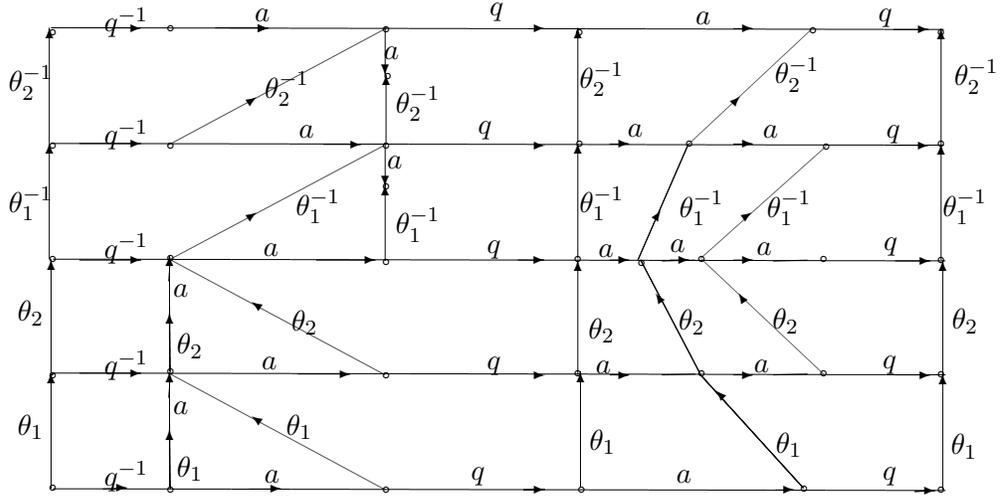

As the name $S$-machine suggests, we can also consider $\sbf$ as a kind of Turing machine with tape letter $a$, state letter $q$ and commands $\theta_1, \theta_2$ (and their inverses). Then we can consider {\em computations}. Say,
\begin{equation}\label{es}q\iv aq aq\tosg{\theta_1} q\iv aq aaq\tosg{\theta_2} q\iv aqaaaq\tosg{\theta_1\iv }q\iv aqaaq\tosg{\theta_2\iv}q\iv aqaq\end{equation} is a {\em reduced} computation of $\sbf$. At the same time if we consider $\sbf$ as a multiple  HNN-extension $S$ of the free group, then this computation corresponds to the van Kampen diagram on Figure \ref{p:1}.

This diagram is called the {\em trapezium} corresponding to the computation (\ref{es}). Three things need to be noticed from this diagram.

1. The trapezium looks like a rectangle with the first word and the last word of the computation on the bottom and top sides. All other words of the computation are on the horizontal paths of the trapezium, and $\theta$'s conjugate each of these words to the next one.

2. The vertical sides of the trapezium are labeled by the same words: the {\em history} of the computation (in the case of (\ref{es}) it is $\theta_1\theta_2\theta_1\iv\theta_2\iv$).

3. The trapezia has three types of {\em bands} (also called in the literature {\em corridors}), i.e. sequences of cells  where each two consecutive cells share an edge with a prescribed  label: horizontal $\theta_i^{\pm 1}$-bands, vertical $q$-bands and $a$-bands. The median lines of these bands serve as "walls" in van Kampen diagrams over $S$-machines, provide necessary rigidity and are crucial for all applications of $S$-machines.

Now let us give some of the main ideas  of the construction and proof of Theorem \ref{t:main} (the actual construction, see Section \ref{am}, is somewhat different and employs different notation).

We start with any Turing machine $\mmm_0$ with one input tape where the input word is written in unary, as a power $\alpha^k$ where $\alpha$ is a tape letter, and non-recursive language of accepted input words.

We convert $\mmm_0$ into an $S$-machine $\ms$ (in the notation below this is $\mtt_3\cup \mtt_4\cup \mtt_5)$. As any $S$-machine, $\ms$ can be viewed as a group, $S$. It has three kinds of generators: $Y$-letters (or tape letters),
$q$-letters (or state letters) and $\theta$-letters (rule letters).

The set $Y$ contains the "input letter" $\alpha$ (as for $\mmm_0$) and several "historic" sets of letters, $Y_1,Y_2,...$. The $S$-machine $\ms$ has the following properties:

\begin{itemize}
\item[(S1)] every configuration of $\ms$ has several "historic sectors" (a sector is a subword that starts and ends with $q$-letters, and all other letters from $Y$);

\item[(S2)] any starting configuration of $\ms$ has no $Y$-letters except the input letter $\alpha$ which can occur in the input sector only (in fact, there are two input sectors but that is not a significant details); if the input sectors contains $\alpha^k$, then the start word is denoted by $W(k)$;

\item[(S3)] there is only one accept configuration of $\ms$, $W_{ac}$. It has no $Y$-letters;

\item[(S4)] if a computation of $\ms$
$$U_1\to U_2\to ...\to U_m$$ is "very long", then there is a computation
$U_1\to ...\to W(k)$ and $U_m\to ..\to W_{ac}$ of length at most a constant multiple of the length of the first word in the computation.
Moreover in the first case the number $k$ is determined by the history of computation;

\item[(S5)] the machine ${\bf \Theta}_3$
inserts the history subwords in all historical sectors;
the computation of the main  machine ${\bf \Theta}_4$ proceeds by executing the rules of the history subwords in non-historic
parts of the configurations;
so if machine ${\bf \Theta}_4$ accepts the input, then ${\bf \Theta}_5$ removes the history from the historic
parts of the configurations, and accepts;

\item[(S6)] $\ms$ takes $W(k)$ to $W_{ac}$ if and only if  $\alpha^k$ is accepted by $\mmm_0$; in particular the language of words $W(k)$ accepted by $\ms$ is not recursive.
\end{itemize}

By \cite{OS06}, the Dehn function of the group $S$ corresponding to $\ms$
is at least $n^2\log n\gg n^2$.
But we prove in the paper that most of the area in van Kampen diagrams of large area over $\ms$ is concentrated in a few standard trapezia which we call {\em big trapezia}. The phenomenon that large part of the area of a van Kampen diagram is concentrated in a few large standard subdiagrams, is interesting and seems to be very common. For example, we proved similar facts for van Kampen diagrams over presentations satisfying the small cancelation condition $C(p)-T(q)$ in the "CAT(0)" case $\frac 1p+\frac 1q=\frac 12$ in \cite{OS-17}. In that case the geometric meaning of existence of large standard subdiagram is very close to a popular topic in CAT(0) geometry: "every quasi-flat in the universal cover of the presentation complex is close to a flat" (see a discussion in \cite{OS-17}). In the case of $S$-machines, we proved similar facts in \cite{OS06} and \cite{O18}, in both cases, as in the present paper, these were crucial steps in the proofs.

Big trapezia over $S$ must correspond to "very long" computations of $\ms$.
Then we augment the $S$-machine $\ms$ by a new $S$-machine $\ms'$ to
obtain an $S$-machine $\mmm$.
As a group $M$ the $S$-machine $\mmm$ is
still a multiple HNN extension of a free group.
The group $G$ is obtained by imposing one relation $W_{ac}^L$ where $L\gg 1$. This relation is called
the hub. The hubs and the disks (that are hubs surrounded by $\theta$-annuli) make the areas of trivial in $M$ words quadratic with respect of the significantly larger \cite{SBR} presentation of $G$ (another important new idea: usually, disks make Dehn functions bigger \cite{SBR}).
 Therefore the presentation of $G$ is highly non-aspherical: the boundaries of the large trapezia can be filled both by diagrams with $\sim n^2\log n$ cells and by diagrams with at most $\sim n^2$ cells.

The new $S$-machine $\mmm$ is obtained by augmenting $\ms$ with two simple $S$-machines $\mtt_1$ and $\mtt_2$
(see Subsection \ref{M6}).
The $S$-machine $\ms'=\mtt_{1}\cup \mtt_2$ starts working with the word $W(0)$ and inserts $\alpha^k$ in the input sector


This augmentation provides us with the property that arbitrary configuration of a ``long computation'' of $\mmm$ can be reached with linear time and space either from $W(0)$  or from the stop configuration of $\mmm$. Afterwards this linearity guaranties quadratic estimates of the areas of both disks and big trapezia over the presentation of $G$. The linearity is achieved by, in particular, adding many so-called history sectors where the history of a computation is non-deterministically written before the actual computation executing that history starts.

In order to connect $\ms'$ with the $S$-machine $\ms$ and obtain the main $S$-machine $\mmm$, we need one rule, denoted $\theta(23)$ which changes the state letters to the start state letters  of $\ms$. However the standard interpretation of $\mmm$ as a group would make the conjugacy problem decidable in the group $M$. So the rule $\theta(23)$ is interpreted in $G$ as turning $L$  words in different alphabets into $L$ input configurations of $\ms$ in the same alphabet (by erasing extra indices). This new ``irregular'' interpretation requires a study of some non-reduced (eligible) computations, i.e., the history of an ``eligible'' computation may contain (many) subwords $\theta(23)\;\theta(23)^{-1}$.

The proof that $G$ has quadratic Dehn function is much harder than the proof of undecidability of the conjugacy problem.  We use several tools developed in \cite{SBR, OS06, O12, O18} and more. As in all our papers where estimates of the Dehn function are produced, we need to consider diagrams with and without hubs separately. This is done in Sections \ref{6} and \ref{midi} respectively. In both cases, one of the main ideas is to assign to the boundary of every van Kampen diagram $\Delta$ over the presentation of $G$ certain numeric invariant $\mu(\Delta)$ (the {\em mixture} on the boundary $\partial\Delta$ from \cite{O18})  which is bounded from above by a quadratic function in terms of the perimeter. We had a somewhat similar numeric invariant called {\em dispersion} in \cite{OS06} but that invariant does not work for diagrams with hubs.

To obtain a quadratic estimate for diagrams $\Delta$ over $M$, we have to consider an artificial $G$-areas instead of areas, and just at the end of this paper we replace the diagrams of quadratic $G$-area over $M$ with diagrams with hubs, having quadratic (usual) areas over $G$. The quadratic upper bound for $G$-area is obtained by induction over the (modified) perimeter $n$ of $\Delta$. We perform surgeries on the diagram, so that each surgery makes the diagram look more "standard" and smaller. Our inductive argument estimates the $G$-area in terms of some linear combination of $n^2$ and  the mixture $\mu(\Delta)$. Although we are not able to choose just one of these two summands for induction, the final upper bound of the $G$-area is $O(n^2)$, because of the aforementioned quadratic estimate of the mixture in terms of $n$.

In the case of diagrams with hubs, we estimate a similar linear combination,
but the inductive parameter is not the (modified) perimeter $n$ but the sum
$\Sigma=n +\sigma(\Delta)$. The invariant  $\sigma(\Delta)=\sigma_{\lambda}(\Delta)$
was invented in \cite{O18}. It is defined by the design formed by maximal
bands of two types in $\Delta$. The important and non-trivial feature of the
$\sigma$-invariant is  the linear inequality $\sigma(\Delta)=O(n)$,
and so the quadratic upper bound of the form $O(\Sigma^2)$ is also
quadratic in terms of the perimeter $n$.

In fact in both cases (over $M$ or over $G$), the proof proceeds by taking a minimal counterexample diagram $\Delta$ and then perform surgeries producing smaller diagrams which
cannot be counterexamples. This provides more and more useful information about $\Delta$, until finally we get a contradiction.

For instance, in Section \ref{midi} where diagrams with hubs are considered, we need to remove one of the disks from the diagram. As in our previous papers (starting with \cite{SBR}, \cite{BORS,O97}), we use hyperbolicity of certain graph associated with hubs (hubs are vertices, $q$-bands connecting hubs are edges), and find a hub connected to the boundary of the whole diagram by almost all bands starting on the hub 2-cell. This gives a subdiagram of $\Delta$ consisting of a subdiagram without hubs, called a clove, and a disk. We would like to remove that subdiagram from $\Delta$ producing a diagram $\Delta'$ with smaller parameter $n'+\sigma(\Delta')$.

A similar task was solved in \cite{SBR}. It is one of the most non-trivial parts of \cite{SBR}.  Using it, we decomposed a diagram in \cite{SBR} into a few disks of small total perimeter, and a diagram without hubs. This was called the {\em snowman decomposition}. But that task is now much harder than in \cite{SBR}. The reason is that in \cite{SBR}, after removing the clove and the disk,  we needed to show that the perimeter of the diagram decreases and the perimeter of the removed disk (only the disk) is linearly bounded by the difference of the perimeters of the old and new diagrams. For the quadratic upper bound this is not enough. We need to get a linear lower bound of the difference in terms of the whole piece that we cut off (the clove and the disk). That can be achieved not always. If not, we get a new information about the disk and the clove and remove the disk together with a certain sub-clove.  The mixture and the $\sigma_\lambda$ invariant help achieve it at the end.

Some estimates used in this paper are very similar to the estimates in \cite{O18}. More precisely for every function $f(n)$ satisfying certain conditions, a finitely presented group $G_f$ with Dehn function $n^sf(n)^3$ (where $s\ge2$) is constructed in \cite{O18}. In particular, if $s=2$ and $f(n)$ is a constant, then $G_f$ has quadratic Dehn function. Although the group $G_f$ in \cite{O18} is very different from the group $G$ in this paper, the underlying $S$-machines have similar enough properties, so that we could use identical and almost identical proofs of several lemmas (which indicates that there is a general theory of $S$-machines for which this paper and \cite{O18} are applications). For the sake of completeness, we include these lemmas here with complete proofs.

\bigskip

{\bf Acknowledgment.} The authors are grateful to the anonymous referee for many helpful comments.


\section{\texorpdfstring{$S$}{S}-machines}

\subsection{\texorpdfstring{$S$}{S}-machines as rewriting systems}\label{SM}

There are several equivalent definitions of \index[g]{S@$S$-machine} $S$-machines (see \cite{S})
We are going to use the following definition which is easily seen to be equivalent to the original definition from \cite{SBR} (essentially the same definition was used in \cite{OS06}):

A \index[g]{S@$S$-machine!hardware of an $S$-machine}"hardware" of an $S$-machine $\bf S$ is a pair $(Y,Q),$ where $Q=\sqcup_{i=0}^nQ_i$ and $Y= \sqcup_{i=1}^n Y_i$ for some $n\ge 1$. Here and below $\sqcup$ denotes the disjoint union of sets.

We always set $Y_{n}=Y_0=\emptyset$ and if $Q_{n}=Q_0$ (i.e., the indices of $Q_i$ are counted $\mod n$), then we say that $\bf S$ is a \index[g]{S@$S$-machine!circular} {\em circular} $S$-machine.

The elements from $Q$ are called \index[g]{S@$S$-machine!state letters of an $S$-machine}{\it state letters}, the elements from $Y$ are \index[g]{S@$S$-machine!tape letters of an $S$-machine}{\it tape letters}. The sets $Q_i$ (resp.
$Y_i$) are called \index[g]{S@$S$-machine!parts of state and tape letters of an $S$-machine}{\em parts} of $Q$ (resp. $Y$).

The {\it language of \index[g]{S@$S$-machine!admissible words of an $S$-machine} admissible words}
consists of reduced words $W$ of the form
\begin{equation}\label{admiss}
q_1u_1q_2\dots u_s q_{s+1},
\end{equation}
where every $q_i$  is a state letter from some part $Q_{j(i)}^{\pm 1}$, $u_i$ are reduced group words in the alphabet of tape letters of the part $Y_{k(i)}$ and for every $i=1,...,s$ one of the following holds:

\begin{itemize}
\item If $q_i$ is from $Q_{j(i)}$ then $q_{i+1}$ is either from $Q_{j(i)+1}$ or is equal to $q_i\iv$; moreover $k(i)=j(i)+1$.

\item If $q_i\in Q_{j(i)}\iv$ then $q_{i+1}$ is either from $Q_{j(i)-1}^{-1}$ or is equal to $q_i\iv$; moreover $k(i)=j(i)$.
\end{itemize}
Every subword $q_iu_iq_{i+1}$ of an admissible word (\ref{admiss}) will be called the \index[g]{S@$S$-machine!admissible words of an $S$-machine!sector of an admissible word} {\em $Q_{j(i)}^{\pm 1}Q_{j(i+1)}^{\pm 1}$-sector} of that word. An admissible word may contain many $Q_{j(i)}^{\pm 1}Q_{j(i+1)}^{\pm 1}$-sectors.

For every word $W$, if we delete all non-$Y^{\pm 1}$ letters from $W$ we get the \index[g]{Y@$Y$-projection of a word} $Y$-projection of the word $W$. The length of the $Y$-projection of $W$ is called the  \index[g]{Y@$Y$-length of a word} $Y$-{\em length} and is denoted by $|W|_Y$.
Usually parts of the set $Q$ of state letters are denoted by capital letters. For example, a part $P$ would consist of
letters $p$ with various indices.

If an admissible word $W$
has the form (\ref{admiss}), $W=q_1u_1q_2u_2...q_s,$
and $q_i\in Q_{j(i)}^{\pm 1},$
$i=1,...,s$, $u_i$  are
group words in tape letters, then we shall say that the \index[g]{S@$S$-machine!admissible words of an $S$-machine!base of an admissible word}{\em base} of $W$ is the word
$Q_{j(1)}^{\pm 1}Q_{j(2)}^{\pm 1}...Q_{j(s)}^{\pm 1}$. Here $Q_i$ are just symbols which denote the corresponding parts of the set of state letters. Note that, by the definition of admissible words, the base is not necessarily a reduced word.

Instead of saying that the parts of the set of state letters of $\sss$ are $Q_0, Q_1, ... , Q_n$ we will write
that the \index[g]{S@$S$-machine!standard base of an $S$-machine}{\em the standard base}
of the $S$-machine is $Q_0...Q_n$.

The \index[g]{S@$S$-machine!software of an $S$-machine}\emph{software} of an $S$-machine with the standard base $Q_0...Q_n$ is a set of \index[g]{S@$S$-machine!rule of an $S$-machine}{\em rules} $\Theta$.  Every $\theta\in \Theta$
is a sequence $[q_0\to a_0q_0'b_0,...,q_n\to a_nq_n'b_n]$ and a subset $Y(\theta)=\sqcup Y_j(\theta)$, where $q_i\in Q_i$,
$a_i$ is a reduced word in the alphabet $Y_{i-1}(\theta)$, $b_i$ is a reduced word in $Y_i(\theta)$, $Y_i(\theta)\subseteq Y_i$, $i=0,...,n$ (recall that $Y_{0}=Y_n=\emptyset$).

Each component $q_i\to a_iq_i'b_i$ is called a \index[g]{S@$S$-machine!rule of an $S$-machine!part of a rule}{\em part} of the rule.
In most cases the sets $Y_j(\theta)$ will be equal to
either $Y_j$ or $\emptyset$. By default $Y_j(\theta)=Y_j$.

Every rule $$\theta=[q_0\to a_0q'_0b_0,...,q_n\to a_nq'_nb_n]$$ has an inverse
$$\theta\iv=[q_0'\to a_0\iv q_0b_0\iv,..., q_n'\to a_n\iv q_nb_n]$$ which is also a rule of $\sss$. It is always the case that
$Y_i(\theta\iv)=Y_i(\theta)$ for every $i$. Thus the set of rules
$\Theta$ of an $S$-machine is divided into two disjoint parts, $\Theta^+$ and
$\Theta^-$ such that for every $\theta\in \Theta^+$, $\theta\iv\in
\Theta^-$ and for every $\theta\in\Theta^-$, $\theta\iv\in\Theta^+$ (in particular $\Theta\iv=\Theta$, that is any $S$-machine is symmetric).

The rules from $\Theta^+$ (resp. $\Theta^-$) are called \index[g]{S@$S$-machine@rule of an $S$-machine!positive or negative}{\em
positive} (resp. {\em negative}).

 To \index[g]{S@$S$-machine!rule of an $S$-machine!application of a rule}apply a rule  $\theta=[q_0\to a_0q'_0b_0,...,q_n\to a_nq'_nb_n]$ as above to an admissible word $p_1u_1p_2u_2...p_s$ (\ref{admiss})
where each $p_i\in Q_{j(i)}^{\pm 1}$ means
\begin{itemize}

\item check if $u_i$ is a word in the alphabet $Y_{j(i)+1}(\theta)$ when $p_i\in Q_{j(i)}$ or if it is a word
    in $Y_{j(i)}(\theta)$ when $p_i\in Q_{j(i)}^{-1}$ ($i=1,\dots,s-1$); and if this property holds,

\item replace each $p_i=q_{j(i)}^{\pm 1}$ by $(a_{j(i)}q'_{j(i)}b_{j(i)})^{\pm 1}$,
\item if the resulting word is not reduced or starts (ends) with $Y$-letters, then reduce the word and trim the first and last $Y$-letters to obtain an admissible word again.
\end{itemize}

For example, applying the rule $[q_1\to a\iv q_1'b, q_2\to cq_2'd]$ to the admissible word $q_1b\iv q_2dq_2\iv q_1\iv$ we first obtain the word
$$a\iv q_1'bb\iv cq_2'ddd\iv (q_2')\iv c\iv b\iv(q_1')\iv a,$$ then after trimming and reducing we obtain $$q_1'cq_2'd (q_2')\iv c\iv b\iv(q_1')\iv.$$

 If a rule $\theta$ is applicable to an admissible word $W$ (i.e., $W$ belongs to the \index[g]{S@$S$-machine!rule of an $S$-machine!domain of a rule}{\em domain} of $\theta$) then we say that $W$ is a \index[g]{S@$S$-machine!$\theta$-admissible word} $\theta$-{\it admissible word} and denote the result of application of $\theta$ to $W$ by $W\cdot \theta$. Hence each rule defines an invertible partial map from the set of configurations to itself, and one can consider an $S$-machine as an inverse semigroup of partial bijections of the set of admissible words.

We call an admissible word with the standard base a  \index[g]{S@$S$-machine!configuration of an $S$-machine} \emph{configuration} of an $S$-machine.

We usually assume that every part $Q_i$ of the set of state letters contains a \index[g]{S@$S$-machine!start state letter of an $S$-machine} {\em start state letter} and an
\index[g]{S@$S$-machine!end state letter of an $S$-machine} {\em end state letter}. Then a configuration is called a \index[g]{S@$S$-machine!start configuration of an $S$-machine}\index[g]{S@$S$-machine!end configuration of an $S$-machine} \emph{start} (\emph{end}) configuration if all state letters in it are start (end) letters. As Turing machines, some $S$-machines are \index[g]{S@$S$-machine!recognizing a language} {\em recognizing a language}. In that case we choose an \index[g]{S@$S$-machine!recognizing a language!input sector of an admissible word of an $S$-machine} {\em input} sector, usually the $Q_0Q_1$-sector, of every configuration. The $Y$-projection of that sector is called the \index[g]{S@$S$-machine!recognizing a language!input of a configuration of an $S$-machine recognizing a language} {\em input} of the configuration. In that case, the end configuration with empty $Y$-projection is called \index[g]{S@$S$-machine!recognizing a language!accept configuration of an $S$-machine recognizing a language} the \emph{accept} configuration. If the $S$-machine (viewed as a semigroup of transformations as above) can take an input configuration with input $u$ to the accept configuration, we say that $u$ is \index[g]{S@$S$-machine!recognizing a language!accepted input word} {\em accepted} by the $S$-machine. We define \index[g]{S@$S$-machine!recognizing a language!accepted configuration of an $S$-machine} {\em accepted} configurations (not necessarily start configurations) similarly.

A \index[g]{S@$S$-machine!computation of an $S$-machine}{\em computation} of  \index[g]{S@$S$-machine!computation of an $S$-machine!length of a computation} {\it length} $t\ge 0$ is a sequence of admissible words $$W_0\stackrel{\theta_1}{\to} \dots\stackrel{\theta_{t-1}}{\to} W_t$$
such that for
every $i=0,..., t-1$ the $S$-machine passes from $W_i$ to $W_{i+1}$ by applying
the rule $\theta_i$ from $\Theta$.  The word $H=\theta_1\dots\theta_{t-1}$ is called the \index[g]{S@$S$-machine!computation of an $S$-machine!history of computation} {\it history}
of the computation, and the word $W_0$ is called \index[g]{S@$S$-machine!$H$-admissible word} $H$-{\it admissible}. Since $W_t$ is determined by $W_0$ and the history $H$, we use notation $W_t=W_0\cdot H$ (if rules
$\theta_j$ are irrelevant, we will omit them in the notation).

A computation is called \index[g]{S@$S$-machine!computation of an $S$-machine!reduced} {\em reduced} if its history is a reduced word.

Note, though, that in this paper, unlike the previous ones, we consider non-reduced computations too because these may correspond to reduced van Kampen diagrams (tra\-pe\-zia) under our present interpretation of $S$-machines in groups.

The {\em space} of a computation $W_0\to\dots\to W_t$ is \index[g]{S@$S$-machine!computation of an $S$-machine!space of a computation} $\max_{i=0}^t ||W_i||$, where $||W_i||$ is the  length of $W_i$.

If for some rule $\theta=[q_0\to a_0q_0'b_0,...,q_n\to a_nq'_nb_n]\in \Theta$ of an $S$-machine $\sss$ the set  $Y_{i+1}(\theta)$ is empty (hence in every admissible word in the domain of $\theta$ every $Q_iQ_{i+1}$-sector has no $Y$-letters) then we say that $\theta$ \index[g]{S@$S$-machine!rule of an $S$-machine!locking a sector} locks the $Q_iQ_{i+1}$-sector. In that case we always assume that $b_i, a_{i+1}$ are empty and we denote the $i$-th part of the rule $q_i\tool a_iq_i'$. If the $Q_iQ_{i+1}$-sector is locked by $\theta$ then we also assume that $a_{i+1}$ is empty too.

\begin{rk} \label{tool}
For the sake of brevity, the substitution $[q_i\tool aq_i', q_{i+1}\to q_{i+1}'b]$
can be written in the form $[q_iq_{i+1}\to aq'_iq'_{i+1}b]$.
Similarly we will allow parts of rules of the form $q_{i}...q_j\to aq_i'...q_j'b$. If the rule locks the $Q_sQ_{s+1}$-sector where $Q_s$ is the part of state letters containing $q_j, q_j'$, then we write $q_{i}...q_j\tool aq_i'...q_j'b$ (in that case $b$ is empty).
\end{rk}

The above definition of $S$-machines resembles the definition of multi-tape Turing machines (see \cite{SBR}).
The main differences are that every state letter of an $S$-machines is blind: it does not "see" tape letters next to it (two state letters can see each other if they stay next to each other). Also $S$-machines are symmetric (every rule has an inverse), can work with words containing negative letters, and words with "non-standard" order of state letters.

It is important that $S$-machines can simulate the work of Turing machines. This non-trivial fact, especially if one tries to get a polynomial time simulation, was first proved in \cite{SBR}. But we do not need a restriction on time, and it would be more convenient for us to use an easier $S$-machine from \cite{OS06}.

Let \index[g]{Turing machine $\mmm_0$}$\mmm_0$ be a deterministic Turing machine accepting a non-recursive language $\lll$ of words in the one-letter alphabet $\{\alpha\}$.

\begin{lemma} \label{S} (\cite{OS06}) There is a recognizing $S$-machine $\mmm_1$ whose language of accepted input words is $\mathcal L$. In every input configuration of $\mmm_1$ there is exactly one input sector, the first sector of the word, and all other sectors are empty of $Y$-letters.
\end{lemma}

We say that two recognizing $S$-machines are \index[g]{equivalent $S$-machines} \emph{equivalent} if they have the same
language of accepted configurations.

We can simplify rules of any $S$-machine in the obvious way.

\begin{lemma} \label{simp}  Every $S$-machine $\sss$ is equivalent to an $S$-machine $\sss'$, where

(*) every
part $q_i\to aq_ib$ of an $S$-rule of $\sss'$ has $||a||\le 1$, $||b||\le 1$, i.e., both words $a$ and $b$ are just letters from $Y^{\pm 1}$ or empty words;

(**) moreover $\sss'$ can be constructed so that for every rule $\theta=[q_0\to a_0q'_0b_0,...,q_n\to a_nq'_nb_n]$ of $\sss'$, we have $\sum_i (||a_i||+||b_i||)\le 1$.
\end{lemma}

For example, a rule $[q\to aq'b]$ is equivalent to the set of two rules $[q\to aq'']$, $[q''\to q'b]$ where $q''$ is a new state letter added to the part containing $q$ and $q'$.

Thus, applying Lemma \ref{S} we will assume that the
$S$-machine \index[g]{S@$S$-machine!M1@$\mmm_1$}$\mmm_1$ satisfies Property (**).

\subsection{Some elementary properties of \texorpdfstring{$S$}{S}-machines}

\label{gpsm}
The base of an admissible word is not always a reduced word. However the following is an immediate corollary of the definition of admissible word.

\begin{lemma}\label{qqiv}
If the $i$-th component of the rule $\theta$ has the form $q_i\tool a_iq_i',$ then the
base of any admissible word in the domain of $\theta$ cannot have subwords $Q_iQ_i\iv$ or $Q_{i+1}^{-1}Q_{i+1}.$
\end{lemma}

In this paper we are often using copies of words.
If $A$ is an alphabet and $W$ is a word involving no letters from $A^{\pm 1}$, then to obtain a \index[g]{copy of a word in a different alphabet}{\em copy}
of $W$ in the alphabet $A$ we substitute letters from $A$ for letters in $W$ so that different letters from $A$
substitute for different letters. Note that if $U'$ and $V'$ are copies of $U$ and $V$ respectively corresponding to
the same substitution, and $U'\equiv V'$, then $U\equiv V,$ where '$\equiv$' means
leter-by-letter equality of words. We also use copies of $S$-machines (defined in the same way).

The following two lemmas also immediately follow from definitions (see details in \cite[Lemmas 2.6,2.7]{O18}).

\begin{lemma}\label{gen} Suppose that the base of an admissible word $W$ is $Q_{i}Q_{i+1}$. Suppose that each rule of
a reduced computation starting with $W\equiv q_iuq_{i+1}$ and ending with $W'\equiv q_i'u'q_{i+1}'$ multiplies the
$Q_iQ_{i+1}$-sector by a letter on the left (resp. right).
And suppose that different rules multiply that sector by
different letters.
Then

(a) the history of computation is a copy
of the reduced form of the word $u'u\iv$ read from right to left
(resp. of the word $u\iv u'$ read from left to right). In particular,
 if $u\equiv u'$, then the computation is empty;

(b) the length of the history $H$ of the computation does not exceed $||u||+||u'||$;

(c) for every configuration $q_i''u''q_{i+1}''$ of the computation, we have
$$||u''||\le \max (||u||, ||u'||).$$
\end{lemma}

\begin{lemma}\label{gen1} Suppose the base of an admissible word $W$ is $Q_{i}Q_{i+1}$. Assume that each rule of a reduced computation starting with $W\equiv q_iuq_{i+1}$ and ending with $W'\equiv q_i'u'q_{i+1}'$ multiplies the
$Q_iQ_{i+1}$-sector by a letter on the left and by a letter from the right.
Suppose different rules multiply that sector by
different letters and the left and right
letters are taken from disjoint alphabets.
Then

(a) for every intermediate configuration $W_j$ of the computation, we have $||W_j||\le \max(||W||, ||W'||)$

(b) the length of the history $H$ of the computation does not exceed $\frac 12 (||u||+||u'||)$.
\end{lemma}

The next statement is Lemma 3.7 from \cite{O12}.

\begin{lemma} \label{gen2}
Suppose the base of an admissible word $W$ of an $S$-machine $\sss$ is $Q_{i}Q_{i}\iv$
(resp., $Q_i\iv Q_i$). Suppose that each rule $\theta$ of a reduced
computation starting with $W\equiv q_iuq_i\iv$ (resp., $q_i\iv uq_i$), where
$u\ne 1$, and ending with $W'\equiv q_i'u'(q_i')\iv$ (resp., $W'\equiv (q_i')\iv
u'q_i')$
has a part $q_i\to a_\theta
q_i'b_\theta,$ where $b_{\theta}$ (resp., $a_{\theta}$) is a letter, and  for different $\theta$-s the $b_\theta$-s
(resp., $a_\theta$-s) are different. Then the history of the computation
has the form $H_1H_2^kH_3,$ where $k\ge0$, $||H_2||\le \min(||u||, ||u'||),$ $||H_1||\le ||u||/2,$
and $||H_3||\le ||u'||/2.$
\end{lemma}

\begin{lemma} \label{wi} Suppose that a reduced computation $W_0\to W_1\to\dots \to W_t$ of an $S$-machine $\sss$ satisfying (*) in Lemma \ref{simp} has a $2$-letter base
and the history of the form  $H \equiv  H_1H_2^k H_3$ ($k\ge 0$). Then for the $Y$-projection $w_i$ of $W_i$ ($i=0,1,\dots,t$) , we have the inequality $$||w_i||\le ||w_0|+ ||w_t||+2||H_1||+ 3||H_2||+2||H_3||$$.
\end{lemma}

\proof
 By (*) we have that the absolute value of $||w_i||-||w_{i-1}||$ is at most $2$
for every $i=1,\dots t$. Therefore for  $i\le ||H_1||$, we have $||w_i||\le||w_0||+2||H_1||$.
Similarly, $||w_i||\le||w_t||+2||H_3||$ for $i\ge t-||H_3||$. It remains to assume that $||H_1||<i< t-||H_3||$.

Denote  the words $w_i$ with $i=||H_1||+j ||H_2||$, by $u_j$, $j=0,1,\dots,k$ and the corresponding words $W_i$ by $U_j$.
Then there exist two words $v_l, v_r$ such that for every $s$ from 1 to $k$, $u_s=v_lu_{s-1}v_r$  in a free group for some $Y$-words $v_l$ and $v_r$ depending on $H_2$. Hence $u_j=v_l^ju_0v_r^j$, where both $v_l$ and $v_r$ have length
at most $||H_2||$ by (*).

By \cite[Lemma 8.1]{OS04},
the length of an arbitrary word $U_j$ then is not greater
than $||v_l||+||v_r||+||U_0||+||U_k||$ provided $0\le j\le k$.

Now we need to estimate the lengths of $W_i$ ($i=|||H_1||,...,t-||H_3||$), such that $w_i$ which are not equal to any $u_j$.
Choose $j$ such that the absolute value of $i -j||H_2||$ does not exceed $||H_2||/2 $. Then the absolute value of $ ||w_i||-||u_j||$ does not exceed $||H_2||$ by (*), and therefore
$||W_i||\le ||v_l||+||v_r||+||U_0||+||U_k||+||H_2||$. Since $||U_0||\le||w_0||+2||H_1||$ and
$||U_k||\le ||w_t||+2||H_3||$, we obtain $$\begin{array}{l}||w_i||\le ||v_l||+||v_r||+||w_0||+||w_t||+2||H_1||+2||H_3||+||H_2||\\
\le ||w_0|+ ||w_t||+2||H_1||+2||H_3||+3||H_2||\end{array}$$ for every $i$, as required.
\endproof

\subsection{The highest parameter principle} \label{param}

In this paper, we estimate length and space of computations of $S$-machines, and also areas and other numerical invariants of van Kampen diagrams. The following constants will be used in the estimates throughout this paper.

\index[g]{parameters used in the paper} \index[g]{H@the highest parameter principle} \index[g]{parameters used in the paper!l@$L_0$ - a number satisfying $c_5\ll L_0\ll L$}

\begin{equation}\label{const} \begin{array}{l} \lambda^{-1}\ll m\ll  N \ll c_0\ll c_1\ll c_2\ll c_3\ll c_4\ll c_5\ll L_0\ll L\ll K\ll \\
J\ll \delta^{-1}\ll c_6\ll c_7\ll
N_1\ll N_2\ll N_3\ll N_4\end{array}
  \end{equation}
where $\ll$ means "much smaller".

For each inequality in this paper involving several of these constants, let $D$ be the biggest constant appearing there. Then the inequality  can be rewritten in the form \begin{center} $D \ge$ some expression involving smaller constants.\end{center} This {\em highest parameter principle} \cite{book} makes the system of inequalities used in this paper consistent.

\section{Auxiliary \texorpdfstring{$S$}{S}-machines and constructions}\label{am}

\subsection{Running state letters} \label{pm}

For every alphabet $Y$ we define a "running state letters" $S$-machine \index[g]{S@$S$-machine!L@$\rhh$}$\rhh(Y)$. We will omit $Y$ if it is obvious or irrelevant. The standard base of $\rhh(Y)$ is $Q^{(1)}PQ^{(2)}$ where $Q^{(1)}=\{q^{(1)}\}$, $P=\{p^{(i)}, i=1,2\}$, $Q^{(2)}=\{q^{(2)}\}$. The state letter $p$ with indices runs from the state letter $q^{(2)}$ to the state letter $q^{(1)}$ and back.  The $S$-machine $\rhh$ will be used to check the "structure" of a configuration (whether the state letters of a configuration are in the appropriate order), and to recognize a computation by its history.

The alphabet of tape letters $Y$ of $\rhh(Y)$ is $Y^{(1)}\sqcup Y^{(2)}$, where $Y^{(2)}$
is a (disjoint) copy of $Y^{(1)}$.
The positive rules of $\rhh$ are defined
as follows.
\begin{itemize}

\item $\zeta^{(1)}(a)=[q^{(1)}\to q^{(1)}, p^{(1)}\to a\iv p^{(1)}a', q^{(2)}\to q^{(2)}]$,
where $a$ is any positive letter from $Y=Y^{(1)}$ and $a'$ is the corresponding letter in the copy $Y^{(2)}$ of $Y^{(1)}$.
\me

{\em Comment.} The state letter $p^{(1)}$ moves left replacing letters $a$ from $Y^{(1)}$ by their copies $a'$
from $Y^{(2)}$.

\me

\item $\zeta^{(12)}=[q^{(1)}p^{(1)}\to  q^{(1)}p^{(2)}, q^{(2)}\to q^{(2)}]$.

\me

{\em Comment.} When $p^{(1)}$ meets $q^{(1)}$, $p^{(1)}$  turns into $p^{(2)}$.

\me

\item $\zeta^{(2)}(a) =[q^{(1)}\to q^{(1)}, p^{(2)}\to ap^{(2)}(a')^{-1}, q^{(2)}\to q^{(2)}]$

{\em Comment.} The state letter $p^{(2)}$ moves right towards $q^{(2)}$ replacing letters $a'$ from $Y^{(2)}$ by their copies $a$
from $Y^{(1)}$.
\end{itemize}

The start (resp. end) state letters of $\rhh$ are $\{q^{(1)},p^{(1)}, q^{(2)}\}$ (resp. $\{q^{(1)},p^{(2)}, q^{(2)}\}$).

\begin{rk}\label{r:rh}Note that each of the rules $(\zeta^j)^{\pm 1}(a)$, ($j=1,2$) either moves the state letter $p$ left or moves
it right, or deletes one letter from left and one letter from right, or insert letters from both sides of itself. In the latter
case, the next rule of a computation must be again $\zeta(j)^{\pm 1}(b)$ for some $b$, and if the computation is
reduced, it again must increase the length of the configuration by two. This observation implies
\end{rk}

\begin{rk} \label{proj} Note that no rule of $\rhh$ changes the projection of a configuration onto the free group with basis $Y^{(1)}$ if the state letters are mapped to $1$ and the letters from $Y^{(2)}$ are mapped to their copies from $Y^{(1)}$.
 This will be later referred to as the
\index[g]{projection argument} {\it projection argument}.
\end{rk}


\begin{lemma}\label{prim}  Let $\ccc\colon W_0\to\dots \to W_t$ be a reduced computation of the $S$-machine $\rhh$ with the standard base. Then

(1) if $|W_i|_Y>|W_{i-1}|_Y$ for some $i=1,\dots,t-1$, then $|W_{i+1}|_Y>|W_i|_Y$;

(2) $|W_i|_Y\le\max(|W_0|_Y, |W_t|_Y)$ for every $i=0,1,\dots t$;

(3) if $W_0\equiv q^{(1)}up^{(1)}q^{(2)}$ and $W_t\equiv q^{(1)}vp^{(2)}q^{(2)}$ for some words $u,v$, then $u\equiv v$,
$|W_i|_Y=|W_0|_Y$ for every $i=0,\dots,t$, $t=2k+1$, where $k=|W_0|_Y$,
and the sector $Q^{(1)}P$ is locked in the transition $W_k\to W_{k+1}$.
Moreover if $W_0$ and $W_t$ have the form  $q^{(1)}up^{(1)}q^{(2)}$ and $q^{(1)}vp^{(2)}q^{(2)}$, then the history $H$ of $\ccc$ is a copy of the word $\bar u\zeta^{(12)} (\bar u')\iv$, where $\bar u$ is the mirror image of $u$ and $\bar u'$ is a copy of $\bar u$. Thus, $W_0, W_t, H$ uniquely determine each other in that case.

(4) if  $W_0\equiv q^{(1)}up^{(1)}q^{(2)}$ and $W_t\equiv q^{(1)}vp^{(1)}q^{(2)}$ for some $u,v$ or
$W_0\equiv q^{(1)}up^{(2)}q^{(2)}$ and $W_t\equiv q^{(1)}vp^{(2)}q^{(2)}$ then $u\equiv v$ and the computation is empty ($t=0$);

(5) if $W_0\equiv q^{(1)}up^{(1)}q^{(2)}$ or $W_0\equiv q^{(1)}p^{(1)}uq^{(2)}$, or $W_0\equiv q^{(1)}up^{(2)}q^{(2)}$, or $W_0\equiv q^{(1)}p^{(2)}uq^{(2)}$ for some word $u$, then $|W_i|_Y\ge |W_0|_Y$ for every $i=0,\dots,t$.
\end{lemma}

\proof For every $i=0,...,t$ let $W_i=q^{(1)}u_ip^{(l_i)}v_iq^{(2)}$ where $u_i$ is a word in $Y$, $v_i$ is a word in $Y'$ (it is easy to check by induction on $i$ that this is true for every $i$).

Suppose that $|W_{i-1}|_Y<|W_{i}|_Y$ for some $i$. That means that the $i$-th rule in the computation is of the form $(\zeta^{(k)}(a))^{\pm 1}$. This rule multiplies $u_{i-1}$ by a letter $a^{\pm 1}$ on the right, and multiplies $v_{i-1}$  by a copy of the inverse of that letter on the left, and these letters do not cancel in $u_{i}$, $v_{i}$.
In particular both $u_{i}$ and $v_{i}$ are not empty. Hence $\zeta^{(12)}$ does not apply to $W_{i}$. Thus the rule in $W_{i}\to W_{i+1}$  is $(\zeta^{(j)}(b))^{\pm 1}$ (with the same $j$) and it multiples $u_{i}=u_{i-1}a$ by $b^{\pm 1}$ on the right and multiples $v_{i}$ by a copy of the inverse of that letter on the left. Since the computation is reduced, $b\ne a\iv$. Therefore $|W_{i+1}|_Y>|W_{i}|_Y$. Continuing in this manner, we establish (1).

To establish (2), we can choose the shortest word $W_j$ in the computation and apply (1) to the computation $W_j\to\dots\to W_t$ and the inverse computation $W_j\to\dots\to W_0$.

Suppose that the assumptions of (3) hold. Then $u\equiv v $ by the projection argument. Since $\zeta^{(12)}$ locks $Q^1P$-sector, the $p$-letter must reach $q^{(1)}$ moving always left to change $p^{(1)}$ to $p^{(2)}$,
 and so $W_k\equiv q^{(1)}p^{(1)}\dots$. If the next rule of the form $\zeta^{(1)}(a)^{\pm 1}$ could increase the length of the configuration, we would obtain a contradiction with Property (1). Since the computation is reduced, the next rule is $\zeta^{(12)}$, and arguing in this way, one uniquely reconstructs the whole computation
in case (3) for given $W_0$ or $W_t$, and vice versa, the history $H$ determines both $u$ and $v$. Propery (4) holds for same reasons.

By the projection argument, we have $|q^{(1)}up^{(1)}q^{(2)}|_Y =||u||\le |W_i|_Y$ if the first assumptions of (5) holds.
The other cases of (5) are similar.
\endproof

The projection argument also immediately gives:

\begin{lemma}\label{ewe} If $W_0\to\dots\to W_t$ is a reduced
computation of $\rhh$ with base $$Q^{(1)}PP^{-1}(Q^{(1)})^{-1} \hbox{ or } (Q^{(2)})^{-1}P^{-1}PQ^{(2)}$$
and $$W_0\equiv q^{(1)}p^{(i)}u(p^{(i)})^{-1}(q^{(1)})^{-1} (i=1,2)$$ or $$W_0\equiv(q^{(2)})^{-1}(p^{(i)})^{-1}v(p^{(i)})q^{(2)} (i=1,2)$$  for some words $u,v$, then $|W_j|_Y\ge |W_0|_Y$ for every $j=0,\dots,t$.
\end{lemma}

\begin{rk} \label{right}
We will also use the right analog \index[g]{S@$S$-machine!R@$\rhr$}$\rhr$ of $\rhh$. The base of $\rhr$ is $Q_1RQ_2$. The state letter $r$ first moves right from $q^{(1)}$ to $q^{(2)}$ and then left. Lemmas "left-right dual" to Lemmas \ref{prim} and  \ref{ewe}
as well as Remark \ref{proj} are true for $\rhr$ as well. \end{rk}

\begin{rk} \label{lrm} For every $m\ge 1$, we will also need the $S$-machine \index[g]{S@$S$-machine!L@$\rhh_m$}$\rhh_m$, that repeats the
work of $\rhh$ $m$ times. That is the $S$-machine $\rhh_m$ runs the state letter $p$ back and forth between $q^{(2)}$ and $q^{(1)}$ $m$ times. Every time $p$ meets $q^{(1)}$ or $q^{(2)}$, the upper index of $p$ increases by $1$ after the application of the rule $\zeta^{(i,i+1)}$ ($i=1,\dots, 2m-1$), so the highest upper index of $p$ is $(2m)$. A precise definition of $\rhh_m$ is obvious and is left to the reader. (Recall that $m$ is one of the system of parameters used in this paper (see Section \ref{param}).)
\end{rk}

\begin{rk}\label{r-prim} The analog of Lemma
\ref{prim} holds for $\rhh_m$. In particular, if $$W_t\equiv q^{(1)}vp^{(2m)}q^{(2)}$$ in the formulaion of (3), then $t=2mk+2m-1$ (the proof is essentially the same and is left to the reader).
\end{rk}

\subsection{Adding history sectors} \label{hs}

We will add new (history) sectors to our $S$-machine $\mmm_1$. If we ignore the new sectors, we get the hardware and the  software of the $S$-machine $\mmm_1$. The new $S$-machine $\mmm_2$ will start with a configuration where in every history sector a copy of the history $H$ of a computation of $\mmm_1$ is written. Then it will execute $H$ on the other (working) sectors simulating the work of $\mmm_1$, while  in the history sector, state letters scan the history, one symbol at a time. Thus if a computation with the standard base starts with a configuration $W$ and ends with configuration $W'$, then the length of the computation does not exceed $||W||+||W'||$.

Here is a precise definition of $\mmm_2$. Recall that the $S$-machine $\mmm_1$ satisfies the
condition (**) of Lemma \ref{simp} and has hardware $(Q,Y)$, where
$Q=\sqcup_{i=0}^n Q_i$, and the set  of rules $\Theta$. The new
$S$-machine \index[g]{S@$S$-machine!M2a@$\mmm_2$} $\mmm_2$ has hardware
$$Q_{0,r}\sqcup Q_{1,\ell}\sqcup Q_{1,r}\sqcup Q_{2,\ell}\sqcup Q_{2,r}\sqcup\dots \sqcup Q_{n,\ell},\;\; Y_h=  Y_1
\sqcup X_1\sqcup Y_2\sqcup\dots \sqcup X_{n-1}\sqcup Y_n$$
where $Q_{i,\ell}$ and $Q_{i,r}$ are (left and right) copies of $Q_i$  and $X_i$ is a disjoint union of two copies of $\Theta^+$, namely $X_{i,\ell}$ and $X_{i,r}$. (The sets $Q_{0,\ell}$, $Q_{n,r}$ 
are empty.) Every letter $q$ from $Q_i$ has two copies $q^{(\ell)}\in Q_{i,\ell}$ and $q^{(r)}\in Q_{i,r}$. By definition, the start (resp. end) state letters  of $\mmm_2$ are copies of the corresponding start (end) state letters of $\mmm_1$.  The $Q_{0,r}Q_{1,\ell}$-sectors are the \index[g]{S@$S$-machine!M2a@$\mmm_2$!input sector of a configuration of $\mmm_2$} {\em input sectors} of configurations of $\mmm_2$.

The positive rules $\theta_h$ of $\mmm_2$ are in one-to-one correspondence with the positive rules $\theta$ of $\mmm_1$. If $\theta=[q_0\to a_0q_0'b_0,...,q_n\to a_nq_n'b_n]$ is a positive rule of $\mmm_1$, then each part $q_i\to a_iq_i'b_i$ is replaced in $\theta_h$ by two parts $$q_{i,\ell}\to a_iq_{i,\ell}'h_{\theta,i}\iv $$ and
$$q_{i,r}\to \overline h_{\theta,i}q_{i,r}' b_i,$$
where $h_{\theta,i}$ (resp., $\overline h_{\theta,i}$) is a copy of $\theta$ in the alphabet $X_{i,\ell}$ (in $X_{i,r}$, respectively).

If $\theta$ is the start (resp. end) rule of $\mmm_1$, then for any word in the domain of $\theta_h$ (resp. $\theta_h\iv$) all $Y$-letters in history sectors are from $\sqcup_i X_{i,\ell}$ (resp. $\sqcup X_{i,r}$).

Thus for every rule $\theta$ of $\mmm_1$, the rule $\theta_h$ of $\mmm_2$ acts in the $Q_{i,r}Q_{i+1,\ell}$-sector in the same way as $\theta$ acts in the $Q_iQ_{i+1}$-sector. In particular, $Y$-letters which can appear in the $Q_{i,r}Q_{i+1,\ell}$-sector of an admissible word in the domain of $\theta_h$ are the same as the $Y$-letters  that can appear in the $Q_iQ_{i+1}$-sector of an admissible word in the domain of $\theta$. Hence if $\theta$ locks $Q_iQ_{i+1}$-sectors, then $\theta_h$ locks $Q_{i,r}Q_{i+1,\ell}$-sectors.

\begin{rk} Note that $\mmm_2$ no longer satisfies Property (**) from Lemma \ref{simp} but it satisfies Property (*) of that Lemma. Property (*) holds for subsequent machines ${\bf M}_3 - {\bf M}_6={\bf M}$ as well.
\label{l:s}
\end{rk}

\begin{rk}
Every computation of the $S$-machine $\mmm_2$ with history $H$ and the standard base coincides with the a computation of $\mmm_1$ whose history is a copy of $H$
if
one observes it only in {\it working} sectors $Q_{i,r}Q_{i+1,l}$.
In the standard base of $\mmm_2$  the \index[g]{S@$S$-machine!M2a@$\mmm_2$!working sector of a configuration of $\mmm_2$} {\em working
sectors} $Q_{i,r}Q_{i+1,\ell}$ alternate with \index[g]{S@$S$-machine!M2a@$\mmm_2$!history sector of a configuration of $\mmm_2$}{\it history} sectors $Q_{i,\ell}Q_{i,r}$. Every positive rule
$\theta_h$ multiplies the content of the history $Q_{i,\ell}Q_{i,r}$-sector  by the corresponding letter $\overline h_{\theta, i}$ from the right
and by letter $h_{\theta,i}^{-1}$ from the left. Thus if the $S$-machine $\mmm_2$ executes the history written in the history sectors, then the history word $H$ in letters from $X_{i,\ell}$ gets rewritten into the  copy of $H$ in letters from $X_{i,r}$. Say, if the copy of the history $H$ was written in a history sector as $h_1h_2h_3$, then during the computation with history $H$ it will transform as follows:

$$h_1h_2h_3\to h_2h_3\overline h_1\to h_3\overline h_1\overline h_2\to \overline h_1\overline h_2\overline h_3.$$

Let \index[g]{S@$S$-machine!M1@$\mmm_1$!I@$I_1(\alpha^k)$ - a start configuration of $\mmm_1$}$I_1(\alpha^k)$ be a start configuration of $\mmm_1$ (i.e.,a configuration in the domain of the start rule of $\mmm_1$) with $\alpha^k$ written in the input sector (all other sectors do not contain $Y$-letters). Then the corresponding start configuration \index[g]{S@$S$-machine!M2a@$\mmm_2$!I@$I_2(\alpha^k,H)$ - a start configuration of $\mmm_2$} $I_2(\alpha^k,H)$ of $\mmm_2$ is obtained by first replacing each state letter $q$ by the product of two corresponding letters $q^{(\ell)}q^{(r)}$, and then inserting a copy of $H$ in the \index[g]{X@$X_{i,\ell}$, a left alphabet} {\em left alphabet} $X_{i,\ell}$ in every history $Q_{i,\ell}Q_{i,r}$-sector. End configurations \index[g]{S@$S$-machine!M2a@$\mmm_2$!A@$A_2(H)$ - an end configuration of $\mmm_2$} $A_2(H)$ of $\mmm_2$ are defined similarly, only the $Y$-letters in the history sectors must be from the \index[g]{X@$X_{i,r}$, a right alphabet} {\em right alphabet} $X_{i,r}$.
\label{r:09}
\end{rk}

\begin{lemma}\label{I2A2} (1) If a word $\alpha^k$ is accepted by the Turing machine $\mmm_0$, then for some word $H$, there is a reduced computation $I_2(\alpha^k,H)\to\dots\to A_2(H)$ of the $S$-machine $\mmm_2$.

(2) If there is a computation $I_2(\alpha^k,H)\to\dots\to A_2(H')$ of  $\mmm_2$, then the word $\alpha^k$ is accepted
by $\mmm_0$ and $H'\equiv H$.
\end{lemma}

\proof (1) The word $\alpha^k$ is accepted by the $S$-machine $\mmm_1$ by Lemma \ref{S}. If $H$ is the history
of the accepting computation of $\mmm_1$, then the computation of $\mmm_2$ with history $H$
starting with $I_2(\alpha^k,H)$ ends with $A_2(H)$ since $\mmm_2$ works as $\mmm_1$ in the working
sectors  and replaces the letters from the left alphabets by the corresponding letters from
the right alphabets in the history sectors.

(2) If $I_2(\alpha^k,H)\cdot H''=A_2(H')$ for some history $H''$ of $\mmm_2$ then the word $\alpha^k$ is accepted by $\mmm_0$ by Lemma \ref{S} and the fact that $\mmm_2$ works as $\mmm_1$ in the working sectors.
Note that both $H$ and $H'$ must be the copies of $H''$, 
because the word $I_2(\alpha^k,H)$ has no letters from
right alphabets, $A_2(H')$ has no letters from left alphabets, and every rule
multiplies the $Y$-projection of every history sector by a letter from $X_{i,\ell}\iv$ (from $X_{i,r}$) on the left
(resp., on the right).
\endproof

The sectors of the form $Q_{i,\ell}Q_{i,\ell}^{-1}$ and $Q_{i,r}^{-1}Q_{i,r}$ (in a non-standard base)
are also called {\it history} sectors. History sectors help obtaining a linear estimate of the
space of every computation $W_0\to\dots\to W_t$ in terms of $||W_0||+||W_t||$.

\begin{lemma}\label{w}  Let $W_0\to\dots\to W_t$ be a reduced computation of $\mmm_2$ with base $Q_{i,\ell}Q_{i,r}$ and history $H$. Assume that all the $Y$-letters of $W_0$ belong to only one
of the alphabets $X_{i,\ell}$ or  $X_{i,r}$.
Then $||H||\le |W_t|_Y$ and $|W_0|_Y\le |W_t|_Y$
\end{lemma}

\proof
Let $W_i = q_iv_iq_i'$, $i=0,...,t$, and assume that $v_0$ has no letters from $X_{i,r}$. Then $v_t=uv_0u'$, where $u$ is a copy of $H^{-1}$ in the alphabet $X_{i,\ell}$ and $u'$ is a copy of $H$
in $X_{i,r}$. So no letter of $u'$ is cancelled in the
product $uv_0u'$, Therefore $|W_t|_Y\ge ||u'||=||H||$ and
$|W_t|_Y\ge |W_0|_Y$.
\endproof

\begin{lemma} \label{9}
For any reduced computation $W_0\to\dots\to W_t$ of $S$-machine $\mmm_2$
with base of length at least $3$, we have $|W_i|_Y\le 9( |W_0|_Y+|W_t|_Y) $ ($0\le i\le t$).
\end{lemma}

\proof Let $Q_{i_1}^{\pm 1} \dots Q_{i_k}^{\pm 1} $ be the base of the computation. We can divide the base into several subwords of length $3$ or $4$, each containing one history sector. Thus we can assume that $k$ is equal to $3$ or $4$ and that the base contains one history sector. Without loss of generality, that history sector is either a $Q_{i,\ell}Q_{i,r}$-sector or a $Q_{i,\ell}Q_{i,\ell}\iv$-sector or a $Q_{i,r}\iv Q_{i,r}$-sector.

Consider two cases.

{\bf 1.} The history sector has the form $Q_{i,\ell}Q_{i,r}$. By Lemma \ref{gen1}, we have $||H||\le \frac 12 (|W_0|_Y+|W_t|_Y)$. It follows from property (*) of Lemma \ref{simp} that $|\;|W_{i+1}|_Y-|W_i|_Y\;|\le 6$ for every
$i$. Therefore $$|W_i|_Y\le \max(|W_0|_Y, |W_t|_Y)+3||H|| \le $$ $$\max(|W_0|_Y, |W_t|_Y)+
\frac 32 (|W_0|_Y +|W_t|_Y)\le \frac 52 (|W_0|_Y +|W_t|_Y)$$

{\bf 2.} The history sector is either a $Q_{i,\ell}Q_{i,\ell}^{-1}$-sector or a $Q_{i,r}^{-1}Q_{i,r}$-sector.
Then one can apply Lemma \ref{gen2} to the history sector and obtain the factorization $H\equiv H_1H_2^cH_3,$ with $c\ge0$, $||H_2||\le \min(||u_0||, ||u_t||),$ $||H_1||\le ||u_0||/2,$
and $||H_3||\le ||u_t||/2,$ where $u_0$ and $u_t$ are the $Y$-projections of the history sectors of $W_0$ and $W_t$, respectively. Since every $W_i$ has at most three sectors,
applying Lemma \ref{wi} to each of them, we obtain:
$$|W_i|_Y\le |W_0|_Y+|W_t|_Y+3(2||H_1||+ 3||H_2||+2||H_3||) \le $$ $$|W_0|_Y+ |W_t|_Y +3|W_0|_Y+ 9\min(|W_0|_Y, |W_t|_Y)+3|W_t|_Y\le 9(|W_0|_Y +|W_t|_Y).$$
\endproof

\begin{lemma} \label{WV} Suppose that a reduced computation $W_0\to\dots\to W_t$  of the $S$-machine $\mmm_2$ starts with an admissible word $W_0$ having no letters
from the alphabets $X_{i,l}$ (resp., from the alphabets $X_{i,r}$) . Assume that the length of its base $B$ is bounded from above by a constant  $N_0$, and $B$ has
a history subword $Q_{i,\ell}Q_{i,r}$. Then there is a constant $c=c(N_0)$ such that
$|W_0|_Y\le c|W_t|_Y$.
\end{lemma}

\proof
 Let $V_0\to\dots\to V_t$
be the restriction of the computation  to
the $Q_{i,\ell}Q_{i,r}$-sector. By Lemma \ref{w}, we have
$t\le |V_t|_Y$ and $|V_0|_Y\le |V_t|_Y$.

It follows from (*) that

$$|W_0|_Y
\le |W_t|_Y+ 2N_0t\le |W_t|_Y+ 2N_0|V_t|_Y\le (2N_0+1)|W_t|_Y$$

It suffices to choose $c=2N_0+1$.
\endproof

\subsection{Adding running state letters}\label{csl}

Our next $S$-machine will be a composition of $\mmm_2$ with $\rhh$ and $\rhr$. The running state letters will control the work of $\mmm_3$.

First we replace every part $Q_i$ of the state letters in the standard base of $\mmm_2$ by three parts $P_iQ_iR_i$ where $P_i, R_i$ contain the running state letters. Thus if $Q_0...Q_s$ is the standard base of $\mmm_2$ then the standard base of  \index[g]{S@$S$-machine!M2b@$\bmm_2$} $\bmm_2$ is
\begin{equation}\label{basec}
P_0Q_0R_0P_1Q_1R_1\dots P_sQ_sR_s,
\end{equation}
where
$P_i$ (resp., $R_i$) contains
copies of running $P$-letters (resp. $R$-letters) of $\rhh$ (resp. $\rhr$), $i=0,\dots, s$.

For every rule $\theta$ of  $\mmm_2$,
its $i$-th part $[q_i\to a_iq_i'b_i]$ is replaced in $\bmm_2$ with
\begin{equation}\label{3p}
[p^{(i)}q_ir^{(i)}\to a_ip^{(i)}q_i'r^{(i)}b_i],  (i=0,\dots, s),
\end{equation}
 where $ p^{(i)}\in P_i, r^{(i)}\in R_i$ do not depend on $\theta$.

 \medskip

{\it Comment.} Thus, the sectors $P_iQ_i$ and $Q_iR_i$ are always locked. Of course, such a modification is useless for solo work of $\mmm_2$. But it will be helpful when one constructs
a composition of $\bmm_2$ with $\rhh$ and $\rhr$ which will be turned on after certain rules of $\bmm_2$ are applied.

If $Q_iQ_{i+1}$-sector is a history sector of $\mmm_2$, then  $Q_iR_i$-, $R_iP_i$-, $P_iQ_i$-sectors are \index[g]{S@$S$-machine!M2b@$\bmm_2$!history sectors of $\bmm_2$} history sectors of $\bmm_2$. Accordingly the $Q_iQ_i\iv$-sectors ($R_iR_i\iv$-sectors, etc.) of admissible words with nonstandard bases will be called history sectors of $\bmm_2$ too. (Alternatively, history sectors of admissible words of $\bmm_2$ are those sectors which can contain letters from left or right alphabets.) The $R_0P_1$-sectors of admissible words are the \index[g]{S@$S$-machine!M2b@$\bmm_2$!input sector of $\bmm_2$} input sectors. The $R_0R_0\iv-$ and $P_1\iv P_1$-sectors are also input sectors of admissible words of $\bmm_2$.

If $B$ is the base of some computation $\ccc$ of $\bmm_2$, and $UV$ is a 2-letter subword of $B$ such that $UV$-sectors of admissible words in $\ccc$ are history (resp. working, input) sectors, then we will call $UV$ a \index[g]{history, working and input subwords of the base of a computation of $\bmm_2$ and $\mmm_3$} history (resp. working, input) subword of $B$.

\subsection{\texorpdfstring{$\mmm_3$}{M3}}\label{M3M5}

The next $S$-machine \index[g]{S@$S$-machine!M3@$\mmm_3$} $\mmm_3$ is the composition of the $S$-machine $\mathbf{\overline M}_2$ with $\rhh$ and $\rhr$.
The $S$-machine $\mmm_3$ has the  input, working and history sectors, i.e. the same base as $\bmm_2$, although the parts of this base have more
state letters than the corresponding parts of $\mathbf{\overline M}_2$. It works as follows. Suppose that $\mmm_3$ starts with a start configuration of $\bmm_2$, a word $\alpha^k$ in the input $R_0P_1$-sector, copies of a history word $H$ in the alphabets $X_{i,\ell}$ in the history sectors, all other sectors empty of $Y$-letters. Then $\mmm_3$ first executes $\rhr$ in all history sectors (moves the running state letter from $R_i$ in the history sectors right and left), then it executes the history $H$ of $\bmm_2$. After that the $Y$-letters in the history sectors are in $X_{i,r}$ and $\mmm_3$ executes copies of $\rhh$ in the history sectors (moves the running state letters left then right). After that $\mmm_3$ executes a copy of $H$ backwards, getting to a copy of the same start configuration of $\bmm_2$, runs $\rhr$, executes a copy of the history $H$ of $\bmm_2$, runs a copy of $\rhh$, etc. It stops after $m$ times running $\rhr, \bmm_2,\rhh, \bmm_2\iv$ and running $\rhr$ one more time.

Thus the $S$-machine $\mmm_3$ is a concatenation of $4m+1$ $S$-machines $\mmm_{3,1}-\mmm_{3,4m+1}$. After one of these $S$-machines terminates, a transition rule changes its end state letters to the start state letters of the next $S$-machine. All these $S$-machines have the same standard bases as $\bmm_2$.

The configuration $I_3(\alpha^k,H)$ of ${\bf M}_3$ is obtained from $I_2(\alpha^k,H)$
by adding the control state letters $r^{(1)}_i$ and $p^{(1)}_i$ according to (\ref{basec}) in Section \ref{csl}.

\medskip

{\bf Set $\mmm_{3,1}$} is a copy of the set of rules of the $S$-machine $\rhr$, with \index[g]{S@$S$-machine!L@$\rhh$!parallel work of $\rhh$ or $\rhr$ in several sectors} \emph{parallel work} in all history sectors, i.e., every subword $Q_{i-1}R_{i-1}P_i$ of the standard base, where $Q_{i-1}Q_i$ is a history sector of $\mmm_2$,
is treated as the base of a copy of $\rhr$, that is $R_{i-1}$ contain the running state letters which run between state letters from $Q_{i-1}$ and $P_i$. Each rule of Set $\mmm_{3,1}$ executes the corresponding rule of $\rhr$ simultaneously in each history sector of $\mmm_2$. The partition of the set of state letters of these copies of $\rhr$ in each history sector is $X_{i,\ell}\sqcup X_{i,r}$ for some $i$ (that is state letters from $R_{i-1}$ first run right replacing letters from $X_{i,\ell}$ by the corresponding letters of $X_{i,r}$ and then run left replacing letters from $X_{i,r}$ by the corresponding letters of $X_{i,\ell}$.

The  transition rule $\chi(1,2)$ changes the state letters to the state letters of start configurations of $\bmm_2$. The admissible words in the domain of $\chi(1,2)^{\pm 1}$
have all $Y$-letters from the left alphabets $X_{i,\ell}$. The rule $\chi(1,2)$ locks all sectors except the history sectors $R_{i-1}P_i$ and the input sector. It does not apply to admissible words containing $Y$-letters from right alphabets.

\medskip

{\bf Set  $\mmm_{3,2}$} is a copy of the set of rules of the $S$-machine $\bmm_2$.

The transition rule $\chi(2,3)$ changes the state letters of the stop configuration of $\bmm_2$ to  their copies in a different alphabet. The admissible words in the domain of $\chi(2,3)^{\pm 1}$
have no $Y$-letters from the left alphabets $X_{i,\ell}$. The rule $\chi(2,3)$ locks all sectors except for the history sectors $R_{i-1}P_i$. It does not apply to admissible words containing $Y$-letters from right alphabets.

\medskip

{\bf  Set $\mmm_{3,3}$} is a copy of the set of rules of the $S$-machine $\rhh$, with parallel work in the same sectors as  $\mmm_{3,1}$ (and the same partition of $Y$-letters in each history sector $X_{i,r}\sqcup X_{i,\ell}$).

The transition rule $\chi(3,4)$ changes the state letters of the stop configuration of $\bmm_2$ to  their copies in a different alphabet. The admissible words in the domain of $\chi(3,4)^{\pm 1}$
have no $Y$-letters from the left alphabets $X_{i,l}$. The rule $\chi(3,4)$ locks all non-history sectors.

\medskip

{\bf Set $\mmm_{3,4}$.} The positive rules of Set $\mmm_{3,4}$ are the copies of the negative rules of the $S$-machine  $\bmm_2$.

The transition rule $\chi(4,5)$ changes the state letters of the start configuration of $\bmm_2$ to  their copies in a different alphabet. The admissible words in the domain of $\chi(4,5)^{\pm 1}$
have no $Y$-letters from the right alphabets $X_{i,r}$. The rule $\chi(4,5)$ locks all  non-history and non-input sectors.

\medskip

{\bf Sets $\mmm_{3,5}, \dots, \mmm_{3,8}$} consist of rules that are  copies of the rules of the Sets $\mmm_{3,1},\dots,$ $\mmm_{3,4}$, respectively.


$\dots$

{\bf Sets $\mmm_{3,4m-3},\dots, \mmm_{3,4m}$} consist of copies of the steps $\mmm_{3,1}, \dots, \mmm_{3,4}$, respectively.

\medskip

{\bf Set $\mmm_{3,4m+1}$} is a  copy of Set $\mmm_{3,1}$. The end configuration for Set $\mmm_{3,4m+1}$,
$A_3(H)$, is obtained from a copy of $A_2(H)$ by inserting the control letters according to (\ref{basec}).

The transition rules $\chi(i,i+1)$ are called \index[g]{S@$S$-machine!M3@$\mmm_3$!Chi@$\chi$-rules of $\mmm_3$} $\chi$-rules.

We say that a configuration $W$ of the $S$-machine $\mmm_3$ is \index[g]{S@$S$-machine!M3@$\mmm_3$!tame configuration of $\mmm_3$} {\it tame} if every $P$- or $R$-letter is next to some $Q$-letter in $W$.


\begin{lemma}\label{Hprim}
Let $\ccc\colon
W_0\to\dots\to W_t$ be
a reduced computation of  $\mmm_3$ consisting of rules of one of the copies of $\rhh$ or $\rhr$ with standard base. Then

(a) $|W_j|_Y\le\max (|W_0|_Y,|W_t|_Y)$ for every configuration $W_j$ of $\ccc$; moreover, $|W_0|_Y\le \dots\le |W_t|_Y$ if $W_0$ is tame;

(b) $ t\le ||W_0||+||W_t||-2$, moreover, $t\le 2||W_t||-2$ if $W_0$ is tame.

\end{lemma}

\proof (a) Let $W_r$ be a shortest word of the computation $\ccc$.
Then either $|W_r|_Y=|W_{r+1}|_Y=\dots= |W_t|_Y$, or
$|W_r|_Y=|W_{r+1}|_Y=\dots= |W_s|_Y<|W_{s+1}|_Y$ for
some $s$. It follows that the number of sectors increasing
their lengths by two at the transition $W_s\to W_{s+1}$ is greater than the number of the sectors decreasing the lengths by
$2$. Now it follows from Lemma \ref{prim} (1) that the lengths of the $Y$-projections will keep increasing: $|W_{s+1}|_Y<|W_{s+2}|_Y<\dots $.
So for every $j\ge r$, we have $|W_j|_Y\le |W_t|_Y$.
Similarly, we have $|W_r|_Y\le |W_0|_Y$ for $j\le r$.
If the word $W_0$ is tame, then it is the shortest configuration by the projection argument.

(b) If the rules do not change the lengths of configurations, then every control letter runs right and left
only one time by Lemma \ref{prim} (4), and the inequality follows.  If $||W_r||<||W_{r+1}||$
for some $r$, then every next transition keeps increasing the
length by Lemma \ref{prim} (1), and so the inequality holds as well.

\endproof

\begin{lemma}\label{M31} Let $\ccc\colon W_0\to\dots\to W_t$ be a reduced computation
of $\mmm_3$. Then for every $i$, there is at most one occurrence of the rules
$\chi(i,i+1)^{\pm 1}$ in the history $H$ of $\ccc$ provided the
base of $\ccc$ has a history $(R_{j-1}P_j)^{\pm 1}$-sector.
\end{lemma}

\proof
Arguing by contradiction, we can assume that
$H=\chi(i,i+1)^{\pm 1}H'\chi(i,i+1)^{\mp 1}$, where
$H'$ is a copy of the history of a computation of either $\rhh$ or $\rhr$
or $\bmm_2$. The two cases $\rhh$ and  $\rhr$ contradict Lemma \ref{prim} (4).
The latter case (namely $\bmm_2$) is also  impossible. Indeed, consider any history subword $(R_{j-1}P_j)^{\pm 1}$ of the base of the computation. Then  the $Y$-projection of the $(R_{j-1}P_j)^{\pm 1}$ -sector of $W_1$ must be a word either in the $X_{j,\ell}$ or in $X_{j,r}$ (depending on the parity of $i$). Without loss of generality assume that it is $X_{j,\ell}$. Then the computation $W_1\to,\dots,\to W_{t-1}$ multiplies the $Y$-projection of the $(R_{j-1}P_j)^{\pm 1}$ -sector of $W_1$ by a word in $X_{j,\ell}$ and a reduced word in $X_{j,r}$. Hence the $(R_{j-1}P_j)^{\pm 1}$ -sector of $W_{t-1}$ contains letters from a right alphabet, hence $W_{t-1}$ cannot be in the domain of $\chi(i,i+1)^{\pm 1}$, a contradiction.
\endproof

\begin{lemma}\label{M3} Let $\ccc\colon W_0\to\dots\to W_t$ be a reduced computation
of $\mmm_3$. Suppose also that the base of $\ccc$ is standard, then

(a) if the history of  $\ccc$
has the form $\chi(i,i+1)H'\chi(i+4,i+5)$, then the word
$W_0$ is a copy of $W_t$;

(b) two subcomputations $\ccc_1$ and $\ccc_2$ of $\ccc$ with histories $\chi(i,i+1)H'\chi(i+4,i+5)$ and $\chi(j,j+1)H''\chi(j+4,j+5)$
have equal lengths; moreover some cyclic permutation of $\ccc_2$
is a copy of $\ccc_1$;

(c) there is a constant \index[g]{parameters used in the paper!c@$c_1$ - the parameter controling the space of a computation of $\mmm_3$ (see Lemma
\ref{M3})} $c_1=c_1(\mmm_3)$ such that $|W_j|_Y\le c_1\max(|W_0|_Y, |W_t|_Y)$ for $j=0,1,\dots, t$;
moreover, $|W_j|_Y\le c_1|W_t|_Y$ if $W_0$ is a tame configuration. (Recall that $c_1$ is one of the parameters from Section \ref{param}.).

\end{lemma}

\proof
(a) Without loss of generality we assume that $i=1$. Consider the projection $H_\chi$ of the history $H$ of $\ccc$ onto the alphabet of $\chi$-rules of $\mmm_3$. By the definition of $\mmm_3$, if $\chi=\chi(j,j+1)^{\pm1}$ is a letter in $H_\chi$, then the next letter in $H_\chi$ is either $\chi\iv$ or $\chi(j-1,j)^{\pm 1}$ or $\chi(j+1,j+2)$. By Lemma \ref{M31}, for the every letter $\chi$, the word $H_\chi$ contains at most one occurrence of $\chi^{\pm 1}$. This implies that $H_\chi\equiv \chi(1,2)\chi(2,3)\chi(3,4)\chi(4,5)\chi(5,6)$.

Therefore the history of $\ccc$ has the form
$$\chi(1,2)H_1\chi(2,3)H_2\chi(3,4)H_3\chi(4,5)H_4\chi(5,6),$$
for some subhistories $H_1$, $H_2, H_3, H_4$ which do not contain $\chi$-rules. By the definition of $\mmm_3$, each $H_i$ is the history of a computation  of a copy of one of the $S$-machines: $\bmm_2,\rhh,\rhr$ (because rules of any two of these mahines have disjoint domains). This implies that $H_1, H_2, H_3, H_4$ are histories of computations of copies of $\bmm_2, \rhh, \bmm_2, \rhr$, respectively.

Let $UV$ be a history 2-letter subword in the base $B$ of the computation $\ccc$. The $Y$-projection $u$ of the $UV$-sector of $W_1$ is a word in a left alphabet, while the $Y$-projection of the $UV$-sector of $W_1\cdot H_1$ is a word in the corresponding right alphabet. Each rule $\theta$ of $H_1$ multiples the $Y$-projection of the $UV$-sector by a letter from the left alphabet on the left and by a letter from the right alphabet on the right. The two letters correspond to the rule $\theta$. Therefore $u$ must be a copy of $H_1$. In particular, this implies that the $Y$-projections of all history sectors of $W_1$ and $W_1\cdot H_1$ are copies of $H_1$.

Applying Lemma \ref{prim} (3) to the subcomputation $W_1\cdot H_1\chi(2,3) \to\dots, W_1\cdot H_1\chi(2,3)H_2$ and considering the history $UV$-sector again, we deduce that $H_2$ is a copy of $${\bar H}_1 \zeta^{(12)}({\bar H}_1')\iv$$ where ${\bar H}_1$ is the mirror image of $H_1$ and ${\bar H}_1'$ is a copy of $H_1$. Moreover $H_2$ is uniquely determined by $W_1\cdot H_1$, hence by $W_1$.

Similar arguments work for the rest of the computation $\ccc$: $H_3$ is a copy of $H_1\iv$ and $H_4$ is a copy of $H_1\zeta^{(12)}H_1'$. This implies (a).

(b) follows from the same argument as (a).

(c) If the history $H$ of $\ccc$ does not have $\chi$-rules, then $\ccc$ is a computation of a copy of one of the $S$-machines $\bmm_2,\rhh, \rhr$ and we can apply Lemmas \ref{Hprim} (b) and \ref{WV}.

Suppose that $H$ contains a $\chi$-rule. Then $H=H_1H_2H_3$ where $H_1, H_3$ do not contain $\chi$-rules, but $H_2$ starts and ends with $\chi$-rules (it is possible that $||H_2||=1$). Let $W_k=W_0\cdot H_1$, $W_s=W_0\cdot H_1H_2=W_t\cdot H_3\iv$. Then $W_k$ is tame being in the domain of a $\chi$-rule. Hence by Lemmas \ref{Hprim} (b) and \ref{WV} for every $i$ between $0$ and $k$ $|W_i|_Y$ does not exceed $c|W_0|_Y$ for some constant $c$. The same argument shows that for $i$ between $s$ and $t$, $|W_i|_Y$ does not exceed $c|W_t|_Y$. The proof of part (a) describes the subcomputation $W_k\to\dots\to W_s$ in detail. This description and Lemma \ref{WV} imply that for $i$ between $k$ and $s$, $|W_i|_Y$ does not exceed a constant times the maximum of $|W_k|_Y$ and $|W_s|_Y$. This implies (c).

\endproof

\begin{lemma}\label{I3A3} (1) If a word $\alpha^k$ is accepted by the Turing machine $\mmm_0$, then for some word $H$, there is a reduced computation $I_3(\alpha^k,H)\to\dots\to A_3(H)$ of the $S$-machine $\mmm_3$.

(2) If there is a computation $\ccc\colon I_3(\alpha^k,H)\to\dots\to A_3(H')$ of  $\mmm_3$, then the word $\alpha^k$ is accepted
by $\mmm_0$ and $H'\equiv H$.
\end{lemma}

\proof (1) is obvious from the definition of $\mmm_3$ (see the informal definition of $\mmm_3$ at the beginning of Section \ref{M3M5}): $H$ is a copy of the history of a computation of $\bmm_2$ accepting $I_2(\alpha^k)$ (which exists by Lemma \ref{I2A2} (1)).

(2) The word $I_3(\alpha^k,H)$ is in the domain of a rule from $\mmm_{3,1}$ while $I_3(H')$ is in the domain of a rule from $\mmm_{3,4m+1}$. For different $i,j$ domains of rules from $\mmm_{3,i}$ and $\mmm_{3,j}$ are disjoint and if rules of sets $\mmm_{3,i}$ and $\mmm_{3,i+1}$ appear in a computation, the computation must also contain the $\chi$-rule $\chi(i,i+1)$. Therefore the projection of the history of $\ccc$ onto the alphabet of $\chi$-rules must contain a subword $\chi(1,2)\chi(2,3)$. Hence $\ccc$ must contain a subcomputation $\ddd$ with history of the form $\chi(1,2)H_1\chi(2,3)$,
where  $H_1$ is the history of a computation of a copy of $\bmm_2$ of the form
$I_2(\alpha^{\ell},H)\to\dots\to A_2(H'')$ for some $\ell, H''$ and the rules in $\ccc$ applied before this $\chi(1,2)$ are from $\mmm_{3,1}$. Since rules of $\mmm_{3,1}$ do not modify the input sector, $k=\ell$. Therefore $\alpha^k$ is accepted by $\bmm_2$. By Lemma \ref{I2A2} then $\alpha^k$ is accepted by $\mmm_0$ and $H''\equiv H$. The fact that $H'\equiv H$ is proved in the same way as in Lemma \ref{I2A2} (2).
\endproof

\subsection{\texorpdfstring{$\mmm_4$}{M4} and \texorpdfstring{$\mmm_5$}{M5}}\label{45}

Let $B_3$ be the standard base of $\mmm_3$ and $B_3'$ be its disjoint copy. By \index[g]{S@$S$-machine!M4@$\mmm_4$} $\mmm_4$ we denote
the $S$-machine with standard base $B_3(B'_3)^{-1}$ and rules $\theta(\mmm_4)=[\theta, \theta]$,
where $\theta\in \Theta$ and $\Theta$ is the set of rules of $\mmm_3$. So the rules of $\Theta(\mmm_4)$ are the same
for $\mmm_3$-part of $\mmm_4$ and for the mirror copy of $\mmm_3$. Therefore we will denote $\Theta(\mmm_4)$ by
$\Theta$ as well.  The sector between the last state letter of $B_3$ and the first state letter of $(B_3')^{-1}$ is locked by any rule from $\Theta$.

The 'mirror' symmetry of the base  will be used in Lemma \ref{001}.

\medskip

The $S$-machine \index[g]{S@$S$-machine!M5@$\mmm_5$}$\mmm_5$ is a circular analog of $\mmm_4$. We add one more base letter $\tt$ to the hardware of $\mmm_4$. So the standard base $B$ of $\mmm_5$ it $\{\tt\}B_3(B_3')^{-1}\{\tt\}$, where the part $\{\tt\}$ has only one letter $\tt$
and the first part $\{\tt\}$ is identified with the last part. For example, $\{\tt\}B_3(B_3')^{-1}\{\tt\}B_3(B_3')^{-1}$ can be a base of an admissible word for $\mmm_5$. Furthermore,  sectors involving $\tt^{\pm 1}$ are
 locked by every rule from $\Theta$. The accordingly modified sets $\mmm_{3,i}$ are denoted by $\mmm_{5,i}$.

In particular, for $\mmm_5$, we have the start and stop words $I_5(\alpha^k,H)$ and $A_5(H)$ similar to the configurations
$I_3(\alpha^k,H)$ and $A_3(H)$, and the following analog of Lemma \ref{I3A3} can be proved in the same way as Lemma \ref{I3A3}.

\begin{lemma}\label{I5A5} (1) If a word $\alpha^k$ is accepted by the Turing machine $\mmm_0$, then for some word $H$, there is a reduced computation of $I_5(\alpha^k,H)\to\dots\to A_5(H)$ of the $S$-machine $\mmm_5$.

(2) If there is a computation $\ccc\colon I_5(\alpha^k,H)\to\dots\to A_5(H')$ of  $\mmm_5$, then the word $\alpha^k$ is accepted
by $\mmm_0$ and $H'\equiv H$.
\end{lemma}

\begin{df}\label{faul}
We call the base of an admissible word of an $S$-machine
\index[g]{S@$S$-machine!admissible words of an $S$-machine!faulty base of an admissible word}
{\it faulty} if

\begin{enumerate}

\item[(1)]{it starts and ends with the same base letter,}

\item[(2)]{only the first and the last letters can occur in the base twice}

\item[(3)]{it is not a reduced word.}

\end{enumerate}

\end{df}

\begin{lemma} \label{nonsta} There is a constant $C=C(\mmm_5)$, such that for every reduced computation $\ccc\colon W_0\to\dots\to W_t$
of $\mmm_5$ with a faulty base and every $j=0,1,\dots,t$, we have $|W_j|_Y\le C\max(|W_0|_Y, |W_t|_Y)$.
\end{lemma}

\proof {\bf Step 1. } One may assume that $|W_r|_Y > \max(|W_0|_Y, |W_t|_Y)$ for every $0<r<t$ since
otherwise it suffices to prove the statement for two shorter computations $W_0\to\dots\to W_r$
and $W_r\to\dots\to W_t$. Since $\chi$-rules do not change the length of configurations, the
history $H$ of $\ccc$ cannot start or end with a $\chi$-rule.

{\bf Step 2.} If the history $H$ of $\ccc$ has no $\chi$-rules, then the statement with  $C\ge 18$ follows
from Lemmas \ref{Hprim} (a), \ref{ewe} and \ref{9}.

{\bf Step 3.} If there is only one $\chi$-rule $\chi$ in $H$, then $H=H'\chi^{\pm 1} H''$, where $H'$ is a copy of the  history of a computation of a copy of $\rhh$ or $\rhr$ and $H''$ is the history of a computation of a copy of $\bmm_2$ (or vice versa). For the computation
$W_r\to\dots\to W_0$ with history $(H')^{-1}$, we have $|W_r|_Y\le |W_0|_Y$ by Lemmas
\ref{Hprim} (a) and \ref{ewe}.
This contradicts the assumption of  Step 1, and so one may assume further that
$H$ has at least two $\chi$-rules.

{\bf Step 4.} The base $B$ of the computation $\ccc$ has no history sectors $PP^{-1}$-, $R^{-1}R$-,  $QQ^{-1}$-, or $Q^{-1}Q$-sectors,
since every $\chi$-rule locks the $PQ$- and $QR$-sectors of the standard base.

The same statement is true for the mirror copies
of the above-mentioned sectors, and this stipulation
works throughout the remaining part of the proof.

{\bf Step 5.}
Assume that the history $H^{\pm 1}$ is of the form $H_1\chi(i-1,i)H_2\chi(i,i+1)H_3$ for
some $i$, where $H_2$ is the history of a computation of a copy of $\bmm_2$. Since $B$ is not reduced, there is a 2-letter subword of the base of the form $U^{\pm1}U^{\mp1}$ (for some part $U$ of the set of state letters).
By Lemma \ref{qqiv}, then this subword must be a history subword of the form $P\iv P$ or $RR\iv$ since
every sector of the standard base of $\mmm_3$, except for history $RP$-sectors is locked either by $\chi(i-1,i)$ or by $\chi(i,i+1)$.

Let us consider the case of $P^{-1}P$ since the second case is similar.
Depending on the parity of $i$ either a prefix $H_3'$ of $H_3$ is the history of a computation of a copy of $\rhh$ or the suffix $H_1'$ of $H_1$ is the history of a computation of a copy of $\rhh$. These two cases are similar so we consider only the first one.

Then between the $P$-letter of the $P\iv P$-sector of an admissible word in the subcomputation of $\ccc$ with the history $H_3'$ and the corresponding $R$-letter in that admissible word, there is always a $Q$-letter or a $P\iv$-letter, hence the $P$-letter never meets the corresponding $R$-letter during that subcomputation and no transition rules rules can apply to any of the admissible words of that subcomputation. Therefore $H_3'=H_3$ and for the subcomputation $\ccc'\colon W_s\to\dots\to W_t$ of $\rhh$ with history  $H_3$  we have $|W_s|_Y\le |W_t|_Y$ by Lemmas
\ref{prim} (1) and \ref{ewe}. This contradicts  Step 1, and so
the assumption made in the beginning of Step 5 was false.

{\bf Step 6.} Assume that there is a history of a subcomputation of $\ccc$ of the form $H_1\chi H_2\chi^{-1}H_3$, where $\chi$ is a $\chi$-rule,
$H_2$ is the history of a computation of a copy of $\bmm_2$. Then we claim that the base of $\ccc$ has no history $P\iv P$- or $RR\iv$-sectors. To prove this, we consider
only the former case since the latter one is similar.

If the subcomputation $\ccc'$ of $\ccc$ with history $H_3$ starts with an admissible word $W$
having in the $P^{-1}P$-sector all $Y$-letters from the right alphabets, then, as in Step 5, $H_3$ corresponds to the work of $\rhh$, which
gives a contradiction as in item 5.

If the $P^{-1}P$-sector of $W$ has all $Y$-letters from the left alphabet, then the subcomputation of $\ccc\iv$ with history $\chi H_2^{-1}$ will conjugate the $Y$-projection of that sector by a non-empty
reduced word from the right alphabet. Therefore in the last admissible word of that subcomputation, there will still be letters from both left and right alphabets, and so it cannot be in the domain of any $\chi$-rule or its inverse, a contradiction.

Together with Step 4, this implies that the base of $\ccc$
has no mutually inverse letters from history sectors staying next to each other.

Since the base is faulty, it must contain an input $P_1\iv P_1$ or $R_0R_1\iv$-sector.
This implies that the base does not contain input $(R_0P_1)^{\pm 1}$-sectors
since
the first and the last letters of the base are equal (say, positive)
and the base has no proper subwords with this property.
In both cases the configuration $W_r$ corresponding to
the transition $\chi\colon W_{r-1}\to W_r$ is the shortest one in $\ccc$ since the $Y$-projection of that word is of the form $\alpha^k$, each rule from $\ccc$ conjugates the $Y$-projection from the input sector, and $\alpha^k$ cannot be
shortened by any conjugation.
This contradicts Step 1.

{\bf Step 7.} It follows from items  2,3, 5 and 6 that $H=H_1\chi H_2 \chi' H_3$, for two $\chi$-rules
(or their inverses). Moreover  $H_2$ is the history of a computation $\ccc_2$ of a copy of $\rhh$ or of
$\rhr$ and $H_1, H_3$ are histories of computations $\ccc_1, \ccc_3$ of copies of $\bmm_2$, i.e.,$H$ has exactly two $\chi$-rules (otherwise $H$ has a subword which is ruled out in the previous steps of the proof).

{\bf Step 8.} We claim that we can assume that the admissible words in the computation $\ccc$ do not have a history
$(PR)^{\pm 1}$-sectors. Indeed, if such a sector exists, then for the subcomputation $\ccc_1\colon W_0\to\dots\to W_r$
with history $H_1\chi $, we have $|W_r|_Y\le c |W_0|$ by Lemma \ref{WV}. A similar estimate
is true for the subcomputation with history $\chi' H_3$ starting with some $W_s$. So in order to prove the inequality from the lemma, it suffices to
apply Step 2 to the three  subcomputations $\ccc_1, \ccc_2. \ccc_3$.

{\bf Step 9.}  Suppose that the base of $\ccc$ contains a history subword of the form $P^{-1}P$.

If the admissible word from $\ccc$ in the domain of $\chi$ has no letters
from the left alphabets, then $H_2$ is the history of a computation of a copy of $\rhh$ and the state $P$-letter will never meet the corresponding state $R$- or $Q$-letter during the computation $\ccc_2$, so an application $\chi'$ is not possible after $\ccc_2$ ends, a contradiction.

Thus we can assume that if the base of $\ccc$ contains a history subword of the form $P\iv P$, then the last admissible word of $\ccc_2$ (which is in the domain of $\chi$) contains letters from the left alphabet.

Similarly, if the base of $\ccc$ contains a history subword of the form $RR\iv$, then the last admissible word in $\ccc_2$ contains letters from the right alphabet. This implies, in particular that the base of $\ccc$ cannot contain both a history subword of the form $P\iv P$, and a history subword $R'(R')\iv$. Without loss of generality,
we will assume that there are no subwords $R'(R')\iv$.


{\bf Step 10.} It follows from Steps 4,8 and 9, that there are no unlocked by $\chi$ history sectors of the base except for $P^{-1}P$-sectors, and  if there is such a sector $UV$,then ${\mathcal C}_2$ is a computation of a copy of $\rhr$. Therefore $UV$
may contain tape letters from a left alphabet,
while every rule $\theta$ of ${\mathcal C}_1^{-1}$ multiplies this
sector from both sides by letters from a right
alphabet. So $\theta$ increases the lengths of every history sectors by 2. The rule $\chi$ locks working sectors (except for the input one), and so  by Lemma \ref{simp} (**), $\theta$ can decrease the lengths of every  working sector at most by
one. Since working sectors alternate with history
ones in any base, we have $||W_r||\le||W_0||$,
contrary to Step 1.



{\bf Step 11.} To complete the proof of the lemma, it remains to assume that there are no history
sectors in the base of $\ccc$. Then  the faulty base of $\ccc$ must contain input subwords  of the form $R_0R_0\iv$ only, because every
$\chi$-rule locks all sectors of  the standard base except for the input and history sectors. Then any admissible word of $\ccc$ from the domain of a  $\chi$-rule in $H$ is the shortest admissible word in $\ccc$  since (as in Step 6) every rule of the computation conjugates
$R_0R_0\iv$-sectors and a word $\alpha^k$
cannot be shortened by any conjugation. The lemma is proved since we can refer to Step  1 again.
\endproof


\section{The main \texorpdfstring{$S$}{S}-machine \texorpdfstring{$\mmm$}{M}}

\subsection{The definition of \texorpdfstring{$\mmm$}{M}} \label{M6}

We use the $S$-machine $\mmm_5$ from Section \ref{45}, $\rhh_m$ from Section \ref{pm} and three more easy $S$-machines to compose
the main circular $S$-machine \index[g]{S@$S$-machine!M@$\mmm$} $\mmm$ needed for this paper. The standard base of $\mmm$ is the same as the standard base of $\mmm_5$, i.e.,$\{\tt\}B_3(B_3')^{-1}$, where $B_3$ has the form (\ref{basec}).
However
we will use ${\tilde Q}_0$ instead of $Q_0$, ${\tilde R}_1$ instead of $R_1$ and so on to denote parts of the set of state letters since $\mmm$ has more
state letters in every part of its hardware.

The rules of $\mmm$ will be partitioned into five sets ($S$-machines) $\mtt_i$ ($i=1,\dots,5$) with transition rules $\theta(i,i+1)$ \index[g]{S@$S$-machine!M@$\mmm$!transition rules $\theta(i,i+1)$} connecting $i$-th and $i+1$-st sets.
The state letters are also disjoint for different sets $\mtt_i$. It will be clear that  ${\tilde Q}_0$ is the disjoint union of 5 disjoint sets including $Q_0$, ${\tilde R}_1$ is the disjoint union of five disjoint sets including $R_1$, etc.

By default, every transition rule $\theta(i,i+1)$ of $\mmm$ locks a sector if this sector
is locked by all rules from $\mtt_i$ or if it is locked by all rules from $\mtt_{i+1}$.
It also changes the end state letters of $\mtt_i$ to the start state letters of $\mtt_{i+1}$.

The \index[g]{start configuration $W_{st}$ of $\bf M$} {\it start configuration} $W_{st}$ of $\mtt$ is $\tilde t b_3(b'_3)^{-1}$,
where $b_3$ and $b'3$ are obtained by replacing every base letter of $B_3$ and $B'_3$
by special start letter. The start rule $\theta_1$ of $\bf M$ changes the letters from $b_3$ and $b'_3$ to
their copies and starts the work of the rules from the set {\bf $\mtt_1$}.

\medskip

 Set {\bf $\mtt_1$} inserts input words in the input sectors. The set contains only one positive rule inserting the letter $\alpha$ in the input sector next to the left of a letter $p$ from ${\tilde P}_1$. It also inserts a copy $\alpha\iv$  next to the right of the corresponding letter $(p')\iv$ (the similar
mirror symmetry is assumed in the definition of all other rules.)
So the positive rule of $\mtt_1$ has the form $$[q_0\tool q_0, r_1\to r_1, p_1\tool \alpha p_1, ..., (p_1')^{-1}\to  (p_1')^{-1}\alpha^{-1}, (r_1')^{-1}\tool (r_1')^{-1}, t\tool t]$$

The rules of $\mtt_1$ do not change state letters, so it has one state letter in each part of its hardware.

The connecting rule $\theta(12)$ changes the state letters of $\mtt_1$ to their copies in a disjoint alphabet. It locks all sectors except for the input sector ${\tilde R}_0{\tilde P}_1$ and the mirror copy of this sector.

\medskip

 Set {\bf $\mtt_2$} is a copy of the $S$-machine $\rhh_m$ working in the input sector and its mirror image in parallel, i.e.,we identify the standard base
of $\rhh_m$ with ${\tilde R}_0 {\tilde P}_1{\tilde Q}_1$. The connecting rule
$\theta(23)$ locks all sectors except for the input sector ${\tilde R}_0{\tilde P}_1$ and its mirror image.
\medskip

 Set {\bf $\mtt_3$} inserts history in the history sectors.
This set of rules is a copy of each of the left alphabets $X_{i,l}$ of
the $S$-machine $\mmm_2$. Every positive rule of $\mtt_3$ inserts a copy of the corresponding positive letter in every
history sector ${\tilde R}_i{\tilde P}_{i+1}$ next to the right of a state letter from ${\tilde R}_i$.

Again, $\mtt_3$ does not change the state letters, so each part of its hardware contains one letter.

The transition rule $\theta(34)$ changes the state letters to their copies from Set $\mmm_{5,1}$ of $\mmm_5$.
It locks all sectors except for the input sectors and the history sectors. The history sectors in admissible words from the domain of $\theta(34)$ have  $Y$-letters from the left
alphabets $X_{i,l}$ of the $S$-machine $\mmm_5$.

\medskip

 Set {\bf $\mtt_4$} is a copy of  the $S$-machine $\mmm_5$. The transition rule
$\theta(45)$ locks all sectors except for history ones. The admissible words in the domain of
$\theta(45)$ have no letters from right alphabets.

\medskip

Set {\bf $\mtt_5$.} The positive rules from $\mtt_5$ simultaneously
erase the letters of the history sectors from the  right  of the state letter from ${\tilde R}_i$. That is, parts of the rules are of the form $r\to ra\iv $ where $r$ is a state letter from ${\tilde R}_i$, $a$ is a letter from the left alphabet of the history sector.

\medskip

Finally the accept rule $\theta_0$ (regarded as a transition rule) from $\mmm$ can be applied when all the sectors are empty, so it locks all the sectors and changes the end state letters of $\mmm_5$ to the corresponding end state letters of $\mmm$.
Thus, the main $S$-machine $\mmm$ has unique accept configuration which we will denote by \index[g]{S@$S$-machine!M@$\mmm$!W@$W_{ac}$, the accept word of $\mmm$} $W_{ac}$.


For every $i=1,2,3,4$, we will sometimes denote $\theta(i,i+1)\iv$ by $\theta(i+1,i)$.

\subsection{Standard computations of \texorpdfstring{$\mmm$}{M}}\label{SC}

We say that the history $H$ of a computation of  $\mmm$ (and the computation itself) is \index[g]{S@$S$-machine!M@$\mmm$!eligible computation of $\mmm$} \index[g]{S@$S$-machine!M@$\mmm$!eligible history of computation of $\mmm$}{\it eligible} if it has no
neighboring mutually inverse letters except possibly for the subwords $\theta(23)\theta(23)^{-1}$.
(The subword $\theta(23)^{-1}\theta(23)$ is not allowed.)
\begin{rk} \label{elig}
Clearly the history $H^{-1}$
is eligible if and only if $H$ is. Every reduced computation is eligible.
\end{rk}

Considering eligible computations instead of just reduced computations is necessary for our interpretation of $\mmm$ in a group.

The history $H$ of an eligible computation of $\mmm$ can be factorized so that every
factor is either a transition rule $\theta(i,i+1)^{\pm 1}$ or a maximal non-empty product of rules of one of the sets $\mtt_1 - \mtt_5$. If, for example, $H=H'H''H'''$,
where $H'$ is a product of rules from $\mtt_2$, $H''$ has only one rule
$\theta(23)$ and $H'''$ is a product of  rules from $\mtt_3$, then we say
that the \index[g]{S@$S$-machine!M@$\mmm$!eligible computation of $\mmm$!step history of a computation of $\mmm$} {\it step history} of the computation is $(2)(23)(3)$.
Thus the step history of a computation is a word in the alphabet $\{(1),(2), (3), (4),(5), (12), (23), (34), (45), (21), (32), (43), (54)\}$, where $(21)$
is used for the rule $\theta(12)^{-1}$ an so on.
For brevity, we can omit some transition symbols, e.g. we may use $(2)(3)$ instead of $(2)(23)(3)$ since the only rule connecting Steps 2 and 3 is $\theta(23)$.


If the step history of a computation consists of only one letter $(i)$, $i=1,\dots,5$, then we call it a \index[g]{S@$S$-machine!M@$\mmm$!one step computation of $\mmm$} {\em one step computation}.
The computations with step histories $(i)(i,i\pm 1)$, $(i\pm 1, i)(i)$  and $(i\pm 1, i)(i)(i,i\pm 1)$ are also considered
as one step computations.
Any eligible one step computation is always reduced by definition.

The step history of any computation cannot contain certain subwords. For example, $(1)(3)$ is not a subword of any step history because domains of rules from $\mtt_1$ and $\mtt_3$ are disjoint. In this subsection, we eliminate some less obvious subwords in step histories of eligible computations.

\begin{lemma} \label{212} If the base of a computation $\ccc$ has at least one history subword $UV$, then there are no reduced computations $\ccc$ of $\mmm$ with step history

(1) $(34)(4)(43)$ or $(54)(4)(45)$,  provided $UV\equiv ({\tilde R}_{i-1}{\tilde P}_i)^{\pm 1}$ for some $i$,

(2) $(23)(3)(32)$.

\end{lemma}
\proof
(1) We consider only the step history $(34)(4)(43)$ since the second  case is similar. Let $W_0$ be the first admissible word of $\ccc$.
Suppose that  the history $H=\theta(34)H'\theta(43)$ of $\ccc$ has $\chi$-letters. By Lemma \ref{M31} each $\chi$ letter $\chi^{\pm 1}$ appears in $H_\chi$ only once. Each $\chi$-rule changes the state letters, and words in the domains of different (positive) $\chi$-rules have different state letters. Therefore $W_0\cdot \theta(34)H'$ has different state letters than $W_0$, hence $W_0\cdot \theta(34)H'$ is not in the domain of $\theta(43)$, a contradiction.

If $H'$ has no $\chi$-letters, then it is a history of $\rhr$, and
we obtain a contradiction with Lemma \ref{prim} (4) (and Remark \ref{r-prim}).

(2) Suppose the step history of $\ccc$ is $(23)(3)(32)$. Since the history sectors are locked by $\theta(23)^{\pm 1}$, the history subwords in the base of $\ccc$ must have the form $(R_{i-1}P_i)^{\pm 1}$ for some $i$. Every rule of $\mtt_3$ inserts a letter next to the left of every $P_i$-letter, different rules insert different letters, same letter for the same rule.
Since at the beginning and at the end of the subcomputation with step history (3) all history sectors are empty of $Y$-letters,
the word inserted during the subcomputation must be freely trivial. That contradicts the assumption that this subcomputation is reduced.
\endproof

By definition,  the rule $\theta(23)$ locks all history sectors of the standard base of $\mmm$
except for the input sector ${\tilde R}_0{\tilde P}_1$ and its mirror copy.
Hence every admissible word in the domain of
$\theta(23)^{-1}$ has the form \index[g]{W@$W(k,k')$ - a word in the domain of $\theta(23)$} $W(k,k')\equiv w_1\alpha^kw_2(\alpha')^{-k'}w_3$,
where $(\alpha')^{-1}$ is the mirror copy of $\alpha$, $k$ and $k'$ are integers, and $w_1, w_2, w_3$ are fixed words in state letters;
$w_1$ starts with $\tt$. Recall that $W_{ac}$ is the accept word of $\mmm$.

\begin{lemma} \label{121}  There are no reduced computations of $\mmm$ with the standard base whose step history is $(12)(2)(21)$  or $(32)(2)(23)$.
\end{lemma}

\proof  Consider only the step history $(12)(2)(12)$. Thus the history $H$ of the computation is $\theta(12)H'\theta(21)\iv$ and $H'$ is a computation of a copy of $\rhh_m$ working in the input sectors of admissible words of $\mmm$. Then applying Lemma \ref{prim} (4) and Remark \ref{r-prim} we can conclude that $H'$ is empty, a contradiction.
\endproof

\begin{lemma} \label{resto} Let
a reduced computation  $\ccc\colon W_0\to\dots\to W_t$
have the history $H$ of
the form (a) $\chi(i-1,i)H'\chi(i,i+1)$ (i.e.,the $S$-machine
works as $\mmm_3$ with step history (4)) or (b) $\zeta^{(i-1,i)} H'\zeta^{(i,i+1)}$ (i.e.,it works as $\rhh_m$
with step history (2)).

 Then the base of the computation $\cal C$ is a reduced word, and
 all configurations of $\ccc$ are uniquely defined by the history
 $H$ and the base of $\ccc$.
Moreover,
$H'$ is the copy of the maximal $Y$-word contained in arbitrary history (resp., input) sector of $W_0$.
\end{lemma}

\proof (a) Every history sector of the standard base is locked either by one of the rules $\chi(i-1,i), \chi(i,i+1)$, or by a rule of $H'$. Every non-history sector of the standard
base is also locked either by $\chi(i-1,i)$ or by $\chi(i,i+1)$. It follows from Lemma \ref{qqiv} that the base of $\cal C$ is a reduced word.
By Lemma \ref{prim} (3), the histories
of the primitive $S$-machines subsequently restore the tape words in all history sectors.
Since one of the rules $\chi(i-1,i), \chi(i,i+1)$ locks all non-history sectors,
Lemma \ref{Hprim} applied to $\ccc$
gives the other statements.

(b) The same proof up to change of the history sectors to the input ones.
\endproof

\begin{lemma}\label{SH} The step history of every eligible computation of $\mmm$ with standard base
either

(A) contains one of the words $(34)(4)(45)$, $(54)(4)(43)$, $(12)(2)(23)$, $(32)(2)(21)$ as a subword or

(B) is a subword of one of the words $$(4)(45)(5)(54)(4), (4)(43)(3)(34)(4), (2)(23)(3)(34)(4),$$ $$(4)(43)(3)(32)(2), (2)(21)(1)(12)(2), (2)(23)(32)(2).$$
\end{lemma}

\proof The statement is obvious if there are neither (2) nor (4) in the step history.
Lemmas \ref{212} (1) (Lemma \ref{121}) implies that if (4) (resp. (2)) is not the first or the last letter in the
step history then it can occur in a  subword of the form $(34)(4)(45)$ or $(45)(4)(34)$
(resp., $(12)(2)(23)$ or $(23)(2)(12)$), i.e.,we have Property (A).

If the first letter in the step history is (2) and Property (A) fails, then the same lemmas give us the longest
possible step histories $(2)(23)(32)(2)$, $(2)(1)(2)$ and $(2)(23)(3)(34)(4)$. The assumption
that the last letter in the step history is (2) adds one more possible longest step history word $(4)(43)(3)(32)(2)$.

Similarly, we may assume that (4) is either the first or the last letter in the step
history and conclude that the step history is a subword of one of the words
$(4)(5)(4)$, $(4)(3)(4)$, $(2)(3)(4)$ and $(4)(3)(2)$ provided Property (A) fails.
\endproof


\begin{lemma} \label{I6A6} (1) If the word $\alpha^k$ is accepted by the Turing machine $\mmm_0$, then there
is a reduced computation of $\mmm$, $W(k,k)\to\dots\to W_{ac}$ whose history has no rules of $\mtt_1$ and $\mtt_2$.

(2) If the history of a computation $\ccc\colon W(k,k)\to\dots\to W_{ac}$ of $\mmm$ has no rules of $\mtt_1$ and $\mtt_2$, then the word $\alpha^k$ is accepted by $\mmm_0$.
\end{lemma}

\proof (1) By Lemma \ref{I5A5}, there is a computation $I_5(a_k,H)\to\dots\to A_5(H)$
of the $S$-machine $\mmm_5$ for some $H$. So we have the corresponding computation
of $\mtt_4$:  $${\mathcal D}\colon I_6(a_k,H)\to\dots\to A_6(H).$$ Now the computation of $\mtt_3$ inserting letters
in history sectors and a computation of $\mtt_5$ erasing these letters extend $\mathcal D$ and provide us with a computation $W(k,k)\to\dots\to I_6(a_k,H)\to\dots\to A_6(H)\to\dots\to W_{ac}$.

(2) By Lemma \ref{212} (1), the step history of $\ccc$ begins with $(3)(4)(5)$, and so
there is a subcomputation of Set 4 of the form $I_5(\alpha^{\ell},H)\to\dots\to A_5(H)$ for some $\ell$ and $H$,
where according to Lemma \ref{I5A5} (2), the word $\alpha^{\ell}$ is accepted by $\mmm_0$.  Since
the computation of Set $\mmm_{3,3}$ does not change the input sector, we have $\ell=k$.
\endproof

\subsection{The first estimates of computations of \texorpdfstring{$\mmm$}{M}}

\begin{lemma}\label{E} Let $\ccc\colon W_0\to\dots\to W_t$ be a computation of $\mmm$ satisfying Property (B) of Lemma \ref{SH} or any computation of $\mmm$ with step history of length at most 2.  Then for some constant $c_2$ (see Section \ref{param})
\index[g]{parameters used in the paper!c@$c_2$ -parameter controlling the space and length of computations of $\mmm$ satisfying Property (B) (see Lemma \ref{E})}

(a) $|W_j|_Y\le c_2\max(|W_0|_Y, |W_t|_Y)$
for $j=0,1,\dots, t$;

(b) $t\le c_2^2(||W_0||+ ||W_t||)$.
\end{lemma}

\proof
(a) If $\ccc$ is a one-step computation and its step history is (1), (3), or (5), then Statement (a) follows from Lemma \ref{gen} (c). For step history (2) (resp. (4) it follows from Lemma \ref{Hprim} (a) (resp., Lemma \ref{M3} (c)).

If there is a transition  rule $\theta(i,i+1)$ of $\mmm$ in the history $H$ of $\ccc$, then $H$ can be decomposed
in at most three factors $H=H_1H_2H_3$, where $H_2$ is a one-step computation of step history (1), (3) or (5), or $H_2=(23)(32)$
and $H_1$, $H_3$, if non-empty, are of step history (2) or (4). Respectively, the computation $\ccc$
is a composition of at most three subcomputations $\ccc_1\colon W_0\to\dots\to W_r$, $\ccc_2\colon W_r\to\dots\to W_s$ and $\ccc_3\colon W_s\to\dots\to W_t$. Now we can bound $|W_r|_Y$ and $|W_s|_Y$
by $c\max(|W_0|_Y, |W_t|_Y)$ applying either Lemma \ref{Hprim} (a) (for  step history (2)) or Lemma \ref{M3} (c) (for step history  (4)) to $\ccc_1$ and $\ccc_3$. The same lemmas applied to
subcomputations $\ccc_1$, $\ccc_2$ and $\ccc_3$ completes the proof since we can assume that $c_2\gg c_1$ (see Section \ref{param}).


(b) It suffices to bound the lengths of at most three one step subcomputations $\ccc'\colon W_j\to\dots\to W_k$, where
$\max(|W_j|_Y,|W_k|_Y)\le c_2\max(|W_0|_Y,|W_t|_Y)$ by (a).
For step history (1), (3) or (5), the history lengths are bounded by Lemma \ref{gen} (b).
For (2), we refer to Lemma \ref{Hprim} (b). The computation
with step history (4) has at most $4m$ $\chi$-rules in the history as follows
from Lemma \ref{M31}. So it has at most $4m+1$ maximal subcomputations of the form $W_l\to\dots\to W_s$, corresponding
to one of the $4m+1$ subsets $\mmm_{3,i}$ of the set of rules of $\mmm_3$, where
$\max(|W_l|_Y,|W_s|)\le c_2\max(|W_0|_Y,|W_t|)$ by part  (a) of the lemma.
Hence we have the same upper bound for $s-l$ by Lemmas \ref{prim} (3)
(if it is a computation of $\rhh$) and \ref{w} (if it is acomputation of $\mmm_2$). This completes the proof of the first inequality since we have $ c_2\gg m$ (Section \ref{param}). 
\endproof

\subsection {Computations of \texorpdfstring{$\mmm$}{M} with faulty bases}\label{fau}

\begin{lemma} \label{nonst} For every eligible computation $\ccc\colon W_0\to\dots\to W_t$
of $\mmm$ with a faulty base and every $j=0,1,\dots,t$, we have $|W_j|_Y\le c_1\max(|W_0|_Y,|W_t|_Y)$.
\end{lemma}

{\bf Step 1.} As in Step 1 of the proof of Lemma \ref{nonsta}, one may assume that
$|W_j|_Y> \max(|W_0|_Y,|W_t|_Y)$ if $1<j<t$ and so the history $H$ of $\ccc$ neither starts nor ends with a transition  rule $\theta(i,i+1)^{\pm 1}$.

{\bf Step 2.} If $\ccc$ is a one step computation and $(i)$ is its step history, then the statement follows from Lemma
\ref{gen} (c) for $i=1,3,5$,  (since $c_1\ge 2$), Lemma \ref{Hprim} (a) for $i=2$ (since $c_1\ge 2$) and Lemma \ref{nonsta} for $i=4$ (since $c_1\ge C$).  Hence one may assume
further that $H$ contains a transition rule $\theta(i,i+1)$ of $\mmm$ or its inverse.

{\bf Step 3.} Assume that $\ccc$ (or the inverse computation) has a transition rule $\theta(23)$, $W_{j+1}=W_j\cdot\theta(23).$ Recall that the $\theta(23)$ does not lock only the input ${\tilde R}_0{\tilde P}_1$-sector  and its mirror copy.  So by Lemma \ref{qqiv}, we should have an input subword ${\tilde R}_0{\tilde R}_0^{-1}$ or ${\tilde P}_1^{-1}{\tilde P}_1$ in the faulty base. Moreover, we must have exactly two such input subwords in the base and no subwords $({\tilde R}_0{\tilde P}_1)^{\pm 1} $ since the first and the last
letters of the base are equal (e.g., positive) and the base has no
proper subwords with this property (see Definition \ref{faul}).

The input sectors of both $W_j$ and $W_{j+1}$ have $Y$-projections of the
form $\alpha^k$, and they are not longer than the corresponding $Y$-words
in the input sectors of any other $W_i$ since $\alpha^k$ cannot be shortened  by
conjugation. It follows that $|W_j|_Y, |W_{j+1}|_Y\le \max(|W_0|_Y,|W_t|_Y)$ contrary to Step 1. Thus, one may assume further that $H$ has no letters $\theta(23)^{\pm 1}$. In particular, $\ccc$ is a reduced computation.

The same argument eliminates letters $\theta(12)^{\pm 1}$ from $H$,
and so the letter (1) from the step history of $\ccc$. Hence one can assume that the step history contains neither (1) nor (2).

{\bf Step 4.} Suppose $H$ (or $H^{-1}$) contains a subhistory $H'\theta(45)$, where $H'$ is a maximal subword of $H$ which is word in $\mtt_4$ (which is a copy of the $S$-machine $\mmm_5$).
By Lemma \ref{qqiv}, the faulty base of the computation $\ccc$ contains one of the history subwords ${\tilde R}_{i-1}{\tilde R}_{i-1}^{-1}$ or ${\tilde P}_i^{-1}{\tilde P}_i$ for some $i$, because all non-history sectors are locked by $\theta(45)$.

Suppose the base of $\ccc$ contains a history subword ${\tilde R}_{j-1}{\tilde R}_{j-1}\iv$ for some $j$. The word $H'$ must have a suffix which is a word in the alphabet of a copy of $\rhr$ working in parallel in the history sectors (see the definition of $\mmm_{3,4m+1}$). The state letters from ${\tilde R}_{j-1}$ in the ${\tilde R}_{j-1}{\tilde R}_{j-1}\iv$-sector will then never meet a letter from either ${\tilde Q}_{j-1}$ or ${\tilde P}_j$. Therefore $H'$ cannot contain the transition rule $\chi(4m,4m+1)^{\pm 1}$ or $\theta(45)\iv$. Thus $H'$ is a prefix of $H$, is a computation of a copy of $\rhr$, and by Lemma \ref{Hprim} (a) applied to the subcomputation of $\ccc\iv$ with history $(H')\iv$, we get a contradiction with Step 1 because admissible words in the domain of $\theta(45)\iv$ is tame.

Suppose the base of $\ccc$ contains  a subword $({\tilde R}_{i-1}{\tilde P}_i)^{\pm 1}$. Then $H$ has no subword $\theta(45)^{-1}H'\theta(45)$ by Lemma \ref{212} (1).
If $H'$ has neither transition rules nor $\chi$-rules, then we have a contradiction by Lemma \ref{Hprim} (a). Hence $H$ has a subword $\chi(4m, 4m+1)H''\theta(45)$, but then by Lemma \ref{prim} (3), $H'$ has a rule locking all the sectors
${\tilde R}_{i-1}{\tilde P}_i$ of the standard base, and we get a contradiction with Lemma \ref{qqiv}.

Finally suppose  all history subwords in the base of $\ccc$  have the form ${\tilde P}_i^{-1}{\tilde P}_i$. Then the rules of a copy of $\rhr$ from $H'$ do not change the history sectors of admissible words in the corresponding subcomputation $\ccc'$ of $\ccc$, hence the lengths of all admissible words in $\ccc'$ stay the same.  Moreover since the state letters in the history sectors do not change during the subcomputation $\ccc$, none of the admissible words in that subcomputation is in the domain of $\chi(4m,4m+1)^{\pm 1}$. Therefore the rules of $H'$ do not change the lengths of admissible words, and { either $H'$ is a prefix of $H$ and} we get a contradiction with Step 1  or we have the subhistory $\theta(45)^{-1}H'\theta(45)$.

In the latter case, we consider the maximal subhistory $H''$ of type 5
following after the rule $\theta(45)$ (or before $\theta(45)^{-1}$).
All the admissible words of the corresponding subcomputation $\cal C''$
have equal lengths since the base has no letters ${\tilde R}_i$.
Arguing in this way we see that the history of $\cal C$ has Steps
4 and 5 only, and all the admissible words in $\cal C$ have equal length,
which proves the inequality of the lemma.

We can conclude that $H$ does not contain  $\theta(45)^{\pm 1}$. By Step 2, (5) is not in the step history of $\ccc$ and the only possible transition rules of $\mmm$ in $H$ are $\theta(34)^{\pm 1}$.

{\bf Step 5.} Assume that there is a subhistory of $H$ of the form $H_1\theta(34) H_2\theta(34)^{-1}H_3$, where
$H_2$ is the history of $\mmm_5$. Then the base of $\ccc$ has no history sectors of the form ${\tilde R}_i{\tilde R}^{-1}_i$ (since, as before, the
machine $\rhr$ starting with $\theta(34)$ would never end with $\chi(12)$).

If there is a history subword ${\tilde R}_{i-1}{\tilde P}_i$ in the faulty base,
then $H_2$ cannot follow by the transition rule $\theta(34)^{-1}$, by Lemma \ref{M31} if $H_2$ contains $\chi$-rules and by Lemma \ref{prim} (4) otherwise, a contradiction.

Thus the base of $\ccc$ has no $\tilde R$-letters from history sectors. It also has no ${\tilde P}_1$-letters from input sectors, because
otherwise the base would contain the letter ${\tilde R}_1$  of the history sector next to the input sector since the sectors ${\tilde P}_1{\tilde Q}_1$ and ${\tilde Q}_1{\tilde R}_1$ are locked by $\theta(34)$.

Thus, all history sectors have the form ${\tilde P}_i^{-1}{\tilde P}_i$ in the faulty base of $\ccc$, and so $H$ cannot have the rule $\chi(1,2)^{\pm 1}$ (for the same reason the rule $\chi(4m,4m+1)$ was eliminated in Step 4). But without $\chi(1,2)^{\pm 1}$,
one cannot get a rule in $H$  changing history sectors ${\tilde P}_i^{-1}{\tilde P}_i$ since the rules of $\mtt_3$ leave such sectors
unchanged. The input sectors ${\tilde R}_0{\tilde R}_0^{-1}$ of the base of $\ccc$ (if any) cannot be shorten by a subcomputation since no conjugation
shortens a power of one letter in a free group.
therefore the rules $\theta(34)^{\pm 1}$ are applied to the
shortest admissible word of $\ccc$, contrary to Step 1.

So our assumption was wrong.

{\bf Step 6.} If there is only one transition rule $\theta(34)$ in $H^{\pm 1}$, then $H^{\pm 1}=H'\theta(34) H''$, where $H''$ is the history of $\mmm_5$.
If $H''$ is the history of a copy of $\rhr$, starting
with an admissible word $W_r$, then $|W_r|_Y\le |W_t|_Y$ by Lemmas
\ref{Hprim} (a) and \ref{ewe}, contrary to Step 1. Otherwise
we have a subhistory $\theta(34)H_0\chi(1,2)$, and by Lemma \ref{prim} (3), there are no history subsectors of the form ${\tilde R}_i{\tilde R}_i^{-1}$ or ${\tilde P}_i^{-1}{\tilde P}_i$ in the base of $\ccc$.
If there is a history sector ${\tilde R}_{i-1}{\tilde P}_i$,
then one can linearly bound $|W_r|_Y$ in terms of $|W_t|_Y$ applying
Lemmas \ref{Hprim} (b) and \ref{WV} several times, namely at most
$4m+1$ times by Lemma \ref{M31}. Since $c_1\gg C, c_1\gg m$ (see Section \ref{param}) one can consider two subcomputations of $\ccc$: $W_0\to\dots\to W_r$ and $W_r\to\dots\to W_t$, and reduce the proof to Step 2.

Thus, one may assume that the base of $\ccc$ has no letters $\tilde P$ and $\tilde R$ from history sectors. This also eliminates the letter
${\tilde P}_1$ of the input sector and gives the inequality $$|W_r|_Y\le \max (|W_0|_Y, |W_t|_Y),$$ contrary to Step 1. Therefore the assumption of Step 6 was wrong.

{\bf Step 7.} It remains to consider the case when $H^{\pm 1}$ is of the form
$$H_1\theta(34)^{-1} H_2\theta(34)H_3,$$ where
$H_2$ is the history of $\mtt_3$, $H_1$ and $H_3$ are histories of $\mtt_4$, and it suffices to repeat the argument of Step 6 with  decomposition of $\mathcal C$ in the product of three subcomputaions, because we did not use there that the subword $H_1\theta(34)$ was absent.

The lemma is proved.
\endproof

\subsection{Space and length of \texorpdfstring{$\mmm$}{M}-computations with standard base}
\label{space}

Let us call a configuration $W$ of $\mmm$ \index[g]{S@$S$-machine!M@$\mmm$!accessible configuration of $\mmm$} {\it accessible} if there is a \index[g]{S@$S$-machine!M@$\mmm$!accessible computation of $\mmm$} $W$-\emph{accessible computation}, i.e., either an
accepting computation starting with $W$ or a  computation
$s_1(\mmm)\to\dots\to W$, where \index[g]{S@$S$-machine!M@$\mmm$!s@$s_1(\mmm)$ - the start configuration of $\mmm$} $s_1(\mmm)$ is the start configuration of $\mmm$ (i.e., the configuration where all state letters are start state letters of $\Theta_1$ and the $Y$-projection is empty).

\begin{lemma} \label{fea} If $W$ is an accessible configuration, then
for a constant $c_3=c_3(\mmm)$, \index[g]{parameters used in the paper!c@$c_3$ - parameter controlling the length of an accessible computations of $\mmm$ whose step history is either a sduffix of $(4)(5)$ or a prefix of $(1)(2)(3)(4)$ (see Lemma \ref{fea})}
there is a $W$-accessible computation $\ccc$ of length at most $c_3||W||$ whose step history is
 either a suffix of $(4)(5)$ or a prefix of $(1)(2)(3)(4)$.
 The $Y$-length of every configuration of $\ccc$ does not
 exceed $c_2|W|_Y$. (Recall that $c_2, c_3$ are parameters in Section \ref{param}.)
 \end{lemma}

 \proof Assume that a $W$-accessible computation $\ccc$ has (4) in its step history
 and its history $H$ has a rule $\chi(i,i+1)$ with $1<i<4m$.
 Since $\ccc$ is accessible, we have by Lemma \ref{M3} (b), a subcomputation  $W_l\to\dots\to W_r$ with history of the form (a) $\chi(i,i+1)H'\chi(i+1,i+2)$ or (b)
 $\chi(i,i+1)^{-1}H'\chi(i-1,i)^{-1}$, where $H'$ is a history of a canonical
 computation of $\mmm_5$. By Lemma \ref{resto} we also conclude
  that every history sector of $W_l$ and of $W_r$ is a copy of $H'$.
 It makes possible to accept $W_r$ using erasing rules of Set 5
 in case (a) or to construct a computation of type $(1)(2)(3)$ starting with $s_1(\mmm)$ and ending with $W_l$ in case (b).

 It follows now from Lemma \ref{M3} that one can choose a accessible computation $\ccc$ having no subhistories of type $(34)(4)(45)$ or $(45)(4)(34)$, and so Set 4 can occur only in the beginning or
 at the end of $H$. In the first case $H$ has to have type $(4)(5)$,
 and the required inequalities follow from Lemma \ref{E} since $c_3\gg c_2$.

 In the second case, the step history ends with $(3)(4)$, and the connection
 $$\theta(34)\colon W_{k-1}\to W_k$$ provides us with copies in all history sectors and in all input sectors since $W_k$ is accessible. Hence one may assume that
 the step history has the form $(1)(2)(3)(4)$. Here $|W_k|_Y\le c_1|W|_Y$ by Lemma \ref{M3} (c). The canonical computation with
 step history $(1)(2)(3)$ does not decrease the lengths of configurations. Now the required estimates follows from Lemma
 \ref{E} for four one-step subcomputations since we chose $c_3$
 after $c_2$.
 \endproof

For any accessible word $W$ we choose an accessible computation $\ccc(W)$
according to Lemma \ref{fea}. \index[g]{S@$S$-machine!M@$\mmm$!C@$\ccc(W)$ the accessible computation corresponding to an accessible word $W$}

\begin{lemma}\label{form}
Let $W_0$ be an accessible  word,  $\ccc\colon  W_0\to\dots\to W_t$ be an eligible computation of $\mmm$ and $H_0$, $H_t$ be the histories of  $\ccc(W_0)$ and $\ccc(W_t)$, respectively.
 Then for some constants $c_4, c_5$ (see Section \ref{param}) either
 \index[g]{parameters used in the paper!c@$c_4, c_5$ - parameters controlling the length of accessible computations of $\mmm$ (see Lemma \ref{form})}

(a) $t\le c_4\max(||W_0||, ||W_t||)$ and $||W_j||\le c_5\max(||W_0||, ||W_t||)$, for every $j=0,\dots,t$ or

(b) $||H_0||+||H_t||\le t/500$ and
the sum of lengths of all subcomputations of $\ccc$ with step histories $(12)(2)(23)$, $(23)(2)(12)$, $(34)(4)(45)$ and $(45)(4)(34)$ is at least $0.99t$.
\end{lemma}

\begin{rk} Using the highest parameter principle (see Section \ref{param}), one can replace $500$ with a much bigger number
and replace $0.99$ with a number which is much closer to $1$. However the chosen values are sufficient for the applications of Lemma \ref{form} in this paper.
\end{rk}

 \proof One may assume that $t>c_4\max(||W_0||, ||W_t||)$,
 because otherwise Property (a) holds for sufficiently large $c_5$
 since an application of every rule can increase the length of
 a configuration by a constant depending on $\mmm$. Hence
 by Lemma \ref{fea}, \\$||H_0||+||H_t||\le 2c_3\max(||W_0||, ||W_t||)\le t/500$.

 The computation $\ccc$ is not a $B$-computation by Lemma \ref{E}
 since $c_2<c_4$. Therefore it is a computation satisfying Property (A) of Lemma \ref{SH},
 and there is a maximal subcomputation $\ccc''\colon W_r\to\dots\to W_s$ starting and ending with subcomputations with step histories (2) or (4), which are listed in part (A) of that lemma. We have $\ccc=\ccc'\ccc''\ccc'''$, where $\ccc'$ and $\ccc'''$ have Property (B).

 Lemma \ref{SH} implies that the subcomputation $\ccc''$ is a product $\ccc_1{\mathcal D}_1\dots \ccc_{k-1}{\mathcal D}_{k-1}\ccc_k$, where $k\ge 1$, every $\ccc_i$ has one of the four step
 histories from item (A) of that lemma, and every ${\mathcal D}_i$
 is a subcomputation having type 1 or 3, or 5, or just empty if
 the history $H(i)$ of $\ccc_i$ ends with $\theta(23)$ and
 $H(i+1)$ starts with $\theta(23)^{-1}$. Let $K(i)$ be the history of ${\mathcal D}_i$.


  Let ${\mathcal D}_i\colon W_x\to\dots\to W_y$. Then on the one hand, $||K_i||\le |V_x|_Y+|V_y|_Y$ by Lemma \ref{gen} (b); here
   $V_x\to\dots\to V_y$ is the restriction ${\mathcal D}_i$  to a sector
    with base of lengths two, where the rules of ${\mathcal D}_i$ insert/delete letters. On the other hand, $||H(i)||\ge 2m|V_x|_Y$,
    as it follows from Remark \ref{r-prim}
    (if $\ccc_i$ has type (2)) and from Lemmas \ref{prim} (3), \ref{M31}, \ref{M3} (a) and the definition of Set 4
    (if $\ccc_i$ has type (4)). Similarly we have $||H(i+1)||\ge 2m|V_y|_Y$, whence $$||H(i)||+||H(i+1)||)/1000 \ge m(|V_x|_Y+|V_y|)/500\ge ||K(i)||$$
    by the choice of $m$. It follows that
    $\sum ||K(i)||\le \sum ||H(i)||/500\le t/500$.

Suppose $||W_r||\le c_2||W_0||$ and $||W_s||\le c_2||W_t||$. Then by Lemma \ref{E} (a), $r\le c_2^2(c_2+1)||W_0||$ and $t-s\le c_2^2(c_2+1)||W_t||$.
So for $\ell'=r$, $\ell'''=t-s$, and large enough $c_4$
(chosen after $c_2$), we have $\max\{\ell',\ell'''\}\le c_4\max\{||W_0||, ||W_t||\}/1000$ and $\ell'+\ell'''\le t/500$. This inequality and the inequality $\sum ||K(i)||\le t/500$ imply that $\sum ||H(i)||>0.99t$, as required.

Suppose now $||W_r||> c_2||W_0||$ or $||W_s||> c_2||W_t||$. As above, it suffices to show that $l'$ and $l'''$ are small in comparison with $t$; we will show that $l'\le t/300$ ($l'''\le t/300$) if $||W_r||> c_2||W_0||$ (resp., if if $||W_s||> c_2||W_t||$).

{\bf Case 1.} The step history of $\cal C''$ starts
with (12)(2)(23). By Lemma \ref{121}, the step history of $\cal C'$ is a suffix of $(2)(21)(1)$. If $W_k\to\dots\to W_r$ is a subcomputation corresponding to Step (1),
then $||W_0||\ge ||W_k||$ by Lemma \ref{Hprim} (1) and $k\le 2m ||W_0||$ by Remark \ref{r-prim} since there are
at most $m$ cycles of the machine ${\bf LR}_m$ at Step (2). Also we have $r-k\le ||W_k||+||W_r||$ since the rules of Step (1)
just insert the same letter $\alpha^{\pm 1}$. Therefore
$\ell'\le (2m+1) (||W_k|| +||W_r||)<2||W_r||$ since $||W_r||> c_2||W_k||$.

Since $||W_r||> c_2||W_0||$, every cycle of the machine
${\bf LR}_m$ has length at least $\kappa ||W_r||$, where $\kappa^{-1}$ is the length of the standard base.
It follows that choosing $m$ large enough, we have $t\ge m\kappa ||W_r|| > 600 ||W_r||\ge 300 \ell'$, as required.

{\bf Case 2.} The step history of $\cal C''$ starts
with (32)(2)(21). By Lemma \ref{121}, the step history of $\cal C'$ is a suffix of $(4)(43)(3)$ or the suffix of $(2)(23)$.

Consider the former option for $\cal C'$. Let the subcomputation $W_0\to\dots \to W_k$ corresponds to the
step history $(4)(34)$ ($k\ge 0$). Then we have $||W_k||\le c_1||W_0|$ by Lemma \ref{M3} (c), and therefore
$||W_k||\le \frac{c_1}{c_2}||W_r||<||W_r||$. However
$W_k\to\dots\to W_r$ is a computation of Step (3), and
$W_r$ has empty historical sectors, which implies
that $||W_r||\le ||W_k||$ giving a contradiction.

If the step history of $\cal C'$ is the suffix of $(2)(23)$, then $||W_r||\le ||W_0||$ by Remark \ref{r-prim}, contrary to the assumption $||W_r||>c_2||W_0||$.

{\bf Case 3.} The step history of $\cal C''$ starts
with (34)(4)(45). By Lemma \ref{212}, the step history of $\cal C'$ is a suffix of $(2)(23)(3)$.
If $W_k\to\dots\to W_r$ is a subcomputation corresponding to Step (3) ($k\ge 0$),
then as in Case 1, $||W_0||\ge ||W_k||$ and $k\le 2m ||W_0||$.
Applying Lemma \ref{gen1} (2) to the historical sectors of Step 3, we have $r-k\le \frac 12 (||W_k||+||W_r||)\le (\frac 12 +\frac {1}{c_2})||W_r||<||W_r||$. Therefore
$\ell' \le (\frac{2m}{c_2}+1)||W_r||< 2||W_r||$.

By Lemma \ref{M3} (b), we have at least $m$ cycles of the machine
${\bf M}_3$ at Step 4, with equal lengths  $\ge\kappa ||W_r||$, where $\kappa^{-1}$ is the length of the standard base.
It follows that $t\ge m\kappa ||W_r|| > 600 ||W_r||\ge 300 \ell'$, as desired.

{\bf Case 4.} The step history of $\cal C''$ starts
with (54)(4)(43). By Lemma \ref{212} (1), the step history of $\cal C'$ is a suffix of $(4)(45)(5)$.
Let the subcomputation $W_0\to\dots\to W_k$ correspond to the
step history $(4)(45)$ ($k\ge 0$). Then, as in Case 2, we have $||W_i||\le c_1||W_0|$  for $i\le k$              and therefore
$||W_k||\le \frac{c_1}{c_2}||W_r||<||W_r||$.

By Lemma \ref{M31}, the subcomputation of $\cal C'$ with step history (4)
is subdivided by subcomputations ${\cal E}_1$, ${\cal E}_2$,... by at most $4m+1$ $\chi$-rules, where each of ${\cal E}_j$-s corresponds either to a work of ${\bf LR}$ (the work of ${\bf RL}$) or
to the work of ${\bf M}_2$. In the former case, the length
of ${\cal E}_j$ does not exceed $2c_1||W_0||$ by Lemma
\ref{Hprim} (b). In the latter case, it does not exceed
 $c_1||W_0||$  by Lemma \ref{gen1} (b) applied to the
 historical sectors. Thus, we have $k\le (8m+2)c_1||W_0||$.

 Applying Lemma \ref{gen1} (2) to the historical sectors of Step 5, we have $r-k\le \frac 12 (||W_k||+||W_r||)\le (\frac 12 +\frac {c_1}{c_2})||W_r||<||W_r||$. Therefore, on the one hand, we obtain
$\ell' \le (\frac{8m+2}{c_2}+1)||W_r||< 2||W_r||$. On the other hand, exactly as in Case 3, we have $t\ge m\kappa ||W_r||$. It follows that $t>300\ell'$ by the choice
of $m$.

Now  the proof is exhaustive by Lemmas \ref{121} and \ref{212}.
\endproof

We call a base $B$ of an eligible computation (and the computation itself)
\label{revolv} {\it revolving} if
$B\equiv x v x$ for some letter $x$ and a word $v$, and $B$  has no proper
 subword of this form.

If $v\equiv v_1zv_2$ for some letter $z$, then the word
$zv_2xv_1z$ is also revolving. One can cyclically permute the
sectors of revolving computation with base $xvx$ and obtain
a uniquely defined computation with the base $zv_2xv_1z$,
which is called a cyclic permutation of the original
computation. The history and lengths of configurations do not change
when one cyclically permutes a computation.

\begin{lemma}\label{narrow} Suppose the base $B$ of an eligible computation $\ccc\colon W_0\to\dots\to W_t$ is revolving. Then one of the following statements hold:

(1) we have inequality $||W_j||\le c_4\max(||W_0||, ||W_t||)$, for every $j=0,\dots,t$  or

(2) we have the following  properties:

(a)the word $xv$ or $v^{-1}x^{-1}$ is a cyclic permutation of the standard base of $\mmm$ and

(b) the corresponding cyclic permutations $W'_0$ and $W'_t$ of the words $W_0$ and $W_t$ are  accessible words, and

(c) 
the step history of $\ccc$ (or of the inverse computation) contains  a subword $(12)(2)(23)$ or
 $(34)(4)(45)$; moreover, the sum of lengths of corresponding subwords of the history is
 at least $0.99 t$
 and

(d) we have $||H'||+||H''||<t$ for the histories $H'$ and $H''$ of $\ccc(W_0)$ and $\ccc(W_t)$.
\end{lemma}

\proof If the computation is faulty, then Property (1) is
given by Lemma \ref{nonst} since $c_4>c_1$. If it is non-faulty, then we have
all sectors of the base in the same order as in the standard base (or its inverse), and we obtain Property (2a). Therefore we may assume now that the base $xv$ is standard and Property (1) does not hold.

If $\ccc$ is a $B$-computation, we obtain a contradiction
with Lemma \ref{E} since $c_4>c_2$. Therefore we assume further that
$\ccc$ is an $A$-computation. So it (or the inverse one) contains
a subcomputation with step history $(12)(2)(23)$ or $(34)(4)(45)$.
In case of $(34)(4)(45)$, we consider the transition $\theta(45)\colon
W_j\to W_{j+1}$. By Lemma \ref{resto}, the words in the history
sectors ${\tilde R}_{i-1}{\tilde P}_i$ are copies of each other. Therefore
they can be simultaneously erases by the rules of Set 5, and so
$W_{j+1}$ and all other configurations are accepted. Similarly
one applies Lemma \ref{resto} in case $(12)(2)(23)$ and concludes
that Property (2b) holds.

Now the second part of (2c) and (d) follow from Lemma \ref{form}.
\endproof

\subsection{Two more properties of standard computations}\label{long}

Here we prove two lemmas needed for the estimates in Subsection \ref{ub}. The first one says (due to Lemma \ref{121} (2)) that if a standard computation $\ccc$
is very long in comparison with the lengths of the first and
the last configuration, then it can be completely restored
if one knows the history of $\ccc$, and the same is true for
the long subcomputations of $\ccc$. This makes the auxiliary
parameter $\sigma_{\lambda}(\Delta)$ useful for some estimates
of areas of diagrams $\Delta$. The second lemma is also helpful
for the proof of Lemma \ref{led} in Subsection \ref{ub}.

\begin{lemma} \label{B} Let $\ccc\colon W_0\to\dots\to W_t$ be a reduced computation with standard base,
where $t\ge  c_4 \max(||W_0||, ||W_t||)$.
Suppose the word $W_0$ is accessible. Then the history of any subcomputation ${\mathcal D}\colon W_r\to\dots\to W_s$ of $\ccc$ (or the inverse for $\mathcal D$) of  length at least $0.4t$ contains a  word
of the form (a) $\chi(i-1,i)H'\chi(i,i+1)$ (i.e.,the $S$-machine
works as $\mmm_3$ at $\Theta_4$) or (b) $\zeta^{i-1,i} H'\zeta^{i,i+1}$ (i.e.,it works as $\rhh_m$
at $\Theta_2$).
\end{lemma}

\proof By Lemma \ref{form},
the sum of lengths of all subcomputations $\ccc'$ of $\ccc$ with step histories $(12)(2)(23)$, $(23)(2)(12)$, $(34)(4)(45)$ and $(45)(4)(34)$ is at least $0.99t$. Therefore ${\mathcal D}$ has to contain a subcomputation ${\mathcal D}'$ of type 2 or 4, which is a subcomputation of some ${\ccc'}$, and $||K'||\ge 0.3||H'||$ for
the histories $K'$ and $H'$ of ${\mathcal D}'$ and ${\ccc'}$, respectively.

It suffices to show that such a subcomputation ${\mathcal D}'$ of a computation $\ccc'$ with step history $(34)(4)(45)$ (with $(12)(2)(23)$) contains a subcomputation of the form (a) (form (b),
resp.) For $\ccc'$ of type $(34)(4)(45)$, this follows from
Lemma \ref{M3} (b) since $m >10$. For $\ccc'$ of type
$(12)(2)(23)$, the same property holds since
the $S$-machine $\rhh_m$ has to repeat the cycles of $\rhh$ $m$ times by Lemma \ref{prim} (3,4).
\endproof

\begin{lemma} \label{pol} Let a reduced computation $\ccc\colon W_0\to\dots\to W_t$ start with an accessible word $W_0$ and have step history of
length 1. Assume that for some index $j$, we have $|W_j|_Y>3|W_0|_Y$.
 Then there is a sector $QQ'$ such that
a state letter  from $Q$ or from $Q'$ inserts an $Y$-letter increasing
the length of this sector after any transition of the subcomputation
$W_j\to\dots\to W_t$.
\end{lemma}

\proof First of all we observe that the $Y$-words in all history sectors
(in all input sectors) of any configuration $W_i$ are copies of each other,
because $W_0$ is accessible. Also the statement is trivial if $t=1$, and so $j=1$ too. Then inducting on $t$,  one can assume that $|W_1|_Y>|W_0|_Y$
since otherwise it suffices to consider the  computation $W_1\to\dots\to W_t$ of length $t-1$.

If we have one of the Sets 1, 3, 5, then inequality $|W_0|_Y<|W_1|_Y$
implies $|W_1|_Y<|W_2|_Y<\dots$
since the second rule cannot be inverse for the first one, and so on,
i.e., we obtain the desired property of any input sector for Set 1
or of any history sector for Sets 3 or 5.

If we have Set 2, then the statement for  any imput sector follows from Lemma \ref{prim} (1)
.

Let  the step history be $(4)$. Recall that the rules of Set 4 are
subdivided in several sets, where each set copies the work
of either $\rhh$ or $\mmm_3$. If a $\bf LR$-rule
of the subcomputation ${\mathcal D}\colon W_0\to\dots\to W_j$
increases the length of a history sector, then we refer to Lemma \ref{prim} (1) as above.  So one may
assume that no $\bf LR$-rules of ${\mathcal D}$ increase the length of history sectors.

Assume now that $\mathcal D$ has an $\mmm_3$-rule
increasing the length of history sectors.
It has to insert a letter from
$X_{i,\ell}$ from the left and a letter from $X_{i,r}$ from the right.
Since the obtained word is not a word over one of these alphabets,
the work of $\mmm_3$ is not over, and the next rule has to increase the length
of the sector again in the same manner since the computation is reduced.
This procedure will repeat until one gets $W_t$. This proves the statement for any
history sector.

It remains to assume that there are no transitions in  $\mathcal D$ increasing the lengths of history sectors and the first
transition $W_0\to W_1$ is provided by a rule $\theta$
of $\mmm_3$. It cannot shorten history sectors (by 2).
Indeed $\theta$ can change the length of neighbor working sectors at most by $1$ (see Lemma \ref{simp} (**)),
which implies $|W_0|_Y\ge |W_1|_Y$, a contradiction. It follows
that no further rules of $\mmm_3$ can shorten history sectors.
Then Lemma \ref{w} implies that all history
sectors in all configurations of $\mathcal D$ have equal lengths.

By Lemma \ref{gen1}
(b) the lengths of the history of the maximal subcomputation ${\mathcal E}\colon W_0\to\dots\to W_s$ of $\mmm_3$ in $\mathcal D$ does not exceed $h$, where
$h$ is the $Y$-length of all history sectors of the configurations from $\mathcal D$.

Every rule of the subcomputation $\mathcal E$ can change the length of any working sector at most by $1$. (See Lemma \ref{simp} (**)). Hence if its
length in $W_0$ is $\ell$, its length in $W_s$ is at most $\ell+h$. It follows that $|W_{s}|_Y\le 3 |W_0|_Y$, because
the working sectors of $\mmm_2$ and its history sectors alternate in the standard base; and the same inequality $|W_r|_Y\le 3|W_0|_Y$ holds for any configuration $W_r$ of $\mathcal E$. Hence $s\ne j$ and the subcomputation $\mathcal E$ is followed in $\mathcal D$ by a subcomputation $\mathcal F$ of $\bf LR$, which does not change the
length of configurations by Lemma \ref{Hprim}.

So $\mathcal F$ has to be followed in $\mathcal D$ by a maximal subcomputation $\mathcal G$
of ${\bf M}_3$ again. Since we have the canonical work of ${\mathcal M}_3$ in history sectors, a prefix of the history of
${\mathcal G}^{-1}$ is a copy of the entire $H({\mathcal E})^{-1}$, where $H({\mathcal E})$ is the history of $\mathcal E$. ($\mathcal G$ cannot be shorter than $\mathcal E$ since otherwise the configuration $W_j$ would
have a copy in $\mathcal E$, whence $|W_j|_Y\le 3|W_0|_Y$, a contradiction.) It follows that a configuration $W_l$ of $\mathcal G$
is a copy of $W_0$, and so $|W_l|_Y=|W_0|_Y$. Since the subcomputation
$W_l\to\dots\to W_j\to\dots\to W_t$ is shorter than $\ccc$, we complete the proof of the lemma inducting on $t$.
\endproof

\section{Groups and diagrams}\label{gd}

\subsection{The groups}\label{MG}

Every $S$-machine can be simulated by a finitely presented group (see \cite{SBR}, \cite{OS04}, \cite{OS06},  etc.). Here we apply a modified construction from \cite{SBR} to the $S$-machine $\mmm.$
To simplify formulas, it is convenient to change the notation. From now on we shall denote by \index[g]{parameters used in the paper!n@$N$ - the length of the standard base of the $S$-machine $\mmm$}$N$ the length of
the standard base of $\mmm$.

Thus the set of state letters is $Q=\sqcup_{i=0}^{N-1}Q_i$ (we set $Q_{N}=Q_0=\{\tt\}$), $Y= \sqcup_{i=1}^{N} Y_i,$ and $\Theta$ is the set of rules of the $S$-machine $\mmm$.

The finite set of generators of the group $M$ \index[g]{M@group $M$}\index[g]{M@group $M$!generators of the group $M$}  consists of {\em $q$-letters},
{\em $Y$-letters} and
{\em $\theta$-letters} defined as follows.

For every letter $q\in Q$ the set of generators of $M$ contains
$L$ copies $q^{(i)}$ of it,  $i=1,\dots, L$, if the letter $q$ occurs in the rules of $\mtt_1$ or $\mtt_2$. \index[g]{parameters used in the paper!l@$L$ - the number of generators $q^{(i)}$ of the group $G$ for each state letter $q$ of $\mmm$, the order of $W_{ac}$ in $G$}(The  number $L$ is one of the parameters from Section \ref{param}.)  Otherwise only the letter $q$ is included
in the generating set of $M$.

For every letter $a\in Y$ the set of generators of $M$ contains $a$
and $L$ copies $a^{(i)}$ of it $i=1,..., L$.

For every $\theta\in \Theta^+$ we have $N$ generators $\theta_0,\dots,\theta_N$ in $M$ (here $\theta_{N}\equiv\theta_0$)
if $\theta$ is a rule of $\Theta_3$ (excluding $\theta(23)$) or $\Theta_4$, or $\Theta_5$. For $\theta$ from $\Theta_1$ or $\Theta_2$ (including $\theta(23)$),
we introduce $LN$ generators $\theta_j^{(i)}$, where $j=0,\dots, N$,
$i=1,\dots, L$ and $\theta_N^{(i)}=\theta_0^{(i+1)}$ (the superscripts are taken modulo $L$).

The \index[g]{M@group $M$!relations of the group $M$} relations of the group $M$ correspond to the rules of the $S$-machine $\mmm$ as follows.
For every rule $\theta=[U_0\to V_0,\dots U_{N}\to V_{N}]\in \Theta^+$ of sets $\Theta_1$ or $\Theta_2$, we have

\begin{equation}\label{rel1}
U_j^{(i)}\theta_{j+1}^{(i)}=\theta_j^{(i)} V_j^{(i)},\,\,\,\, \qquad \theta_j^{(i)} a^{(i)}=a^{(i)}\theta_j^{(i)}, \,\,\,\, j=0,...,N,\;\;
i=1,\dots L,
\end{equation}
for all $a\in Y_j(\theta)$, where $U_j^{(i)}$ and $V_j^{(i)}$ are obtained from $U_j$ and $V_j$
by addiing the superscript ${i}$ to every letter.

For $\theta=\theta(23)$, we introduce relations

\begin{equation}\label{rel11}
U_j^{(i)}\theta_{j+1}^{(i)}=\theta_j^{(i)} V_j,\,\,\,\, \qquad a^{(i)}\theta_j^{(i)}=\theta_j^{(i)} a,
\end{equation}
i.e.,the superscripts are erased in the words $U_j^{(i)}$ and in the $Y$-letters after an application of (\ref{rel11}).

For every rule $\theta=[U_0\to V_0,\dots U_{N}\to V_{N}]\in \Theta^+$ from $\Theta_3$ or $\Theta_4$, or $\Theta_5$, we define

\begin{equation}\label{rel111}
U_j\theta_{j+1}=\theta_j V_j,\,\,\,\, \qquad a\theta_j=\theta_j a
\end{equation}

The first type of relations (\ref{rel1} - \ref{rel111}) will be
called $(\theta,q)$-{\em relations}, the second type - \label{thetaar}
$(\theta,a)$-{\em relations}.

Finally, the required \label{groupG} group $G$ is given by the generators and
relations of the group $M$ and by two more additional
relations, namely the \label{hubr} {\it hub}-relations

\begin{equation}\label{rel3}
W_{st}^{(1)}\dots W_{st}^{(L)}=1\;\; and \;\; (W_{ac})^L=1,
\end{equation}
where the word $W_{st}^{(i)}$ is a copy with superscript $(i)$ of the start word $W_{st}$ (of length $N$) of the $S$-machine $\mmm$
and $W_{ac}$
is  the accept word of $\mmm$.

\begin{rk} The main difference of the construction of $M$ and the groups based  on $S$-machines with hubs from our previous papers \cite{SBR, OS04, OS06, O18} and others, is that relations (\ref{rel1}) are defined differently for different rules of the $S$-machine. We also use two hub relations instead of just one, although
it is easy to see that one hub relation follows from the other (and other relations).

Note also that, as usual, $M$ is a multiple HNN extension  of the free group generated by all $Y$- and $q$-letters, because by Tietze transformations using $(\theta,q)$-relations, all $\theta$-letters, except for one for every rule
$\theta$, can be eliminated.
\end{rk}

\subsection{Van Kampen diagrams}\label{md}

Recall that a \index[g]{van Kampen diagram} van Kampen {\it diagram} $\Delta $ over a presentation
$P=\langle A | \rrr\rangle$ (or just over the group $P$)
is a finite oriented connected and simply--connected planar 2--complex endowed with a
\index[g]{van Kampen diagram!labeling function} {\em labeling function} $\Lab \colon E(\Delta )\to A^{\pm 1}$, where $E(\Delta
) $ denotes the set of oriented edges of $\Delta $, such that $\Lab
(e^{-1})\equiv \Lab (e)^{-1}$. Given a \index[g]{van Kampen diagram!cell} {\em cell} (that is a 2-cell) $\Pi $ of $\Delta $,
we denote by $\partial \Pi$ the boundary of $\Pi $; similarly, \index[g]{van Kampen diagram!boundary $\partial(\Delta)$}
$\partial \Delta $ denotes the boundary of $\Delta $. The labels of
$\partial \Pi $ and $\partial \Delta $ are defined up to cyclic
permutations. An additional requirement is that the label of any
cell $\Pi $ of $\Delta $ is equal to (a cyclic permutation of) a
word $R^{\pm 1}$, where $R\in \rrr$. The label and the \index[g]{combinatorial length of a path} combinatorial length $||\bf p||$ of
a path $\bf p$ are defined as for Cayley graphs.

The van Kampen Lemma \cite{LS,book, Sbook} states that a word $W$ over the alphabet $A^{\pm 1}$
represents the identity in the group $P$ if and only
if there exists a diagram $\Delta
$ over $P$ such that
$\Lab (\partial \Delta )\equiv W,$ in particular, the combinatorial perimeter $||\partial\Delta||$ of $\Delta$ equals $||W||.$
(\cite{LS}, Ch. 5, Theorem 1.1; our formulation is closer to
Lemma 11.1 of \cite{book}, see also \cite[Section 5.1]{Sbook}). The word $W$ representing $1$ in $P$ is freely equal
to a product of conjugates to the words from $R^{\pm 1}$. The minimal number
of factors in such products is called the \index[g]{area of a word} {\em area} of the word $W.$ The \index[g]{van Kampen diagram!area} {\it area}
of a diagram $\Delta$ is the number of cells in it.
The proof of the van Kampen Lemma \cite{book, Sbook} shows that  $\area(W)$ is equal
to the area of a van Kampen diagram having
the smallest number of cells among all van Kampen diagrams with  boundary label $\Lab (\partial \Delta )\equiv W.$

We will study diagrams over the group presentations of $M$ and $G$. The edges labeled by state
letters ( = $q$-{\it letters}) will be called
\label{qedge} $q$-{\it edges}, the edges labeled by tape
letters (= $Y$-{\it letters}) will be called \label{Yedge} $Y$-{\it edges}, and the edges labeled by
$\theta$-letters are \label{thedge} $\theta$-{\it edges}.

We denote by $|{\bf p}|_Y$ (by $|{\bf p}|_{\theta}$, by
$|{\bf p}|_q$)
the \label{alength} $Y$-{\it length} (resp., the \label{thlength} $\theta$-{\it length}, the \label{qlength} $q$-length) of a path/word $\bf p,$ i.e., the number of
$Y$-edges/letters (the number of $\theta$-edges/letters, the number of $q$-edges/letters) in $\bf p.$

The cells corresponding
to relations (\ref{rel3})  are called \index[g]{hub}{\it hubs}, the cells corresponding
to $(\theta,q)$-relations are called \index[g]{t@$(\theta,q)$-cell}$(\theta,q)$-{\it cells},
and the cells are called \index[g]{t@$(\theta,a)$-cell}$(\theta,a)$-{\it cells} if they correspond to $(\theta,a)$-relations.

A Van Kampen diagram is \index[g]{van Kampen diagram!reduced}{\em reduced}, if
it does not contain two cells (= closed $2$-cells) that have a
common edge $e$ such that the boundary labels of these two cells are
equal if one reads them starting with $e$
(if such pairs of cells exist, they can be removed to obtain a  diagram of smaller area and with the same boundary label).

\subsubsection{The superscript shift of a van Kampen diagram over $M$ or $G$}

\begin{rk}\label{ss} If one changes all superscripts of the generators of  $M$ or $G$ by adding the same integer $k$: $(i)\to (i+k)$
(modulo $L$) in all letters having a superscript, then  one obtains the relations again, as it is clear from
formulas (\ref{rel1} - \ref{rel3}). Therefore similar change
$\Delta\to\Delta^{(+k)}$ of the edge labels transforms a (reduced)
diagram $\Delta$ to a (reduced) diagram $\Delta^{(+k)}$. Let us
call such a transformation \index[g]{superscript shift or $k$-shift}{\it superscript shift} (or $k$-{\it shift}) of $\Delta$.
\end{rk}

\subsubsection{Bands}
To study (van Kampen) diagrams
over the group $G$ we shall use their simpler subdiagrams such as bands and trapezia, as in \cite{O97},  \cite{SBR}, \cite{BORS}, etc.
 Here we repeat one more necessary definition.

\index[g]{band}\begin{df}Let $\mathcal Z$ be a subset of the set of letters in the set of generators of the group $M$. A
$\mathcal Z$-band $\bb$ is a sequence of cells $\pi_1,...,\pi_n$ in a reduced \vk
diagram $\Delta$ such that

\begin{itemize}
\item Every two consecutive cells $\pi_i$ and $\pi_{i+1}$ in this
sequence have a common boundary edge ${\mathbf e}_i$ labeled by a letter from ${\mathcal Z}^{\pm 1}$.
\item Each cell $\pi_i$, $i=1,...,n$ has exactly two $\mathcal Z$-edges in the boundary $\partial \pi_i$,
${\mathbf e}_{i-1}^{-1}$ and ${\mathbf e}_i$ (i.e.,edges labeled by a letter from ${\mathcal Z}^{\pm 1}$) with the requirement that either
both $\Lab(e_{i-1})$ and $\Lab(e_i)$ are positive letters or both
are negative ones.

\item If $n=0$, then $\bb$ is just a $\mathcal Z$-edge.
\end{itemize}
\end{df}

The counter-clockwise boundary of the subdiagram formed by the
cells $\pi_1,...,\pi_n$ of $\bb$ has the factorization ${\mathbf e}\iv {\mathbf q}_1{\mathbf f} {\mathbf q}_2\iv$
where ${\mathbf e}={\mathbf e}_0$ is a $\mathcal Z$-edge of $\pi_1$ and ${\mathbf f}={\mathbf e}_n$ is an $\mathcal Z$-edge of
$\pi_n$. We call ${\mathbf q}_1$ the {\em bottom} of $\bb$ and ${\mathbf q}_2$ the
{\em top} of $\bb$, denoted $\bott(\bb)$ and $\topp(\bb)$.

Top/bottom paths and their inverses are also called the {\em
sides} of the band. \index[g]{band!sides} The $\mathcal Z$-edges ${\mathbf e}$
and ${\mathbf f}$ are called the {\em start} \index[g]{band!start and end edges} and {\em end} edges of the
band. If $n\ge 1$ but ${\mathbf e}={\mathbf f},$ then the $\mathcal Z$-band is called a  $\mathcal Z$-{\it annulus} \index[g]{band!annulus}.

If $\bb$ is a ${\mathcal Z}$-band with ${\mathcal Z}$-edges $e_1,...,e_n$ (in that order), then we can form a broken line connecting midpoints of the consecutive edges $e_1, ..., e_n$ and laying inside the union of the cells from $\bb$ which will be called the {\em median} of $\bb$. \index[g]{band!median}

We will consider \index[g]{band!q@$q$-band} $q$-{\it bands}, where $\mathcal Z$ is one of the sets $Q_i$ of state letters
for the $S$-machine $\mmm$, \index[g]{band!t@$\theta$-band}
$\theta$-{\it bands} for every $\theta\in\Theta$, and $Y$-{\it bands}\index[g]{band!a@$Y$-band}, where
${\mathcal Z}=\{a,a^{(1)},\dots, a^{(L)}\}\subseteq Y$.
 The convention is that $Y$-bands do not
contain $(\theta,q)$-cells, and so they consist of $(\theta,a)$-cells  only.

\begin{lemma}\label{tba} Let ${\mathbf e}\iv {\mathbf q}_1{\mathbf f} {\mathbf q}_2\iv$ be the boundary of a $\theta$-band $\mathcal B$ with bottom ${\mathbf q}_1$
and top ${\mathbf q}_2$ in a reduced diagram.

(1) If the start and the end edges ${\mathbf e}$
and ${\mathbf f}$ have different labels, then $\mathcal B$ has $(\theta,q)$-cells.

(2) For every $(\theta,q)$-cell $\pi_i$ of $\mathcal B$, one of its
boundary $q$-edges belongs in ${\mathbf q}_1$ and another one
belongs in ${\mathbf q}_2$.
\end{lemma}

\proof (1) If every cell $\pi_i$ of $\mathcal B$ is a $(\theta,a)$-cell, then
both $\theta$-edges of the boundary $\partial\pi_i$ have equal labels, as it follows from the definition of $(\theta,a)$-relations.
Then the definition of band implies that $\Lab({\mathbf e})=\Lab({\mathbf f})$, a contradiction.

(2) Proving by contradiction, we have that that $\pi_i$ and $\pi_j$ ($i\ne j$) share a boundary $q$-edge $\bf g$. We may assume that
the difference $j-i>0$ is minimal, and so the  subband
formed by $\pi_{i+1},\dots,\pi_{j-1}$ has no $(\theta,q)$-cells.
It folows from (1) that $\pi_i$ and $\pi_j$ have the same boundary
labels if one read then starting with $\Lab(\bf g)$, contrary to the
assumption that the diagram is reduced.

\endproof

\begin{rk} \label{tb} To construct the top (or bottom) path of a band $\mathcal B$, at the beginning
one can just form a product ${\mathbf x}_1\dots {\mathbf x}_n$ of the top paths ${\mathbf x}_i$-s of the cells $\pi_1,\dots,\pi_n$ (where each $\pi_i$ is a $\mathcal Z$-bands of length $1$).
No  $\theta$-letter is being canceled in the word
$W\equiv\Lab({\mathbf x}_1)\dots\Lab({\mathbf x}_n)$ if $\mathcal B$ is  a $q$- or $Y$-band since
otherwise two neighbor cells of the band would make the diagram non-reduced. By Lemma \ref{tba} (2), there are no cancellations of
$q$-letters of
$W$ if $\mathcal B$ is  a $\theta$-band.

If $\mathcal B$ is a $\theta$-band then no cancellations of $q$-letters
are possible in $W$ by Lemma \ref{tba} (2).
The {\it trimmed} top/bottom label of $\mathcal B$
are the maximal subwords of the top/bottom labels starting and ending
with $q$-letters. \index[g]{band!top path $\topp(\bb)$} \index[g]{band!bottom path $\bott(\bb)$}

However a few cancellations of $Y$-letters  are possible in $W.$ (This can happen if one of $\pi_i, \pi_{i+1}$
is a $(\theta,q)$-cell and another one is a $(\theta,a)$-cell.) We will always assume
that the top/bottom label of a $\theta$-band is a reduced form of the word $W$.
This property  is easy to achieve: by folding edges
with the same labels having the same initial vertex, one can make
the boundary label of a subdiagram in a \vk diagram reduced (e.g., see \cite{book} or
\cite{SBR}).
\end{rk}

We shall call a $\mathcal Z$-band \index[g]{band!maximal band}{\em maximal} if it is not contained in
any other $\mathcal Z$-band.
Counting the number of maximal $\mathcal Z$-bands
 in a diagram we will not distinguish the bands with boundaries
 ${\mathbf e}\iv {\mathbf q}_1{\mathbf f} {\mathbf q}_2\iv$ and ${\mathbf f} {\mathbf q}_2\iv {\mathbf e}\iv {\mathbf q}_1,$ and
 so every $\mathcal Z$-edge belongs to a unique maximal $\mathcal Z$-band.

We say that a ${\mathcal Z}_1$-band and a ${\mathcal Z}_2$-band \index[g]{band!crossing bands}{\em cross} if
they have a common cell and ${\mathcal Z}_1\cap {\mathcal Z}_2=\emptyset.$

Sometimes we specify the types of bands as follows.
A $q$-band corresponding to one
letter $Q$ of the base is called a \index[g]{Q@$Q$-band}$Q$-band. For example, we will consider \index[g]{band!t@$\tt$-band}$\tt$-{\it band} corresponding to the part $\{\tt\}$.

Our previous papers (see \cite{SBR}, \cite{BORS}, etc.) contain the proof of the next lemma in a more general setting.  The difference
caused by different simulation of the $S$-machine $\mmm$ by defining
relations of $M$ does not affect the validity of the proof since the proof uses the properties
mentioned in Lemma \ref{tba} and Remark \ref{tb}. To convince
the reader, below we recall the proof of one of the following claims.

\begin{lemma}\label{NoAnnul}
A reduced van Kampen diagram $\Delta$ over $M$ has no
$q$-annuli, no $\theta$-annuli, and no  $Y$-annuli.
Every $\theta$-band of $\Delta$ shares at most one cell with any
$q$-band and with any $Y$-band.
\end{lemma}

\proof We will prove only the property that a $\theta$-band
$\mathcal T$ and a $q$-band $\mathcal Q$ cannot cross each other two
times. Taking a minimal counter-example, one assumes that
these bands have exactly two common cells $\pi$ and $\pi'$, and $\Delta$ has no cells outside the region bounded by $\mathcal T$ and $\mathcal Q$. Then $\mathcal Q$ has exactly two cells since otherwise a
maximal $\theta$-band starting with a cell $\pi''$ of $\mathcal Q$,
where $\pi''\notin \{\pi, \pi'\}$, has to end on $\mathcal Q$,
bounding with a part of $\mathcal T$ a smaller counter-example.

\begin{figure}[ht]
\unitlength=1.00mm
\special{em:linewidth 0.4pt}
\linethickness{0.4pt}
\begin{picture}(113.00,45.00)
\put(98.67,38.67){\line(0,1){6.33}}
\put(98.67,45.00){\line(1,0){14.33}}
\put(113.00,45.00){\line(0,-1){16.67}}
\put(113.00,28.33){\line(-3,-4){3.67}}
\put(109.33,23.33){\line(-1,0){7.33}}
\put(101.67,23.00){\line(-3,-4){3.33}}
\put(98.33,18.67){\vector(0,-1){9.67}}
\put(98.33,9.00){\line(0,-1){7.67}}
\put(98.33,1.33){\line(1,0){14.67}}
\put(113.00,1.33){\line(0,1){17.33}}
\put(113.00,18.67){\line(-4,5){3.67}}
\put(105.67,26.67){\makebox(0,0)[cc]{$q$}}
\put(100.00,33.00){\makebox(0,0)[cc]{$\theta$}}
\put(100.67,14.33){\makebox(0,0)[cc]{$\theta'$}}
\bezier{800}(98.33,38.67)(-0.33,26.00)(98.67,9.00)
\bezier{568}(98.33,28.33)(27.67,26.00)(98.33,19.00)
\put(70.00,22.33){\vector(-3,-1){11.33}}
\put(82.67,20.33){\vector(-2,-3){4.67}}
\put(98.50,28.17){\vector(0,1){10.50}}
\put(90.83,28.00){\vector(-1,4){2.33}}
\put(82.00,27.67){\vector(-1,2){4.00}}
\put(73.33,27.00){\vector(-2,3){4.50}}
\put(64.50,25.83){\vector(-4,1){11.83}}
\put(101.67,23.33){\line(-2,3){3.33}}
\end{picture}
\
\begin{center}
\caption{A $Q$-band intersects a $\theta$-band twice.}
\end{center}

\end{figure}
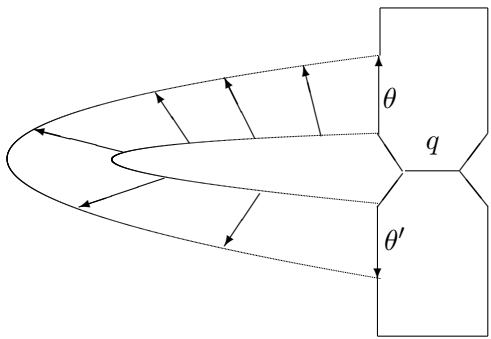

Thus, the boundaries of $\pi$ and $\pi'$ share a $q$-edge.

For the similar  reason, $\mathcal T$ has no $(\theta,q)$-cells
except for $\pi$ and $\pi'$, and by Lemma \ref{tba} (1),
these cells have the same pairs of $\theta$-edges in the boundaries.
This makes the diagram non-reduced, a contradiction.
\endproof

If $W\equiv x_1...x_n$ is a word in an alphabet $X$, $X'$ is another
alphabet, and $\phi\colon X\to X'\cup\{1\}$ (where $1$ is the empty
word) is a map, then $\phi(W)\equiv\phi(x_1)...\phi(x_n)$ is called the
\label{projectw}{\em projection} of $W$ onto $X'$. We shall consider the
projections of words in the generators of $M$ onto
$\Theta$ (all $\theta$-letters map to the
corresponding element of $\Theta$,
all other letters map to $1$), and the projection onto the
alphabet $\{Q_0\sqcup \dots \sqcup Q_{N-1}\}$ (every
$q$-letter maps to the corresponding $Q_i$, all other
letters map to $1$).

\begin{df}\label{dfsides}
{\rm  The projection of the label
of a side of a $q$-band onto the alphabet $\Theta$ is
called the {\em history} of the band. \index[g]{band!history of a $q$-band} The step history of this projection
is the \index[g]{band!step history of a $q$-band}{\it step history} of the $q$-band. The projection of the label
of a side of a $\theta$-band onto the alphabet $\{Q_0,...,Q_{N-1}\}$
is called the {\em base} of the band, i.e., the base of a $\theta$-band
is equal to the base of the label of its top or bottom}\index[g]{band!base of a $\theta$-band}
\end{df}
As in the case of words, we will  use representatives of
$Q_j$-s in base words.

If $W$ is a word in the generators of $M$, then by $W^{\emptyset}$
we denote the projection of this word onto the alphabet of the
$S$-machine $\mmm$, we obtain this projection after deleting all
superscripts in the letters of $W$. In particular, $W^{\emptyset}\equiv W$, if there are no superscripts in the letters of $W$.

We call a word $W$ in $q$-generators and $Y$-generators {\it permissible} \index[g]{permissible word} if the
word $W^{\emptyset}$ is admissible, and the letters of any 2-letter subword of $W$
have equal superscripts (if any), except for the subwords $(q\tt)^{\pm 1}$,
where the letter  $q$ has some superscript $(i)$ and $q^{\emptyset} \in Q_{N-1}$; in this case the superscript of the letter $\tt$ must be $(i+1)$
(modulo $L$).

\begin{rk} \label{perad} It follows from the definition that if $V$ is $\theta$-admissible for a rule $\theta$
of $\{\theta(23)\iv \}\cup  \mtt_3\cup \{\theta(34)\}\cup \mtt_4\cup\{\theta(45)\}\cup \mtt_5$, then there is exactly one permissible
word $W$ such that $W^{\emptyset}\equiv V$, namely, $W\equiv V$.
If $\theta$ is a rule of $\mtt_1\cup\{\theta(12)\}\cup \mtt_2\cup\{\theta(23)\}$, then the permissible word $W$
with property $W^{\emptyset}\equiv V$ exists and it is uniquely defined if one choose arbitrary
superscript for the first letter (or for any particular letter) of $W$.
\end{rk}

\begin{lemma} \label{perm} (1) The trimmed bottom and top labels $W_1$ and $W_2$ of any reduced $\theta$-band $\mathcal T$ containing at least one $(\theta,q)-cell$ are permissible
and $W_2^{\emptyset}\equiv W_1^{\emptyset}\cdot\theta$.

(2) If $W$ is a $\theta$-admissible word, then for a
permissible word $W_1$ such that $W_1^{\emptyset}\equiv W$ (given by Remark \ref{perad})
one can construct a reduced $\theta$-band with the trimmed bottom
label $W_1$ and the trimmed top label $W_2$, where $W_2^{\emptyset}\equiv
W_1^{\emptyset}\cdot\theta$.

\end{lemma}

\proof (1) It follows from Lemma \ref{tba} (2) that $W_1\equiv q_1^{\pm 1}u_1q_2^{\pm 1}\dots u_k q_{k+1}^{\pm 1}$, where $q_j^{\pm 1}$ and $q_{j+1}^{\pm 1}$ are
the labels of $q$-edges of some cells $\pi(j)$ and $\pi(j+1)$
such that the subband connecting these cells has no $(\theta, q)$-cells. Therefore by Lemma \ref{tba} (1), all the $\theta$-edges
between $\pi(j)$ and $\pi(j+1)$ have the same labels. It follows
from the list of $(\theta,a)$-relations that all $Y$-letters
of the word $u_j$ have to belong to the same subalphabet.
In particular, if we have the subword $q_ju_jq_{j+1}$, then the projection of this subword is a subword of $W_1^{\emptyset}$
satisfying the first condition from the definition of admissible
word. Similarly one obtains other conditions if $q_j$ or/and $q_{j+1}$ occur in $W_1$ with exponent $-1$.
Hence the word $W_1^{\emptyset}$ (and $W_2^{\emptyset}$) are
admissible, and the words $W_1, W_2$ are permissible since again
the condition on $2$-letter subwords follows from Lemma \ref{tba} and the relations (\ref{rel1} - \ref{rel111}).

If ${\mathbf x}={\mathbf x}_1\dots {\mathbf x}_n$ ( ${\mathbf y}={\mathbf y}_1\dots {\mathbf y}_n$ ) is the product of the top paths ${\mathbf x}_i$-s (bottom  paths ${\mathbf y}_i$-s) of the all cells $\pi_1,\dots,\pi_n$ of $\mathcal T$, as in Remark \ref{tb}, then the transition from the trimmed label of $\bf x$ to the trimmed label of $\bf y$ with erased superscripts, is the application
of $\theta$, as it follows from relations (\ref{rel1} - \ref{rel111}). Since by definition, the application of $\theta$ automatically implies possible cancellations, we have $W_2^{\emptyset}\equiv W_1^{\emptyset}\cdot\theta$ for the reduced words $W_1$ and $W_2$, as required.

Since $W$ is $\theta$-admissible, there is an equality $W'\equiv W\cdot \theta$. Therefore we can simulate the application of $\theta$
to every letter of $W$ as follows. We draw a path $\bf p={\mathbf e}_1\dots {\mathbf e}_n$ labeled by $W_1$ and attach a cell $\pi_i$  corresponding to one of the defining relations of $M$ to every edge ${\mathbf e}_i$ of $\bf p$ from the left. Since the word $W_1$ is permissible, the $\theta$-edges started with the common vertex
of $\pi_i$ and $\pi_{i+1}$ must have equal labels, and so these
two edges can be identified. Finally, we obtain a required $\theta$-band. It is reduced diagram since the permissible word $W_1$
is reduced.
\endproof

\subsubsection{Trapezia}

\begin{df}\label{dftrap}
{\rm Let $\Delta$ be a reduced  diagram over $M$,
which has  boundary path of the form ${\mathbf p}_1\iv {\mathbf q}_1{\mathbf p}_2{\mathbf q}_2\iv,$ where
${\mathbf p}_1$ and ${\mathbf p}_2$ are sides of $q$-bands, and
${\mathbf q}_1$, ${\mathbf q}_2$ are maximal parts of the sides of
$\theta$-bands such that $\Lab({\mathbf q}_1)$, $\Lab({\mathbf q}_2)$ start and end
with $q$-letters.

\begin{figure}[ht]
\begin{center}
\includegraphics[width=1.0\textwidth]{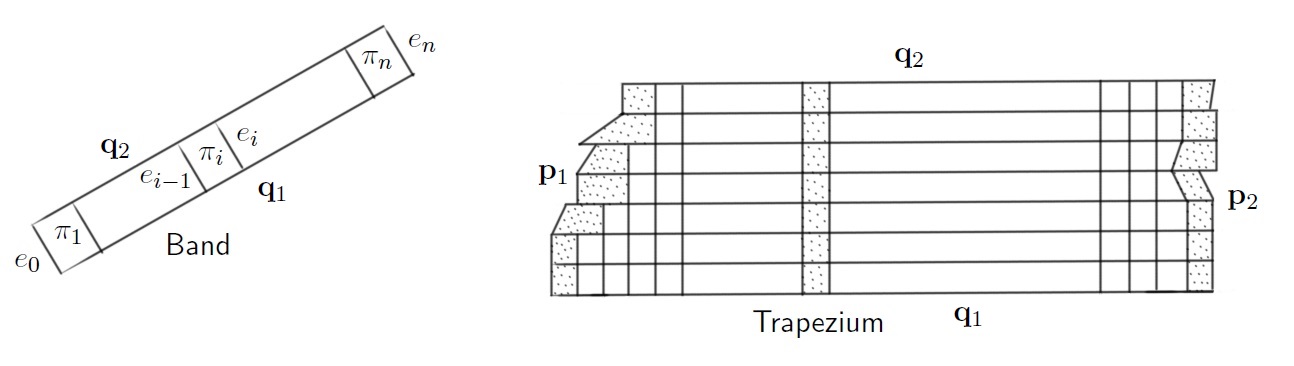}
\end{center}
\caption{Band and Trapezium}\label{bt}
\end{figure}

Then $\Delta$ is called a \index[g]{trapezium}{\em trapezium}. The path ${\mathbf q}_1$ is
called the \index[g]{trapezium!bottom}{\em bottom}, the path ${\mathbf q}_2$ is called the \index[g]{trapezium!top}{\em top} of
the trapezium, the paths ${\mathbf p}_1$ and ${\mathbf p}_2$ are called the \index[g]{trapezium!left and right sides}{\em left
and right sides} of the trapezium. The history \index[g]{trapezium!step history} \index[g]{trapezium!history}(step history) of the $q$-band
whose side is ${\mathbf p}_2$ is called the {\em history} (resp., step history) of the trapezium;
the length of the history is called the \index[g]{trapezium!height}{\em height}  of the
trapezium. The base of $\Lab ({\mathbf q}_1)$ is called the \index[g]{trapezium!base}{\em base} of the
trapezium.}
\end{df}

\begin{rk} Notice that the top (bottom) side of a
$\theta$-band $\ttt$ does not necessarily coincides with the top
(bottom) side ${\mathbf q}_2$ (side ${\mathbf q}_1$) of the corresponding trapezium of height $1$, and ${\mathbf q}_2$
(${\mathbf q}_1$) is
obtained from $\topp(\ttt)$ (resp. $\bott(\ttt)$) by trimming the
first and the last $Y$-edges
if these paths start and/or end with $Y$-edges.
We shall denote the
\index[g]{band!t@$\theta$-band!trimmed}{\it trimmed} \index[g]{band!t@$\theta$-band!trimmed bottom and top paths: $\tbott$, $\ttopp$} top and bottom sides of $\ttt$ by $\ttopp(\ttt)$ and
$\tbott(\ttt)$. By definition, for arbitrary $\theta$-band $\mathcal T,$ $\ttopp(\mathcal T)$
is obtained by such a trimming only if $\mathcal T$ starts and/or ends with a
$(\theta,q)$-cell; otherwise $\ttopp(\mathcal T)=\topp(\mathcal T).$
The definition of $\tbott(\mathcal T)$ is similar.
\end{rk}

By Lemma \ref{NoAnnul}, any trapezium $\Delta$ of height $h\ge 1$
can be decomposed into $\theta$-bands $\ttt_1,...,\ttt_h$ connecting
the left and the right sides of the trapezium.


\begin{lemma}\label{simul} (1) Let $\Delta$ be a trapezium
with history $H\equiv\theta(1)\dots\theta(d)$ ($d\ge 1$).
Assume that $\Delta$
has consecutive maximal $\theta$-bands  ${\mathcal T}_1,\dots
{\mathcal T}_d$, and the words
$U_j$
and $V_j$
are the  trimmed bottom and the
trimmed top labels of ${\mathcal T}_j,$ ($j=1,\dots,d$).
Then $H$ is an eligible word, $U_j$, $V_j$ are
permissible
words,
$$V_1^{\emptyset}\equiv U_1^{\emptyset}\cdot \theta(1),\;\; U_2\equiv V_1,\;\;\dots,\;\; U_d \equiv V_{d-1},\;\; V_d^{\emptyset}\equiv U_d^{\emptyset}\cdot \theta(d)$$

Furthemore, if the first and the last $q$-letters of the word $U_j$ or of the word $ V_j$ have some superscripts $(i)$ and $(i')$, then the difference $i'-i$ (modulo $L$) does not depend on on the choice of $U_j$ or $V_j$.

(2) For every eligible computation $U\to\dots\to U\cdot H \equiv V$ of $\mmm$
with $||H||=d\ge 1$
there exists
a trapezium $\Delta$ with bottom label $U_1$ (given by Remark \ref{perad}) such that $U_1^{\emptyset}\equiv U$, top label $V_d$ such that $V_d^{\emptyset}\equiv V$, and with history $H.$

\end{lemma}

\proof (1) The trimmed top side of one of the bands $\ttt_j$ is the same as trimmed bottom side of $\ttt_{j+1}$ ($j=1,\dots, d-1$), and the
equalities $U_2\equiv V_1,\dots, U_d\equiv V_{d-1}$ follow.
The equalities $V_j^{\emptyset}\equiv U_j^{\emptyset}\cdot\theta(j)$ ($j=1,\dots d)$ are given
by Lemma \ref{perm} (1). By the same lemma the words $U_j$ and $V_j$
are permissible.

Assume that there is a cancellation: $\theta(i+1)\equiv \theta(i)^{-1}$. Since $\Delta$ is a reduced diagram, any pair
of $(\theta,q)$-cells $\pi\in {\mathcal T}_i$ and $\pi'\in {\mathcal T}_{i+1}$ with a common $q$-edge $\bf e$ are not cancellable. Hence
the relations given by these cells are not uniquely defined
by the $q$-letter $\Lab({\mathbf e})$ and the history letter $\theta(i)$. It follows
from the list of defining relations (\ref{rel1} - \ref{rel111}) that $\Lab ({\mathbf e})$ has no superscripts while other labels of the boundary edges
of these two cells do have superscripts. Thus, these relations
are in the list (\ref{rel11}) and $\theta(i)\equiv \theta(23)$,
which prove that the history $H$ is eligible.

Since by Lemma \ref{NoAnnul} every maximal $q$-band of $\Delta$ connects
the top and the bottom of $\Delta$, it suffices to prove the last claim
under assumption that the base of $\Delta$ is a word $Q^{\pm 1}(Q')^{\pm 1}$of length $2$. Then by definition of permissible
word, $i'-i=0$, except for the base $Q_{N-1}Q_N$ (or the inverse one) with $i'-i=1$ modulo $L$ (resp., $i'-i=-1$ modulo $L$). Since all the words $U_j$ and $V_j$ have equal bases, the
last statement of (1) is proved.

(2) We can obtain the $\theta(1)$-band ${\mathcal T}_1$ by Lemma \ref{perm} (2). By induction, there is a trapezium $\Delta'$ of
height $d-1$ with bottom label $U_2\equiv U_1$ an top label $V$
such that $U_2^{\emptyset}\equiv U_1^{\emptyset}\cdot \theta(1)$ and $V_d^{\emptyset}\equiv V$, such that the union $\Delta$ of ${\mathcal T}_1$  and $\Delta'$ has history $H$. If $\Delta$ is not reduced
then we have a pair of cancellable cells $\pi\in {\mathcal T}_1$
and $\pi'\in{\mathcal T}_2$. Then as in item (1) we conlude that
$\theta(1)\equiv \theta(23)$, and so the top $\bf q$ of ${\mathcal T}_1$
has no superscript in the boundary label. Therefore one can replace
$\Delta'$ with its subscript shift $(\Delta')^{+1}$ in $\Delta$.
After such a modification, $\Delta$ becomes a reduced diagram
since for any pair cells $\pi$ and $\pi'$ with common boundary
edge from $\bf q$, the other edges have now different superscripts
in their labels. Since $V_d^{\emptyset} $ does not change under the
superscript shift, the lemma is proved.
\endproof

\subsubsection{Big and standard trapezia}

Using Lemma \ref{simul}, one can immediately derive properties of trapezia from the properties of computations obtained earlier.

If $H'\equiv \theta(i)\dots\theta(j)$ is a subword of the history $H$
from Lemma \ref{simul} (1), then the bands ${\mathcal T}_i,\dots, {\mathcal T}_j$ form a subtrapezium
$\Delta'$ of the trapezium $\Delta.$ This subtrapezium is uniquely defined by the
subword $H'$ (more precisely, by the occurrence of $H'$ in the word $\theta_1\dots\theta_d$), and $\Delta'$ is called the \index[g]{trapezium!H@$H'$-part of a trapezium where $H'$ is a subhistory} $H'$-{\it part} of $\Delta.$
\medskip

\begin{df}\label{d:standard}
We say that a trapezium $\Delta$ is \index[g]{trapezium!standard} {\it standard} if the base of $\Delta$
  is the standard base $\bf B$ of $\mmm$ or ${\mathbf B}^{-1}$, and the
  history of $\Delta$ (or the inverse one) contains one of the words (a) $\chi(i-1,i)H'\chi(i,i+1)$ (i.e.,the $S$-machine
works as $\mtt_4$) or (b) $\zeta^{i-1,i} H'\zeta^{i,i+1}$ (i.e.,it works as $\mtt_2$).
\end{df}

\begin{df} \label{bigt} We say that a trapezium $\Gamma$ is $\it big$ if
\index[g]{trapezium!big}

(1) the base of $\Gamma$ or the inverse word has the form $xvx$, where $xv$ a cyclic permutation of the $L$-th power of the standard base;

(2) the diagram $\Gamma$ contains a standard trapezium.

\end{df}

\begin{lemma}\label{or}

Let $\Delta$ be a trapezium  whose base is
$xvx$, where $x$ occurs in $v$ exactly $L-1$ times and other letters occur $<L$ times each. Then either $\Delta$ is big or the length of a side of every $\theta$-band of $\Delta$ does not exceed $c_5(||W||+||W'||)$, where $W, W'$ are the labels of its top and bottom, respectively.
\end{lemma}

\proof The diagram $\Delta$ is covered by $L$  subtrapezia $\Gamma_i$ ($i=1,\dots,L$) with bases $xv_ix$.

Assume that the the step history of $\Delta$ (or inverse step history) contains one of the subwords $\chi(i-1,i)H'\chi(i,i+1)$ or (b) $\zeta^{i-1,i} H'\zeta^{i,i+1}$. Then
by Lemma \ref{resto}  (and \ref{simul}), the base of $\Delta$ has the form $(xu)^Lx$,
where $xu$ is a cyclic permutation of the standard base (or the inverse one).
Since $\Delta$ contains a standard subtrapezia, it is is big.

Now, under the assumption that the step history has no subwords mentioned in the previous paragraph, it suffices to bound the
the length of a side of every $\theta$-band of arbitrary $\Gamma_i$ by $ \le c_4(||V|_Y+||V'||)$, where $V$ and $V'$ are the labels of the top and the bottom of $\Gamma_i$.

Assume that the word $xv_ix$ has a proper subword $yuy$, where
$u$ has no letters $y$, and any other letter occurs in $u$ at most once.
Then the word $yuy$ is faulty since $v_i$ has no letters $x$. By Lemma
\ref{nonst}, we have $|U_j|_Y\le c_1\max(|U_0|_Y,|U_t|_Y)$ for every
configuration $U_j$ of the computation given by Lemma \ref{simul} (1) restricted to the base $yuy$. Since $c_4>c_1$, it suffices to obtain
the desired estimate for the computation whose base is obtained
by deleting the subword $yu$ from $xv_ix$. Hence inducting on the
length of the base of $\Gamma_i$, one may assume that it has no proper
subwords $yuy$, and so the base of $\Gamma_i$ is revolving.
Now the required upper estimate for $\Gamma_i$ follows from Lemma \ref{narrow} (see (1) and (2c) there).
\endproof

\section{Diagrams without hubs} \label{6}

\subsection{A modified length function} \label{lf}

Let us modify the length function on the group words in $q$-, $Y$- and $\theta$-letters, and paths. The
standard length of a word (a path) will be called its
\index[g]{combinatorial length of a word (path)}{\em
combinatorial length}. From now on we use the word {\em 'length'} for
the modified length.

\begin{df}\label{d:ml}
We set the length of every
$q$-letter equal to 1, and the length of every $Y$-letter equal a small
enough number $\delta$ \index[g]{parameters used in the paper!d@$\delta$ - the length of $Y$-letters in the groups $M$ and $G$, $\delta\iv$ is between $J$ and $c_6$} given in the list of parameters (\ref{param}).

We also set to 1 the length of every word of length $\le 2$ which
contains exactly one $\theta$-letter and no $q$-letters (such
words are called $(\theta,Y)$-{\em syllables}). The length of a
decomposition of an arbitrary word into a product of letters and \index[g]{modified length function}
$(\theta,Y)$-syllables is the sum of the lengths of the factors.

The length $|w|$ of a word $w$ is the smallest length of such
decompositions. {\em The length} $|\bf p|$ {\em of a path} in a diagram is the
length of its label. The {\em perimeter} $|\partial\Delta|$ of a
van Kampen diagram over $G$ is similarly defined by cyclic decompositions
of the boundary $\partial\Delta$.
\end{df}

The next statement follows
from the property of $(\theta,q)$-relations and their cyclic permutations: the
subword between two $q$-letters
in an arbitrary $(\theta,q)$-relation is a syllable. This, in turn, follows
from Property (*) of the $S$-machine $\mmm_2$ (see Remark \ref{l:s}).

\begin{lemma} \label{ochev}
Let $\s$ be a path in a diagram $\Delta$ having $c$ $\theta$-edges
and $d$ $Y$-edges. Then

(a) $ |\s|\ge \max(c, c+(d-c)\delta)$;

(b) $|\s|=c$ if $\s$ is a top or a bottom of a $q$-band.

(c) For any product ${\s=\s_1\s_2}$ of two paths in a  diagram, we
have


\begin{equation}
|\s_1|+|\s_2|\ge |\s|\ge |\s_1|+|\s_2|-\delta \label{delta}
\end{equation}

(d) Let $\mathcal T$ be a $\theta$-band with base
of length $l_b$.  Let $l_Y$ be the number of $Y$-edges in the top
path ${\topp}(\mathcal T)$. Then the length of $\mathcal T$ (i.e., the number of cells in $\mathcal T$) is between
$l_Y-l_b$ and $l_Y+3l_b$.

\end{lemma}

\subsubsection{Rim bands}

Let ${\mathbf e}\iv {\mathbf q}_1{\mathbf f} {\mathbf q}_2\iv$ be the standard factorization of the boundary of a $\theta$-band.
If the path $({\mathbf e}\iv {\mathbf q}_1{\mathbf f})^{\pm 1} $ or the path $({\mathbf f} {\mathbf q}_2\iv {\mathbf e}\iv)^{\pm 1}$
 is the subpath of the boundary path of $\Delta$ then the band is called
 a \index[g]{band!rim band}
{\it rim} band of $\Delta.$

From now on we shall fix a constant $K$
\begin{equation}\label{kk}K>2K_0=4LN\end{equation}

The following basic facts  will allow us to remove short enough rim bands from van Kampen diagrams (see Lemma  \ref{nori} below).

\begin{lemma} \label{rim} Let $\Delta$ be a van Kampen diagram
whose rim $\theta$-band $\ttt$ has base with at most \index[g]{parameters used in the paper!k@$K$ - the maximal length of the base of a rim $\theta$-band which can be removed from a diagram, it is between $L$ and $J$, $K>2K_0=4LN$ (see (\ref{kk} and Lemma \ref{rim})}$K$ letters.
Denote by
$\Delta'$ the subdiagram $\Delta\backslash \ttt$. Then
$|\partial\Delta|-|\partial\Delta'|>1$.
\end{lemma}

{\proof Let ${\bf s}=\topp(\ttt)$ and
${\bf s}\subset\partial\Delta$. Note that
the difference between the number of $Y$-edges in
${\bf s'}=\bott(\ttt)$ the number of $Y$-edges in $s$ cannot be greater than
$2K$, because every $(\theta,q)$-relator has at most two $Y$-letters by Property (*) and the commutativity relations do not increase the number of $Y$-letters. Hence ${\bf |s'|-|s|}\le 4LN\delta$. However, $\Delta'$ is obtained
by cutting off $\ttt$ along ${\bf s'}$, and its boundary contains two
$\theta$-edges fewer than $\Delta$. Hence we have ${\bf |s_0|-|s'_0|}\ge
2-2\delta$ for the complements ${\bf s}_0$ and ${\bf s'}_0$
of $s$ and $s'$, respectively, in the boundaries $\partial\Delta$
and $\partial\Delta'$. Finally,
$$|\partial\Delta|-|\partial\Delta'|\ge
2-2\delta-2K\delta-4\delta>1 $$ by (\ref{param}),
(\ref{delta}) and the highest parameter principle .
\endproof}

\begin{df}

We call a base word $w$ \label{tightb} {\it tight} if \index[g]{S@$S$-machine!admissible words of an $S$-machine!base of an admissible word!tight}

(1) for some letter $x$ the word $w$  has the form $uxvx$,
where the letter $x$ does not occur in $u$ and $x$ occurs in $v$ exactly $L-1$ times,

(2) every proper prefix $w'$  of $w$
does not
satisfy property (1).
\end{df}

\begin{lemma}\label{width}
If a base $w$ of a $\theta$-band has no tight prefixes, then
$||w||\le K_0$, where $K_0=2LN$.
\end{lemma}

{\proof The hub base includes
every base letter $L$ times. Hence every word in this group alphabet of length $\ge K_0+1$ includes one of the letters $L+1$ times. \endproof}

\bigskip

\subsubsection{Combs}

\index[g]{comb}
\begin{df} \index[g]{comb}We say that a reduced diagram $\Gamma$ is a {\em comb} if it has a
maximal $q$-band $\mathcal Q$ (\index[g]{comb!handle of a comb}the {\em handle} of the comb), such that

\begin{enumerate}
\item[$(C_1)$]$\bott(\mathcal Q)$ is a part of $\partial \Gamma$, and every maximal
$\theta$-band of $\Gamma$ ends at a cell in $\mathcal Q$.
\end{enumerate}

If in addition the following properties hold:
\begin{enumerate}
\item[$(C_2)$] one of the maximal $\theta$-bands $\mathcal T$ in
$\Gamma$ has a tight base (if one reads the base towards the handle) and

\item[$(C_3)$] the other maximal
$\theta$-bands in $\Gamma$ have tight  bases or bases
without tight prefixes
\end{enumerate}
then the comb is called {\em tight}.
\index[g]{comb!tight comb}

The
number of cells in the handle $\mathcal Q$ is the \index[g]{comb!height of a comb} {\em height} of the
comb, and the
maximal length of the bases of the $\theta$-bands of a
comb is called the {\em basic width} of the comb.
\index[g]{comb!basic width}
\end{df}

\begin{figure}[ht]
\begin{center}
\includegraphics[width=0.4\textwidth]{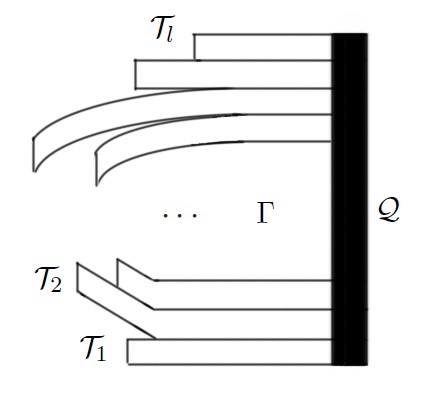}
\end{center}
\caption{A comb} \label{Pic6}
\end{figure}

Notice that every trapezium is a comb.

\begin{lemma} \label{comb} (\cite{OS06}, Lemma 4.10) Let $l$ and $b$ be the length and
the basic width of a comb $\Gamma$ and let ${\mathcal T}_1,\dots {\mathcal T}_l$ be
consecutive $\theta$-bands of $\Gamma$ (as in Figure \ref{Pic6}). We
can assume that $\bott({\mathcal T}_1)$ and ${\topp}({\mathcal T}_l)$ are contained in
$\partial\Gamma$. Denote by $\nu_a=|\partial\Gamma|_Y$ the number
of $Y$-edges in the boundary of $\Gamma$, and by $\nu_a'$ the
number of $Y$-edges on ${\bott}({\mathcal T}_1)$. Then $\nu_a + 2lb\ge
2\nu_a'$, and the area of $\Gamma$ does not exceed $c_0bl^2+2\nu_a
l$ for some constant $c_0=c_0(\mmm)$ \index[g]{parameters used in the paper!c@$c_0$ - the parameter controlling the area of a comb (see Lemma \ref{comb})}. (Recall that $c_0$ is one of the parameters from Section \ref{param}.)
\end{lemma} $\Box$

\begin{rk} The inequality with $\nu'_a$ looks stronger in
Lemma \ref{comb} than in \cite{OS06} due to the new restriction (*) from Lemma
\ref{simp}.
\end{rk}

\begin{df}
We say that a subdiagram $\Gamma$ of a diagram $\Delta$ is a
\index[g]{comb!subcom of a comb}{\it subcomb} of $\Delta$ if $\Gamma$ is a comb, the handle of $\Gamma$ divides $\Delta$ in two parts, and $\Gamma$ is one of these parts.
\end{df}

\begin{lemma} \label{est'takaya}[Compare with Lemma 4.9 of \cite{OS06}]
Let $\Delta$ be a reduced diagram
over $G$
with non-zero area, where
every rim $\theta$-band has base
of length at least $K$.
Assume that

(1) $\Delta$ is a diagram over the group $M$ or

(2) $\Delta$ has a subcomb
   of basic width at least $K_0$.

   Then there exists a maximal $q$-band $\mathcal Q$ dividing
$\Delta$ in two parts, where one of the parts is a tight subcomb with handle $\mathcal Q$.
\end{lemma}

{\proof Let $\ttt_0$ be a rim band of $\Delta$ (fig.\ref{Pic7}). Its base $w$ is of
length at least $K$, and therefore $w$ has disjoint prefix and
suffix of lengths $K_0$ since $K>2K_0$ by (\ref{kk}). The prefix of
this base word must have its own tight subprefix $w_1$, by
Lemma \ref{width} and the definition of tight words.
A $q$-edge of $\ttt_0$ corresponding to the last $q$-letter of
$w_1$ is the start edge of a maximal $q$-band $\qq'$ which bounds a
subdiagram $\Gamma'$ containing a band $\ttt$ (a subband of
$\ttt_0$) satisfying property ($C_2$). It is useful to note that a
minimal suffix $w_2$ of $w$, such that $w_2\iv$ is tight, allows us
to construct another band $\qq''$ and a subdiagram $\Gamma''$ which
satisfies ($C_2$) and has no cells in common with $\Gamma'$.

\begin{figure}[ht]
\begin{center}
\includegraphics[width=0.7\textwidth]{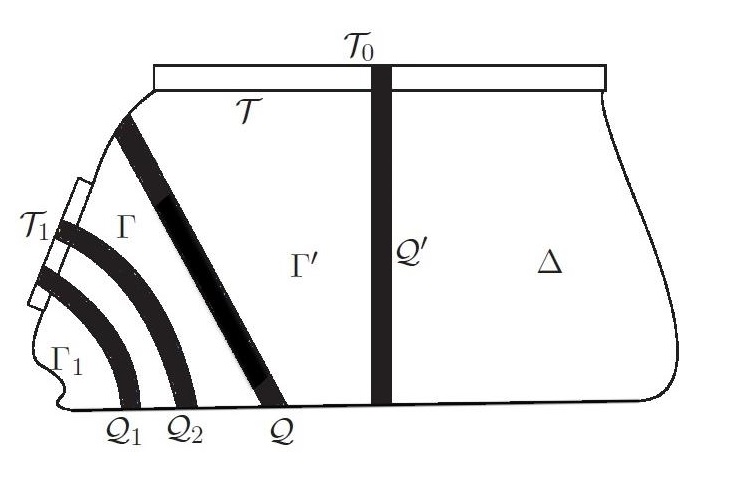}
\end{center}
\caption{Lemma \ref{est'takaya}.} \label{Pic7}
\end{figure}

Thus, there are $\qq$ and $\Gamma$ satisfying ($C_2$). Let us choose
such a pair with minimal $\area(\Gamma)$. Assume that there is a
$\theta$-band in $\Gamma$ which does not cross $\qq$. Then there
must exist a rim $\theta$-band $\ttt_1$ which does not cross $\qq$ in
$\Gamma$. Hence one can apply the construction from the previous
paragraph to $\ttt_1$ and construct two bands $\qq_1$ and $\qq_2$
and two disjoint subdiagrams $\Gamma_1$ and $\Gamma_2$ satisfying
the requirement ($C_2$) for $\Gamma$. Since $\Gamma_1$ and
$\Gamma_2$ are disjoint, one of them, say $\Gamma_1$, is inside
$\Gamma$. But the area of $\Gamma_1$ is smaller than the area of
$\Gamma$, and we come to a contradiction. Hence $\Gamma$ is a comb
and condition ($C_1$) is satisfied.

Assume that the base of a maximal $\theta$-band $\ttt$ of $\Gamma$
has a tight proper prefix (we may
assume that $\ttt$ terminates on $\qq$), and again one obtain a
$q$-band $\qq'$ in $\Gamma$, which provides us with a smaller
subdiagram $\Gamma'$ of $\Delta$, satisfying ($C_2$), a
contradiction. Hence $\Gamma$ satisfies property ($C_3$) as well.

(2) The proof is shorter since a comb is given in the very beginning. \endproof}


We will also need the definition of a derivative subcomb from \cite{O12}.

\begin{df} \label{d:dersubcomb}  \index[g]{comb!derivative subcomb of a comb} If $\Gamma$  is a comb with handle $\cal C$  and $\cal B$  is another maximal q-band in $\Gamma$,
then $\cal B$ cuts up $\Gamma$ in two parts, where the part that does not contain $\cal C$ is a comb $\Gamma_0$
with handle $\cal B$. It follows from the definition of comb, that every maximal $\theta$-band of $\Gamma$
crossing $\cal B$ connects $\cal B$ with $\cal C$. If $\cal B$  and $\cal C$ can be connected by a $\theta$-band
containing no $(\theta; q)$-cells, then $\Gamma_0$ is called the \emph{derivative subcomb} of $\Gamma$. Note that no
maximal $\theta$-band of $\Gamma$ can cross the handles of two derivative subcombs.

\end{df}

\subsection{The mixture}\label{mix}

We will need a numerical parameter associated  with van Kampen diagrams introduced in \cite{O12}, it was called \emph{mixture}.

Let $O$ be a circle with two-colored (black and white)
finite set of points (or vertices) on it. We call $O$ a \index[g]{necklace} {\it necklace} with black and white \index[g]{necklace!beads} {\it beads} on it.

Assume that there are $n$ white beads and $n'$ black ones on $O$. We define sets \index[g]{necklace!sets ${\bf P}_j$}
${\bf P}_j$ of
ordered pairs of distinct white beads as follows. A pair $(o_1,o_2)$
($o_1\ne o_2$) belongs to the set ${\bf P}_j$ if the simple arc of $O$
drawn from $o_1$ to $o_2$ in the clockwise direction has at least $j$
black beads. We denote by \index[g]{necklace!m@$\mu_J(O)$}$\mu_J(O)$ the sum $\sum_{j=1}^J \card(
 {\bf P}_j)$ (the $J$-{\it mixture}\index[g]{necklace!j@$J$-mixture} of $O$). Below similar sets for
another necklace $O'$ are denoted by ${\mathbf P'}_j$. \index[g]{parameters used in the paper!j@$J$ - the parameter of the mixture of a van Kampen diargam over $G$, it is between $K$ and $\delta\iv$}. In this subsection, $J\ge
1$, but later on it will be a fixed large enough number $J$ from the list (\ref{const}).

\begin{lemma}\label{mixture} (\cite{O12}, Lemma 6.1) (a) $\mu_J(O)\le J(n^2-n)$.

 (b) Suppose a
necklace $O'$ is obtained from $O$ after removal of a   white bead
$v$. Then \\ $
\card({\mathbf P'}_j) \le \card({\bf
P}_j)$ for every $j$, and $
\mu_J(O')\le \mu_J(O).$

(c) Suppose a necklace $O'$ is obtained from $O$ after removal of a
black bead $v$. Then  $\card({\mathbf P'}_j) \le \card({\bf P}_j)$ for
every $j,$ and $\mu_J(O')\le \mu_J(O).$

(d) Assume that there are three black beads $v_1, v_2, v_3$ of a necklace
$O,$ such that the clockwise arc $v_1 - v_3$ contains $v_2$ and has
at most $J$ black beads (excluding $v_1$ and $v_3$), and the arcs
$v_1-v_2$ and $v_2-v_3$ have $m_1$ and $m_2$ white beads,
respectively. If $O'$ is obtained from $O$ by removal of $v_2$, then
$\mu_J(O')\le\mu_J(O)-m_1m_2.$
\end{lemma}

\begin{df}
For any diagram $\Delta$ over $G$, we introduce the following invariant
$\mu(\Delta)=\mu_{J} (\partial\Delta)$
depending on the boundary of $\Delta$ only.
 To define it, we
consider the boundary $\partial(\Delta),$ as a {\it necklace}, i.e.,
we consider a circle $O$ with $||\partial\Delta||$ edges labeled as the
boundary path of $\Delta.$ By definition, the
white beads are the mid-points  of the $\theta$-edges of
$O$ and black beads are the mid-points of the $q$-edges
$O$.

The \index[g]{mixture $\mu(\Delta)$ of a diagram}
${\it mixture}$
of $\Delta$ is $\mu(\Delta)=\mu_J(O).$
\end{df}

\subsection{Quadratic upper bound for quasi-areas of  diagrams over \texorpdfstring{$M$}{M}}\label{qub}

\subsubsection{The $G$-area of a diagram over $M$}

The Dehn function of the group $M$ is super-quadratic (in fact by \cite{OS06} it is at least $n^2\log n$ because $M$ is a mulltiple HNN extension of a free group and has undecidable conjugacy problem). However we
are going to obtain a quadratic Dehn function of $G$, and
first we want to bound the areas of the words vanishing in $M$
with respect to the presentation of $G$. For this goal we artificially introduce the concept of $G$-area, as in \cite{O18}. The $G$-area of a big trapezia can be much smaller than the real area of it in $M$. This concept will be justified at the end of this paper, where some big trapezia are replaced by diagrams with hubs whose areas do not exceed the $G$-area of the trapezia.

\begin{df}\label{abt} \index[g]{G@$G$-area $\area_G(\Gamma)$ of a big trapezium $\Gamma$}The $G$-area $\area_G(\Gamma)$ of a big trapezium $\Gamma$ is, by definition, the minimum of the half of its area (i.e., the number of cells) and the product $$c_5h(||\topp(\Gamma)||+||\bott(\Gamma)||),$$ where $h$ is the height of $\Gamma$ and $c_5$ is one of the parameters from (\ref{const}).

To define the $G$-area of a diagram $\Delta$ over $M$, we consider a family $\bf S$ of big subtrapezia (i.e.,subdiagrams, which are big trapezia) and single
cells of $\Delta$ such that every cell of $\Delta$ belongs to
a member $\Sigma$ of this family, and if a cell $\Pi$ belongs to different $\Sigma_1$ and $\Sigma_2$ from $\bf S$, then both $\Sigma_1$ and $\Sigma_2$ are big subtrapezia of $\Delta$ with bases $xv_1x$, $xv_2x$, and $\Pi$ is a $(\theta,x)$-cell.
(In the later case, the intersection $\Sigma_1\cap\Sigma_2$
must be an $x$-band.)
There is
such a family 'covering' $\Delta$, e.g., just the family of all cells of $\Delta$.

The $G$-area of $\bf S$ is the sum of $G$-areas of all big
trapezia from $\bf S$ plus the number of single cells from $\bf S$
(i.e.,the $G$-area of a cell $\Pi$ is $\area_G(\Pi)=1$). Finally,
the \label{areaGd} $G$-{\it area} $\area_G(\Delta)$ is the minimum of the $G$-areas
of all "coverings" $\bf S$ as above.

\end{df}

It follows from the Definition \ref{abt} that $\area_G(\Delta)\le \area(\Delta)$ since the $G$-area of a big trapezium does not
exceed a half of its area and no cell belongs to three big trapezia of a covering.


\begin{lemma} \label{GA}  Let $\Delta$ be a reduced diagram, and suppose every cell $\pi$ of $\Delta$ belongs in one of subdiagrams $\Delta_1,\dots,\Delta_m$, where any intersection $\Delta_i\cap\Delta_j$ either has no cells or it is a $q$-band.
Then $\area_G(\Delta)\le \sum_{i=1}^m \area_G(\Delta_i)$.
\end{lemma}

{
\proof Consider the families ${\bf S}_1,\dots, {\bf S}_m$ given
by the definition of $G$-areas for the diagrams $\Delta_1,\dots,\Delta_m$. Then the family ${\bf S}={\bf S}_1\cup\dots\cup {\bf S}_m$
'covers' the entire $\Delta$ according to the above definition. This implies the required inequality for $G$-areas,
\endproof}

\subsubsection{Combs of a potential counterexample}

In this section we show that for some constants \index[g]{parameters used in the paper!n@$N_1, N_2$ - parameters controling the $G$-area of a van Kampen diagram over $M$ in terms of its perimeter and the mixture} $N_1, N_2$ the $G$-area of any reduced diagram
$\Delta$ over $M$ with perimeter $n$ does not exceed $N_2n^2+N_1\mu(\Delta)$.

Using the quadratic upper bound for
$\mu(\Delta)$ from Lemma \ref{mixture} (a), one then deduces that the $G$-area is bounded by $N'n^2$ for some constant $N'$.

Roughly speaking, we are doing the
following. We use induction on the perimeter of the diagram. First
we remove rim $\theta$-bands (those with one side and both ends on the boundary of
the diagram) with short bases. This operation decreases the
perimeter and preserves the sign of $$N_2n^2+N_1\mu(\Delta)-\area_G(\Delta),$$ so we can assume that the diagram
does not have rim $\theta$-bands. Then we use Lemma \ref{est'takaya} and
find a tight comb inside the diagram with a handle $\ccc$. We also
find a long enough $q$-band $\ccc'$ that is close to $\ccc$. We use
a surgery which amounts to removing a part of the diagram between
$\ccc'$ and $\ccc$ and then gluing the two remaining parts of
$\Delta$ together. The main difficulty is to show that, as a result
of this surgery, the perimeter decreases and the
mixture changes in such a way that the expression $$N_2n^2+N_1\mu(\Delta)-\area_G(\Delta)$$ does not change its sign.

In the
proof, we need to consider several cases depending on the shape of
the subdiagram between $\ccc'$ and $\ccc$. Note that neither
$N_2n^2$ nor $N_1\mu(\Delta)$ nor $\area_G(\Delta)$ alone behave in
the appropriate way as a result of the surgery, but the expression
$$N_2n^2+N_1\mu(\Delta)-\area_G(\Delta)$$ behaves as needed.

\medskip

%

Arguing by contradiction in the remaining part of this section, we consider a {\bf counter-example}
$\Delta$ with minimal perimeter $n$, so that

\begin{equation}\label{counter}
\area_G(\Delta) > N_2n^2+N_1\mu(\Delta)
\end{equation}

Of course, the  $G$-area of $\Delta$ is
positive, and, by Lemma \ref{NoAnnul}, we have at least 2
$\theta$-edges on the boundary $\partial\Delta$, so $n\ge 2$.

\begin{lemma} \label{notwo}  (1) The diagram $\Delta$ has no two disjoint subcombs $\Gamma_1$ and $\Gamma_2$ of basic widths at most $K$ with handles ${\mathcal B}_1$ and ${\mathcal B}_2$ such that some ends of these handles are
connected by a subpath $\bf x$ of the boundary path of $\Delta$
with $|{\bf x}|_q\le N$.

(2) The boundary of every subcomb $\Gamma$ with basic width $s\le K$ has $2s$ $q$-edges.
\end{lemma}

{\proof We will prove the Statements (1) and (2) simultaneously. We use induction on $A=\area(\Gamma_1)+\area(\Gamma_2)$ for Statement (1) and induction on $A=\area(\Gamma)$ for Statement (2).  Suppose that our diagram $\Delta$ is also a counterexample for Statement (1) or (2) with minimal possible $A$.

\begin{figure}[ht]
  \centering
  \includegraphics[width=1.0\textwidth]{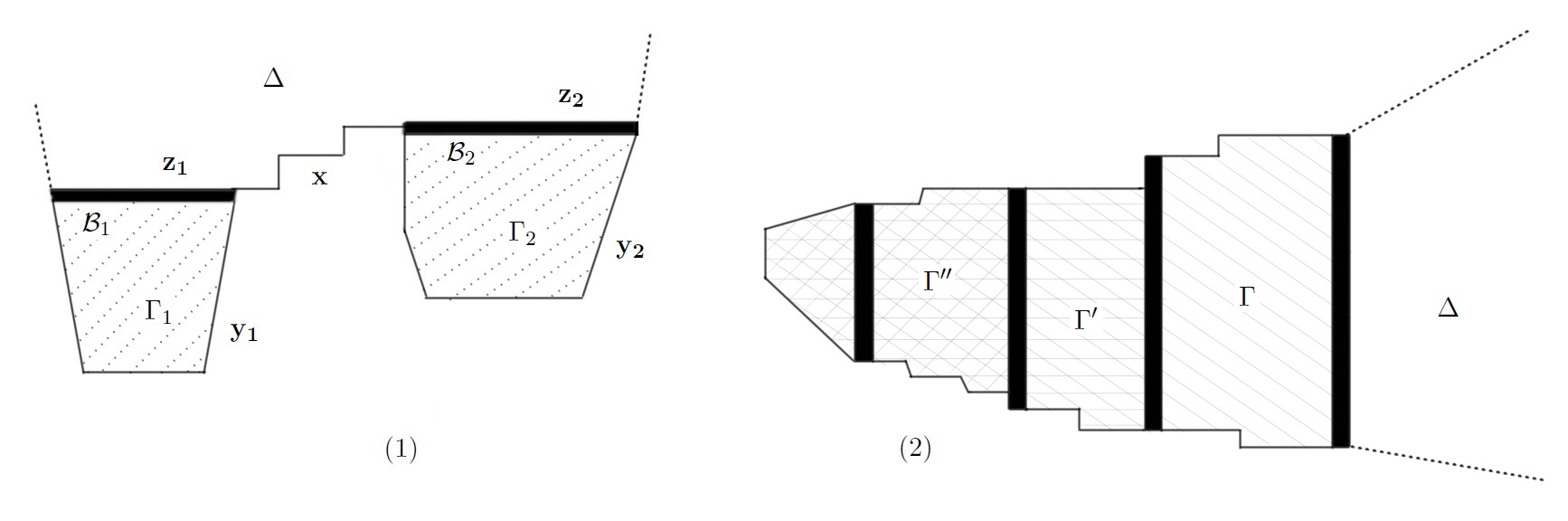}
  \caption{Lemma \ref{notwo}}\label{Pic8}
\end{figure}

Suppose that $\Delta$ is a counterexample to (1). Since the area of $\Gamma_i$
($i=1,2$) is smaller than $A$, we may use Statement (2) for $\Gamma_i$, and so we have at most $2K$ $q$-edges in $\partial\Gamma_i$.

Let $h_1$ and $h_2$ be the lengths of the handles ${\mathcal B}_1$ and ${\mathcal B}_2$ of $\Gamma_1$ and $\Gamma_2$, resp.
Without loss of generality, we assume that $h_1\le h_2$.
Denote by ${\bf y}_i{\bf z}_i$ the boundaries of $\Gamma_i$ ($i=1,2$), where ${\bf z}_i$ is the part of $\partial\Delta$
and ${\bf y}_i$ is the side of the handle of $\Gamma_i$ (so ${\bf y}_1{\bf x}{\bf y}_2$ is the part of the boundary path of $\Delta$, see Figure \ref{Pic8} (1)). Then each of the $\theta$-edges $\bf e$
of ${\bf y}_1$ is separated in $\partial\Delta$ from every $\theta$-edge $\bf f$ of ${\bf y}_2$ by
less than $4K+N < J$ $q$-edges.  Hence every such pair $({\bf e,f})$ (or the pair of white beads on these edges) makes a contribution to $\mu(\Delta)$.

Let $\Delta'$ be the diagram obtained by deleting the subdiagram
$\Gamma_1$ from $\Delta$. When passing from $\partial\Delta$ to $\partial\Delta'$,
  one replaces the $\theta$-edges (black beads) from ${\bf z}_1$ by the $\theta$-edge
  of ${\bf y}_1$ (black bead) belonging to the same maximal  $\theta$-band. The same is true for white beads.

  But each of the $h_1h_2$ pairs in the corresponding set $P'$ of white beads is
  separated in $\partial\Delta'$ by a smaller number of black
  beads than for the pair defined by $\Delta$.
Indeed, since the handle of $\Gamma_1$ is removed when one replaces $\partial\Delta$ by $\partial\Delta'$, two black beads at the ends of this handle are removed, and  therefore
 \begin{equation}\label{mudd}
 \mu (\Delta)-\mu (\Delta')\ge h_1h_2
 \end{equation}
  by Lemma \ref{mixture} (d).

  Let $\nu_a$ be the number of $Y$-edges in $\partial\Gamma_1$.
  It follows from Lemma \ref{comb} that the area, and so the $G$-area of $\Gamma_1$, does not exceed $
J(h_1)^2+2\nu_a h_1$ since $J > c_0K$.


  Since the boundary of $\Delta'$ has at
least two $q$-edges fewer than $\Delta$ and $|{\bf z}_1|=h_1\le |{\bf y}_1|$, we have $|\partial
\Delta'|\le |\partial\Delta|-2$. Moreover, we have from Lemma
\ref{ochev} (a) and Lemma \ref{NoAnnul} that
\begin{equation}\label{deriv}
|\partial \Delta|-|\partial\Delta'|\ge\gamma =\max(2,
\delta(\nu_a-2h_1))
\end{equation}
because the top/the bottom path of ${\mathcal B}_1$  has at most $h_1$
$Y$-edges.

Since $\Delta$ is a counter-example to (\ref{counter}) with minimal perimeter, $\Delta'$ is not
a counter-example by (\ref{deriv}), and so  the $G$-area of $\Delta'$ does not exceed

$$N_2|\partial\Delta'|^2+ N_1\mu(\Delta')\le  N_2(n-\gamma)^2+ N_1\mu(\Delta')$$

Hence by inequality (\ref{mudd}), we have

$$\area_G(\Delta')\le N_2(n-\gamma)^2+ N_1\mu(\Delta)-N_1h_1h_2$$

Adding the $G$-area of $\Gamma_1$
we  see that the $G$-area of $\Delta$ does
not exceed $$N_2n^2-N_2\gamma n +N_1\mu(\Delta)
-N_1h_1h_2+Jh_1^2+2\nu_a h_1.$$
Since $h_1\le h_2$, this will contradict
inequality (\ref{counter}) when we prove that

\begin{equation}\label{nado1}
- N_2\gamma n -N_1h_1^2+Jh_1^2+2\nu_a h_1<0
\end{equation}

If $\nu_a\le 4h_1$, then inequality (\ref{nado1}) follows from
the inequalities $\gamma\ge 2$ and
\begin{equation}\label{param3i}
N_1\ge J+ 8
\end{equation}

 Assume that $\nu_a> 4h_1$. Then by (\ref{deriv}), we have
$\gamma\ge \frac12 \delta \nu_a$ and so
\begin{equation}\label{first}
N_2\gamma n \ge \frac12 \delta \nu_aN_2 n>2\nu_a h_1
\end{equation}
because $n\ge 2h_1$ by Lemma \ref{NoAnnul} and
\begin{equation}\label{param4}
N_2> 2\delta^{-1}.
\end{equation}
Note that $N_1h_1^2>Jh_1^2$ by (\ref{param3i}), and this inequality
together with (\ref{first}) imply
inequality (\ref{nado1}).

(2) If there are at least two derivative subcombs of $\Gamma$, then one can
find two of them satisfying the assumptions of Statement (1).

Indeed, the derivative subcombs of $\Gamma$ are ordered linearly in a natural way  (as they are connected with the handle of $\Gamma$ by $\theta$-bands).
Consider two neighbor derivative subcombs $\Gamma_1$, $\Gamma_2$. The handle of $\Gamma_i$ are intersected by two collections of $\theta$-bands $\ccc_1, \ccc_2$ which connect these handles with the handle of $\Gamma$ (by Definition \ref{d:dersubcomb}). The  maximal $\theta$-bands that intersect the handle of $\Gamma$ and are between the two collections $\ccc_1, \ccc_2$ do not intersect any derivative combs, hence they do not intersect  $q$-bands 
except for the handle of $\Gamma$. Therefore the handles of $\Gamma_1$ and $\Gamma_2$ are connected by a
subpath $x$ of $\partial\Delta$ with no $q$-edges, so $|x|_q=0<N$.

We deduce that $\area(\Gamma_1)+\area(\Gamma_2)<\area(\Gamma)=A$, a contradiction. Therefore
there is a most one derivative subcomb $\Gamma'$ in $\Gamma$ (Figure \ref{Pic8} (2)). In turn, $\Gamma'$
has at most one derivative subcomb $\Gamma''$, and so on. It follows that
there are no maximal $q$-bands in $\Gamma$ except for the handles of
$\Gamma', \Gamma'',\dots $. Since the basic width of $\Gamma$ is $s$, we have
$s$ maximal $q$-bands in $\Gamma$, and the lemma is proved.
\endproof}


\begin{lemma} \label{twocombs}  There is no pair of subcombs $\Gamma$ and $\Gamma'$
in $\Delta$ with handles $\mathcal X$ and $\mathcal X'$ of length $\ell$ and $\ell'$ such that $\Gamma'$
is a subcomb of $\Gamma$, the basic width of $\Gamma$ does not exceed $K_0$ and
$\ell'\le \ell/2$.
\end{lemma}

{\proof
Proving by contradiction, one can choose $\Gamma'$ so that $\ell'$ is minimal for
all subcombs in $\Gamma$ and so $\Gamma'$  has no
proper subcombs, i.e. its basic width is $1$ (fig. \ref{Pic9}).
It follows from Lemma \ref{comb}
that for $\nu_Y' =|\partial\Gamma'|_Y$, we have
\begin{equation}\label{G'}
\area_G(\Gamma')\le\area(\Gamma')\le c_0(\ell')^2+2\nu_Y' \ell'
\end{equation}

\begin{figure}[ht]
\begin{center}
\includegraphics[width=0.9\textwidth]{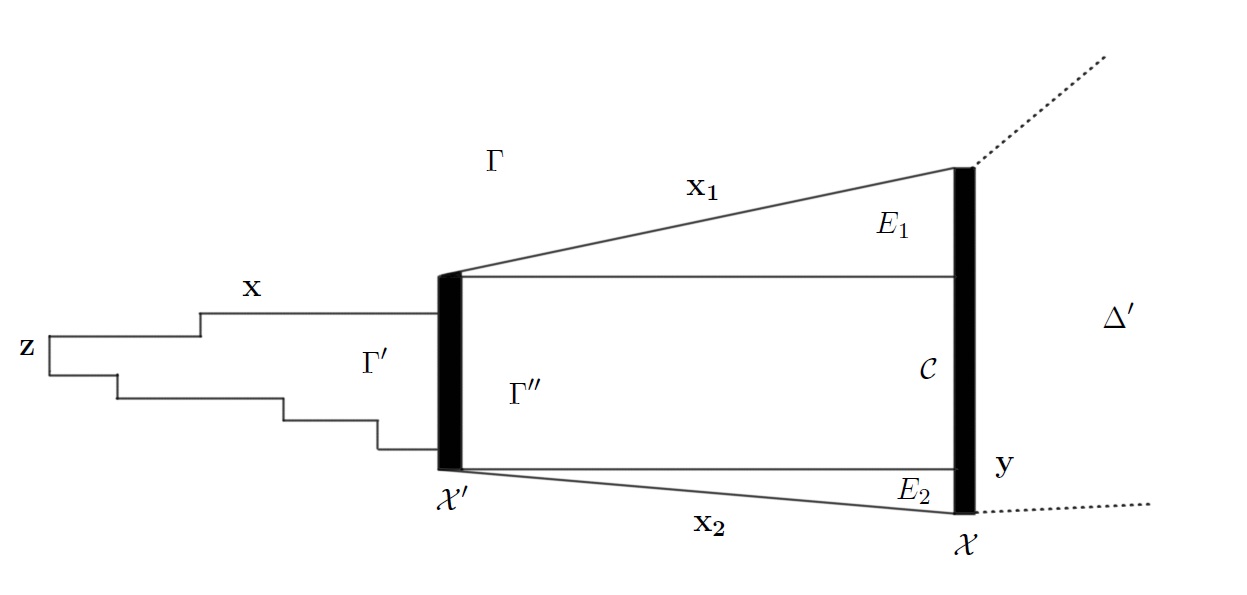}
\end{center}
\caption{Lemma \ref{twocombs}}\label{Pic9}
\end{figure}

Let $\Delta'$ be the diagram obtained after removing the subdiagram
$\Gamma'$ from $\Delta$. The following inequality is the analog of (\ref{deriv})
(where $h_1$ is replaced by $\ell'$)
\begin{equation}\label{eq679}
|\partial \Delta|-|\partial\Delta'|\ge\gamma =\max(2,
\delta(\nu_Y'-2\ell'))
\end{equation}

The $q$-band $\mathcal X$ contains a subband $\mathcal C$ of length
$\ell'$. Moreover one can choose $\mathcal C$ so that all maximal $\theta$-bands of
$\Gamma$  crossing the handle $\mathcal X'$ of $\Gamma'$, start from $\mathcal C$.
These $\theta$-bands form a comb $\Gamma''$ contained in $\Gamma$, and in turn,
$\Gamma''$ contains $\Gamma'$. The two parts of the
complement ${\mathcal X}\backslash{\mathcal C}$ are the handles of two subcombs $E_1$ and $E_2$
formed by maximal $\theta$-bands of $\Gamma$, which do not cross $\mathcal X'$.
Let the length of these two handles be $\ell_1$ and $\ell_2$, respectively, and so we have $\ell_1+\ell_2=\ell-\ell'>\ell'$. ($E_1$ or $E_2$ can be empty;
then $\ell_1$ or $\ell_2$ equals $0$.)

It will be convenient to assume that $\Gamma$ is drawn from the left of the vertical handle $\mathcal X$. Denote by ${\bf yz}$ the boundary path of 
of $\Gamma$, where ${\bf y}$ is the right side of the band $\mathcal X$. Thus, there are $\ell_1$ (resp., $\ell_2$)
$\theta$-edges on the common subpath ${\bf x}_1$ (subpath ${\bf x}_2$) of ${\bf z}$ and $\partial E_1$ (and $\partial E_2$).

By Lemma \ref{notwo} (2),
the path ${\bf z}$ contains at most $2K_0$ $q$-edges, because the basic width of $\Gamma$ is at most $K_0$.

Consider the factorization ${\bf z=x}_2{\bf xx}_1$, where ${\bf x}$ is a subpath
of $\partial \Gamma'$.
It follows that between every white bead on ${\bf x}_1$ (i.e. the middle point of the $\theta$-edges on ${\bf x}_1$) and a white bead on $\bf x$ we
have at most $2K_0$ black beads (i.e. the middle points of the
$q$-edges of the path $\bf x$). Since $J$ is greater than $2K_0$,
every pair of white beads, where one bead belongs in $\bf x$
and another one belongs in ${\bf x}_1$ (or, similarly,  in ${\bf x}_2$) contributes $1$ to $\mu(\Delta)$. Let $P$ denote the set
of such pairs. By the definition of $E_1$
and $E_2$, we have $\card(P) =\ell'(\ell_1+\ell_2)=\ell'(\ell-\ell')>(\ell')^2$.

  When passing from $\partial\Delta$ to $\partial\Delta'$,
  one replaces the left-most $\theta$-edges of every maximal $\theta$-band from $\Gamma'$ with the right-most $\theta$-edges
  lying on the right side of ${\mathcal X}'$. The same is true for white beads.
  But each of the $\ell'(\ell-\ell')$ pairs in the corresponding set $P'$  of white beads is
  separated in $\partial\Delta'$ by smaller number of black
  beads since the $q$-band $\mathcal X'$ is removed. Therefore every
  pair from $P'$ gives less by $1$ contribution to the mixture, as it follows from the definition of mixture. Hence $\mu (\Delta)-\mu (\Delta')\ge \ell'(\ell-\ell')\ge(\ell')^2$.
  This inequality and inequality (\ref{eq679}) imply that
$$\area_G(\Delta')\le N_2|\partial\Delta'|^2 +N_1\mu(\Delta')\le N_2(n-\gamma)^2+ N_1\mu(\Delta)-N_1(\ell')^2,$$
because the perimeter of $\Delta'$ is less than the perimeter of the minimal counter-example $\Delta$.
Adding the estimate of $G$-area of $\Gamma'$ (\ref{G'})
we see that

$$\area_G(\Delta)\le N_2 n^2 +N_1\mu(\Delta)- N_2\gamma n
-N_1(\ell')^2+c_0(\ell')^2+2\nu_Y' \ell'.$$
This will contradict the fact that $\Delta$ is a counterexample of (\ref{counter}) when we prove that
\begin{equation}\label{nado}
- N_2\gamma n -N_1(\ell')^2+c_0(\ell')^2+2\nu_Y' \ell'<0,
\end{equation}
  Consider two cases.

(a) Let $\nu_Y'\le 4\ell'$. Then inequality (\ref{nado}) follows from
the inequalities $\gamma\ge 2$ and
 $$N_1\ge c_0+ 8.$$

(b) Assume that $\nu_Y'> 4\ell'$. Then by (\ref{eq679}) we have
$\gamma\ge \frac12 \delta \nu_Y'$ and so
\begin{equation}\label{second}
N_2\gamma n\ge\frac12 \delta \nu_Y' N_2n  >2 \nu_Y' \ell'
\end{equation}
by (\ref{param4}) since $n\ge 2\ell\ge 4\ell'$ by Lemma \ref{NoAnnul}.

Also we have  $N_1(\ell')^2>c_0(\ell')^2$, which together with (\ref{second}) implies
 (\ref{nado}).

Thus, the lemma is proved by contradiction.
\endproof}


\subsubsection{Removing rim $\theta$-bands}



Recall that $K>2K_0=4LN$.

\begin{lemma} \label{nori} $\Delta$ has no rim $\theta$-band
whose base has  $s\le K$ letters.
\end{lemma}

\begin{figure}[ht]
\begin{center}
\includegraphics[width=0.5\textwidth]{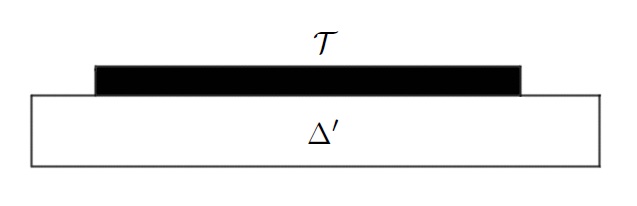}
\end{center}
\caption{Rim $\theta$-band}\label{Pic10}
\end{figure}

{
\proof Assume by contradiction that such a rim $\theta$-band $\mathcal T$ exists,
and ${\bf top}(\mathcal T)$ belongs in
$\partial(\Delta)$ (fig.\ref{Pic10}). When deleting $\mathcal T$, we obtain, by Lemma
\ref{rim}, a diagram $\Delta'$ with $ |\partial\Delta'|\le
n-1$. Since ${\bf top}(\mathcal T)$ lies on $\partial\Delta$, we have
from the definition of the length , that the
number of $Y$-edges in ${\bf top}(\mathcal T)$ is less than
$\delta^{-1}(n-s)$. By Lemma \ref{ochev}, the length of $\mathcal T$ is at
most $3s+\delta^{-1}(n-s)< \delta^{-1}n$. Thus, by applying the inductive
hypothesis to $\Delta'$, we have that $G$-area of $\Delta$ is not
greater than $N_2(n-1)^2 +N_1\mu(\Delta)
+\delta^{-1}n$ because $\mu(\Delta')\le\mu(\Delta)$ by Lemma
\ref{mixture} (b). But the first term of this sum does not exceed $N_2 n^2-N_2n$
and so the entire sum is bounded by
$N_2 n^2 +N_1\mu(\Delta)$
 provided

\begin{equation}\label{param5}
N_2\ge\delta^{-1}.
\end{equation}

This contradicts the choice of $\Delta$, and the lemma is proved.
\endproof}

\subsubsection{The quadratic estimate}

The next lemma is one of the main ingredient in this section.

\begin{lemma}\label{main} The $G$-area of a reduced diagram $\Delta$ over $\mmm$ does not exceed  $N_2n^2 +N_1\mu(\Delta)$, where
$n=|\partial\Delta|$.
\end{lemma}

{\it Proof.} We continue studying the hypothetical  counter-example $\Delta$ of minimal  possible perimeter. By Lemma \ref{nori}, now we can apply Lemma \ref{est'takaya} (1). By that
lemma, there exists a tight subcomb $\Gamma\subset\Delta$. Let $\mathcal T$ be a
$\theta$-band of $\Gamma$ with a tight base.

\begin{figure}[ht]
\begin{center}
\includegraphics[width=1.0\textwidth]{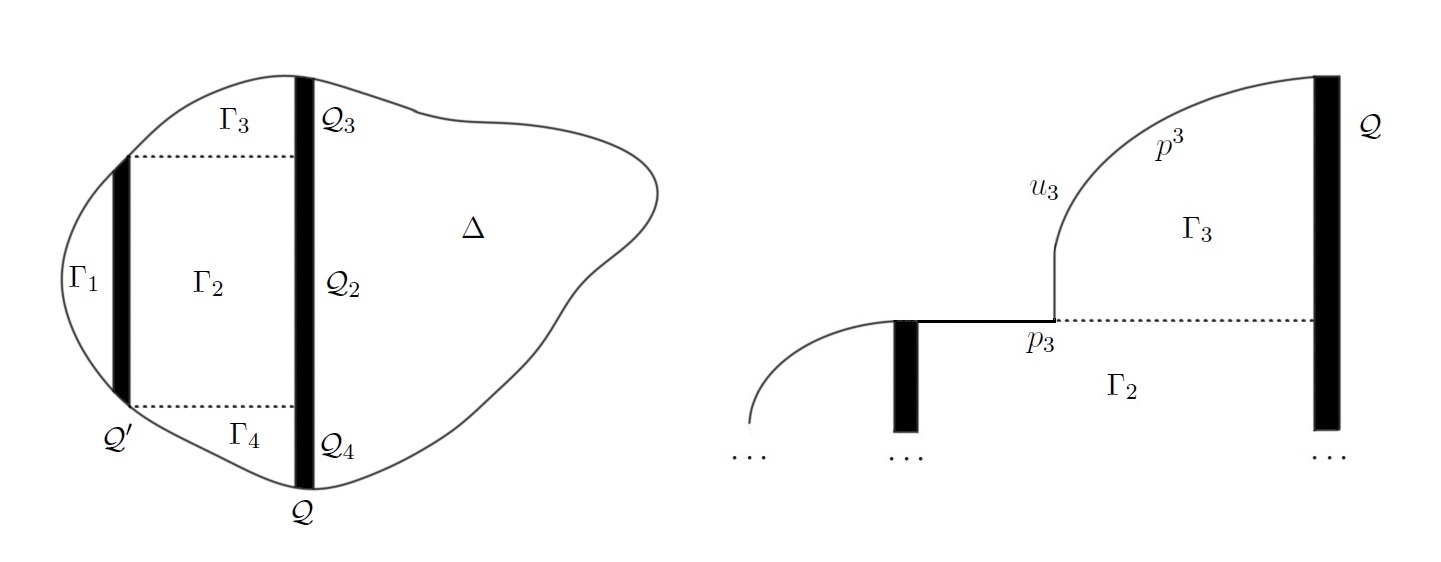}
\end{center}
\caption{Lemma \ref{main}.} \label{Pic18}
\end{figure}

The basic width of $\Gamma$ is less than $K_0$ by Lemma \ref{width}.
Since the base of $\Gamma$ is tight, it is equal to $uxvx$ for
some $x$, where the last occurrence of $x$ corresponds to the handle $\mathcal Q$ of $\Gamma$,
 the word $u$ does not contain $x$, and $v$ has exactly $L-1$  occurrences of $x$.
Let $\mathcal Q'$ be the maximal
$x$-band of $\Gamma$ crossing $\mathcal T$ at the cell corresponding to
the first occurrence of $x$ in $uxvx$ (fig. \ref{Pic18} (a)).

We consider the smallest subdiagram $\Gamma'$ of $\Delta$ containing
all the $\theta$-bands of $\Gamma$ crossing the $x$-band $\mathcal Q'$. It
is a comb with handle ${\mathcal Q}_2\subset{\mathcal Q}$. The comb $\Gamma'$ is
covered by a trapezium $\Gamma_2$ placed between $\mathcal Q'$ and $\mathcal Q$,
and a comb $\Gamma_1$ with handle $\mathcal Q'$. The band $\mathcal Q'$
belongs to both $\Gamma_1$ and $\Gamma_2$. The remaining part of
$\Gamma$ is a disjoint union of two combs $\Gamma_3$ and $\Gamma_4$
whose handles ${\mathcal Q}_3$ and ${\mathcal Q}_4$ contain the cells of $\mathcal Q$ that do
not belong to the trapezium $\Gamma_2$. The handle of $\Gamma$ is
the composition of handles ${\mathcal Q}_3$, ${\mathcal Q}_2$, ${\mathcal Q}_4$ of $\Gamma_3$,
$\Gamma'$ and $\Gamma_4$ in that order.

Let the lengths of ${\mathcal Q}_3$ and ${\mathcal Q}_4$ be $\ell_3$ and $\ell_4$,
respectively. Let $\ell'$ be the length of the handle of $\Gamma'$.
Then by Lemma \ref{twocombs}, we have
\begin{equation}\label{ll'}
\ell'>\ell/2 \;\;\; and \;\;\;\ell=\ell'+\ell_3+\ell_4
\end{equation}

For $i\in \{3,4\}$ and $\nu_i=|\partial\Gamma_i|_Y$, Lemma
\ref{comb} and the highest parameter principle (\ref{const}) give inequalities
\begin{equation}\label{gamma34}
A_i\le J\ell_i^2+2\nu_i \ell_i,
\end{equation}
where $A_i$ is the $G$-area of $\Gamma_i$. (We take into account that $G$-area cannot exceed area.)

Let
${\bf p}_3, {\bf p}_4$ be the top and the bottom of the trapezium $\Gamma_2$.
Here ${\bf p}_3^{-1}$ (resp. ${\bf p}_4^{-1}$) shares some initial edges with $\partial
\Gamma_3$ (with $\partial \Gamma_4$), the rest of these paths belong
to the boundary of $\Delta$. We denote by $d_3$ the number of
$Y$-edges of ${\bf p}_3$ and by $d'_3$ the number of the $Y$-edges of ${\bf p}_3$ which do
not belong to $\Gamma_3$. Similarly, we introduce $d_4$ and $d'_4$.

Let $A_2$ be the $G$-area of $\Gamma_2$. Then by Lemma \ref{or} and the definition of the $G$-area for big trapezia (if $\Gamma_2$ is big),
we have

\begin{equation}\label{gamma2}
A_2\le c_5\ell'(d_3+d_4+2K)<J\ell'(d_3+d_4+1)
\end{equation}
 because the basic width of $\Gamma_2$ is less than $K$ and $J>2Kc_5$ by (\ref{const}).

Recall that the first and the last base letters of the base of the trapezium $\Gamma_2$
are equal to $x$. So for every maximal $\theta$-band $\cal T$, the first and the
last $(\theta,q)$-cells have equal boundary labels up to some superscript shift
$+k$ (if there are superscripts in these labels). However $k$ does not depend on
the choice of $\cal T$ by the last statement of Lemma \ref{simul} (1).
Therefore the whole ${\cal Q'}^{(+k)}$  is a copy of ${\cal Q}_2$, and so there is
a superscript shift $\Gamma_1^{(+k)}$ of the entire comb $\Gamma_1$ such
that the handle $({\mathcal Q}')^{(+k)}$ of $\Gamma_1^{(+k)}$ is a copy of ${\mathcal Q}_2$.

This makes the following
surgery possible. The diagram $\Delta$ is covered by two
subdiagrams: $\Gamma$ and another subdiagram $\Delta_1$, having only
the band ${\cal Q}_2$ in common. We construct a new auxiliary diagram by
attaching
$\Gamma_1^{(+k)}$ to $\Delta_1\cup \cal Q$
    with identification of the band
$({\mathcal Q}')^{(+k)}$ of $\Gamma_1^{(+k)}$  and the band ${\mathcal Q}_2$.
We denote the constructed diagram by $\Delta_0$.

Note that $\area_G(\Gamma_1^{(+k)})= \area_G(\Gamma)$  and $\Delta_0$ is a reduced diagram because every pair of
its cells having a common edge, has a copy either in $\Gamma_1$ or
in $\Delta_1\cup \mathcal Q$. Now we need the following claim.

\begin{lemma} \label{A0}The $G$-area $A_0$ of $\Delta_0$ is at least the
sum of the $G$-areas of $\Gamma_1$ and $\Delta_1$ minus $\ell'$.
\end{lemma}

\proof Consider a minimal covering $\bf S$ of $\Delta_0$ from Definition \ref{abt}of $G$-area, and assume that there is a big trapezium $E\in \bf S$, such that neither $\Gamma_1^{(+k)}$ nor $\Delta_1$ contains it. Then $E$ has a base $ywy$, where
$(yw)^{\pm 1}$ is a cyclic permutation of the $L$-th power of the standard base,
and the first $y$-band of $E$ is in $\Gamma_1^{(+k)}$, but it is not a subband of $\mathcal Q'$.

Since the history $H$ of the big trapezium $E$ is a subhistory
of the history of $\Gamma_2$,
and $H$ uniquely determines the base starting with given letter by Lemma  \ref{resto}, we conclude that $\Gamma_2$ is a big trapezium itself, and therefore $(xv)^{\pm 1}$ is an $L$-th power
of the standard base. Since the first $y$ occurs in uxvx
before the first $x$ it follows that we have the $(L+1)-th$
occurrence of $y$ before the last occurrence of $x$ in the word $uxvx$. But this contradicts the definition of tight comb $\Gamma$.

Hence every big trapezium from $\bf S$ entirely belongs
either in $\Gamma_1^{(+k)}$ or in $\Delta_1$. Therefore one can obtain
'coverings' $\bf S'$ and $\bf S''$ of these two diagrams
if (1) every $\Sigma$ from $\bf S$ is assigned either to $\bf S'$ or to $\bf S''$ and then (2) one add at most $\ell'$ single cells since the common band ${\mathcal Q}'$ in $\Delta_0$ should
be covered twice in disjoint diagrams $\Gamma_1^{(+k)}$ and $\Delta_1$.
These construction complete the proof of the lemma.
\endproof

Let us continue the proof of Lemma \ref{main}.

By Lemma \ref{GA}, the $G$-area of $\Delta$ does not exceed the sum
of $G$-areas of the five subdiagrams
$\Gamma_1$, $\Gamma_2$, $\Gamma_3$, $\Gamma_4$ and $\Delta_1$. But
the direct estimate of each of these values is not efficient.
Therefore we will use Lemma \ref{A0} to bound the $G$-area of the
auxiliary diagram $\Delta_0$ built of two pieces $\Gamma_1$ and $\Delta_1$.

It follows from our constructions and lemmas \ref{GA}, \ref{A0},  that
\begin{equation}\label{AGD}
 \area_G(\Delta)\le A_2+A_3+A_4+A_0+l'
 \end{equation}

Let ${\bf p}^3$ be the segment of the boundary $\partial\Gamma_3$ that
joins $\mathcal Q$ and $\Gamma_2$ along the boundary of $\Delta$ (fig. \ref{Pic18} (b)). It
follows from the definition of $d_3$, $d'_3$, $\ell_3$ and $\nu_3$,
that the number of $Y$-edges lying on ${\bf p}^3$ is at least $\nu_3
-(d_3-d'_3)-\ell_3$.

Let ${\bf u}_3$ be the part of $\partial\Delta$ that contains ${\bf p}^3$ and
connects $\mathcal Q$ with $\mathcal Q'$. It has $l_3$ $\theta$-edges. Hence we
have, by Lemma \ref{ochev}, that
$$|{\bf u}_3|\ge \max(\ell_3, \ell_3+\delta(|p^3|_Y-\ell_3))\ge
\max(\ell_3, \ell_3+\delta(\nu_3 -(d_3-d'_3)-2\ell_3)).$$ Since ${\bf u}_3$ includes a subpath of length $d'_3$ having no $\theta$-edges, we also have by Lemma \ref{ochev} (c) that $|{\bf u}_3|\ge
\ell_3+\delta(d'_3-1)$.

One can similarly define ${\bf p}^4$ and ${\bf u}_4$ for $\Gamma_4$. When
passing from $\partial\Delta$ to $\partial\Delta_0$ we replace the
end edges of $\mathcal Q'$, ${\bf u}_3$ and ${\bf u}_4$ by two subpaths of
$\partial\mathcal Q$ having lengths $\ell_3$ and $\ell_4$. Let
$n_0=|\partial\Delta_0|$. Then it follows from the previous paragraph
that

\begin{equation}\label{nbezn0}
n-n_0\ge2+\delta(\max(0,d'_3-1, \nu_3-(d_3-d'_3)-2\ell_3)+\max(0,
d'_4-1,\nu_4-(d_4-d'_4)-2\ell_4))
\end{equation}

In particular, $n_0\le n-2$. By the inductive hypothesis,
\begin{equation}\label{A01}
A_0\le N_2 n_0^2 +N_1 \mu(\Delta_0)
\end{equation}

We note that the mixture $\mu(\Delta_0)$ of $\Delta_0$ is not
greater than $\mu(\Delta)- \ell'(\ell-\ell')$ . Indeed, by Lemma \ref{twocombs} (2), one can use the same trick as in Lemma \ref{twocombs} as follows.
For every pair of white beads, where
one bead corresponds to a $\theta$-band of $\Gamma_2$ and
another one to a $\theta$-band of $\Gamma_3$ or $\Gamma_4$,
the contribution of this pair to $\mu(\Delta_0)$ is less than the contribution to $\Delta$. It remains to count
the number of such pairs: $\ell'(\ell_3+\ell_4)=\ell(\ell-\ell')$.

Therefore, by inequality (\ref{A01}), the $G$-area of $\Delta$ is not
greater than
\begin{equation}\label{ADelta}
N_2 n^2 +N_1\mu(\Delta)- N_2 n(n-n_0) -N_1
\ell'(\ell-\ell')+A_2+A_3+A_4+\ell'
\end{equation}

In view of inequalities (\ref{gamma2}), (\ref{gamma34}) for the terms $A_2$, $A_3$ and $A_4$, to
obtain the desired contradiction with (\ref{counter}), it suffices to prove that
\begin{equation}\label{tsel}
N_2 n(n-n_0)+ N_1 \ell'(\ell-\ell')\ge  Jl'(d_3+d_4+1)+J(\ell_3^2+\ell_4^2)+
2\nu_3 \ell_3 + 2\nu_4\ell_4+\ell'
\end{equation}



First we can choose $N_1$ big enough so that
$N_1 \ell'(\ell-\ell')/3\ge J(\ell_3+\ell_4)^2\ge J(\ell_3^2+\ell_4^2)$.
Indeed, by (\ref{ll'}), we obtain
$\frac{N_1}{3} \ell'(\ell-\ell')\ge \frac{N_1}{3}
(\ell_3+\ell_4)(\ell_3+\ell_4)$, so it is enough to assume that
\begin{equation}
\label{param6} N_1>3J.
\end{equation}

We also have that
\begin{equation}\label{a4}
\frac{N_2}2 n(n-n_0)\ge Jl' +\ell'
\end{equation}
because $n-n_0\ge 2$, $n\ge 2\ell'$ and $N_2\ge J$ by (\ref{param6}).

It remains to prove that \begin{equation}\label{t1} \frac{N_2}2
n(n-n_0)+\frac{2N_1}{3}\ell'(\ell-\ell')>
J\ell'(d_3+d_4)+2\nu_3\ell_3+2\nu_4\ell_4.
\end{equation}

We assume without loss of generality that $\nu_3\ge\nu_4$, and
consider two cases.

\medskip

(a) Suppose $\nu_3\le 2J(\ell-\ell')$.

Since $d_i\le \nu_i+d'_i$ for $i=3,4$,  by inequality
(\ref{nbezn0}), we have $$d_3+d_4\le
\nu_3+\nu_4+d_3'+d_4'<4J(\ell-\ell')+\delta^{-1}(n-n_0)+2-2\delta\iv<
4J(\ell-\ell')+\delta\iv(n-n_0).
$$

Therefore

\begin{equation}\label{a1}
\frac{N_1}{3} \ell'(\ell-\ell')+\frac{N_2}2 n(n-n_0)\ge
4J^2\ell'(\ell-\ell')+J\delta\iv(n-n_0)l'>
  Jl'(d_3+d_4)
\end{equation}
since we can assume by (\ref{const}) that

\begin{equation}\label{param7}
N_1> 12J^2,\qquad N_2/2> J\delta^{-1}.
\end{equation}

We also have  by (\ref{ll'}):
\begin{equation}\label{a3}
\frac{N_1}{3} \ell'(\ell-\ell')\ge \frac{N_1}{3}(\ell_3+\ell_4)(\ell_3+\ell_4)\ge
\frac{N_1}{3}\frac{\nu_3+\nu_4}{4J}(\ell_3+\ell_4)>
2\nu_3\ell_3+2\nu_4\ell_4
\end{equation}
since we can assume by (\ref{const}) that
\begin{equation}\label{param8}
N_1>24J.
\end{equation}

 The sum of inequalities (\ref{a1}) and (\ref{a3}) gives us the
desired inequality (\ref{t1}).

\medskip

(b) Assume now that $\nu_3>2J(\ell-\ell')$. Then, applying Lemma
\ref{comb} to the comb $\Gamma_3$, we obtain
\begin{equation}\label{dd'}
d_3-d'_3<\frac12\nu_3+K_0l_3\le\frac56\nu_3
\end{equation}
since $\ell_3\le
\ell-\ell'<\frac{\nu_3}{2J}$ and
\begin{equation}\label{param9} J>3K_0.\end{equation} We also have
$d_4-d'_4<\frac12\nu_4+K_0\ell_4\le\frac56 \nu_3$. These two
inequalities and inequality (\ref{nbezn0}) lead to
\begin{equation}\label{d34}
d_3+d_4\le \frac53\nu_3+\delta^{-1}(n-n_0)
\end{equation}

It follows from (\ref{dd'}) that
$$\nu_3-(d_3-d'_3)-2l_3\ge \frac16 \nu_3 -
\frac{2}{2J}\nu_3\ge\frac17\nu_3,$$ since $\ell_3\le
\ell-\ell'<\frac{\nu_3}{2J}$ and $J>42$ by (\ref{const}). Therefore,
by (\ref{nbezn0}),
\begin{equation}\label{raznitsa}
n-n_0\ge \frac17 \delta\nu_3.
\end{equation}
Thus, by (\ref{d34}),

\begin{equation}
\label{d10} d_3+d_4<13 \delta\iv(n-n_0).
\end{equation}

Since $2\ell'<n$ and $n-n_0\ge 2$, inequality (\ref{d10}) implies

\begin{equation}\label{b1}
\frac{N_2}3 n(n-n_0)> Jl'(d_3+d_4)
\end{equation}
because we can assume that \begin{equation}\label{param10}
N_2\gg J\delta\iv
\end{equation}
($N_2>21J\delta\iv$ is enough).

Inequalities (\ref{raznitsa}), (\ref{param10}),
$\nu_3\ge\nu_4$, and $4(\ell_3+\ell_4)\le n$ give us
\begin{equation}\label{b2}
\frac{N_2}6 n(n-n_0)\ge \frac72 J \delta^{-1}(n-n_0)n\ge
2\nu_3(\ell_3+\ell_4)\ge 2\nu_3\ell_3+2\nu_4 \ell_4
\end{equation}

The inequality (\ref{t1}) follows now from inequalities (\ref{b1}),
and (\ref{b2}). $\Box$

\section{Minimal diagrams over \texorpdfstring{$G$}{G}}\label{midi}


\subsection{Diagrams with hubs}

Given a reduced diagram $\Delta$ over the group $G,$ the maximal $q$-bands start and end either
on the boundary $\partial\Delta$ or on the boundaries of hubs. Therefore one can construct a planar graph whose vertices
are the hubs of this diagram plus one improper vertex outside $\Delta,$ and the
edges are the maximal $\tt$-bands of $\Delta.$

\subsubsection{Eliminating pairs of hubs connected by two $\tt$-bands}\label{elim}

Let us consider two hubs $\Pi_1$ and $\Pi_2$ in a reduced diagram,
connected by two neighbor
$\tt$-bands
$\ccc$ and  $\ccc'$, and
there are no other hubs between these $\tt$-bands.
By Lemma \ref{NoAnnul}, these bands, together with parts of
\begin{figure}[ht]
\begin{center}
\includegraphics[width=0.8\textwidth]{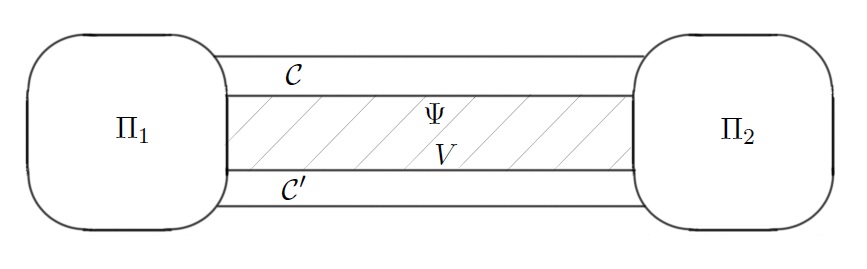}
\end{center}
\caption{Cancellation of two hubs}\label{2h}
\end{figure}
$\partial\Pi_1$ and $\partial\Pi_2,$ bound either a subdiagram
having no cells,
or  a trapezium $\Psi$ of height $\ge 1$ (fig. \ref{2h}).

The former case is impossible.  Indeed, in this case
the hubs have to correspond to the same hub relation since the relations
(\ref{rel3}) have no common letters. Hence the diagram is not reduced since a cyclic permutation of a hub relation starting with a fixed copy of the letter $\tt$ is unique.

We want to show that the latter case is not possible either if the diagram $\Delta$ is
chosen with minimal number of hubs among the diagrams with the same boundary label.

Indeed, by Lemma \ref{perm} (1), the $\tt$-band $\ccc'$ is a
$k$-shift of $\ccc$ In fact, $k=\pm 1$ since the superscripts
of the letters in $W_{st}^L$ change by one after every $\tt$-letter. One may assume that $k=1$.  So if we construct a $1$-shift $\Psi_2$
of $\Psi_1=\Psi$, then  the first maximal $\tt$-band of $\Psi_2$
is a copy of $\ccc'$ (the second $\tt$-band in $\Psi_1$).
Similarly one can construct $\Psi_3=\Psi_2^{(+1)}=\Psi_1^{(+2)},\dots,
\Psi_L=\Psi_1^{(+L)}$. Let us separately construct an auxiliary
diagram $\Delta_1$ consequently attaching the bottoms of $\Psi_1, \Psi_2,\dots,\Psi_L$ to $\Pi_1$ and identifying the second $\tt$-band
of $\Psi_i$ with the first $\tt$-band of $\Psi_{i+1}$ (indices modulo $L$). This is possible since the $L$-shift of any diagram is equal
to itself. Now we can attach $\Pi_2$ to the tops of $\Psi_i$-s in $\Delta_1$ and
obtain a spherical diagram $\Delta_2$. The diagram $\Delta_2$ contains a copy of the subdiagram $\Gamma$ of $\Delta$ formed by
$\Pi_1$, $\Pi_2$ and $\Psi$. Hence the boundary label of $\Gamma$
is equal to the boundary label of the complement $\Gamma'$ of
(the copy of) the subdiagram $\Gamma$ in $\Delta_2$. Thus, one can replace $\Gamma$ with $\Gamma'$ in $\Delta$ decreasing the number
of hubs.

\subsubsection{Disks}\label{di}

\begin{df} \label{dw} A permissible word $V$ is called a \index[g]{disk word} {\it disk word} if $V^{\emptyset}\equiv W^L$ for
some accessible word $W$. The cyclic permutations of $W$ and $W^{-1}$ are also disk words by definition.
\end{df}

\begin{lemma}\label{trivial} Every disk word $V$ is equal to $1$ in the group $G$.
\end{lemma}

\proof Assume there is an eligible  computation $W_{st}\to\dots\to W$, where  $V^{\emptyset}\equiv W^L$. Then the computation
$W_{st}^L\to\dots\to W^L$ with the same history is eligible too. By Lemma \ref{simul} (2),
one can construct a trapezium  $\Delta$ with bottom label $W_{st}^{(1)}\dots W_{st}^{(L)}$ and top label $V'$ such that
$(V')^{\emptyset}\equiv V^{\emptyset}$, and so $V'$ is a cyclic permutation of the word $V$.
The two sides of $\Delta$ have equal labels since the $L$-shift preserves superscripts.
So one can identify these sides and attach the obtained annulus to the hub cell
labeled by $W_{st}^{(1)}\dots W_{st}^{(L)}$. Since $V'$ is the boundary
label of the obtained disk diagram, we have $V'=1$ in $G$, and so $V=1$, as required.
If there there is an eligible computation $W\to\dots\to W_{ac}$, then the proof is similar
with bottom label of $\Delta$ equal to $W_{ac}^L$.
\endproof

\begin{rk}\label{dr} In fact, for the disk word $W$, we have built a van Kampen
diagram using one hub and $L$ trapezia corresponding to an accessible computation
for $W$.
\end{rk}

We will increase the set of relations of $G$ by adding
the (infinite) set of {\it disk relations} $V$ , one for
every disk word $V$. So we will consider diagrams
with \label{diskr} {\it disks}, where every disk cell (or just {\it disk}) is labeled by
such a word  $V$. (In particular, a hub is a disk.)

If two disks are connected by two $\tt$-bands and there are no other disks between these $\tt$-bands, then one can reduce the number of disks in the diagram. To achieve this,
it suffices to apply
the trick exploited for a pair of hubs in Subsection \ref{elim}.

\begin{df}\label{minimald}
We will call a reduced diagram $\Delta$ {\it minimal} if
\index[g]{minimal diagram over $G$}

(1) the number of disks is minimal for all diagrams with the same boundary label as $\Delta$ and

(2) $\Delta$ has minimal number of $(\theta,t)$-cells among
the diagrams with the same boundary label and with minimal number of  disks.

Clearly, a subdiagram of a minimal diagram is minimal itself.

\end{df}

Thus, no two disks of a minimal diagram are connected by  two $\tt$-bands,
such that the subdiagram bounded by them contains no other disks.
This property makes the disk graph of a reduced diagram
hyperbolic in the sense that the degree $L$ of every proper vertex (=disk) is high ($L\gg 1$) and there are no multiple edges.
Below we give a more precise formulation (proved for diagrams with such a disk graph, in particular,
in \cite{SBR}, Lemma 11.4 and in  \cite{O97}, Lemma 3.2).

\begin{figure}[ht]
\begin{center}
\includegraphics[width=0.7\textwidth]{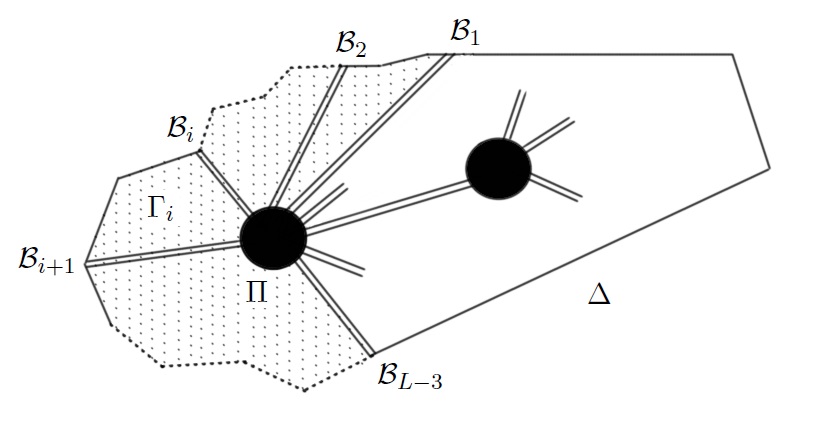}
\end{center}
\caption{Lemma \ref{extdisk}}\label{extd}
\end{figure}

\begin{lemma} \label{extdisk} If a minimal diagram  contains a least one disk,
then there is a disk $\Pi$ in $\Delta$ such that $L-3$ consecutive maximal $\tt$-bands ${\mathcal B}_1,\dots
{\mathcal B}_{L-3} $ start on $\partial\Pi$ , end on the boundary $\partial\Delta$, and for any $i\in [1,L-4]$,
there are no disks in the subdiagram $\Gamma_i$ bounded by ${\mathcal B}_i$, ${\mathcal B}_{i+1},$ $\partial\Pi,$ and $\partial\Delta$ (fig. \ref{extd}).
\end{lemma}

A maximal $q$-band starting on a disk of  a diagram is called a \index[g]{spoke}
{\it spoke}.



\subsubsection{The band moving transformation}

Recall the following band moving transformation for diagrams
with disks, exploited earlier in \cite{O97}, \cite{SBR}. Assume there is a disk $\Pi$ and a $\theta$-band $\mathcal T$ subsequently crossing some spokes ${\mathcal B}_1,\dots, {\mathcal B}_k$ which start (say, counter-clockwise) from $\Pi$.
Assume that $k\ge 2$ and there are no other cells between $\Pi$
and the bottom of $\mathcal T$, and so there is  a subdiagram $\Gamma$
formed by $\Pi$ and $\mathcal T$.

\begin{figure}[ht]
\begin{center}
\includegraphics[width=0.8\textwidth]{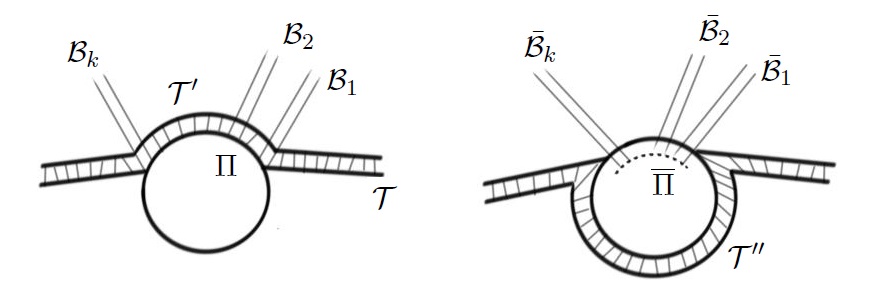}
\end{center}
\caption{The band moving transformation of a $\theta$-band and a disk}\label{tran}
\end{figure}

We describe the \index[g]{band moving transformation} {\it band moving transformation} (see, e.g., \cite{SBR}) as follows.
By Lemma \ref{perm} (1), for some $s$, we have a word $$V\equiv (\tilde t^{(s)} W)(\tilde t^{(s)}W)^{(+1)}\dots (\tilde t^{(s)}W)^{(+(k-2))}(\tilde t^{s})^{(+(k-1))}$$
(or  $V^{-1}\equiv (\tilde t^{(s)} W)(\tilde t^{(s)}W)^{(+1)}\dots (\tilde t^{(s)}W)^{(+(k-2))}   (\tilde t^{s})^{+(k-1))}$)
written on the top of the subband ${\mathcal T'}$ of
$\mathcal T$, that starts on ${\mathcal B}_1$ and ends on ${\mathcal B}_k$. (There are no superscripts in $V$ if $V$ is $\theta$-admissible word for a rule $\theta\in \Theta_3 - \Theta_5$.) The bottom ${\mathbf q}_2$ of $\mathcal T'$ is the subpath
of the boundary path ${\mathbf q}_2{\mathbf q}_3$ of $\Pi$ (fig. \ref{tran}), its label is a part of a disk word, and so is $V$ by Lemma \ref{perm}.


Therefore one can construct a new disk $\overline{\Pi}$
with boundary label $$(\tilde t^{(1)}W)(\tilde t^{(1)}W)^{(+1)}\dots (\tilde t^{(1)}W)^{(+(L-1))}$$ and boundary ${\mathbf s}_1{\mathbf s}_2$,
where $\Lab({\mathbf s}_1)\equiv V$. Also one construct
an auxiliary band ${\mathcal T}''$ with top label $$(W^{-1}(\tilde t^{(s)})^{-1})^{(+(L-1))}\dots (W^{-1}(\tilde t^{(s)})^{-1})^{(+k)}(W^{-1})^{(+(k-1))},$$  and attach it
to ${\mathbf s}^{-1}_2$, which has the same label. Finally
we replace the subband $\mathcal T'$ by ${\mathcal T''}$
(and make cancellations in the new $\theta$-band $\overline{\mathcal T}$ if any appear). The new diagram $\overline\Gamma$
formed by $\overline{\Pi}$ and $\overline{\mathcal T}$ has the same
boundary label as $\Gamma$.

\begin{rk}\label{less} After the band moving, the first
$(\theta,t)$-cells of $\tt$-spokes ${\mathcal B}_1,\dots, {\mathcal B}_k$ are removed and the total number of common $(\theta,t)$-cells of the new spokes $\overline{\mathcal B}_1,\dots, \overline{\mathcal B}_k$ of $\overline\Pi$ and $\overline{\mathcal T}$ is less
than the number of common $(\theta,t)$-cells of ${\mathcal B}_1,\dots, {\mathcal B}_k$ and $\mathcal T$ at least by $k$.
In particular, if $k>L-k$, then
the number of $(\theta,t)$-cells in $\overline\Gamma$ is less than the number of $(\theta,t)$-cells in $\Gamma$.
This observation implies
\end{rk}

\begin{lemma} \label{withd} Let $\Delta$ be a minimal diagram.

(1) Assume that a $\theta$-band ${\mathcal T}_0$ crosses $k$
$\tt$-spokes ${\mathcal B}_1,\dots, {\mathcal B}_k$ starting on a disk $\Pi$, and there are no
disks in the subdiagram $\Delta_0$, bounded by these spokes, by ${\mathcal T}_0$ and by $\Pi$. Then $k\le L/2$.

(2) Assume that there are two disjoint $\theta$-bands $\mathcal T$ and ${\mathcal S}$ whose bottom paths are parts of the boundary of a disk $\Pi$ and these bands correspond to the same rule $\theta$ (if their histories are read towards the disk), and $\theta\ne \theta(23)$.
Suppose $\mathcal T$ crosses $k\ge 2$  $\tt$-spokes starting
on $\partial\Pi$ and ${\mathcal S}$ crosses $\ell\ge 2$ $\tt$-spokes
starting on $\partial\Pi$. Then $k+\ell\le L/2$.

(3) $\Delta$ contains no $\theta$-annuli.

(4) A $\theta$-band
cannot cross a maximal $q$-band (in particular, a spoke) twice.
\end{lemma}
\proof
(1) Since every cell, except for disks, belongs to a
maximal $\theta$-band, it follows from Lemma \ref{NoAnnul} that there is a $\theta$-band $\mathcal T$
such that $\mathcal T$ crosses all ${\mathcal B}_1,\dots, {\mathcal B}_k$ and $\Delta_0$ has no cells between $\mathcal T$ and $\Pi$. If $k>L/2$, then by Remark \ref{less}, the band moving $\mathcal T$ around  $\Pi$
would decrease the number of $(\theta, t)$-cells in $\Delta$, a contradiction, since $\Delta$ is a minimal diagram.

(2) As above, let us move the band $\mathcal T$ aroud $\Pi$. This operation
removes $k$ $(\theta,t)$-cells but add $L-k$ new $(\theta,t)$-cells in $\overline{\mathcal T}$.
However $\ell$ $(\theta,t)$-cells of $\mathcal S$ and $\ell$ $(\theta,t)$-cells of $\overline{\mathcal T}$
will form mirror pairs, because for $\theta\ne\theta(23)$, the boundary label of a $(\theta, q)$-cell $\pi$, considered as
a $\theta$-band, is uniquely determined by the history $\theta$ and the label of the top
$q$-edge of $\pi$. So after cancellations one will have at most $L-k-2\ell$ new
$(\theta,t)$-cells. This number is less than $k$ if $k+\ell>L/2$ contrary to the minimality
of the original diagram. Therefore $k+\ell\le L/2$.

(3) Proving by contradiction, consider the
subdiagram $\Delta'$ bounded by a $\theta$-annulus.
It has to contain disks by Lemma \ref{NoAnnul}. Hence
it must contain spokes ${\mathcal B}_1,\dots, {\mathcal B}_{L-3}$ introduced in Lemma \ref{extdisk}. But this
contradits to item (1) of the lemma since $L-3>L/2$.

(4) The argument of item (3) works if there is
a subdiagram $\Delta'$ of $\Delta$ bounded by an $q$-band and a $\theta$-band.
\endproof

The band moving will be used for removing disks from quasi-trapezia.

\subsubsection{Quasi-trapezia}

\begin{df}A \index[g]{quasi-trapezium} {\it quasi-trapezium} is the same as trapezium (Definition \ref{dftrap}), but may contain disks. (So a quasi-trapezium without disks is a trapezium.)
\end{df}

\begin{lemma} \label{qt} Let a minimal diagram $\Gamma$ be a quasi-trapezium with standard factorization of the boundary as ${\mathbf p}_1^{-1}{\bf q}_1{\bf p}_2{\bf q}_2^{-1} $. Then there is a
diagram $\Gamma'$ such that

(1) the boundary of $\Gamma'$ is $$({\bf p}'_1)^{-1}{\bf q}'_1{\bf p}'_2({\bf q}'_2)^{-1},$$ where $$\Lab({\mathbf p'}_j)\equiv\Lab({\mathbf p}_j)$$ and $$\Lab({\mathbf q'}_j)\equiv\Lab({\mathbf q}_j)$$ for $j=1,2$;

(2) the numbers of hubs and
$(\theta, q)$-cells in $\Gamma'$ are the same as in $\Gamma$;

(3) the vertices $({\mathbf p'}_1)_-$ and $({\mathbf p'}_2)_-$ (the vertices $({\mathbf p}_1')_+$ and $({\mathbf p'}_2)_+$) are connected by
a simple path ${\mathbf s}_1$ (by ${\mathbf s}_2$, resp.)
such that we have three subdiagrams $\Gamma_1,\Gamma_2,\Gamma_3$ of $\Gamma'$,
where $\Gamma_2$ is a trapezium
with standard factorization of the boundary ${\bf p'}_1^{-1}{\bf s}_1{\bf p}'_2{\bf s}_2^{-1} $ and all cells
of the subdiagrams $\Gamma_1$ and $\Gamma_3$
with boundaries $\bf q'_1s_1^{-1}$ and $\bf s_2 (q'_2)^{-1}$ are disks;

(4) All maximal $\theta$-bands of $\Gamma$ and all
maximal $\theta$-bands of $\Gamma_2$ have the same number
of $(\theta,t)$-cells (equal for $\Gamma$ and $\Gamma_2$) .
\end{lemma}

\proof Every maximal $\theta$-band of $\Gamma$ must
connect an edge of $\bf p_1$ with an edge of ${\mathbf p}_2$; this follows from Lemma \ref{withd} (3). Hence
we can enumerate these bands from bottom to top:
${\mathcal T}_1,\dots,{\mathcal T}_h$, where $h=|{\mathbf p}_1|=|{\mathbf p}_2|$.

If $\Gamma$ has a disk, then by Lemma
\ref{extdisk}, there is a disk $\Pi$ such that at
least $L-3$  $\tt$-spokes of it end on ${\mathbf q}_1$ and ${\mathbf q}_2$, and there are no disks between the spokes
ending on ${\mathbf q}_1$ (and on ${\mathbf q}_2$). By Lemma \ref{withd} (2), at least $L-3 -L/2\ge 2$ of these spokes must end on  ${\mathbf q}_1$ (resp., on ${\mathbf q}_2$).

If $\Pi$ lies
between ${\mathcal T}_j$ and ${\mathcal T}_{j+1}$,
then the number of its $\tt$-spokes crossing
${\mathcal T}_j$ (crossing ${\mathcal T}_{j+1}$) is at least $2$.
So one can move each of the two $\theta$-bands around $\Pi$. So we can move the disk toward $\bf q_1$ (or toward $\bf q_2$) until the disk
is removed from the quasi-trapezium. (We use the property that
if $k$ $\tt$-spokes ${\mathcal B}_1,\dots,{\mathcal B}_k$ of $\Pi$ end on $\bf q_1$, then after moving bands  toward $\bf q_1$, we again have $k$
$\tt$-spokes $\overline{\mathcal B}_1,\dots,\overline{\mathcal B}_k$ of $\overline\Pi$ ending
on $\bf q_1$. - See the notation of Remark \ref{less}.)

No pair ${\mathcal T}_j$ and ${\mathcal T}_{j+1}$ corresponds to two mutual inverse letters $\theta\theta^{-1}$ of the history if $\theta\ne\theta(23)$. This follows from Lemma \ref{simul} (1) if there are no disks between these $\theta$-bands. If there is a disk, then this is
impossible too by Lemma \ref{qt} (2) since one could choose a disk
$\Pi$ as in the previous paragraph. So the projection of the
label of ${\mathbf p}_1$  on the history is eligible.

Let us choose $i$ such that the number $m$ of $(\theta,t)$-cells in ${\mathcal T}_i$ is minimal.
It follows that $\Gamma$ has at least $hm$ $(\theta,t)$-cells.

If the disk $\Pi$ lies above ${\mathcal T}_i$,
we will move it upwards using the band moving transformation. So after
a number of iterations all such (modified) disks
will be placed above the $\theta$-band number $h$
and form the subdiagram $\Gamma_1$. Similarly
we can form $\Gamma_3$ moving other disks downwards.

In the resulting diagram
$\Gamma_2$ lying between $\Gamma_1$ and $\Gamma_3$,
every $\theta$-band is reduced by the definition of band moving. The neighbor maximal $\theta$-band
of $\Gamma_2$ cannot be mirror copies of each other since the labels of ${\mathbf p}_1$ and
${\mathbf p}'_1$ are equal and $\Lab ({\mathbf p}_1) $ is a reduced word by Remark \ref{tb}.
It follows that the diagram $\Gamma_2$ (without disks)
is a reduced diagram, and so it is a trapezium of height
$h$.

The $\theta$-band ${\mathcal T}_i$ did not participate
in the series of band moving transformations above. Therefore it is a maximal $\theta$-band of $\Gamma_2$. Hence the trapezium $\Gamma_2$ contains exactly $mh$ $(\theta,t)$-cells, which does not
exceed the number of $(\theta,t)$-cells in $\Gamma$. In fact these two numbers are equal since $\Gamma$ is
a minimal diagram. So every maximal $\theta$-band
of $\Gamma$ and every maximal $\theta$-band of $\Gamma_2$ has $m$ $(\theta,t)$-cells.

\endproof

\subsubsection{Shafts}

We say that a history word $H$ is \index[g]{standard history} {\it standard} if there is a standard trapezium with history $H$.

\begin{df}\label{shf}
Suppose we have a disk $\Pi$ with boundary label $V$, $V^{\emptyset}\equiv(\tilde t W)^L$, and $\cal B$ be a $\tilde t$-spoke starting on $\Pi$.
Suppose there is a subband $\cal C$ of $\cal B$, which also starts on $\Pi$ and has a standard history $H$, for which the word $\tilde t W$ is $H$-admissible.
Then we call the $\tt$-band $\cal C$ a {\it shaft}.\index[g]{shaft}

For a constant \index[g]{parameters used in the paper!l@$\lambda$ - the parameter of $\lambda$-shafts (see Definition \ref{shf})}$\lambda\in [0;1/2)$ we also define
a stronger concept of \index[g]{shaft!l@$\lambda$-shaft} $\lambda$-shaft at $\Pi$ as follows.
A shaft $\ccc$ with history $H$ is a $\lambda$-{\em shaft} if
for every factorization of  the history $H\equiv H_1H_2H_3$, where $||H_1||+||H_3||<\lambda ||H||$, the middle part $H_2$ is still a standard history.
(So a shaft is a $0$-shaft).
\end{df}

\begin{lemma} \label{str} Let $\Pi$ be a disk in a minimal diagram
$\Delta$ and $\ccc$ be a $\lambda$-shaft at
$\Pi$ with history $H$. Then $\ccc$ has no factorizations  $\ccc=
\ccc_1\ccc_2\ccc_3$ such that

(a) the sum of lengths of $\ccc_1$ and $\ccc_3$ do not exceed $\lambda ||H||$ and

(b) $\Delta$ has a quasi-trapezium $\Gamma$
such that top (or bottom) label of $\Gamma$ has $L+1$ occurrences of $\tt$-letters and $\ccc_2$ starts on  the bottom and ends on the top of $\Gamma$.
\end{lemma}

\proof Proving by contradiction, we first replace $\Gamma$ by a trapezium $\Gamma'$ according to Lemma \ref{qt}. The transpositions used for this goal affect neither
$\Pi$ nor $\cal C$ since $\cal C$ crosses all the
maximal $\theta$-bands of $\Gamma$. Also one can replace $\Gamma'$ by a trapezium with shorter base and so we assume that the base of it starts and ends with letter $\tilde t$.

\begin{figure}[ht]
\begin{center}
\includegraphics[width=0.9\textwidth]{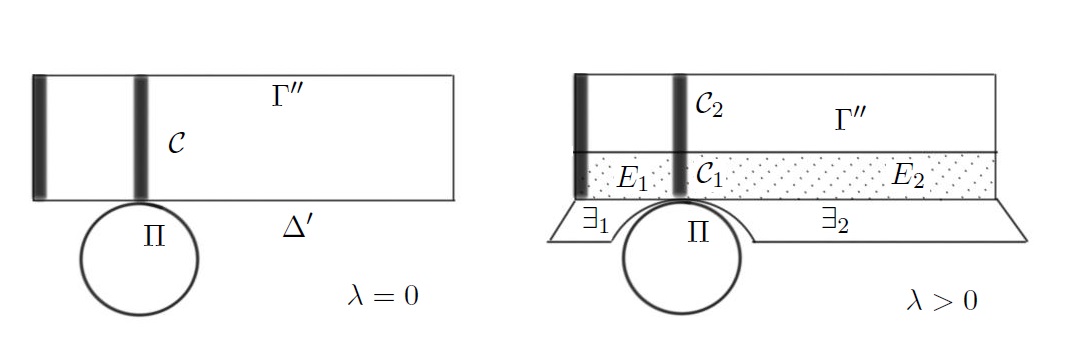}
\end{center}
\caption{Lemma \ref{str}.}\label{shft}
\end{figure}

For the beginning,  we assume that $\ccc$ is a shaft
(i.e.,$\lambda=0$). Then it follows from the definition
of shaft and Lemma \ref{resto} that ${\bf bot} (\Gamma')$ is labeled by a word $Vt$ such that
$V^{\emptyset}\equiv (tW)^L$, where the word $tW$ has standard base.
Now it follows from Remark \ref{perad} and Lemma \ref{simul}that $V$ is the boundary label of $\Pi$.
One can remove  the last maximal
$\tt$-band from $\Gamma'$ and obtain a subtrapezium
$\Gamma''$ whose bottom label coincides with
the label of $\partial\Pi$ (up to cyclic permutation),
and $\partial\Gamma''$ shares a $\tilde t$-edge with $\partial\Pi$ (fig.\ref{shft} with $\lambda=0$). It follows that
the subdiagram $\Delta'=\Pi\cup\Gamma''$ has boundary label freely equal to $\Lab({\bf top}(\Gamma''))$.
However $\Lab({\bf top}(\Gamma'')\equiv V'$,
where $(V')^{\emptyset}=V^{\emptyset}\cdot H$ by Lemma \ref{simul}, and so
there is a disk $\Pi'$ with boundary label $V'$.
Therefore the subdiagram $\Delta'$ can be replaced
by a single disk. So we decrease the number of $(\theta,t)$-cells contrary to the minimality
of $\Delta$.

Now we consider the general case, where ${\cal C}=
{\cal C}_1{\cal C}_2{\cal C}_3$. As above, we  replace
$\Gamma$ by a trapezium $\Gamma'$ and obtain a trapezium $\Gamma''$
after removing of one $\tt$-band in $\Gamma'$.
To obtain a contradiction, it suffices to consider
the diagram $\Delta'=\Pi\cup {\cal C}_1{\cal C}_2\cup\Gamma''$ (forgetting of the complement of $\Delta'$ in $\Delta$) and find another diagram $\Delta''$ with one disk and fewer $(\theta,t)$-cells such that
$\Lab(\partial\Delta'')=\Lab(\partial\Delta')$ in the free group.

Since both histories $H$ and $H_2$ (and so $H_1H_2$)
are standard, one can enlarge $\Gamma''$ and construct
a trapezium $\Gamma'''$ with history $H_1H_2$.
(The added parts $E_1$ and $E_2$  are dashed in figure \ref{shft} with $\lambda >0$).
Note that we add $<\lambda ||H||L$ new $(\theta,t)$-cells since every maximal $\theta$-band of $\Gamma'''$
has $L$ such cells. As in case $\lambda=0$, this trapezium $\Gamma'''$ and the disk $\Pi$ can be replaced
by one disk $\Pi'$. However to obtain the boundary label equal to $\Lab(\partial\Delta')$, we should
attach the mirror copies $\exists_1$ and $\exists_2$ of $E_1$ and $E_2$ to $\Pi'$.
The obtained diagram $\Delta''$ has at most
$\lambda||H_1||L$ $(\theta,t)$-cells, while $\Delta'$
has at least $||H_2||L\ge (1-\lambda)||H||$ $(\theta,t)$-cells. Since $\lambda<1-\lambda$,
we have the desired contradiction.
\endproof

\subsubsection{Designs}

As in \cite{O18}, we are going to use {\em designs}.

\begin{figure}[ht]
\begin{center}
\includegraphics[width=0.45\textwidth]{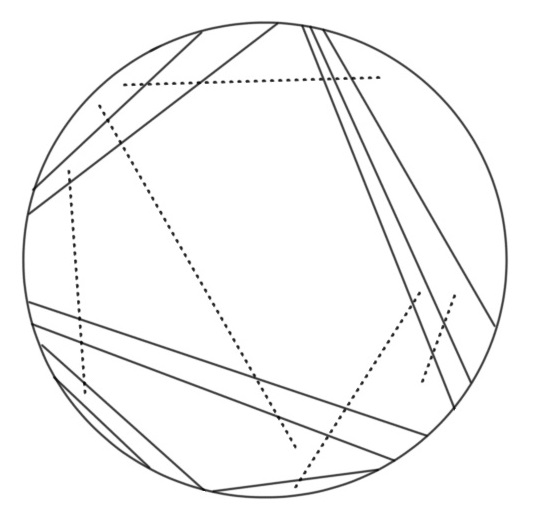}
\end{center}
\caption{Design}\label{desig}
\end{figure}

Let $\mathcal D$ be the Euclidean unit disk and $\mathbf T$ be a finite set of disjoint \index[g]{design!chords} {\em  chords} (solid lines in fig. \ref{desig}) and $\mathbf Q$ a finite set of disjoint simple curves in $\mathcal D$ (dotted lines in fig. \ref{desig}). We assume that a curve is a non-oriented broken line, i.e., it is built from finitely many finite line segments. To distinguish
the elements from $\mathbf T$ and $\mathbf Q$, we will say that the elements of $\mathbf Q$ are \index[g]{design!arcs} {\it arcs}.

We shall assume that the arcs belong to the open disk $D^o$, an arc may cross a  chord transversally  at most once, and the intersection point cannot coincide with one of the two ends of an arc.

Under these assumptions, we shall say that the pair
$(\mathbf T, Q)$ is a \index[g]{design} {\it design}. The number of elements in $\mathbf T$ and $\mathbf Q$ are denoted by $\#{\mathbf T}$ and $\#{\mathbf Q}$.

By definition, the \index[g]{design!length of an arc} {\it length} $|C|$ of an arc $C$
is the number of the chords crossing $C$. The term
{\it subarc} will be used in the natural way. Oviously
one has $|D|\le |C|$ if $D$ is a subarc of an arc $C$.

We say that an arc $C_1$ is \index[g]{design!parallel arcs}  {\it parallel} to an arc $C_2$
and write $C_1\parallel C_2$ if every chord (from $\bf T$) crossing $C_1$ also crosses $C_2$. So the relation
$\parallel$ is transitive (it is not necessarily symmetric). For example, the arc of length $2$ is parallel to
the arc of length $5$ in fig. \ref{desig}.

\begin{df} \label{propP} Given  $\lambda\in (0;1)$ and an integer
$n\ge 1$, the property $P(\lambda,n)$ of a design
says that for any $n$ different ars $C_1,\dots, C_n$,
there exist no subarcs $D_1,\dots, D_n$, respectively,
such that $|D_i|>(1-\lambda)|C_i|$ for every $i=1,\dots, n$ and $D_1\parallel D_2\parallel\dots\parallel D_n$.
\end{df}

By definition, the length $\ell(\bf Q)$ of the set
of arcs $\bf Q$ is defined by the equality
\begin{equation}
\ell({\bf Q})= \sum_{C\in \bf Q} |C|
\end{equation}

The number of chords will be denoted by $\#\bf T$. Here is the main statement about designs from \cite{O18}.

\begin{theorem}[Theorem 8.2 \cite{O18}] \label{design} There is a constant $c=c(\lambda,n)$
such that for any design $(\bf T,Q)$ with property
$P(\lambda,n)$, we have

\begin{equation}\label{QT}
\ell({\bf Q})\le c(\#\bf T)
\end{equation}
\end{theorem}

\subsubsection{Designs and the $\sigma_\lambda$ invariant}

Let $\lambda\in [0,1/2)$. For every
$\tt$-spoke $\mathcal B$ of a minimal diagram $\Delta$, we choose the
$\lambda$-shaft of maximal length in it (if $\mathcal B$ contains a $\lambda$-shaft). It starts on the boundary of a disk, and it is a unique maximal shaft in $\mathcal B$ if the spoke connects
the disk and the boundary $\partial\Delta$. If $\mathcal B$
connects two disks $\Pi_1$ and $\Pi_2$, then there can
be two maximal $\lambda$-shafts: at $\Pi_1$ and  at $\Pi_2$.
We denote by \index[g]{s@$\sigma_{\lambda}(\Delta)$ - the $\sigma_\lambda$-invraiant of a diagram} $\sigma_{\lambda}(\Delta)$ the sum of lengths of all maximal $\lambda$-shafts in the spokes of $\Delta$.

\begin{lemma} \label{clam} There is a constant $c=c(\lambda)$ such that
$\sigma_{\lambda}(\Delta)\le c |\partial\Delta|$ for
every minimal diagram $\Delta$ over the group $G$.
\end{lemma}

\proof Let us associate the following design with $\Delta$. We say that the median lines of the maximal $\theta$-bands are the chords and the median
lines of the maximal $\lambda$-shafts are the arcs.
Here we use two disjoint median lines if two maximal
$\lambda$-shafts share a $(\theta,\tt)$-cell. By Lemma \ref{withd}
(3), (4), we indeed obtain a design.

Observe that the length $|C|$ of an arc is the
number of cells in the $\lambda$-shaft and $\#{\bf T}\le |\partial\Delta|/2$ since every maximal $\theta$-band
has two $\theta$-edges on $\partial\Delta$.

Thus, by Theorem \ref{design},  it suffices to show
that the constructed design satisfies the condition
$P(\lambda,n)$, where $n$ does not depend on $\Delta$.

Let $n=2L+1$. If the property $P(\lambda,n)$ does not hold, then we have $n$ maximal $\lambda$-shafts
${\mathcal C}_1,\dots, {\mathcal C}_n$ and a subband $\mathcal D$ of
${\mathcal C}_1$, such that $|{\cal D}|>(1-\lambda)|{\cal C}_1|$, and every maximal $\theta$-band crossing $\cal D$ must cross each of
${\mathcal C}_2,\dots,{\mathcal C}_n$. (Here $|\mathcal B|$ is the length of a $\tt$-band $\cal B$.) It follows that  each of these $\theta$-band crosses at least $L+1$
maximal $\tt$-bands. (See Lemma \ref{withd} (3,4). Here we take
into account that the same $\tt$-spoke can generate two
arcs in the design.) Hence using the $\lambda$-shaft
${\mathcal C}_1$ one can construct a quasi-trapezium of height
$|\mathcal D|$, which contradicts Lemma
\ref{str}.
\endproof

\subsection{Upper bound for \texorpdfstring{$G$}{G}-areas of diagrams over the group \texorpdfstring{$G$}{G}. }\label{ub}

\subsubsection{The area of a disk is quadratic}
By definition, the \label{areaGdi} $G$-{\it area of a disk} $\Pi$ is just the minimum of areas of the diagrams over the presentation (\ref{rel1}) through (\ref{rel3}) of $G$ having the same
boundary label as $\Pi$.

\begin{lemma} \label{disk} There is a constant \index[g]{parameters used in the paper!c@$c_6 $ - the parameter controlling the area of a disk in terms of its perimeter (see Lemma \ref{disk})}$c_6$ such that both  area and the $G$-area of any disk does not exceed
$c_6|\partial\Pi|^2$.
\end{lemma}

\proof By Remark \ref{dr}, a disk with boundary label $V$ can be built of one hub and $L$ trapezia corresponding to an accessible computation $\ccc$
for $W$, where $W^L\equiv V^{\emptyset}$.  By Lemma \ref{fea}, the length of $\ccc$ can be bound by $c_2||W||$ and the length of every configuration of $\ccc$ does not exceed $c_1||W||$
Hence by Lemma \ref{ochev}, the area and the $G$-area of the disk is bounded by $c_6|\partial\Pi|^2$ since
the constant $c_6$ can be chosen after $c_1, c_2$ and $\delta$.
\endproof

By definition, the \index[g]{G@$G$-area} $G$-area of a minimal diagram $\Delta$ over $G$ is the sum of $G$-areas of its disks plus the $G$-area of the complement $\Gamma$. For the complement, as in subsection \ref{qub}, we consider a family $\bf S$ of big subtrapezia and single
cells of $\Gamma$ such that every cell of $\Gamma$ belongs to
a member $\Sigma$ of this family, and if a cell $\Pi$ belongs to different $\Sigma_1$ and $\Sigma_2$ from $\bf S$, then both $\Sigma_1$ and $\Sigma_2$ are big subtrapezia of $\Gamma$ with bases $xv_1x$, $xv_2x$, and $\Pi$ is an $(\theta,x)$-cell.)
Hence the statement of Lemma \ref{GA} holds for minimal diagrams
over $G$ as well.

\subsubsection{Weakly minimal diagrams.}

We want to prove that for big enough constant $N$, $\area_G(\Delta)\le N n^2$  for every minimal diagram $\Delta$, which will imply in Subsection \ref{end1} that the boundary label  of $\Delta$
has quadratic area with respect to the finite presentation of $G$. However to prove this property by induction, we have to consider a larger class of diagrams, called {\it weakly minimal} .

Let $\cal C$ be a cutting $q$-band of a reduced diagram $\Delta$ with disks, i.e. it starts and ends on $\partial\Delta$ and cut up the diagram. We call $\cal C$ a \index[g]{stem band}{\it stem} band, if it either
is a rim band of $\Delta$ or both components of $\Delta\backslash \cal C$ contain disks.
The (unique) maximal subdiagram  of $\Delta$, where every cutting $q$-band is a stem, is called
\index[g]{stem of a van Kampen diagram}{\it the stem}  $\Delta^*$ of $\Delta$. It is obtained by removing all \index[g]{crown cell}{\it crown} cells from
$\Delta$,  where a cell $\pi$ is called {\it crown}, if it belongs to the component $\Gamma$ defined by a cutting $q$-band $\cal B$, where $\Gamma$ contains no disks and $\pi$ is not in $\cal B$. In particular, all the disks and $q$-spokes of $\Delta$ belong to the stem $\Delta^*$.
The stem of a diagram without disks is empty.

\begin{df} \label{wminimald} A reduced diagram $\Delta$ (with disks) is called {\em weakly minimal}
\index[g]{weakly minimal diagram}
if the stem
$\Delta^{*}$ is a minimal diagram \footnote{Unfortunately this definition was missed in \cite{O18}; it can be found in the arXiv version of that paper.}.
\end{df}

\begin{lemma} \label{wmin} (a) If $\Delta_1$ is a subdiagram of weakly minimal diagram $\Delta$,
then $\Delta_1$ is weakly minimal and $\Delta_1^*\subset\Delta^*$;

(b) under the same assumption, we have $\sigma_{\lambda}(\Delta_1^*)\le \sigma_{\lambda}(\Delta^*)$;

(c) There is a constant $c=c(\lambda)$ such that
$\sigma_{\lambda}(\Delta^*)\le c |\partial\Delta|$ for
every weakly minimal diagram $\Delta$ over the group $G$;

(d) If a diagram $\Delta$ has a cutting $q$-band $\cal C$ and two components
$\Delta_1$ and $\Delta_2$ of the complement of $\cal C$ such that $\Delta_1\cup\cal C$
is a reduced diagram without disks and ${\cal C}\cup \Delta_2$ is a weakly minimal
diagram, then $\Delta$ is weakly minimal itself;

(e) a weakly minimal diagram $\Delta$ contains no $\theta$-annuli, and.
a $\theta$-band cannot cross a $q$-band of $\Delta$ twice.

\end{lemma}

\proof (a) Every crown cell $\pi$ of $\Delta$ belonging is $\Delta_1$ is crown in $\Delta_1$
since the cutting $q$-band $\cal B$ separating $\pi$ from all the disks of $\Delta$
separates (itself or the subbands of $\cal B$ in the intersection of $\cal B$ and  $\Delta_1$)
$\pi$ from $\Delta_1^*$. Therefore we have $\Delta_1^*\subset\Delta^*$, and so  $\Delta_1^*$
is minimal being a subdiagram of a minimal diagram.

(b) Now it follows from the definition of shaft, that every $\lambda$-shaft of $\Delta_1^*$ is
a $\lambda$-shaft in $\Delta^*$, which implies  inequality $\sigma_{\lambda}(\Delta_1^*)\le \sigma_{\lambda}(\Delta^*)$.

(c) If a cutting  $q$-band $\cal C$ of a reduced diagram $\Delta$ gives a decomposition $\Delta=\Gamma_1\cup{\cal C}\cup \Gamma_2$, where $\Delta_1=\Gamma_1\cup{\cal C}$ has no disks,
then every maximal $\theta$-band starting in the subdiagram $\Delta_1$ with $\cal C$
cannot ends on $\partial\Gamma_1$ by Lemma \ref{NoAnnul}. Hence $|\partial\Delta_2|\le|\partial\Delta|$ by Lemma \ref{ochev}. So removing subdiagrams as
$\Gamma_1$ from $\Delta$, we obtain by induction that $|\partial\Delta^*|\le
|\partial\Delta|$. Now the property (c) follows from Lemma  \ref{clam} applied to the minimal
subdiagram $\Delta^*$.

(d) The diagram $\Delta$ is reduced since both $\Delta_1\cup\cal C$ and $\Delta_2\cup\cal C$
are reduced subdiagrams sharing the cutting band $\cal C$. Since $\Delta_1$ has no disks,
we have $\Delta^*=(\Delta_2\cup\cal C)^*$ by the definition of stem. Therefore
the stem $\Delta^*$ is a minimal diagram and $\Delta$ is weakly minimal.

(e) The statement follows from Lemma \ref{withd} (3, 4) if the bands belong to the
stem $\Delta^*$. By the same reason, a $\theta$-band cannot cross a rim $q$-band
of $\Delta^*$ twice. It remains to assume that the bands belong to the crown of $\Delta$,
and in this case, the statement follows from Lemma \ref{NoAnnul} since the crown
is a union of disjoint reduced subdiagrams over the group $M$.
\endproof

\begin{rk} The statement (d) of Lemma \ref{wmin} fails if one replaces the words
``weakly minimal'' with ``minimal''.
\end{rk}

We  will prove that for large enough parameters \index[g]{parameters used in the paper!n@$N_3$,$N_4$ - parameters controlling the $G$-area of a diagram $\Delta$ over $G$ in terms of the perimeter, the mixture, and $\sigma_{\lambda}(\Delta^{*})$ where $\Delta^{*}$  is the stem of $\Delta$} $N_3$
and $N_4$, $\area_G(\Delta)\le N_4 (n+\sigma_{\lambda}(\Delta^{*}))^2
+N_3\mu(\Delta)$ for every weakly minimal diagram $\Delta$
with perimeter $n$. For this aim, we will argue by contradiction in this section and study a
weakly minimal {\bf counter-example} $\Delta$ satisfying  the opposite inequality
\begin{equation}\label{ce}
\area_G(\Delta)> N_4 (n+\sigma_{\lambda}(\Delta^{*}))^2
+N_3\mu(\Delta)
\end{equation}
with minimal possible sum $n+\sigma_{\lambda}(\Delta^*)$.

\subsubsection{Getting rid of rim bands with short base}

\begin{lemma} \label{norim} The diagram $\Delta$ has no rim
$\theta$-bands with base of length at most $K$.
\end{lemma}

\proof  The proof of Lemma
\ref{nori} works for the weakly minimal counter-example over $G$. It suffices
to replace $N_2$ and $N_1$ with $N_4$ and $N_3$, resp., replace $n$ with $n+\sigma_{\lambda}(\Delta^*)$, and notice that the subdiagram $(\Delta')^*$ is weakly minimal and $\sigma_{\lambda}((\Delta')^*)\le \sigma_{\lambda}(\Delta^*)$ by Lemma \ref{wmin}
(a,b).
\endproof

\subsubsection{The cloves}

By Lemma \ref{main},
$\Delta$ has at least one disk. Taking into account that all disks and their spokes belong
 to the stem $\Delta^*$, we can apply Lemma \ref{extdisk} to the weakly minimal diagram $\Delta^*$ and fix a disk $\Pi$ in $\Delta$ such that $L-3$ consecutive maximal $\tt$-bands ${\mathcal B}_1,\dots
{\mathcal B}_{L-3} $ start on $\partial\Pi$, end on the boundary $\partial\Delta$ , and for any $i\in [1,L-4]$,
there are no disks in the subdiagram  bounded by ${\mathcal B}_i$, ${\mathcal B}_{i+1},$ $\partial\Pi,$ and $\partial\Delta.$
(See fig. \ref{extd}.)

\index[g]{clove $\Psi=cl(\pi,{\mathcal B}_1,{\mathcal B}_{L-3})$ of the minimal counterexample from Section \ref{midi}}
We denote by
 $\Psi=cl(\Pi,{\mathcal B}_1,{\mathcal B}_{L-3})$ the subdiagram without disks bounded by the spokes ${\mathcal B}_1$, ${\mathcal B}_{L-3}$
(and including them) and by subpaths of the boundaries of $\Delta$
and $\Pi,$ and call this subdiagram a \index[g]{clove} {\it clove}. Similarly one can define the cloves
$\Psi_{ij}=cl(\Pi,{\mathcal B}_i,{\mathcal B}_j)$ if $1\le i<j\le L-3$.

\subsubsection{A clove cannot contain "wide" subcombs}

Below we use the following
analog of Lemma \ref{notwo} (with identical proof):

\begin{lemma} \label{notw} (1) The counter-example $\Delta$ has no two disjoint subcombs $\Gamma_1$ and $\Gamma_2$ in $\Psi$ with basic widths at most $K$ and handles ${\cal C}_1$ and ${\cal C}_2$ such that some ends of these handles are
connected by a subpath ${\bf x}$ of the boundary path of $\Delta$ with $|{\bf x}|_q\le N$.

(2) The boundary of every subcomb  $\Gamma$ of $\Delta$ with basic width $s\le K$ has $2s$ $q$-edges provided $\Gamma\subset \Psi$.
\end{lemma}
$\Box$

\begin{lemma} \label{nocomb} The clove $\Psi=cl(\Pi,{\mathcal B}_1,{\mathcal B}_{L-3})$ has no subcombs of basic width at least $K_0$.
\end{lemma}

\proof The proof is similar to the proof of Lemma \ref{main}. Proving by contradiction, we may assume that there is a tight
subcomb $\Gamma$ by Lemma \ref{est'takaya} (2). Then we can use Lemma \ref{notw} (which is the analog of Lemma \ref{notwo}) and can  repeat the proofs of the statements of Lemmas \ref{twocombs} - \ref{main} to obtain a contradiction with the minimality of the counter-example $\Delta$. Some modifications are needed in the proof of Lemma \ref{main} only.
Namely, considering the weakly minimal diagram $\Delta$ over the presentation of $G$ and the subcomb $\Gamma$, we should now
replace $N_2$ and $N_1$ with $N_4$ and $N_3$, replace $n$ with $n+\sigma_{\lambda}(\Delta^*)$, and notice that the value of $\sigma_{\lambda}$ does nor increase when we pass from $\Delta$ to a subdiagram by Lemma \ref{wmin} (b). We should use Lemma \ref{wmin} (e) instead of Lemma
\ref{NoAnnul} used in the proofs of Lemmas \ref{notwo} - \ref{A0}. The diagram $\Delta_0$
is weakly minimal because it is constructed from the reduced diagram $\Gamma_1^{(+k)}\cup\cal Q$ over $M$ and the weakly minimal diagram $\Delta_1\cup\cal Q$ according to the assumption of Lemma \ref{wmin} (d).
\endproof




\subsubsection{$\theta$-bands in a clove}

\begin{lemma}\label{psi1}
(1) Every maximal
$\theta$-band of $\Psi$ crosses either ${\mathcal B}_1$ or ${\mathcal B}_{L-1}$.
\index[g]{clove $\Psi=cl(\pi,{\mathcal B}_1,{\mathcal B}_{L-3})$ of the minimal counterexample from Section \ref{midi}!r@$r$: the $\theta$-bands of $\Psi$
crossing ${\mathcal B}_{L-3}$
do not cross ${\mathcal B}_r$, and the $\theta$-bands of $\Psi$ crossing ${\mathcal B}_1$
do not cross ${\mathcal B}_{r+1}$ }

(2) There exists $r$, $L/2-3\le r \le L/2$,   such that the $\theta$-bands of $\Psi$
crossing ${\mathcal B}_{L-3}$
do not cross ${\mathcal B}_r$, and the $\theta$-bands of $\Psi$ crossing ${\mathcal B}_1$
do not cross ${\mathcal B}_{r+1}$.
\end{lemma}

{\proof (1) If the claim were wrong, then one could find a rim $\theta$-band $\mathcal T$ in $\Psi$, which
crosses neither ${\mathcal B}_1$ nor ${\mathcal B}_{L-3}$. By Lemma
\ref{norim}, the basic width of $\mathcal T$ is greater than $K$.
Since (1) a disk has $LN$ spokes, (2) no $q$-band of $\Psi$ intersects $\mathcal T$ twice by Lemma
  \ref{NoAnnul}, (3) $\mathcal T$ has at least $K$  $q$-cells, and (4) $K>2K_0+LN$, there exists a maximal $q$-band $\mathcal C'$
 such that a subdiagram $\Gamma'$ separated from $\Psi$ by $\mathcal C'$ contains no edges of the spokes
  of $\Pi$ and the part of $\mathcal T$ belonging to $\Gamma'$ has at least $K_0$  $q$-cells (fig. \ref{lempsi}).

\begin{figure}[ht]
\begin{center}
\includegraphics[width=1.0\textwidth]{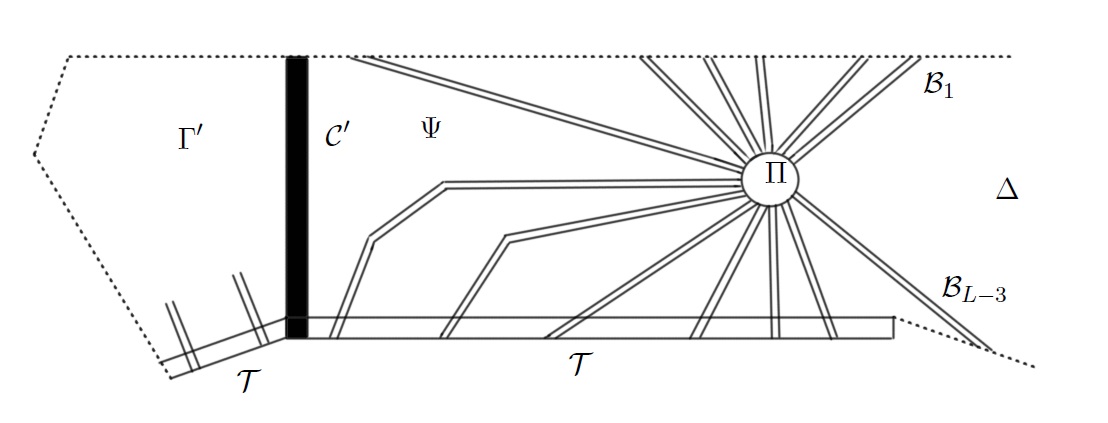}
\end{center}
\caption{Lemma \ref{psi1}}\label{lempsi}
\end{figure}

  If $\Gamma'$ is not a comb, and so a maximal $\theta$-band of it does not cross $\mathcal C',$
  then $\Gamma'$ must contain another rim band $\mathcal T'$ having at least $K$
  $q$-cells. This makes possible to find a subdiagram $\Gamma''$ of
  $\Gamma'$ such that a part of $\mathcal T'$ is a rim band of $\Gamma''$ containing at least $K_0$
  $q$-cells, and $\Gamma''$ does not contain $\mathcal C'$.
  Since $\area(\Gamma')>\area(\Gamma'')>\dots$ , such a procedure must stop. Hence, for some $i$, we
  obtain a subcomb $\Gamma^{(i)}$ of basic width $\ge K_0$,
  contrary to Lemma \ref{nocomb}.

  (2) Assume there is a maximal $\theta$-band $\mathcal T$ of $\Psi$ crossing
  the spoke ${\mathcal B}_1$.   Then assume that $\mathcal T$ is the closest
  to the disk $\Pi$, i.e. the intersection of $\mathcal T$ and ${\mathcal B}_1$ is
  the first cell of the spoke ${\mathcal B}_1$. If ${\cal B}_1,\dots, {\cal B}_r$ are all the spokes crossed by $\cal T$, then $r\le L/2$ by Lemma \ref{width}, which is applicable
  here since all the spokes belong to the stem $\Delta^*$, which is a minimal diagram.
  Since the band $\mathcal T$ does not cross the spoke ${\mathcal B}_{r+1}$, no other
  $\theta$-band of $\Psi$ crossing ${\mathcal B}_1$ can cross ${\mathcal B}_{r+1}$.
  and no $\theta$-band crossing the spoke ${\mathcal B}_{L-3}$ can cross ${\mathcal B}_r$.
  The same argument shows that $r+1\ge L/2 -2$ if there is a $\theta$-band of $\Psi$
  crossing the spoke ${\mathcal B}_{L-3}$.
  \endproof}

\index[g]{clove $\Psi=cl(\pi,{\mathcal B}_1,{\mathcal B}_{L-3})$ of the minimal counterexample from Section \ref{midi}!p@${\mathbf p}(\Psi)$ the common subpath of $\partial\Psi$ and
$\partial\Delta$ starting with the $\tt$-edge of ${\mathcal B}_1$ and ending
with the $\tt$-edge of ${\mathcal B}_{L-3}$}
  For the clove $\Psi=cl(\pi,{\mathcal B}_1,{\mathcal B}_{L-3})$ in  $\Delta$,
we denote by ${\mathbf p}(\Psi)$ the common subpath of $\partial\Psi$ and
$\partial\Delta$ starting with the $\tt$-edge of ${\mathcal B}_1$ and ending
with the $\tt$-edge of ${\mathcal B}_{L-3}.$ Similarly we define the (outer)
path ${\mathbf p}_{ij}={\mathbf p}(\Psi_{ij})$ for every smaller clove $\Psi_{ij}$.
\index[g]{clove $\Psi=cl(\pi,{\mathcal B}_1,{\mathcal B}_{L-3})$ of the minimal counterexample from Section \ref{midi}!p@${\mathbf p}_{ij}(\Psi)$ the common subpath of  $\partial\Psi$ and
$\partial\Delta$ starting with the $\tt$-edge of ${\mathcal B}_i$ and ending
with the $\tt$-edge of ${\mathcal B}_{j}$}

\subsubsection{The clove $\Psi$ and related subdiagrams.}

  \begin{lemma} \label{2K0} Every path ${\mathbf p}_{i,i+1}$ ($i=1,\dots, L-4$) has fewer than $3K_0$ $q$-edges.
\end{lemma}

{ \proof Let a maximal $q$-band $\mathcal C$ of $\Psi$ start on ${\bf p}_{i,i+1}$ and suppose it does not end on $\Pi$. Then is has to end on ${\bf p}_{i,i+1}$ too.
If $\Gamma$ is the subdiagram
(without disks) separated by $\mathcal C$, then every maximal $\theta$-band $\cal T$ of $\Gamma$ has to
cross the $q$-band $\mathcal C$ since the extension of $\cal T$ in $\Psi$ must cross either ${\mathcal B}_1$
or ${\mathcal B}_{L-3}$ by Lemma \ref{psi1}. Therefore $\Gamma$ is a comb with handle $\mathcal C$.

Consider the $q$-bands of this kind defining maximal subcombs $\Gamma_1,\Gamma_2,\dots\Gamma_k$
in $\Psi_{i,i+1}$. The basic width of each of them is smaller than $K_0$ by Lemma \ref{nocomb}.
Therefore $k\le 1$ since otherwise one can get two subcombs contradicting to Lemma \ref{notw} (1), because there are at most $N+1$ maximal $q$-bands
starting on $\partial\Pi$ in $\Psi_{i,i+1}$.
By Lemma \ref{notw} (2),
such a subcomb has at most $2K_0$ $q$-edges in the boundary. Hence there are at most $2K_0+N<3K_0$ $q$-edges in the path ${\bf p}_{i,i+1}$.
\endproof}

\begin{figure}[ht]
\begin{center}
\includegraphics[width=0.8\textwidth]{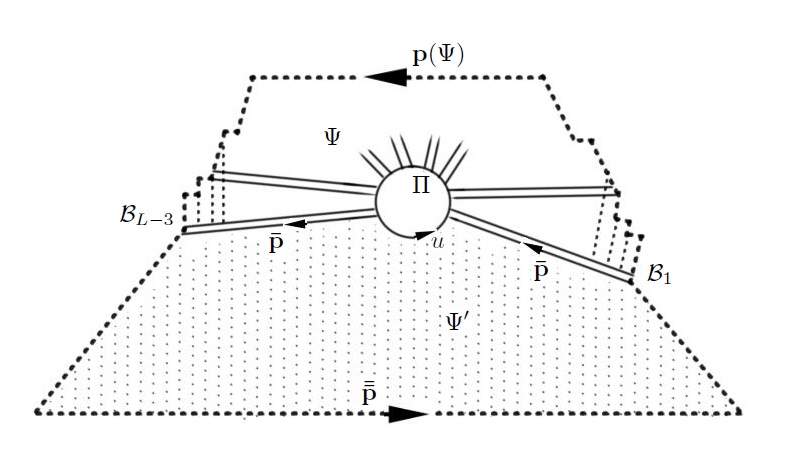}
\end{center}
\caption{Boundaries of $\Psi$ and $\Psi'$}\label{bou}
\end{figure}

\index[g]{clove $\Psi=cl(\pi,{\mathcal B}_1,{\mathcal B}_{L-3})$ of the minimal counterexample from Section \ref{midi}!d@$\overline\Delta$ the subdiagram formed by $\Pi$ and $\Psi$}
We denote
by $\overline\Delta$ the subdiagram formed by $\Pi$ and $\Psi$, and
denote by $\bf \overline p$ the path $\topp({\mathcal B}_1) {\mathbf u}^{-1}
\bott({\mathcal B}_{L-3})^{-1},$ where $\bf u$ is a subpath of $\partial\Pi,$ such that
$\bf \overline p$ separates $\overline\Delta$ from the remaining subdiagram $\Psi'$
of $\Delta$ (fig. \ref{bou}).

\index[g]{clove $\Psi=cl(\pi,{\mathcal B}_1,{\mathcal B}_{L-3})$ of the minimal counterexample from Section \ref{midi}!d@$\overline\Delta_{ij}$}
\index[g]{clove $\Psi=cl(\pi,{\mathcal B}_1,{\mathcal B}_{L-3})$ of the minimal counterexample from Section \ref{midi}!p@${\bf\overline p}_{i,j}=\topp({\mathcal B}_i) {\mathbf u}_{ij}^{-1}
\bott({\mathcal B})_{j}^{-1},$ where ${\mathbf u}_{ij}$ is a subpath of $\partial\Pi$}
Similarly we define subdiagrams  $\overline\Delta_{ij}$, paths
${\bf\overline p}_{i,j}={\topp}({\mathcal B}_i) {\mathbf u}_{ij}^{-1}
\bott({\mathcal B})_{j}^{-1},$ where ${\mathbf u}_{ij}$ is a subpath of $\partial\Pi$, and
the subdiagram $\Psi'_{ij}$.

We denote by $H_1,\dots, H_{L-3}$ the  histories  of the spokes ${\mathcal B}_1,\dots,
{\mathcal B}_{L-3}$ (read starting from the disk $\Pi$) and by $h_1,\dots,h_{L-3}$ their lengths, i.e., the numbers of $(\theta,t)$-cells. By Lemma \ref{psi1},
these lengths non-increase and then non-decrease as follows:

\begin{equation}\label{downup}
h_1\ge h_2\ge\dots \ge h_r;\;\; h_{r+1}\le\dots \le h_{L-3}\;\;(L/2-3\le r \le L/2),
\end{equation}
and therefore $H_{i+1}$ is a prefix of $H_{i}$ ($H_j$ is a prefix $H_{j+1}$) for $i=1,\dots, r-1$ (resp., for $j=r+1,\dots, L-4$).

Recall that by Definition \ref{dw} the boundary label of $\partial\Pi$ is a disk word $V$, where $V^{\emptyset}\equiv W^L$
and $W$ is an accessible word.

\begin{lemma} \label{ppbar} We have the following inequalities
$$ |\mathbf{\overline p}_{ij}|\le h_i+h_j+ (L-j+i)|W|-1$$
and, if $i\le r$ and $j\ge r+1$, then $$|{\mathbf p}_{ij}|\ge |{\mathbf p}_{ij}|_{\theta}+ |{\mathbf p}_{ij}|_q\ge h_i+h_j+(j-i)N+1$$
\end{lemma}

{\proof The first inequality follows from Lemma \ref{ochev} (b) since the
path ${\bf u}_{ij}$ has $L-j+i-1$ $\tilde t$-edges. To prove the second inequality,
we observe that the path ${\bf p}_{ij}$ has $(j-i)N+1$ $q$-edges and it has
$h_i+h_j$ $\theta$-edges by Lemma \ref{psi1}.
\endproof}

\begin{lemma} \label{muJ} If $j-i>L/2$, then we have $$\mu(\Delta)-\mu(\Psi'_{ij})>
-2Jn(h_i+h_j)\ge -2Jn|{\mathbf p}_{ij}|$$
\end{lemma}

{\proof
The number of $q$-edges in the path ${\bf\bar p}_{ij}$ (or in the path ${\bf u}_{ij}$)
does not exceed the similar number for ${\bf p}_{ij}$ provided $j-i\ge L/2$. Therefore any two white beads $o, o'$ of the necklace
on $\partial\Delta$, provided they both  do not belong to ${\bf p}_{ij},$ are separated by at least the same
number of black beads in the necklace for $\Delta$ as in the necklace for $\Psi'_{ij}$  (either the clockwise arc $o-o'$ includes ${\bf p}_{ij}$ or not). So such a pair  contributes
to $\mu(\Delta)$ at least the amount it contributes to $\mu(\Psi'_{ij})$.
Thus, to estimate $\mu(\Delta)- \mu(\Psi'_{ij})$ from below, it suffices to consider
the contribution to $\mu(\Psi')$ for the pairs $o, o', $ where one of the two beads
lies on ${\bf p}_{ij}$. The number of such (unordered) pairs is bounded by $n(h_{i}+h_{j})$,
 because it follows from Lemma \ref{psi1} (1) that every maximal $\theta$-band starting on ${\bf p}_{ij}$ has to cross either ${\cal B}_i$ or ${\cal B}_j$, i.e. $|{\bf p}_{ij}|_{\theta}\le h_i+h_j$. Taking into account the definition of $\mu $ for diagrams and inequalities (\ref{downup}), we get the required statement.
\endproof}

\begin{lemma} \label{epsi} If $j-i>L/2$, then  the following inequality holds: $|{\mathbf p}_{ij}|<(1+\varepsilon)|\mathbf{\overline p}_{ij}|$, where $\varepsilon=N_4^{-\frac 12}$.
Moreover, we have $|{\mathbf p}_{ij}|+\sigma_{\lambda}(\overline\Delta_{ij}^*)<(1+\varepsilon)|\mathbf{\overline p}_{ij}|$.
\end{lemma}

{
\proof It suffices to prove the second statement. Let $d$ be the difference $$|{\bf p}_{ij}|+\sigma_{\lambda}(\bar\Delta_{ij}^*)-|{\bf \bar p}_{ij}|$$ and
assume, by contradiction, that $d\ge\varepsilon |{\bf \bar p}_{ij}|$. Then $$d\ge
|{\bf p}_{ij}|+\sigma_{\lambda}(\bar\Delta_{ij}^*)-\varepsilon^{-1}d, $$
whence
\begin{equation}\label{de}
d\ge (1+\varepsilon^{-1})^{-1}(|{\bf p}_{ij}|+\sigma_{\lambda}(\bar\Delta_{ij}^*))
\ge\frac{\varepsilon}{2} (|{\bf p}_{ij}|+\sigma_{\lambda}(\bar\Delta_{ij}^*))\ge \frac{\varepsilon y}{2},
\end{equation}
where by definition, $y=|{\bf p}_{ij}|+\sigma_{\lambda}(\bar\Delta_{ij}^*)$.

We have
\begin{equation}\label{d}
(|\partial\Delta|+\sigma_{\lambda}(\Delta^*))-(|\partial\Psi'_{ij}|+\sigma_{\lambda}((\Psi'_{ij})^*))\ge d >0,
\end{equation}
because $$|\partial\Delta|-|\partial\Psi'_{ij}|\ge |{\bf p}_{ij}|-|{\bf \bar p}_{ij}|$$
and by Lemma \ref{wmin} (a) $$\sigma_{\lambda}(\bar\Delta_{ij}^*)+\sigma_{\lambda}((\bar\Psi'_{ij})^*)\le \sigma_{\lambda}(\Delta^*)$$ since $\Psi'_{ij}$ and $\bar\Delta_{ij}$ have no common spokes. Therefore  for $$x=n+\sigma_{\lambda}(\Delta^*),$$ we obtain from  the weak minimality of the counter-example $\Delta$, that  $\Psi'_{ij}$ is not a counter-example. Hence using inequality (\ref{d}), we obtain

$$\area_G(\Psi'_{ij})\le N_4 (x-d)^2+N_3\mu(\Psi'_{ij})\le N_4 x^2- N_4 x d +N_3\mu(\Psi'_{ij})$$
By Lemma \ref{muJ}, this implies

\begin{equation}\label{DP}
\area_G(\Psi'_{ij})\le N_4 x^2- N_4 x d +N_3\mu(\Delta)
+ 2N_3Jn|{\bf  p}_{ij}|
\le N_4 x^2
+N_3\mu(\Delta)- N_4 x d+ 2N_3Jny
\end{equation}
By Lemma \ref{ppbar}, we have $|{\bf\bar p}_{ij}|< |{\bf p}_{ij}|+|\partial\Pi|$, and so the perimeter  $|\partial\Psi_{ij}|$ is less
than $2|{\bf p}_{ij}|+|\partial\Pi|$. Since $|\partial\Pi|\le
L|{\bf\bar p}_{ij}|,$ we obtain:
\begin{equation}\label{Lp2}
|\partial\Psi_{ij}|<(2+L)|{\bf p}_{ij}|\le (L+2)y
\end{equation}

By the inequality (\ref{Lp2}) and  Lemma \ref{main}, we have
\begin{equation}\label{aGP}
\area_G(\Psi_{ij})\le N_2(2+L)^2y^2+N_1\mu(\Psi_{ij})\le N_2(J+1)(2+L)^2y^2,
\end{equation}
where the second inequality follows from Lemma \ref{mixture} (a) since $N_2>N_1$.

By Lemma \ref{disk} and (\ref{Lp2}), the $G$-area of $\Pi$ does not exceed
$c_6|\partial\Pi|^2\le c_6 (L+2)^2 y^2$, and so there is
a constant $c_7=c_7(L)$ such that $\area_G(\Pi)\le c_7 y^2$.

This estimate and (\ref{aGP}) give the inequality $$\area_G(\bar\Delta_{ij})\le N_2(J+1)(2+L)^2y^2+c_7 y^2,$$ and we
obtain with (\ref{DP}) that
$$\area_G(\Delta)\le \area_G(\Psi_{ij}')+\area_G(\bar\Delta_{ij})\le $$ $$ N_4 x^2
+N_3\mu(\Delta) -N_4 x d+ 2N_3Jny
+N_2(J+1)(2+L)^2y^2+c_7 y^2$$

To obtain the desired contradiction with (\ref{ce}), it suffices to show that here, the number $T=N_4 x d/3$ is greater than each of the last three summands.
Recall that $x\ge n$, $d>\varepsilon y/2$ by (\ref{de}), $\varepsilon=N_4^{-1/2}$, and so $T>2N_3Jny$ if $N_4$ is large enough in comparison with $N_3$
and other constant chosen earlier. Also we have $T> N_2(J+1)(2+L)^2y^2$, because $$x=n+\sigma_{\lambda}(\Delta^*)>
|{\bf p}_{ij}|+\sigma_{\lambda}(\bar\Delta_{ij}^*)=y$$ by Lemma \ref{wmin} (a), and so
$xd>x\varepsilon y/2\ge \varepsilon y^2/2$.
Finally, $T>c_7 y^2$ since $$xd>x\varepsilon y/2\ge y^2\varepsilon/2$$
\endproof}

 \index[g]{clove $\Psi=cl(\pi,{\mathcal B}_1,{\mathcal B}_{L-3})$ of the minimal counterexample from Section \ref{midi}!q@${\mathbf q}_{i,i+1}$: a shortest path homotopic
to ${\mathbf p}_{i, i+1}$ in the subdiagram $\Psi_{ij}$, such that the first and the last $\tt$-edges of ${\mathbf q}_{i,i+1}$ coincide with the first  and the last $\tt$-edges of ${\mathbf p}_{i,i+1}$}
For every path ${\mathbf p}_{i,i+1}$ we will fix a shortest path ${\mathbf q}_{i,i+1}$ homotopic
to ${\mathbf p}_{i, i+1}$ in the subdiagram $\Psi_{i,i+1}$, such that the first and the last $\tt$-edges of ${\mathbf q}_{i,i+1}$ coincide with the first  and the last $\tt$-edges of ${\mathbf p}_{i,i+1}$. For $j>i+1$ the path
\index[g]{clove $\Psi=cl(\pi,{\mathcal B}_1,{\mathcal B}_{L-3})$ of the minimal counterexample from Section \ref{midi}!q@${\mathbf q}_{i,j}={\mathbf q}_{i,i+1},\dots {\mathbf q}_{j-1,j} $ if $j>i+1$}
 ${\mathbf q}_{i,j}$
is formed by ${\mathbf q}_{i,i+1},\dots, {\mathbf q}_{j-1,j} $.

\begin{lemma}\label{qq}
If $i\le r$ and $j\ge r+1$, then $$|{\mathbf q}_{ij}|\ge |{\mathbf q}_{ij}|_{\theta}+ |{\mathbf q}_{ij}|_q\ge h_i+h_j+(j-i)N+1$$
\end{lemma}

{ The proof is similar to the second part of Lemma
\ref{ppbar}.}

\medskip

Let
\index[g]{clove $\Psi=cl(\pi,{\mathcal B}_1,{\mathcal B}_{L-3})$ of the minimal counterexample from Section \ref{midi}!p@$\Psi_{ij}^0$ (resp. $\Psi^0$, $\Delta^0$) is the subdiagram of $\Psi_{i, j}$ (resp. of $\Psi$, of $\Delta$) obtained after replacing the subpath ${\mathbf p}_{ij}$ (of ${\bf p=p(\Psi)}$ ) by ${\mathbf q}_{ij}$ (resp. by ${\mathbf q}={\mathbf q}_{1,L-3}$) in the boundary}
$\Psi_{ij}^0$ (let $\Psi^0$, $\Delta^0$) be the subdiagram of $\Psi_{ij}$ (of $\Psi$, of $\Delta$) obtained after replacement of the subpath ${\mathbf p}_{ij}$ (of $\bf p$ ) by ${\mathbf q}_{ij}$ (by ${\mathbf q}={\mathbf q}_{1,L-3}$, resp.) in the boundary. 

\begin{lemma} \label{noq} (1) The subdiagram  $\Psi_{i,j}^0$ has no maximal $q$-bands except for the $q$-spokes starting from $\partial\Pi$.

(2) Every $\theta$-band of $\Psi_{i,i+1}^0$ ($i=1,\dots,L-4$) is crossed
by the path ${\mathbf q}_{i,i+1}$ at most once.
\end{lemma}

\proof (1) Assume there is a $q$-band $\mathcal Q$ of $\Psi_{ij}^0$ starting and ending on ${\bf q}_{ij}$. Then $j=i+1$ and ${\bf q}_{i,i+1}=\bf uevfw$, where $\mathcal Q$ starts with the $q$-edge $\bf e$ and ends with the $q$-edge $\bf f$. Suppose that $\mathcal Q$ has length $\ell$. Then $|{\bf v}|\ge \ell$
since every maximal $\theta$-band of $\Psi_{i,i+1}^0$ crossing $\mathcal Q$ has to end on the
subpath  $\bf v$. So one has $|{\bf evf}|\ge \ell+2$, and replacing the subpath ${\bf evf}$
by a side of $\mathcal Q$ of length $\ell$ one replaces the path ${\bf q}_{i,i+1}$ with
a shorter homotopic path by Lemma \ref{ochev}. This contradicts the choice
of ${\bf q}_{i,i+1}$, and so statement (1) is proved.

(2) Assume there is a $\theta$-band $\mathcal T$ of $\Psi_{i, i+1}^0$ starting and ending on ${\bf q}_{i,i+1}$. Then ${\bf q}_{i,i+1}=\bf uevfw$, where $\mathcal T$ starts with the $\theta$-edge $\bf e$ and ends with the $\theta$-edge $\bf f$. Moreover, one can chose $\cal T$ such that
$v$ is a side of this $\theta$--band. By Statement (1) the band $\cal T$ has less than $N$
$(\theta, q)$-cells. Therefore if $v'$ is another side of $\cal T$, we have
$|v'|_Y-|v|_Y \le 2N$. It follows from the definition of length in Subsection \ref{lf}
that $|evf|-|v'|\ge 2 - 2\delta N>1+2\delta$. Therefore, by Lemma \ref{ochev} (c), replacing the subpath $evf$ with $v'$ we decrease the length of  ${\bf q}_{i,i+1}$ at least by $1$,
a contradiction.
\endproof

It follows from Lemma \ref{psi1} that between the spokes ${\mathcal B}_j$ and ${\mathcal B}_{j+1}$ ($1\le j\le r-1$), \index[g]{clove $\Psi=cl(\pi,{\mathcal B}_1,{\mathcal B}_{L-3})$ of the minimal counterexample from Section \ref{midi}!g@$\Gamma_j$ subtrapezia between ${\mathcal B}_{j+1}$ and ${\mathcal B}_j$}
\index[g]{clove $\Psi=cl(\pi,{\mathcal B}_1,{\mathcal B}_{L-3})$ of the minimal counterexample from Section \ref{midi}!h@$h_{j+1}$ the hight of $\Gamma_j$}
there is a trapezium
$\Gamma_j$ of height $h_{j+1}$ with the side $\tt$-bands . Similarly, we have trapezia $\Gamma_j$ for $r+1\le j\le L-4$.
By Lemma \ref{noq} (2),  every trapezium $\Gamma_j$ is contained in both $\Psi_{j,j+1}$ and $\Psi_{j,j+1}^0$.

\index[g]{clove $\Psi=cl(\pi,{\mathcal B}_1,{\mathcal B}_{L-3})$ of the minimal counterexample from Section \ref{midi}!y@${\mathbf y}_j$ the bottom path of $\Gamma_j$}
 \index[g]{clove $\Psi=cl(\pi,{\mathcal B}_1,{\mathcal B}_{L-3})$ of the minimal counterexample from Section \ref{midi}!z@${\mathbf z}_j$ the top path of $\Gamma_j$}

The bottom paths ${\mathbf y}_j$ of all trapezia $\Gamma_j$ are contained in $\partial\Pi$ and have the same label $W\tt$. We will use  ${\mathbf z}_j$ for the top paths of these trapezia.  Since $\Gamma_j$ and $\Gamma_{j-1}$ ($2\le j\le r-1$) have the same bottom labels and the history $H_j$ is a prefix of $H_{j-1}$,
by Lemma \ref{simul}, $h_j$ different $\theta$-bands of $\Gamma_{j-1}$ form the copy $\Gamma'_{j}$ of the trapezium $\Gamma_j$ (more precisely, a copy of a superscript shift $\Gamma_j^{(+(\pm 1))}$) with top and bottom paths ${\mathbf z}'_j$ and ${\mathbf y}'_j={\mathbf y}_{j-1}$.

\begin{figure}[ht]
\begin{center}
\includegraphics[width=1.0\textwidth]{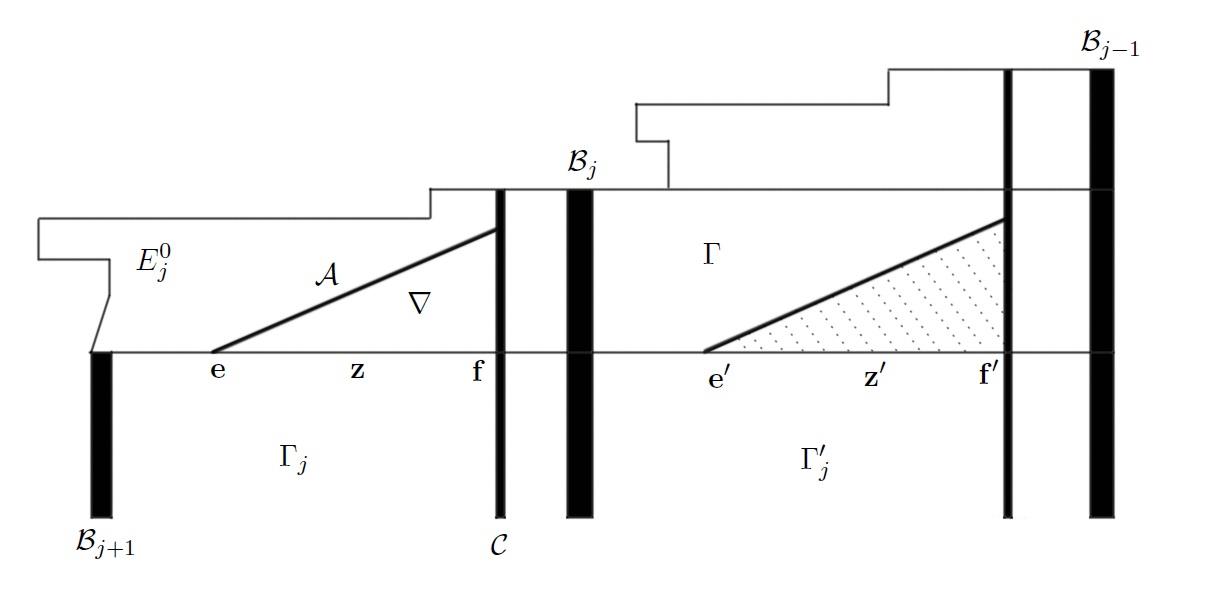}
\end{center}
\caption{Lemma \ref{copy}}\label{Pic15}
\end{figure}

We denote by  $E_j$
\index[g]{clove $\Psi=cl(\pi,{\mathcal B}_1,{\mathcal B}_{L-3})$ of the minimal counterexample from Section \ref{midi}!e@$E_j$ (resp. $E_j^0$): the comb formed by the maximal $\theta$-bands of $\Psi_{j,j+1}$
(resp.of $\Psi_{j,j+1}^0$)}
(by $E_j^0$ )
the comb formed by the maximal $\theta$-bands of $\Psi_{j,j+1}$
(of $\Psi_{j,j+1}^0$, respectively)
crossing the $\tt$-spoke ${\mathcal B}_j$  but not crossing ${\mathcal B}_{j+1}$ ($1\le j\le r-1$, see fig. \ref{Pic15}). Its handle $\ccc_j$ of
height $h_j-h_{j+1}$ is contained in ${\mathcal B}_j$. The boundary $\partial E_j$  (resp., $\partial E_j^0$) consists of the side of this handle, the path ${\mathbf z}_j$ and the path ${\mathbf p}_{j,j+1}$
(the path ${\mathbf q}_{j,j+1}$, respectively).

Assume that a maximal $Y$-band $\mathcal A$ of $E_j^0$ ($2\le j\le r-1$) starts on the path ${\mathbf z}_j$ and ends
on a side $Y$-edge of a maximal $q$-band $\ccc$ of $E_j^0$. Then $\mathcal A$, a part of $\ccc$ and a part ${\mathbf z}$ of ${\mathbf z}_j$
bound a comb $\nabla$.

\begin{lemma} \label{copy} There is a copy of the comb $\nabla$ in the trapezium $\Gamma=\Gamma_{j-1}\backslash\Gamma'_j$. It is a superscript shift of $\nabla$.
\end{lemma}

{\proof The subpath ${\bf z}$ of ${\bf z}_j$ starts with an $Y$-edge $\bf e$ and ends with a $q$-edge $\bf f$.
There is a copy ${\bf z}'$ of $\bf z$ in ${\bf z}'_j$ starting with ${\bf e'}$ and ending with $\bf f'$. Note that the $\theta$-cells $\pi$ and $\pi'$ attached to $\bf f$ and to $\bf f'$ in $\nabla$ and in $\Gamma$ are copies of each other up to superscript shift, since they correspond to the same letter of the history. Now moving from $\bf f$ to $\bf e$, we see that
the whole maximal $\theta$-band ${\mathcal T}_1$ of $\nabla$ containing $\pi$ has a copy in $\Gamma$. Similarly
we obtain a copy of the next maximal $\theta$-band ${\mathcal T}_2$ of $\nabla$, and so on.
\endproof}

\subsubsection{Bounding the number of $Y$-bands in a sector of a clove}

\begin{lemma} \label{le6} At most $N$ $Y$-bands starting on the path ${\mathbf y}_j$ can end on
a $(\theta,q)$-cells of the same $\theta$-band. This property holds
for the $Y$-bands starting on ${\mathbf z}_j$ too.
\end{lemma}

{\proof We will prove the second claim only since the proof of the first
one is similar. Assume that the $Y$-bands ${\mathcal A}_1,\dots, {\mathcal A}_s$
start from ${\bf z}_j$ and end on some $(\theta,q)$-cells of a $\theta$-band $\mathcal T$. Let ${\cal T}_0$ be the minimal subband of $\cal T$, where
the $Y$-bands ${\mathcal A}_2,\dots, {\mathcal A}_{s-1}$ end and ${\bf \bar z}_j$ be the minimal subpath
of ${\bf z}_j$, where they start. Then by Lemma \ref{NoAnnul}, every maximal $q$-band starting
on ${\bf\bar z}_j$ has to cross the band ${\mathcal T}_0$ and vice versa. Hence the base of ${\mathcal T}_0$
is a subbase of the standard base (or of its inverse). Since every rule of $\bf M$ can
change at most $N-2$ $Y$-letters  in a word with standard base, all $(\theta,q)$-cells of ${\mathcal T}_0$ have at most $N-2$ $Y$-edges,
and the statement of the lemma follows.
\endproof}

Without loss of generality, we assume that
\begin{equation}\label{L0}
h=h_{L_0+1}\ge h_{L-L_0-3}.
\end{equation}
(Recall that $L_0$ is one of the parameters used in the paper, a number between $c_5$ and $L$, Section \ref{param}.)


\subsubsection{Estimating the sizes of trapezia $\Gamma_j$}

Recall that the integer $r$ was defined in Lemma (\ref{psi1}) (2).

\begin{lemma} \label{02} If $h\le L_0^2|W|_Y$, then the number of trapezia $\Gamma_j$ with the
properties $|{\mathbf z}_j|_Y\ge |W|_Y/c_5N$ for $j\in [L_0+1,r-1]$ or $j\in [r+1, L-L_0-5]$,
is smaller than $L/5$.
\end{lemma}

{\proof
Consider $\Gamma_j$ as in the assumption of the lemma with $j\in [L_0+1,r-1]$.
The subcomb $E_j^0$ has at most $N$ maximal $q$-bands by Lemma \ref{noq}.
So there are at most $N$ maximal $Y$-bands starting on ${\bf z}_j$ and ending
on each of the $\theta$-bands of $E_j^0$. If $g_j$ is the length of the handle of $E_j^0$ for an index $j$ from the set $S= [L_0+1,r-1]\cup [r+1, L-L_0-5]$, then $\sum_{j\in S} g_i\le 2h$.
Hence at most $2hN$ maximal $Y$-bands starting on all ${\bf z}_j$-s, $j\in S$ (denote this
set of $Y$-bands by ${\bf A}$), end on some $(\theta,q)$-cells.

Proving by contradiction, we
have at least $L|W|_Y/5c_5N$  $Y$-bands in ${\bf A}$.
Hence at least $L|W|_Y/5c_5N-2hN$ bands from
$\bf A$ end on the subpaths ${\bf q}_{j,j+1}$ for $j\in S$. Since the path ${\bf q}_{j,j+1}$
has at most $2h$ $\theta$-edges by Lemma \ref{noq}. Therefore   by Lemma \ref{ochev},
at least $L|W|_Y/5c_5N-2hN-2h$ $Y$-edges contribute in the length of this
path. It follows from Lemma
\ref{qq} that
$$|{\bf p}_{L_0+1, L-L_0-5}|\ge |{\bf q}_{L_0+1, L-L_0-5}|\ge h_{L_0+1}+h_{L-L_0-5}
+LN/2+ \delta(L|W|_Y/5c_5N-2hN-2h)$$
\begin{equation}\label{pp}
\ge h_{L_0+1}+h_{L-L_0-5}
+LN/2+ \delta L|W|_Y/10c_5N
\end{equation}
since  $2hN+2h\le 3L_0^2N|W|_Y$ by the assumption of the lemma, which is less than   $L_0^3|W|_Y/10c_5N\le L|W|_Y/10c_5N$ because
$L_0 \ll L$ (see Section \ref{param}).

Also by Lemma \ref{ppbar}, we have
$$|{\bf \bar p}_{L_0+1, L-L_0-5}|\le  h_{L_0+1}+h_{L-L_0-5}
+3L_0N+ 3L_0\delta|W|_Y$$
\begin{equation}\label{barbar}
\le h_{L_0+1}+h_{L-L_0-5}
+3L_0N+ \delta L|W|_Y/20c_5N,
\end{equation}
because by Section \ref{param}, we have $3L_0< L/20c_5N$.
The inequalities (\ref{pp}, \ref{barbar}) give us

\begin{equation}\label{pmpb}
|{\bf p}_{L_0+1, L-L_0-5}|-|{\bf\bar p}_{L_0+1, L-L_0-5}|\ge LN/3+\delta L|W|_Y/20c_5N
\end{equation}
because $L
\gg L_0$. Since $h_{L_0+1}+h_{L-L_0-5}\le 2h\le 2L_0^2|W|_Y<L|W|_Y$, it follows from (\ref{barbar}))
that $$|{\bf\bar p}_{L_0+1, L-L_0-5}|<L|W|_Y +3L_0N+ \delta L|W|_Y/20c_5N,$$
which implies, together with (\ref{pmpb}), that
\begin{equation}\label{40}
\frac{|{\bf p}_{L_0+1, L-L_0-5}|-|{\bf\bar p}_{L_0+1, L-L_0-5}|}{|{\bf\bar p}_{L_0+1, L-L_0-5}|}\ge\min\big(\frac{3L_0N}{LN/3} , \frac{\delta L|W|_Y/20c_5N}{2L|W|_Y}\big)>\delta/40c_5N
\end{equation}
Finally, for the right-hand side, we have $\delta/40c_5N>\varepsilon=N_4^{-1/2}$ by the choice of $N_4$ and the inequality (\ref{40}) implies
$$\frac{|{\bf p}_{L_0+1, L-L_0-5}|}{|{\bf\bar p}_{L_0+1, L-L_0-5}|}> 1+\varepsilon$$ which
 contradicts Lemma \ref{epsi}. The lemma is proved.
\endproof}

\begin{lemma} \label{hh} If $h\le L_0^2|W|_Y$, then the histories $H_1$ and $H_{L-3}$ have different first letters unless all these letters are equal to $\theta(23)^{-1}$.
\end{lemma}
\proof Let $\mathcal T$ and $\mathcal S$ be the maximal $\theta$-bands of $\Psi$
crossing ${\mathcal B}_1$ and ${\mathcal B}_{L-3}$, respectively, and the closest to the disk $\Pi$. Suppose they cross spokes number $k$ and $\ell$
of $\Pi$, respectively. Note that $\Gamma_j$ has height zero if it is not crossed either by $\mathcal T$ or by  $\mathcal S$, and then
$|{\bf z}_j|_Y=|W|_Y$. Therefore by Lemma \ref{02}, $k+\ell>L-L/5-3L_0>2L/3$,
and also $k,\ell\ge 2$ since $L/2-3\le r \le L/2$. It follows from  Lemma \ref{withd} (2) (applied to $\Delta^*$) that the first letters of $H_1$ and
$H_{L-3}$ are different.
\endproof

\begin{lemma} \label{Wh0} If $h\le L_0^2|W|_Y$, then

\begin{equation}\label{WY}
 |W|_Y>\frac{LN}{4L_0}
\end{equation}

\end{lemma}
\proof Assume that $|W|_Y\le  LN/4L_0$. By Lemma \ref{ppbar} for $i=L_0+1$ and $j=L-L_0-3$, we have
$|{\bf p}_{i.j}|\ge h_i+h_j+(L-3L_0)N$ and $|\mathbf{\overline p}_{ij}|\le h_i+h_j+3L_0(N+|W|_Y)$, whence
\begin{equation}\label{Wh1}
|{\bf p}_{i.j}|-|\mathbf{\overline p}_{ij}|\ge (L-6L_0)N -3L_0|W|_Y
> (L-6L_0)N - \frac 34 LN > LN/5,
\end{equation}
because $L>>L_0$. It follows from inequalities (\ref{downup}, \ref{L0}) that $h_i+h_j\le 2h$. Hence
\begin{equation}\label{Wh2}
|\mathbf{\overline p}_{ij}|\le 2h +3L_0(N+LN/4L_0) \le 2L_0^2\frac{LN}{4L_0}+ LN< L_0LN
\end{equation}
Inequalities (\ref{Wh1} and \ref{Wh2}) imply
$$\frac{|{\bf p}_{i.j}|-|\mathbf{\overline p}_{ij}|}{|\mathbf{\overline p}_{ij}|}>\frac{1}{5L_0}>\varepsilon$$
since $\varepsilon=N_4^{-1/2}$, which contradicts  the statement of Lemma \ref{epsi}.
\endproof

\begin{lemma}\label{arG} We have $h> L_0^2|W|_Y$.
\end{lemma}
\proof Proving by contradiction, we have inequality (\ref{WY}) from Lemma \ref{Wh0}.

By Lemma \ref{02}, there are at least $L-L/5-3L_0> 0.7L$ trapezia
$\Gamma_j$ with $|{\mathbf z}_j|_Y<|W|_Y/c_5N$, and so one can choose two such trapesia
$\Gamma_k$ and $\Gamma_{\ell}$ such that  $k<r$, $\ell\ge r+1$ and $\ell-k>0.6L$.
Since $H_{k+1}$ (resp. $H_{\ell}$) is a prefix of $H_1$ (of $H_{L-3}$), it follows from
Lemma \ref{hh} that the first letters of $H_{k+1}$ and $H_{\ell}$ are different unless
they are equal to $\theta(23)^{-1}$.

Since the bottom paths of $\Gamma_k$ and $\Gamma_{\ell}$ (which belong to $\partial\Delta$) have the same label,
up to a superscript shift,
one can construct an auxiliary trapezium $E$ identifying the bottom of a copy of $\Gamma_k$
and the bottom of a mirror copy of $\Gamma_{\ell}$. The history of $E$ is $H_{\ell}^{-1}H_{k+1}$, which is an eligible word if the first letters of $H_k$ and $H_{\ell}$ are different.

If both first letters are $\theta(23)^{-1}$, then the word $H_{\ell}^{-1}H_{k+1}$ is also eligible
by definition. If the bottom $\theta$-bands of $\Gamma_k$ and $\Gamma_l$ are just copies of each
other then the above constructed diagram $E$ is not reduced. However one can modify
the construction replacing $\Gamma_k$ by an auxiliary superscript shift $\Gamma_k^{(+1)}$. By the
definition of relations (\ref{rel11}), the bottom labels of $\Gamma_k^{(+1)}$, $\Gamma_k$ and
$\Gamma_l$ are all equal, but the top labels of the first $\theta$-bands of  $\Gamma_k^{(+1)}$
and $\Gamma_l$ are not mirror copies of each other (they differ by $1$-shift), and so the
diagram $E$ obtained by identifying the bottom path of a copy of $\Gamma_k^{(+1)}$
and the bottom path of a mirror copy of $\Gamma_{\ell}$ is reduced, i.e.,we can obtain the trapezium
$E$ in any case.

The top $W_0$ and the bottom $W_t$ of $E$  have $Y$-lengths less than $|W|_Y/c_5N$. Without loss of generality, one may assume that
$h_{k+1}\ge h_{\ell}$, and so $h_{k+1}\ge t/2$, where $t$ is the
height of $E$.

Note that the difference of $Y$-lengths $|W|_Y-|W|_Y/c_5N>|W|_Y/2$, and so
\begin{equation}\label{2N}
h_{k+1}, h_{\ell}> |W|_Y/2N
\end{equation}
 since the difference of $Y$-lengths for the top and the bottom of every maximal $\theta$-band of $E$ does not exceed $N$. Therefore  by (\ref{WY}), we obtain inequality
 \begin{equation}\label{tg}
 t>\frac{|W|_Y}{N}\ge \frac{L}{4L_0}
 \end{equation}

 If $|W_0|_Y=|W_t|_Y=0$, then $||W_0||=||W_t||=N$, and so \\ $\max (||W_0||,||W_t||)<\frac{L}{4c_4L_0}<t/c_4$ by the choice of $L$ and (\ref{tg}). If \\ $\max (||W_0||,||W_t||)\ge 1$, then $$\max (||W_0||,||W_t||)\le N+1+\max (|W_0|_Y,|W_t|_Y)<N+1+\frac{|W|_Y}{c_5N}< \frac{2|W|_Y}{c_5N}$$
 by inequality (\ref{WY}) since $\frac{L}{4L_0c_5}>N+1$ by the choice of $L$. It follows from the choice of $c_5$ and (\ref{tg}) that
 $\max (||W_0||,||W_t||)<\frac{2|W|_Y}{c_5N}< \frac{|W|_Y}{c_4N}<t/c_4$.
Therefore in both  cases, the computation corresponding $E$ satisfies the assumption of Lemma \ref{B}.

So for every factorization $H'H''H'''$ of the
history of $\Gamma_k$, where $||H'||+||H''||\le \lambda ||H'H''H'''||$, we have $||H''||>0.4t$, since we can choose $\lambda< 1/5$ in (\ref{param}).
Therefore by Lemma \ref{B}, the spoke ${\mathcal B}_{k+1}$
is a $\lambda$-shaft.

Using Lemma \ref{ppbar}, we obtain:
\begin{equation}\label{pb}
|{\mathbf  p}_{k+1, \ell}|+\sigma_{\lambda}(\overline \Delta_{k+1,\ell}^*) \ge h_{k+1}+h_{\ell}+0.6 LN+h_{k+1}
\end{equation}

By inequality (\ref{2N}), we have  $\delta L|W|_Y\le 2LN\delta h_{k+1}<h_{k+1}$ by the choice of $\delta$. This inequality and  Lemma \ref{ppbar} provide us with

\begin{equation}\label{bp}
|\mathbf{\overline p}_{k+1,\ell}|\le h_{k+1}+h_{\ell}+0.4 LN+0.4L\delta|W|_Y\le h_{k+1}+h_{\ell}+h_{k+1}/2
\end{equation}

The right-hand side of the inequality (\ref{pb}) divided  by the right-hand side of (\ref{bp}) is
greater than $1.1$
(because $h_{k+1}\ge h_{\ell}$),
which contradicts Lemma \ref{epsi}. Thus, the lemma is proved.
\endproof

\begin{lemma} \label{hi} We have $h_i> \delta^{-1}$ for every $i=1,\dots, L_0$.
\end{lemma}

{\proof By inequalities (\ref{L0}) and (\ref{downup}), we have $h_i\ge h_{L-L_0-3}$.
Proving by contradiction, we obtain $|W|_Y<h_i
\le \delta^{-1}$ for some $i=1,\dots, L_0$
by Lemma \ref{arG}.
Then $$|{\bf \bar p}_{i, L-L_0-3}|< h_{i}+h_{L-L_0-3}+3L_0(N+\delta^{-1} \delta)\le
h_{i}+h_{L-L_0-3}+4L_0N$$ by Lemma \ref{ppbar}, and the inequality
$|{\bf p}_{i, L-L_0-3}|\ge h_{i}+h_{L-L_0-3}+LN/2$. Since $h_{i}+h_{L-L_0-3}\le 2\delta^{-1}$
and $4L_0N< LN/4$, we see that $\frac{|{\bf p}_{i, L-L_0-3}|}{|{\bf \bar p}_{i, L-L_0-3}|}> 1+\delta>1+\varepsilon$ which contradicts Lemma \ref{epsi}.
\endproof.}

\subsubsection{Bounding shafts in a clove and corollaries of the bound}

\begin{lemma} \label{shaft} None of the spokes ${\mathcal B}_1,...,{\mathcal B}_{L_0}$ contains
a $\lambda$-shaft at $\Pi$ of  length at least $\delta h$.
\end{lemma}

{\proof On the one hand, by Lemmas \ref{ppbar} and \ref{arG},
\begin{equation}\label{pbp}
|{\bf \bar p}_{L_0+1, L-L_0-3}|< h_{L_0+1}+h_{L-L_0-3}+3L_0(N+\delta|W|_Y)<h_{L_0+1}+h_{L-L_0-3}+3L_0(N+\delta L_0^{-2}h).
\end{equation}
On the other hand, by Lemma \ref{ppbar},
\begin{equation}\label{bpb}
|{\bf p}_{L_0+1,L-L_0-3}|> h_{L_0+1}+h_{L-L_0-3}+(L-3L_0)N.
\end{equation}

If the statement of the lemma were wrong, then we would have $\sigma_{\lambda}(\bar\Delta^*)\ge \delta h$, and  inequalities (\ref{pbp}) and (\ref{bpb}) would imply that
$$|{\bf p}_{L_0+1,L-L_0-3}|-|{\bf \bar p}_{L_0+1, L-L_0-3}|+\sigma_{\lambda}(\bar\Delta^*)\ge (L-6L_0)N - 3L_0^{-1}\delta h +\delta h\ge LN/2 +\delta h/2 $$
The right-hand side of the last inequality divided  by the right-hand side of (\ref{pbp}) is greater than $\varepsilon=N_4^{-\frac12}$,
because $h\ge h_{L_0+1}, h_{L-L_0-3}$,
which contradicts  Lemma \ref{epsi}. Thus, the lemma is proved.
\endproof}

\begin{lemma}\label{zgh}
For every $j\in [1,L_0-1]$, we have
$|{\mathbf z}_{j}|_Y> h_{j+1}/c_5$.
\end{lemma}

{\proof If $|{\bf z}_{j}|_Y\le h_{j+1}/c_5$, then $$||{\bf z}_{j}||\le |{\bf z}_{j}|_Y+N+1
\le 2h_{j+1}/c_5 \le h_{j+1}/c_4$$ since by (\ref{downup}) and Lemma \ref{hi}, we have $h_{j+1}/c_5\ge h/c_5\ge \delta^{-1}/c_5 >N+1$. Similarly by Lemma \ref{arG}, $$||{\bf y}_{j}||\le | W_{j}|_Y+N+1\le N+1 + h_{j+1}/L_0^2\le  2h_{j+1}/L_0^2< h_{j+1}/c_4$$
since $N+1 < \delta^{-1}/L_0^2$ by Section \ref{param}.

Thus, the computation ${\mathcal C}: \;W_0\to\dots\to W_t$
corresponding to the trapezium $\Gamma_j$ satisfies the assumption of Lemma \ref{B},
since $t=h_{j+1}$.
Hence ${\mathcal B}_{j+1}$ is a $\lambda$-shaft by Lemma \ref{B} since $\lambda<1/2.$
We obtain a contradiction with Lemma \ref{shaft} since $\delta h\le h\le h_{j+1}$,
and the lemma is proved.
\endproof}

\begin{lemma} \label{2000N} For every $j\in [1,L_0-1]$, we have $h_{j+1}< (1-\frac{1}{10c_5N})h_j$.
\end{lemma}

{\proof By Lemma \ref{zgh}, we have $|{\bf z}_{j}|_Y\ge h_{j+1}/c_5$.
Let us assume that $h_{j+1}\ge (1-\frac{1}{10c_5N})h_j$, that is the handle ${\mathcal C}_j$ of $E_j$ has height at most $h_j/10c_5N$.
By Lemma \ref{le6}, at most $h_{j}/10c_5$ maximal $Y$-bands of $E_j$ starting on ${\bf z}_j$ can end on
the $(\theta,q)$-cells of $E_j$. Hence at least $$|{\bf z}_j|_Y-h_{j}/10c_5\ge |{\bf z}_j|_Y-2h_{j+1}/10c_5\ge h_{j+1}/c_5-h_{j+1}/5c_5= 0.8h_{j+1}/c_5> 0.7h_{j}/c_5$$ of them have to end on the path ${\bf p}_{j,j+1}$.

The path ${\bf p}_{j,j+1}$ has at most $h_j-h_{j+1}\le \frac{h_j}{10c_5N}\;\;$ $\theta$-edges by Lemma \ref{psi1}. Hence by Lemma \ref{ochev},
$$|{\bf p}_{j,j+1}|\ge h_j-h_{j+1}+\delta(0.7h_{j}/c_5-h_j/10c_5N)\ge h_j-h_{j+1}+0.6\delta h_{j}/c_5.$$
By Lemma \ref{ppbar}, the path ${\bf p}_{j+1, L-L_0-3}$ has length at least $2LN/3+h_{j+1} + h_{L-L_0-3}$
and therefore ,
$$|{\bf p}_{j, L-L_0-3}|\ge |{\bf p}_{j,j+1}|+ |{\bf p}_{j+1, L-L_0-3}|-1> LN/2+h_j + h_{L-L_0-3}+0.6\delta h_{j}/c_5.$$
On the other hand by Lemma \ref{ppbar}, we have
$$|{\bf\bar p}_{j, L-L_0-3}|\le h_j + h_{L-L_0-3}+ 3NL_0+3L_0\delta|W|_Y\le  h_j + h_{L-L_0-3}+ 3NL_0+3L_0^{-1}\delta h_{j+1}$$
by Lemma \ref{arG} and inequality $h\le h_{j+1}$. Hence
$\frac{|{\bf p}_{j, L-L_0-3}|}{|{\bf\bar p}_{j, L-L_0-3}|}\ge (1+\delta/10c_5)$
since $h_{L-L_0-3}\le h_{L_0+1}\le h_{j+1}\le h_j$ and $L_0\gg c_5$. We have a contradiction with Lemma \ref{epsi} since
$\delta/10c_5>\varepsilon$.
The lemma is proved by contradiction.
\endproof

The proof of the next lemma is similar.}

\begin{lemma} \label{zlh} For every $j\in [2,L_0-1]$.
we have $|{\mathbf z}_j|_Y \le 2Nh_j$,
\end{lemma}

{\proof Assume that $|{\bf z}_j|_Y \ge 2Nh_j$.
By Lemma \ref{le6}, at most $Nh_j$ maximal $Y$-bands of $E_j$ starting on ${\bf z}_j$ can end on
the $(\theta,q)$-cells of $E_j$. Hence at least $|{\bf z}_j|_Y-Nh_j\ge Nh_j$ of them have to end on the path ${\bf p}_{j,j+1}$.
The path ${\bf p}_{j,j+1}$ has at most $h_j$ $\theta$-edges. Hence by Lemma \ref{ochev},
$$|{\bf p}_{j,j+1}|\ge h_j-h_{j+1}+\delta(Nh_j-h_j)= h_j-h_{j+1}+\delta(N-1)h_j$$ and therefore by Lemma \ref{ppbar},
$$|{\bf p}_{j, L-L_0-3}|\ge LN/2+h_j + h_{L-L_0-3}+\delta(N-1)h_j.$$
On the other hand by Lemmas \ref{ppbar} and \ref{arG}, we have
$$|{\bf\bar p}_{j, L-L_0-3}|\le h_j + h_{L-L_0-3}+ 3NL_0+3L_0\delta|W|_Y\le  h_j + h_{L-L_0-3}+ 3NL_0+\frac{3\delta h_j}{L_0}$$ because $h\le h_j$. Since $h_j \ge h\ge  h_{L-L_0-3} $, we have
$\frac{|{\bf p}_{j, L-L_0-3}|}{|{\bf\bar p}_{j, L-L_0-3}|}\ge (1+\varepsilon)$,
a contradiction by Lemma \ref{epsi}.
\endproof}

\subsubsection{Certain subtrapezia with one step history do not exist in the clove}

\begin{lemma} \label{001} There is no $i\in [2,L_0-3]$ such that the histories
$H_{i-1}=H_iH'=H_{i+1}H''H'=H_{i+2}H'''H''H'$ and the computation $\ccc$ with history $H_i$ corresponding to the trapezium $\Gamma_{i-1}$ satisfy the following condition:

(*) The history $H'''H''H'$ has only one step, and for the subcomputation $\mathcal D$
with this history,
there is a sector $Q'Q$ such that
a state letter  from $Q$ or from $Q'$ inserts a letter increasing
the length of this sector after every transition of $\mathcal D$.
\end{lemma}

{\proof

Recall that the standard base of $\bf M$  is the product of the standard base $B$ of ${\bf M}_4$
and its inverse copy $(B')^{-1}$, and letter $\tt$.
Due to the mirror symmetry of the standard base, we have mirror symmetry for any accessible computation, in particular, we have it for $\mathcal C$ and $\mathcal D$. Therefore proving by contradiction,
we may assume that the $Y$-letters are inserted from the left of $Q$.

Let  $\mathcal Q$ be the maximal $q$-spoke of the subdiagram
$E_i^0\subset \Gamma_i$ corresponding to the base letter $Q$.
If ${\mathcal Q'}$ is the neighbor from the left $q$-spoke for $\mathcal Q$
(the spokes are directed from the disk $\Pi$),
then the subpath $\bf x$ of ${\bf z}_i$ between these two $q$-spokes has at least
$h_{i+1}-h_{i+2}=||H'''||$ $Y$-letters. Indeed, $\Gamma_i$ contains
a copy $\Gamma'_{i+1}$ of $\Gamma_{i+1}$, the bottom of the trapezium
$\Gamma_i\backslash\Gamma'_{i+1}$ is the copy ${\bf z}'_{i+1}$ of ${\bf z}_{i+1}$
and the top of it iz ${\bf z}_i$,
and so the subcomputation with history $H'''$ has already increased
the length of the $Q'Q$-sector. Thus,  by lemmas \ref{2000N}, \ref{arG} and the choice of $L_0>100c_5N$, we have
\begin{equation}\label{xa}
|{\bf x}|_Y\ge  h_{i+1}-h_{i+2}\ge \frac{1}{10c_5N}h_{i+1}\ge 10L_0 |W|_Y
\end{equation}

Note that an $Y$-band $\mathcal A$ starting on ${\bf x} $ cannot end on a $(\theta,q)$-cell from $\mathcal Q$.
Indeed, otherwise by Lemma \ref{copy}, there is a copy of this configuration in
the diagram $\Gamma_{i-1}$, i.e. the copy of $\mathcal A$ ends on the copy of $\mathcal Q$ which contradicts the assumption that the rules of computation with history $H'''H''H'$ do not delete $Y$-letters.

Let us consider the comb bounded by $\mathcal Q$, ${\mathcal Q}'$, $\bf x$
and the boundary path of $\Delta^0$ (without the cells from ${\mathcal Q}'$). If the lengths of the parts of $\mathcal Q$ and ${\cal Q}'$ bounding this comb are
$s$ and $s'$, respectively, then there are $|{\bf x}|+s$ maximal
$Y$-bands starting on $\bf x$ and $\mathcal Q$ and ending either on
${\mathcal Q}'$ or on $\partial\Delta^0$ since the comb has no
maximal $q$-bands by Lemma \ref{noq}. At most $s'<s$ of these
$Y$-bands can end on ${\mathcal Q}'$. Therefore at least $|{\bf x}|+s-s'$
of them end on the segment of the boundary path of $\Delta^0$ lying between the ends of ${\mathcal Q}'$ and $\mathcal Q$.

Since by Lemma \ref{noq} (2),
this segment has $s-s'$ $\theta$-edges, its length is at least $s-s'+\delta|{\bf x}|_Y$ by Lemma \ref{ochev}.
This inequality and inequality (\ref{xa})  imply

$$|{\bf p}_{i,L-L_0-3}|\ge |{\bf q}_{i,L-L_0-3}|\ge |{\bf q}_{i,L-L_0-3}|_q+|{\bf q}_{i,L-L_0-3}|_{\theta}+\frac{\delta}{10c_5N}h_{i+1},$$

and so by Lemma \ref{qq}, we have
$$|{\bf p}_{i,L-L_0-3}| \ge LN/2+h_{i}+h_{L-L_0-3}+\frac{\delta}{10c_5N}h_{i+1}$$ $$ \ge LN/2+h_{i}+h_{L-L_0-3}+10\delta L_0|W|_Y$$
Therefore by Lemma \ref{ppbar}, we obtain
\begin{equation}\label{pgp}
|{\bf p}_{i,L-L_0-3}|-\frac{7\delta}{100c_5N}h_{i+1}> 3L_0N+h_{i}+h_{L-L_0-3}+3\delta L_0|W|_Y\ge |{\bf \bar p}_{i,L-L_0-3}|,
\end{equation}
  Since $\Delta$ is a minimal counter-example, it follows from (\ref{pgp})
and Lemma \ref{wmin} (a,b) that the subdiagram $\Psi'_{i,L-L_0-3}$ (whose boundary path is
obtained from $\partial\Delta$ by replacing the subpath ${\bf p}_{i,L-L_0-3}$ with ${\bf \bar p}_{i,L-L_0-3}$) is weakly minimal but it is not a counter-example. Therefore  we obtain from
(\ref{pgp}) and Lemma \ref{wmin} (a,b):
$$\area_G(\Psi'_{i,L-L_0-3})\le N_4(|\Psi'_{i,L-L_0-3}|+\sigma_{\lambda}((\Psi'_{i,L-L_0-3})^*))^2+N_3\mu(\Psi'_{i,L-L_0-3})$$
$$\le N_4(n+\sigma_{\lambda}(\Delta^*)-\frac{7\delta}{100c_5N}h_{i+1})^2+N_3\mu(\Psi'_{i,L-L_0-3})$$
\begin{equation}\label{areai}
\le N_4(n+\sigma_{\lambda}(\Delta^*))^2-N_4\frac{7\delta n}{100c_5N}h_{i+1}+N_3\mu(\Psi'_{i,L-L_0-3})
\end{equation}

By Lemma \ref{arG}, $|W|_Y \le L_0^{-2}h_i$, and by Lemma \ref{hi}, $h_i>\delta^{-1}>100L_0N$,  whence $$|{\bf \bar p}_{i,L-L_0-3}|\le 2h_i+3L_0N+3\delta L_0|W|_Y\le(2+0.03+\frac{3\delta}{L_0})h_i \le 2.1h_i$$
by Lemma \ref{ppbar}, because $|{\bf \bar p}_{i,L-L_0-3}|\le |{\bf \bar p}_{i,L-L_0-3}|_q+|{\bf \bar p}_{i,L-L_0-3}|_{\theta}+\delta|{\bf \bar p}_{i,L-L_0-3}|_Y$.
This estimate and Lemma \ref{epsi},
give us
\begin{equation}\label{22}
|{\bf  p}_{i,L-L_0-3}|\le (1+\varepsilon)|{\bf \bar p}_{i,L-L_0-3}|<2.2h_i
\end{equation}

Note that $|\Psi_{i,L-L_0-3}|\le |{\bf p}_{i,L-L_0-3}|+|{\bf \bar p}_{i,L-L_0-3}|\le 2|{\bf p}_{i,L-L_0-3}|\le 5 h_i$ by inequalities (\ref{pgp}, \ref{22})). Hence by Lemmas \ref{main}, we have for the disk-free subdiagram  $\Psi_{i,L-L_0-3}$:
\begin{equation}\label{arP}
\area_G(\Psi_{i,L-L_0-3})\le
N_2|\Psi_{i,L-L_0-3})|^2+N_1\mu(\Psi_{i, L-L_0-3})\le 25N_2h_i^2+N_1\mu(\Psi_{i, L-L_0-3})
\end{equation}
Since by Lemma \ref{mixture} (a), $\mu(\Psi_{i, L-L_0-3})\le J|\Psi_{i, L-L_0-3}|^2<25Jh_i^2$,  it follows from (\ref{arP}) that
\begin{equation}\label{areaP}
\area_G(\Psi_{i,L-L_0-3})\le 25N_2h_i^2+25 N_1Jh_i^2\le 30N_2h_i^2
\end{equation}
since $N_2> 5N_1J$.

By Lemma \ref{disk}, the $G$-area of $\Pi$ is bounded by $c_6|\partial\Pi|^2$.
Inequalities (\ref{pgp}) and (\ref{22}) imply the inequality
$|\partial\Pi|<L|{\bf\bar p}_{i,L-L_0-3}|<L|{\bf p}_{i,L-L_0-3}|<3Lh_i$.
Therefore one may assume that the constant $c_7$ is chosen so that
\begin{equation}\label{c7}
\area_G(\Pi)< c_6|\partial\Pi|^2< c_7h_i^2
\end{equation}
It follows from (\ref{areaP}) and (\ref{c7}) that
\begin{equation}\label{Dij}
\area_G(\bar\Delta_{i,L-L_0-3})\le 30N_2h_i^2+ c_7h_i^2
\end{equation}
 Summing inequalities (\ref{Dij} and \ref{areai}), we have $$\area_G(\Delta)\le \area_G(\Psi'_{i, L-L_0-3})+\area_G(\bar\Delta_{i, L-L_0-3})\le $$
\begin{equation}\label{arD}
\le N_4(n+\sigma_{\lambda}(\Delta^*))^2-N_4\frac{7\delta n}{100c_5N}h_{i+1}+N_3\mu(\Psi'_{i,L-L_0-3})+ 30N_2h_i^2+ c_7h_i^2
\end{equation}

\begin{figure}[ht]
\begin{center}
\includegraphics[width=0.8\textwidth]{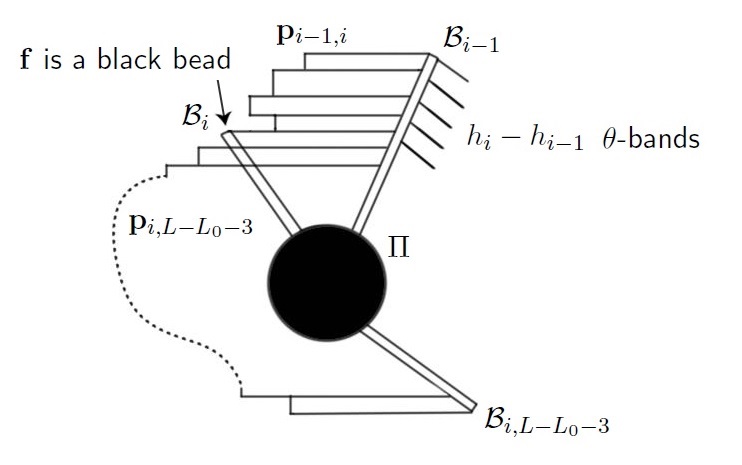}
\end{center}
\caption{$\mu(\Psi'_{i+1,L-L_0-3})-\mu(\Psi'_{i,L-L_0-3})$}\label{Pic17}
\end{figure}

Now we need to estimate the difference $\mu(\Psi'_{i+1,L-L_0-3})-\mu(\Psi'_{i,L-L_0-3})$. Observe that by Lemma \ref{psi1}, the common $q$-edge $\bf f$ of the spoke
${\mathcal B}_{i}$ and $\partial\Delta$
 separates at least
$h_{i-1}-h_{i}=m_1\;\;\;\;\theta$-edges of the path ${\bf p}_{i-1,i}$
and $m_2$ ones lying on ${\bf p}_{i,L-L_0-3}$,
where $m_2=h_{i}+h_{i,L-L_0-3}$ by Lemma \ref{psi1} (2)(see fig. \ref{Pic17}).
Since the number of $q$-edges of ${\bf p}={\bf p}(\Psi)$ is less than $3K_0L<J$ by Lemma \ref{2K0},  one
decreases $\mu(\Psi'_{i+1,L-L_0-3})$ at least by $m_1m_2$ when erasing the black
bead on $\bf f$ in the necklace on $\partial\Psi'_{i+1,L-L_0-3}$ by Lemma \ref{mixture} (d,b,c).
(The white
beads of the subpath ${\bf p}_{i,i+1}$ will be moved to the side of ${\cal B}_i$
along $\theta$-bands when one replaces $\partial\Psi'_{i+1,L-L_0-3}$ with the
boundary $\partial\Psi'_{i,L-L_0-3}$ of smaller diagram.)
Hence
$$\mu(\Psi'_{i+1,L-L_0-3}))-\mu(\Psi'_{i,L-L_0-3})\ge m_1m_2$$ $$=(h_{i-1}-h_{i})(h_{i}+h_{L-L_0-3})\ge\frac{1}{10c_5N}h_{i-1}(h_{i}+h_{L-L_0-3})$$
by Lemma \ref{2000N}.
This inequality and Lemma \ref{muJ} applied to $\Psi_{i+1, L-L_0-3}$, imply
$$\mu(\Delta)-\mu(\Psi'_{i,L-L_0-3})=(\mu(\Delta)-\mu(\Psi'_{i+1,L-L_0-3}))+(\mu(\Psi'_{i+1,L-L_0-3})
-\mu(\Psi'_{i,L-L_0-3}))$$ $$\ge -2Jn(h_{i+1}+h_{L-L_0-3})+\frac{1}{10c_5N}h_{i-1}(h_{i}+h_{L-L_0-3})$$
Note that $(h_{i+1}+h_{L-L_0-3})\le 2h_{i+1}$ by (\ref{downup}) and (\ref{L0}). Hence
\begin{equation}\label{mumu}
N_3\mu(\Delta)-N_3\mu(\Psi'_{i,L-L_0-3})\ge -4N_3Jn h_{i+1}+\frac{N_3}{10c_5N}h_{i-1}(h_{i}+h_{L-L_0-3})
\end{equation}
It follows from inequalities \ref{arD} and \ref{mumu} that $$\area_G(\Delta)\le N_4(n+\sigma_{\lambda}(\Delta^*))^2+N_3\mu(\Delta)-N_4\frac{7\delta n}{100c_5N}h_{i+1}-
$$ $$-\frac{N_3}{10c_5N}h_{i-1}(h_{i}+h_{L-L_0-3})+4N_3Jnh_{i+1}+ 30N_2h_i^2+c_7h_i^2$$

Here we come to a contradiction with (\ref{ce}) obtaining inequality $\area_G(\Delta)\le N_4(n+\sigma_{\lambda}(\Delta^*))^2+N_3\mu(\Delta)$, because by the choice of  parameters, $$N_4\frac{7\delta}{100c_5N}>4N_3J,\;\;\;
 \frac{N_3}{10c_5N}>30N_2+c_7\;\;\; and\;\;\; h_{i-1}\ge h_i$$
\endproof}

\subsubsection{A clove with a disk can be removed}

\begin{lemma} \label{led} There exists no counter-example $\Delta$ (see (\ref{ce})), and therefore $\area_G(\Delta)\le N_4(n+\sigma_\lambda(\Delta^*))^2+N_3\mu(\Delta)$ for any weakly minimal diagram $\Delta$ with $|\partial\Delta|=n$.
\end{lemma}
\proof Recall that when proving by contradiction we obtained in Lemma \ref{2000N} that
\begin{equation}\label{j1j}
h_{j+1}< (1-\frac{1}{10c_5N})h_j\; (j=1,\dots,L_0-1),
\end{equation}
and by lemmas \ref{zgh} and \ref{zlh}, we have inequalities
\begin{equation}\label{jk}
|{\mathbf z}_j|_Y\ge h_{j+1}/c_5 \;(j=1,\dots,L_0-1)\;\;\; and \;\;\;|{\mathbf z}_k|_Y\le 2Nh_k \;(k=2,\dots L_0-1).
 \end{equation}

 One can choose an integer
$\rho=\rho(\mmm)$ depending on $c_5$ and $N$ (and so on the $S$-machine $\mmm$ only) so that $(1-\frac{1}{10c_5N})^{\rho}<\frac{1}{6Nc_5}$, and so by (\ref{j1j}, \ref{jk}), we obtain that $h_{j+1}>6Nc_5h_k$
if $k-j-1\ge \rho $. Together with (\ref{j1j}, \ref{jk}, this implies inequalities
$$|{\mathbf z}_j|_Y\ge h_{j+1}/c_5 > 6Nh_k >3|{\mathbf z}_k|\;\; if \;\; k-j-1\ge \rho $$

If $L_0$ is large enough, say $L_0>2000\rho$, one can obtain $1000$ indices
$j_1<j_2<\dots<j_{1000}<L_0$ such that for $i=2,\dots, 1000$, one obtains inequalities
$j_i-j_{i-1}-2>\rho$, and so

\begin{equation}\label{zh}
|{\mathbf z}_{j_{i-1}}|>3|{\mathbf z}_{j_i}|\;\;and\;\; h_{j_{i-1}}\ge h_{j_{i-1}+1}> 6c_5N h_{j_i}
\end{equation}

Let $\ccc\colon\;W\equiv W_0\to\dots\to W_t$  be the computation corresponding to the trapezium $\Gamma_{j_2}$. Since it contains the copy $\Gamma_{{j_2}+1}'$ of $\Gamma_{{j_2}+1}$, which in turn
contains a copy of $\Gamma_{{j_2}+2}$ and so on, we have
some configurations $W(k)$ in $\ccc$ ($k=1,\dots, 999$), that are the labels of some ${\mathbf z}_{i_k}$ (but without superscripts) and $|W(k+1)|_Y>3|W(k)|_Y$ for $k=1,\dots,998$.
If for some $k$ we obtained one-step
subcomputation $W(k)\to\dots\to W(k+4)$, then the statement of Lemma \ref{pol} would give
a subcomputation $W(k+1)\to\dots\to W(k+4)$
contradicting   the statement of Lemma \ref{001}. Hence no five consecutive words $W(k)$-s are configuration
of a one-step subcomputation, and so the number of steps in $W(1)\to\dots \to W(999)$ is at least $100$.

It follows now from Lemma \ref{SH} that the step
history of $\Gamma_{j_2}\backslash \Gamma$, where $\Gamma$ is the copy of $\Gamma_{L_0}$
in $\Gamma_{j_2}$,  has a subword $(34)(4)(45)$ or $(54)(4)(43)$, or $(12)(2)(23)$, or $(32)(2)(21)$.

Let us consider the case $(34)(4)(45)$ (or $(45)(4)(34)$).
Then the
history $H_{j_2+1}$ of $\Gamma_{j_2}$ can be decomposed as $H'H''H'''$, where $H''$ has form $\chi(i-1,i)H_0\chi(i,i+1)$ (the $S$-machine works as ${\mathbf M_3}$) and $||H'||\ge h$ since the height of $\Gamma$ is at least $h$. Moreover, by Lemma \ref{M3} (b), one can choose $i$
so that
$||H'||\ge ||H''||$
since the number of cycles $m$ is large enough.

Since $h_{j_1+1} > 2h_{j_2}$ by (\ref{zh}), the history $H_{j_1+1}$
of $\Gamma_{j_1}$ has a prefix $H'H''H^*$, where $||H^*||=||H'||\ge ||H''||$, and so the $\tt$-spoke ${\mathcal B}_{j_1+1}$ has a $\tt$-subband $\ccc$ starting with $\partial\Pi$
and having the history $H'H''H^*$.

For any factorization $\ccc=\ccc_1\ccc_2\ccc_3$
with $||\ccc_1||+||\ccc_2||\le||\ccc||/3$, the history
of $\ccc_2$ contains the subhistory $H''$, since $||H^*||=||H'||\ge ||H''||$. It follows that $\ccc$ is
a $\lambda$-shaft, because $H''=\chi(i-1,i)H_0\chi(i,i+1)$ and $\lambda<1/3$. The shaft has length
at least $||H'||\ge h$ contrary to Lemma \ref{shaft}.

The case of $(12)(2)(23)$ (of $(23)(2)(12)$) is similar but $H''=\zeta^{i-1,i}H_0\zeta^{i,i+1}$
(the $S$-machine works as $\rhh_m$ and the cycles of $\rhh_m$ have equal lengths by Lemma \ref{prim} (3)).
We come to the final contradiction in this section.
\endproof

\section{Proof of Theorem \ref{t:main}}\label{end}


\subsection{The Dehn function of the group \texorpdfstring{$G$}{G}}\label{end1}

\begin{lemma} \label{big} For every big trapezium $\Delta$, there is a diagram
$\tilde\Delta$ over the finite presentation (\ref{rel1}) - (\ref{rel3}) of $G$ with the same boundary label, such that
the area of $\tilde\Delta$ does not exceed $2\area_G(\Delta)$.
\end{lemma}

\proof Consider the computation $\ccc\colon\; V_0\to\dots\to V_t$ corresponding to $\Delta$
by Lemma \ref{simul}, i.e. $t=h$. According to Definition \ref{abt},
one may assume that $\area_G(\Delta)=c_5h(||V_0||+||V_t||)$ since
otherwise $\tilde\Delta=\Delta$.

$\Delta$ is then covered by $L$ trapezia $\Delta_1,\dots,\Delta_L$
with base $xvx$, where $xv$ (or the inverse word) is a cyclic shift ot the standard base
of $\mmm$. By Lemmas \ref{resto} and \ref{NoAnnul}, all $\Delta_1,\dots,\Delta_L$ are superscript shifts of each other.
Let us apply Lemma \ref{narrow} to any of them, say to $\Delta_1$,
whose top and bottom have labels $W_0$ and $W_t$. If we have Property (1) of that lemma, then the area of $\Delta_1$ does not exceed
$c_4h(||V_0||+||V_t||)$ since every maximal $\theta$-band of $\Delta_1$ has at most $c_4(||V_0||+||V_t||)$ cells in this case.
Hence area of $\Delta$ does not exceed $$Lc_4h(||W_0||+||W_t||)\le 2c_4h(||V_0||+||V_t||)<c_5h(||V_0||+||V_t||)=\area_G(\Delta),$$  i.e.
$\tilde\Delta=\Delta$ in this case too.

Hence one may assume that Property (2) of Lemma \ref{narrow} holds
for $\Delta_1$. By that Lemma, items (b,d), the corresponding cyclic permutations $(W'_0)^{\emptyset}$ and $(W'_t)^{\emptyset}$
are accessible, and so removing the last letters $x$ from $V_0$ and $V_t$ we obtain disk words $V'_0$ and $V'_t$.
For the histories $H'$ and $H''$ of $\ccc((W'_0)^{\emptyset})$ and
$\ccc((W'_t)^{\emptyset})$, Lemma \ref{narrow} gives inequality $||H'||+||H''||\le t$.

Denote by $\Delta_-$ the diagram $\Delta$ without one maximal rim
$x$-band. So $\Delta_-$ has the boundary ${\mathbf p}_1{\mathbf q}_1
{\mathbf p}_2^{-1}{\mathbf q}_2^{-1}$, where $\Lab({\mathbf p}_1)$ and
$\Lab({\mathbf p}_2)$ are disk words
and $\Lab({\bf q}_1)\equiv\Lab({\bf q}_2)$ since the first and the last maximal
two $x$-bands of $\Delta$ are $L$-shifts of each other by Lemma \ref{simul} (1).

If we attach disks $\Pi_1$ and $\Pi_2$ (of radius $\le t$ each) along their boundaries
to the top and the bottom of $\Delta_-$, we obtain a diagram,
whose boundary label is trivial in the free group. Hence there
is a diagram $E$ with two disks whose boundary label is equal
to the boundary label of $\Delta_-$, and the area is less than
$\le 3c_2 t (||V_0'||+||V_t'||)$ by Lemma \ref{fea}. If we attach one $x$-band
of length $t$ to $E$, we construct the required diagram
$\tilde\Delta$ of area at most $$\le 3c_1 t (||V_0||+||V_t||)< c_5h(||V(1)||+||V(2)||)=\area_G(\Delta)$$
\endproof

\begin{lemma} \label{quadr} The Dehn function $d(n)$ of the group $G$ is $O(n^2)$.
\end{lemma}
\proof To obtain the quadratic upper bound for $d(n)$ (with respect to
the finite presentation of $G$ given in Section \ref{gd}),
 it suffices, for every word $W$ vanishing in $G$ with $||W||\le n$, to find a diagram over $G$ of area $O(n^2)$ with boundary label $W$. Since $|W|\le ||W||$, van Kampen's lemma and Lemma \ref{arG} provide us with a minimal diagram $\Delta$ such that
 $\area_G(\Delta)\le N_4(n+\sigma_\lambda(\Delta^*))^2+N_3\mu(\Delta)$ for some constants $N_3$ and $N_4$ depending on the presentation of $G$.
 By Lemmas \ref{wmin} (c), $\sigma_\lambda(\Delta^*)\le cn$, and by Lemma \ref{mixture} (a) and the definition of $\mu(\Delta)$, we have $\mu(\Delta)\le Jn^2$, Thus, we conclude that $\area_G(\Delta)\le C_0n^2$ for some constant $C_0$.

 Recall that in the definition of $G$-area, the subdiagrams, which are big trapezia $\Gamma, \Gamma', \dots,$ can have common cells in their rim $q$-bands only.  By Lemma \ref{big}, any big trapezia $\Gamma$ from this list with top path ${\mathbf p}_1$ and bottom path ${\mathbf p}_2$ can be replaced
 by a diagram $\tilde \Gamma$ with (combinatorial) area at most $2\area_G(\Gamma)$ over the finite presentation (\ref{rel111}, \ref{rel3}).
 When we replace all big trapezia $\Gamma, \Gamma', \dots,$ in this way, we add $q$-bands for the possible intersection of
 big trapezia, but for every $\Gamma$ of height $h$, we add at most $2h$ new cells. So the area of the modified diagram $E$ is at most
 $3\area_G(\Delta)\le 3C_0n^2$. Hence a required diagram is found for the given word $W$.
\endproof

\subsection{The conjugacy problem in \texorpdfstring{$G$}{G}}

Recall that the rule $\theta(23)$ locks all sectors of the standard base of $\mmm$
except for the input sector ${\tilde R}_0{\tilde P}_1$ and its mirror copy. Hence every
$\theta(23)^{-1}$-admissible word has the form $W(k,k')\equiv w_1\alpha^kw_2(\alpha')^{-k'}w_3$,
where $k$ and $k'$ are integers and $w_1, w_2, w_3$ are fixed word in state letters;
$w_1$ starts with $\tt$.

\begin{lemma} \label{conj} A word $W(k,k)$ is a conjugate of the word $W_{ac}$
in the group $G$ (and in the group $M$) if and only if the input $\alpha^k$
is accepted by the Turing machine $\mmm_0$.
\end{lemma}

\proof Let the Turing machine $\mmm_0$ accept $\alpha^k$.
Then by Lemma \ref{I6A6}, we have an accepting computation $\ccc$ of $\mmm$ starting with $W(k,k)$ and ending with $W_{ac}$. By Lemma \ref{simul}, one can construct a corresponding trapezium $\Delta$.
Since the computation $\ccc$ uses neither the rules of Step 1, nor the rules of Step 2, nor the rules $\theta(23)^{\pm 1}$, the labels of the  edges of $\Delta$ have no superscripts. Hence the bottom path of $\Delta$ is labeled
by $W(k,k)$, the top label is $W_{ac}$ and the sides of $\Delta$ have equal labels since the $S$-machine
$\mmm$ have cyclic standard base. It follows from van Kampen
Lemma that the words  $W(k,k)$ and $W_{ac}$ are conjugate in the group $M$, as required.

\medskip

For the converse statement, we assume that the words $W(k,k)$ and $W_{ac}$ are conjugate in $G$.
Recall that the definition of annular diagram $\Delta$ over a group $G$ is similar to the definition of van Kampen diagram, but the complement of $\Delta$ in the plane has two
connected components. So $\Delta$ has two boundary components.
By the van Kampen-Schupp lemma (see \cite{LS}, Lemma 5.2 or \cite{book}, Lemma 11.2) there is an  annular diagram $\Delta$ whose boundary components ${\mathbf p}_1$ and ${\mathbf p}_2$ have clockwise labels $W(k,k)$ and $W_{ac}$. As for van Kampen diagrams (see Subsection \ref{di}), one may assume that $\Delta$ is a minimal diagram and there are no two disks in $\Delta$ connected
by two $\tt$-spokes $\cal B$ and $\cal C$ provided there are neither disks nor boundary components
of $\Delta$ between $\cal B$ and $\cal C$. This property makes the disk graph of $\Delta$
hyperbolic as in Subsection \ref{di}: if $\Delta$ has a disk, then there is a disk with
at least $L/2$ $\tilde t$-spokes ending on $\partial\Delta$ (see Corollary 10.1 in \cite{book}).

However each of ${\mathbf p}_1$ and ${\mathbf p}_2$ has only one $\tt$-edge, and it follows
that $\Delta$ has no disks since $L/2>2$. Hence a unique maximal $\tt$-band $\mathcal B$ of $\Delta$ has to connect
these $\tt$-edges. Cutting $\Delta$ along a side $\bf q$ of $\mathcal B$, we obtain a reduced van
Kampen diagram $\Gamma$ over the group $M$. Its boundary path is ${\mathbf p}_1{\mathbf q}{\mathbf p}_2^{-1}{\mathbf q}'^{-1}$,
where $\Lab({\mathbf q}')\equiv\Lab({\mathbf q})$. The maximal $\theta$-bands of $\Gamma$ connect
${\mathbf q}$ and ${\mathbf q}'$ since they cannot cross a $q$-band twice by Lemma \ref{NoAnnul}.
Hence $\Gamma$ is a trapezium with top path ${\mathbf p}_1$ and bottom path ${\mathbf p}_2$. The base
of $\Gamma$ is standard since the top/bottom labels have standard base.

The equality $\Lab({\mathbf q}')\equiv\Lab({\mathbf q})$ implies that the side edges have
no superscripts because $\Lab({\mathbf q}')$ has to be a $\pm 1$-shift of $\Lab({\mathbf q})$.
It follows from Lemma \ref{simul} and the definition of $(\theta,q)$-relations that $\Gamma$ corresponds
to a reduced computation $\ccc\colon W(k,k)\to\dots\to W_{ac}$  having no rules of Steps 1,2 and no $\theta(23)^{\pm 1}$. Therefore the word $\alpha^k$ is accepted by $\mmm_0$ by Lemma \ref{I6A6} (2).
\endproof

{\bf Proof of Theorem \ref{t:main}.} Since the language accepted by the Turing machine $\mmm_0$ is non-recursive, the conjugacy problem is undecidable for the group $G$
by Lemma \ref{conj}. The Dehn function of $G$ is at most quadratic by Lemma \ref{quadr}.
To obtain a lower quadratic estimate, it suffices to see that if a $\theta$-letter
$\theta$ and a $Y$-letter $a$ commute,  then by Lemmas \ref{extdisk} and \ref{NoAnnul}, the area of the word $a^n\theta^na^{-n}\theta^{-n}$ is equal to $n^2$ (or to use \cite{Bow}: every non-hyperbolic finitely presented group has at least quadratic Dehn function). The theorem
is proved. $\Box$

\addcontentsline{toc}{section}{Index}

\printindex[g]

\vskip .5 in

\noindent Alexander Yu. Ol'shanskii,\\ Department of Mathematics,
Vanderbilt University
and\\ Department of
Higher Algebra, MEHMAT,
 Moscow State University
 \\ alexander.olshanskiy@vanderbilt.edu

\vskip .2 in

\noindent Mark Sapir, \\ Department of Mathematics,
Vanderbilt University \\ m.sapir@vanderbilt.edu

\end{document}